\documentclass[a4paper,11pt,twoside,openright]{book}

\usepackage{amsmath}
\usepackage{amsthm}
\usepackage{amssymb}
\usepackage{amsfonts}
\usepackage{mathrsfs}
\usepackage{color}
\usepackage[utf8]{inputenc}
\usepackage{enumitem}
\usepackage{multicol}
\usepackage[english]{babel}
\usepackage[babel]{csquotes}
\usepackage{mathtools}
\usepackage{rotating}
\usepackage{latexsym}
\usepackage{epsfig}
\usepackage{verbatim}
\usepackage{thmtools,thm-restate}
\usepackage{listings} 
\usepackage[dvipsnames]{xcolor}
\usepackage[toctitles]{titlesec}
\usepackage{floatrow}
\usepackage{pgf,tikz}
\usetikzlibrary{arrows, plotmarks,calc,intersections,petri,topaths,babel,shapes.geometric}
\usetikzlibrary{decorations}

\usepackage{makecell,booktabs,amsxtra}

\usepackage{adjustbox}
\usepackage{ifthen}
\usepackage[position=top]{subfig}
\usepackage{emptypage}
\usepackage{pgfplots}
\pgfplotsset{xticklabels={}, yticklabels={} width=7cm,compat=1.8}
\usepackage{pgfplotstable}
\pgfmathsetseed{1138} 
\pgfplotstableset{ 
	create on use/x/.style={create col/expr={42+2*\pgfplotstablerow}},
	create on use/y/.style={create col/expr={(0.6*\thisrow{x}+130)+5*rand}}
}
\pgfplotstablenew[columns={x,y}]{30}\loadedtable

\usepackage[ruled]{algorithm2e}

\SetAlFnt{\small}
\SetAlCapFnt{\small}
\SetAlCapNameFnt{\small}
\SetAlCapHSkip{0pt}
\IncMargin{-\parindent}

\usepackage{url}
\usepackage[square,numbers,sort]{natbib}

\usepackage{pgfplots,pgfplotstable}
\pgfplotsset{compat=1.8}

\usepackage[bookmarks=false,colorlinks=true,linkcolor=blue,urlcolor=blue,citecolor=red,pdfborder={0 0 0.5}]{hyperref}

\usepackage{etoolbox}

\apptocmd{\thebibliography}{\csname phantomsection\endcsname\addcontentsline{toc}{chapter}{\bibname}}{}{}

\apptocmd{\listoffigures}{\csname phantomsection\endcsname\addcontentsline{toc}{chapter}{List of Figures}}{}{}
\apptocmd{\listofalgorithms}{\csname phantomsection\endcsname\addcontentsline{toc}{chapter}{List of Algorithms}}{}{}

\theoremstyle{plain}
\newtheorem{theorem}{Theorem}[chapter]
\newtheorem{lemma}[theorem]{Lemma}

\newtheorem{corollary}[theorem]{Corollary}
\newtheorem{conjecture}[theorem]{Conjecture}

\declaretheoremstyle[
notefont=\bfseries, notebraces={}{},
bodyfont=\normalfont\itshape,
headformat=\NAME \NOTE
]{nopar}
\declaretheorem[style=nopar,name=Corollary]{corollary*}
\declaretheorem[style=nopar,name=Theorem]{theorem*}

\newtheorem{manualtheoreminner}{Theorem}
\newenvironment{manualtheorem}[1]{%
  \renewcommand\themanualtheoreminner{#1}%
  \manualtheoreminner
}{\endmanualtheoreminner}

\theoremstyle{definition}
\newtheorem{definition}[theorem]{Definition}
\newtheorem{remark}[theorem]{Observation}
\newtheorem{observation}[theorem]{Observation}
\newtheorem{question}[theorem]{Problem}

\newcount \mycount
\newcount \mycountb

\renewcommand{\b}[1]{\boldsymbol{#1}}
\DeclareMathOperator{\N}{\mathbb N}
\DeclareMathOperator{\Z}{\mathbb Z}

\DeclareMathOperator{\R}{\mathbb R}
\DeclareMathOperator{\csp}{CSP}
\newcommand{\cocsp}{\overline{\operatorname{CSP}}}
\DeclareMathOperator{\dom}{Dom}

\DeclareMathOperator{\A}{{\mathbb A}}
\DeclareMathOperator{\B}{{\mathbb B}}
\DeclareMathOperator{\C}{{\mathbb C}}
\DeclareMathOperator{\D}{{\mathbb D}}
\DeclareMathOperator{\F}{{\mathbb F}}
\DeclareMathOperator{\K}{{\mathbb K}}
\DeclareMathOperator{\T}{{\mathbb T}}
\renewcommand{\S}{{\mathbb S}}
\DeclareMathOperator{\G}{{\mathbb G}}
\DeclareMathOperator{\Hb}{{\mathbb H}}
\renewcommand{\H}{\Hb}
\DeclareMathOperator{\Idemp}{{\operatorname{\mathbb Idem}}}
\renewcommand{\P}{\mathbb P}
\DeclareMathOperator{\Ord}{\operatorname{{\mathbb Ord}}}
\DeclareMathOperator{\OrdDG}{\T_4^{\setminus(03)}}
\DeclareMathOperator{\Zig}{\mathbb Zig}
\DeclareMathOperator{\ZigCyc}{Cyc\mathbb Zig}
\DeclareMathOperator{\Instance}{\mathbb I}
\DeclareMathOperator{\Aut}{Aut}
\DeclareMathOperator{\End}{End}
\DeclareMathOperator{\Emb}{Emb}
\DeclareMathOperator{\Hom}{Hom}
\DeclareMathOperator{\Pol}{Pol}
\DeclareMathOperator{\Ind}{Ind}
\DeclareMathOperator{\lcm}{lcm}

\DeclareMathOperator{\AC}{AC}

\DeclareMathOperator{\Maltsev}{{\operatorname{Malt}}}
\DeclareMathOperator{\Majority}{{\operatorname{Maj}}}
\DeclareMathOperator{\Siggers}{{\operatorname{Sigg}}}
\newcommand{\KK}[1]{\operatorname{KK}(#1)}
\newcommand{\HMcK}[1]{\operatorname{HMcK}(#1)}
\newcommand{\HM}[1]{\operatorname{HM}(#1)}
\newcommand{\GFS}[1]{\operatorname{Ele}(#1)}
\newcommand{\NN}[1]{{\operatorname{NN}(#1)}}
\newcommand{\J}[1]{{\operatorname{J}(#1)}}
\newcommand{\TS}[1]{{\operatorname{TS}(#1)}}
\newcommand{\NU}[1]{{\operatorname{NU}(#1)}}
\newcommand{\WNU}[1]{{\operatorname{WNU}(#1)}}
\newcommand{\Const}{\operatorname{Const}}

\newcommand{\tominor}{\stackrel{\operatorname{minor}}{\to}}

\providecommand{\dotdiv}{
  \mathbin{
    \vphantom{+}
    \text{
      \mathsurround=0pt 
      \protect\ooalign{
        \noalign{\kern-.45ex}
        \hidewidth$\smash{\cdot}$\hidewidth\cr 
        \noalign{\kern.45ex}
        $-$\cr 
      }%
    }%
  }%
}
\providecommand{\cupdot}{
  \mathbin{
    \vphantom{+}
    \text{
      \mathsurround=0pt 
      \protect\ooalign{
        \noalign{\kern-.4ex}
        \hidewidth$\smash{\cdot}$\hidewidth\cr 
        \noalign{\kern.4ex}
        $\cup$\cr 
      }%
    }%
  }%
}
\DeclareMathOperator{\median}{median}

\newcommand{\Flat}{\operatorname{rad}}

\newcommand{\PL}{\ensuremath{\operatorname{PCL}}}
\newcommand{\Downset}{{\Gamma}}

\newcommand{\FD}{\ensuremath{\mathcal F_D(\omega)}}

\newcommand{\homeq}{\mathrel{\substack{\textstyle\rightarrow\\[-0.4ex]
                      \textstyle\leftarrow}}}
\newcommand{\ppeq}{\mathrel{\equiv}}
\newcommand{\ppleq}{\mathrel{\leq}}
\newcommand{\ppstrictlyLess}{<}
\newcommand{\ppgeq}{\mathrel{\geq}}
\newcommand{\ppstrictlyGreater}{>}
\DeclareMathOperator{\mPCL}{M_{\operatorname{PCL}}}
\DeclareMathOperator{\mCL}{M_{\operatorname{CL}}}
\newcommand{\MPCL}{\mathfrak M_{\operatorname{PCL}}}
\newcommand{\MCL}{\mathfrak M_{\operatorname{CL}}}
\DeclareMathOperator{\PCPoset}{\mathfrak P_{\operatorname{UPC}}}
\DeclareMathOperator{\SDPoset}{\mathfrak P_{\operatorname{UC}}}
\DeclareMathOperator{\PPPoset}{\mathfrak P_{\operatorname{fin}}}
\DeclareMathOperator{\DGPoset}{\mathfrak P_{\operatorname{Digraphs}}}
\newcommand{\shiftTuple}[2]{{#1}_{#2}}
\newcommand{\Cyc}[1]{\C_{#1}}

\DeclareMathOperator{\AND}{\wedge}
\DeclareMathOperator{\bigAND}{\bigwedge}

\newcommand\xqed[1]{%
  \leavevmode\unskip\penalty9999 \hbox{}\nobreak\hfill
  \quad\hbox{#1}}
\newcommand\eoe{\xqed{$\triangle$}}

\def\ifinlinemath#1#2{#2}
\def\mathshift{$}
\catcode`$=13
\def${\ifinner\expandafter\mathshift\else\expandafter\myshift\fi}
\def\myshift#1${{\def\ifinlinemath##1##2{##1}\raisebox{0ex}[0ex][0ex]{\mathshift#1\mathshift}}}
\AtBeginDocument{\catcode`$=13}

\theoremstyle{definition}
\newtheorem{examplex}[theorem]{Example}
\newenvironment{example}
{\pushQED{\qed}\examplex}
{\popQED\endexamplex}

\theoremstyle{definition}

\makeatletter
\newcommand\footnoteref[1]{\protected@xdef\@thefnmark{\ref{#1}}\@footnotemark}
\makeatother

\makeatletter
\newcommand*{\rom}[1]{\expandafter\@slowromancap\romannumeral #1@}
\makeatother

\makeatletter
\def\blfootnote{\xdef\@thefnmark{}\@footnotetext}
\makeatother

\usepackage{setspace}

\usepackage{epigraph}
\newcommand{\edges}[1]{E({#1})}

\DeclareMathOperator{\dist}{dist}
\DeclareMathOperator{\levelfunction}{\height}
\newcommand{\level}[1]{\levelfunction({#1})}
\DeclareMathOperator{\id}{id}

\DeclareMathOperator{\KMM}{KMM}

\newcommand{\nCk}[3]{\C^{#1,#2}_{#3}}
\DeclareMathOperator{\height}{\operatorname{lvl}}

\newcommand{\OriP}[1]{\P(#1)}
\newcommand{\OriC}[1]{\C(#1)}

\newcommand{\prog}[1]{\mathcal #1}
\newcommand{\TLinP}{\ensuremath{{\mathrm{3Lin}}_p}}
\newcommand{\ModL}[1]{\ensuremath{{\mathrm{Mod}}_{#1}{\mathrm{L}}}}
\newcommand{\NModL}[1]{\ensuremath{{\mathrm{NMod}}_{#1}{\mathrm{L}}}}
\DeclareMathOperator{\HornSAT}{Horn-3SAT}
\DeclareMathOperator{\Pclass}{P}
\DeclareMathOperator{\NP}{NP}
\DeclareMathOperator{\NL}{NL}
\DeclareMathOperator{\Lclass}{L}

\usepackage{siunitx}
\newcommand{\codelink}{\url{https://github.com/WhatDothLife/TheSmallestHardTrees}}

\newcommand{\tikzKThree}{
\tikz{
\node at (-30:0.2) [circle, fill, scale=0.3] (0) {};
\node at (90:0.2) [circle, fill, scale=0.3] (1) {};
\node at (210:0.2) [circle, fill, scale=0.3] (2) {};
\path 
    (0) edge (1)
    (1) edge (2)
    (2) edge (0);
}}

\newcommand{\drawDot}{
\tikz[baseline=-0.8mm]{
\clip(-0.15,-0.053) rectangle (0.15,0.35);
\node[circle, fill, scale=0.3] at (0,0) (0) {};
}}

\newcommand{\drawEdge}{
\tikz[baseline=-1mm]{
\clip(-0.15,-0.053) rectangle (0.65,0.35);
\node[circle, fill, scale=0.3] at (0,0) (0) {};
\node[circle, fill, scale=0.3] at (0.5,0) (1) {};
\path (0) edge[->,>=stealth'] (1);
}}

\newcommand{\drawUndirectedEdge}{
\tikz[baseline=-1mm]{
\clip(-0.15,-0.053) rectangle (0.65,0.35);
\node[circle, fill, scale=0.3] at (0,0) (0) {};
\node[circle, fill, scale=0.3] at (0.5,0) (1) {};
\path (0) edge (1);
}}

\newcommand{\drawLoop}{
\tikz{
\clip(-0.15,-0.053) rectangle (0.15,0.35);
\node[circle, fill, scale=0.3] at (0,0) (0) {};
\path (0) edge[out=120, in=60, looseness=20,->,>=stealth'] (0);
}}

\newcommand{\drawKThree}{
\tikz{
\node at (-30:0.2) [circle, fill, scale=0.3] (0) {};
\node at (90:0.2) [circle, fill, scale=0.3] (1) {};
\node at (210:0.2) [circle, fill, scale=0.3] (2) {};
\path 
    (0) edge (1)
    (1) edge (2)
    (2) edge (0);
}}

\newcounter{anzahl}
\def\nodeDist{0.4}

\newcommand{\tikzpathFreeVariables}[1]{
\begin{tikzpicture}[scale=0.5]
\setcounter{anzahl}{0}

\foreach \y/\c [count=\xi from 1] in #1 {
    \ifthenelse{\c=0}{
    \node[var-f] (\xi) at (\nodeDist*\xi,\nodeDist*\y) {};}{
    \node[var-b] (\xi) at (\nodeDist*\xi,\nodeDist*\y) {};}
    
    \setcounter{anzahl}{\xi}
    }
    
\foreach \y [count=\yi from 1,count=\yii from 2] in {2,...,\arabic{anzahl}} 
    \path (\yi) edge (\yii);
\end{tikzpicture}
}

\newcommand{\toEdge}{\mathbin
{\begin{tikzpicture}[baseline =-1mm]
    \path[>=stealth',->] (0,0) edge (0.37,0);
\end{tikzpicture}}}
\newcommand{\fromEdge}{\mathbin
{\begin{tikzpicture}[baseline =-1mm]
    \path[>=stealth',<-] (0,0) edge (0.37,0);
\end{tikzpicture}}}
\tikzstyle{var-b}=[circle,fill,draw=white,inner sep=0pt,minimum size=3.2pt]
\tikzstyle{bullet}=[circle,fill,draw=white,inner sep=0pt,minimum size=3.2pt]
\tikzstyle{var-f}=[circle,draw,inner sep=0pt,minimum size=3pt]

\usetikzlibrary{arrows,arrows.meta}
\usetikzlibrary{bending,calc}
\newcommand{\cycle}[2]{
\def \n {#1}
\def \radius {{#1*0.035+0.15}}
\def \margin {3.8/(#1*0.035+0.15)} 
\node[scale={0.8+#1/10}] at #2 {#1};

\foreach \s in {1,...,\n}
{
  \draw #2+({360/\n * (\s - 1)}:\radius) circle (1pt);
  \draw[>=stealth',arrows=-{>[bend]}] #2+({360/\n * (\s - 1)+\margin}:\radius) 
    arc ({360/\n * (\s - 1)+\margin}:{360/\n * (\s)-\margin}:\radius);
}}

\newcommand{\cycleEmpty}[3]{
\def \n {#1}
\def \radius {{#1*0.035+0.15}}
\def \margin {3.8/(#1*0.035+0.15)} 

\foreach \s in {1,...,\n}
{
  \draw #2+({90+360/\n * (\s - 1)}:\radius) circle (1pt);
}
\foreach \s in #3
{
    \draw[>=stealth',arrows=-{>[bend]}] #2+({90+360/\n * (\s - 1)+\margin}:\radius) 
    arc ({90+360/\n * (\s - 1)+\margin}:{90+360/\n * (\s)-\margin}:\radius);
}
}

\begin{document}
\frontmatter

\begin{titlepage}
\begin{center} \centering

\vspace*{.5cm}
\begin{spacing}{2}
{\huge \bfseries 
Digraphs modulo primitive positive constructability}
\end{spacing}
 
  \vspace{1cm}
    \makebox[0.7\textwidth][s]{\LARGE \textbf{D I S S E R T A T I O N}}

    \vspace{1cm}
   {\LARGE zur Erlangung des akademischen Grades}

    \vspace{0.5cm}
   \LARGE {Doctor of Philosophy\\
    (Ph.\ D.)}

    \vspace{0.3cm}
     {\normalsize vorgelegt}

    \vspace{0.3cm}
   \LARGE{ dem Bereich Mathematik und Naturwissenschaften \\der Technischen Universit\"at Dresden}

    \vspace{0.3cm}
{\normalsize von}

    \vspace{0.3cm}
     {\LARGE   Dipl.-Inf.\ Florian Starke}

\vspace{1.5cm}
\Large{ Eingereicht am 22.\ Januar 2024\\
  Verteidigt am 18.\ April 2024
 }

      \vspace{.7cm}
\large {Gutachter: \\
	\normalsize { Prof.\ Dr.rer.nat.\ Manuel Bodirsky \\
		Prof.\ Ph.D.\  Víctor Dalmau}}

\vspace{1cm}
   \normalsize {Die Dissertation wurde in der Zeit von Oktober 2018 bis\\ Januar 2024 am Institut f\"ur Algebra angefertigt.}

\end{center}

\end{titlepage}


\cleardoublepage\thispagestyle{empty}

\noindent%
This work is licensed under the Creative Commons Attribution-ShareAlike
4.0 International
License. To view a copy of this license, visit
\url{http://creativecommons.org/licenses/by-sa/4.0/deed.en_GB}.
\cleardoublepage

\thispagestyle{empty}
\vskip20\baselineskip

\epigraph{``But Uli, it is just a PhD.”
\vskip \baselineskip
 --- \small{\textup{Unknown author}{}}
}

\cleardoublepage
\chapter*{Acknowledgements}
Where to even start? 
There are so many people that have supported me over  the last years to whom I am deeply grateful. First and foremost I want to thank my supervisor Manuel Bodirsky, who one beautiful summer day drove by me on his red bike (admittedly, I am not sure whether he already drove his red bike back then) hopped off and asked me whether I wanted to become his PhD student. Thank you Manuel for the endless discussions of which I only understood half for the first couple of years. 

I am deeply grateful to Albert, my academic brother, with whom I shared an office for many years and who introduced me to the Poset. We also worked on our first big project together: The cycles paper. Dear Albert, I hope that in the future we can work together again to bring more light into dark places (of the Poset).
I also want to thank my other office mates Erkko, Marga, and Sebastian. Two of whom also became co-authors of mine.  

Thank you Ulrike, Antje, and Christian, who where my first contacts to the world of university mathematics and with whom I now teach new generations of students.
Thank you Uli, for all our very motivating conversations. Although my brain shuts down whenever I hear monoidal category I always left our meetings energized and motivated to do more mathematics. Thank you for giving me the opportunity to teach with you in Ghana, it was a wonderful time.  
Thanks to Frau Sturm, who encouraged me to dare to try out giving exercise classes. I am sorry to tell you that I 
still have not written my Master Thesis and probably never will.

Finally, I want to thank my fiancée and my friends and family who have supported me over all those years.

\newpage



\begin{minipage}{.476\textwidth}
	\includegraphics[width=.7\textwidth]{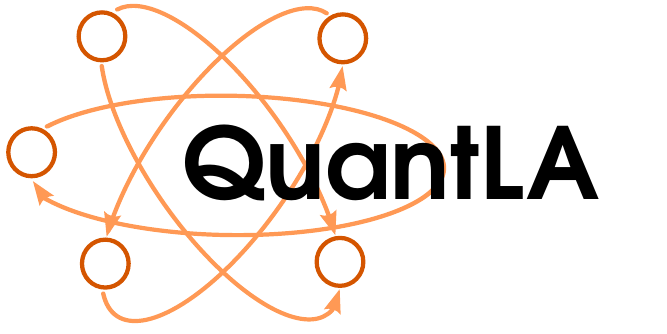}
\end{minipage}\begin{minipage}{.476\textwidth}
	\includegraphics[width=\textwidth]{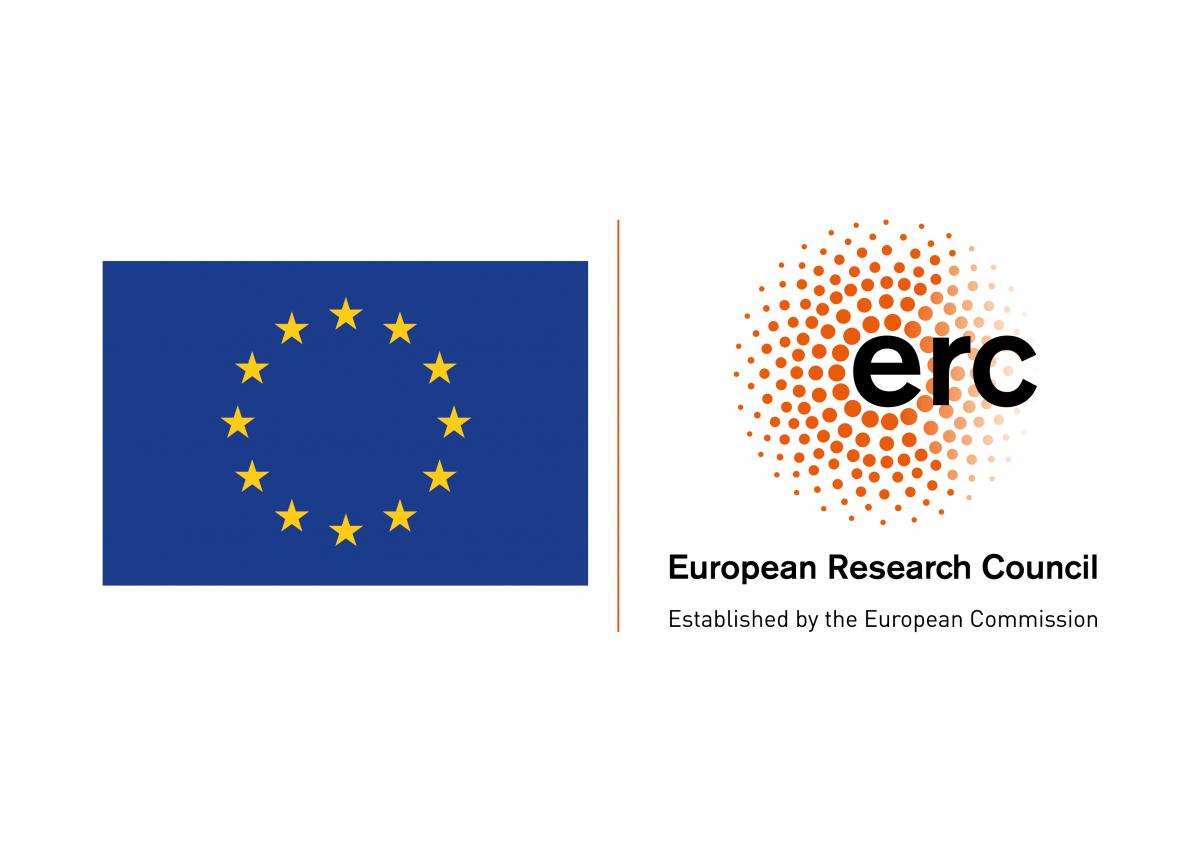}
\end{minipage}

The author was supported by DFG Graduiertenkolleg 1763 (QuantLA) and the European Union  (ERC, POCOCOP, 101071674) and  (grant agreement No. 681988, ``CSP-Infinity").
Views and opinions expressed are however those of the author only and do not necessarily reflect those of the European Union or the European Research Council Executive Agency. Neither the European Union nor the granting authority can be held responsible for them.

\clearpage
\tableofcontents



\mainmatter

\chapter*{Introduction}
\markboth{INTRODUCTION}{INTRODUCTION}
\phantomsection
\addcontentsline{toc}{chapter}{Introduction}

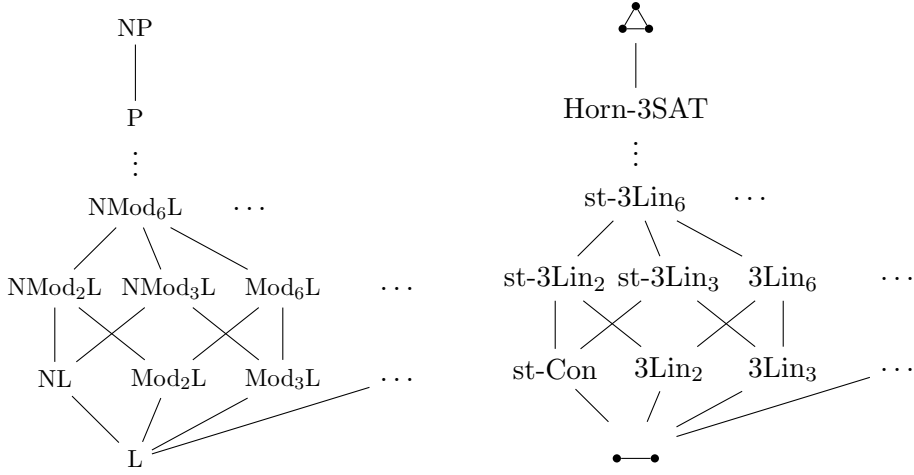
\begin{figure}
    \centering
    \begin{tikzpicture}[node distance = 1.5cm]
    \node[scale=0.85] (NP) {NP};
    
    \node[below of = NP,yshift=3mm,scale=0.85] (P) {P};
    
    \node[below of = P,yshift=10mm]  {$\vdots$};
    
    \node[below of = P,yshift=3mm,scale=0.85] (NMod23) {$\NModL{6}$};
    \node[right of = NMod23] (Mod235) {$\dots$};
    
    \node[below left of = NMod23,scale=0.85] (NMod2) {$\NModL{2}$};
    \node[right of = NMod2,scale=0.85] (NMod3) {$\NModL3$};
    \node[right of = NMod3,scale=0.85] (Mod23) {$\ModL{6}$};
    \node[right of = Mod23] (Mod25) {$\dots$};

    \node[below of = NMod2,yshift=3mm,scale=0.85] (NL) {NL};
    \node[right of = NL,scale=0.85] (Mod2) {$\ModL2$};
    \node[right of = Mod2,scale=0.85] (Mod3) {$\ModL3$};
    \node[right of = Mod3] (Mod5) {$\dots$};
    
    \node[below right of = NL,scale=0.85] (L) {L};
    
    \path 
        (P) edge (NP)
        
        (NL) edge (NMod2)
        (NL) edge (NMod3)
        (Mod2) edge (NMod2)
        (Mod2) edge (Mod23)
        (Mod3) edge (NMod3)
        (Mod3) edge (Mod23)
        
        (NMod2) edge (NMod23)
        (NMod3) edge (NMod23)
        (Mod23) edge (NMod23)
        
        
        (NL) edge (L)
        (Mod2) edge (L)
        (Mod3) edge (L)
        (Mod5) edge (L)
        ;
    \end{tikzpicture}\hspace{8mm}
    \begin{tikzpicture}[node distance = 1.5cm]
    \node (NP) {\drawKThree};
    \node[below of = NP,yshift=3mm] (P) {$\HornSAT$};

    \node[below of = P,yshift=10mm]  {$\vdots$};
    
    \node[below of = P,yshift=3mm] (NMod23) {st-3Lin$_{6}$};
    \node[right of = NMod23] (Mod235) {$\dots$};
    
    \node[below left of = NMod23] (NMod2) {st-3Lin$_{2}$};
    \node[right of = NMod2] (NMod3) {st-3Lin$_3$};
    \node[right of = NMod3] (Mod23) {3Lin$_{6}$};
    \node[right of = Mod23] (Mod25) {$\dots$};

    \node[below of = NMod2,yshift=3mm] (NL) {st-Con};
    \node[right of = NL] (Mod2) {3Lin$_2$};
    \node[right of = Mod2] (Mod3) {3Lin$_3$};
    \node[right of = Mod3] (Mod5) {$\dots$};
    
    \node[below right of = NL] (L) {\drawUndirectedEdge};
    
    \path 
        (P) edge (NP)
        
        (NL) edge (NMod2)
        (NL) edge (NMod3)
        (Mod2) edge (NMod2)
        (Mod2) edge (Mod23)
        (Mod3) edge (NMod3)
        (Mod3) edge (Mod23)
        
        (NMod2) edge (NMod23)
        (NMod3) edge (NMod23)
        (Mod23) edge (NMod23)
        
        
        (NL) edge (L)
        (Mod2) edge (L)
        (Mod3) edge (L)
        (Mod5) edge (L)
        ;
    \end{tikzpicture}
    \caption{Complexity classes from the Larose-Tesson Conjecture ordered by inclusion (left) and structures whose CSPs are complete for the respective complexity class (right).}
    \label{fig:complexityZoo}
\end{figure}

One of the most fundamental questions in theoretical computer science is whether P equals NP, i.e., whether every decision problem for which a solution can be verified in polynomial time (NP) can also be solved in polynomial time (P). Ladner's theorem states that if P is not equal to NP then there are problems in NP that are neither in P nor NP-hard. However, the NP-intermediate  problems constructed in Ladner's proof are very artificial. 
In this thesis we will focus on a class of decision problems that is  tame enough so that it does not contain NP-intermediate problems, while still being large enough to remain interesting, namely the class of \emph{constraint satisfaction problems} (CSPs).
This class consists of decision problems which can be expressed by a set of constraints on a set of variables. More formally, a CSP over a structure $\A$ with a finite relational signature is the problem of deciding whether a given finite structure has a homomorphism into $\A$. CSPs capture many natural and practical problems, such as the  reachability problem for digraphs, the 3-colorability problem for graphs, and the satisfiability problem for Boolean formulas. In 1999, Feder and Vardi conjectured that there are no NP-intermediate CSPs over finite structures, i.e., every CSP over a finite structure is in P or NP-complete.  This conjecture has been proven independently by Bulatov and by Zhuk in 2017~\cite{ZhukFVConjecture,BulatovFVConjecture}. 
Interestingly, this dichotomy has a logical counterpart: the finite structures with NP-complete CSP are precisely those that \emph{primitively positively construct} (pp-construct) $\tikzKThree$, the complete graph on three vertices (unless $\text{P}= \text{NP}$). 
One connection between pp-constructability and the complexity of CSPs used for this dichotomy is the following: Whenever a finite structure $\A$ pp-constructs a finite structure $\B$, the CSP of $\B$ has a log-space reduction to the CSP of $\A$~\cite{wonderland}.
It turns out that pp-constructability is a quasi-order on the class of all finite structures. 
Hence, the poset of CSPs ordered by log-space reductions is a coarsening of the poset of finite structures ordered by  pp-constructability (called \emph{pp-constructability poset} or $\PPPoset$ for short). Understanding the poset $\PPPoset$ can also lead to a better understanding of CSPs within sublasses of P, like L and NL. 
One long-standing open problem in this area is a conjecture by Larose and Tesson which states that every CSP over a finite structure is in L, NL-hard, or Mod$_p$L-hard for some prime $p$~\cite{LaroseTesson}.
The hardness is always considered under log-space reductions. Again, the conjecture comes with a logical counterpart: 
\begin{conjecture}[\cite{LaroseTesson}]\label{conj:LaroseTesson}
Let $\A$ be a finite structure. Then the following are true.
\begin{enumerate}
    \item If $\A$ can neither pp-construct $\HornSAT$ nor $\TLinP$ for any prime $p$, then $\csp(\A)$ is in $\NL$.
    \item If $\A$ can neither pp-construct $\operatorname{st-Con}$ nor $\TLinP$ for any prime $p$, then $\csp(\A)$ is in $\Lclass$.
\end{enumerate}
    
\end{conjecture}
This conjecture implies in particular that every CSP over a finite template in NL that is not NL-hard is in L, i.e., there are no L-NL intermediate CSPs over finite templates. It also implies that there are no L-$\ModL p$ intermediate CSPs over finite templates.
However, if $\Pclass$, $\NL$, $\ModL2$, and $\operatorname{NC}$ are distinct, then
the CSP of the disjoint union of st-Con and $\operatorname{3Lin}_2$ is $\NL$-hard, not in $\NL$, in $\operatorname{NC}$, and not $\Pclass$-hard. Hence it is an $\NL$-$\Pclass$ intermediate problem. 
Our next goal is to give a conjecture for the full picture of the complexity of CSPs over finite templates. To do this, we first need a few more definitions.
Let $n\in\N$, $n>1$.
Assume that st-Con and $\operatorname{3Lin}_n$ have distinct domains and distinct signatures. 
Define the structure $\operatorname{st-3Lin}_n$ as the disjoint union of st-Con and $\operatorname{3Lin}_n$. 
Define the complexity classes $\ModL n$ and $\NModL n$ thought their complete problems $\operatorname{st-3Lin}_n$ and  $\operatorname{3Lin}_n$, respectively. Define $\ModL1\coloneqq\Lclass$ and $\NModL 1\coloneqq\NL$.  
It is open whether the following stronger version of the Larose-Tesson Conjecture is true.
\begin{question}
Let $\A$ be a finite structure and let $C$ be the set of all primes $p$ such that $\A$ can pp-construct $\TLinP$. Then the following are true:
\begin{enumerate}
    \item If $|C|=\infty$, then $\A$ can pp-construct $\drawKThree$ and $\csp(\A)$ is $\NP$-hard.
    \item If $|C|<\infty$ and $\A$ can pp-construct $\HornSAT$, then $\csp(\A)$ is $\Pclass$-complete.
    \item If $|C|<\infty$ and $\A$ cannot pp-construct $\HornSAT$, then let $n$ be the product of all primes in $C$;
    \begin{enumerate}
        \item if $\A$ can pp-construct $\operatorname{st-Con}$, then $\csp(\A)$ is $\NModL n$-complete.
        \item if $\A$ cannot pp-construct $\operatorname{st-Con}$, then $\csp(\A)$ is $\ModL n$-complete.
    \end{enumerate}
    
\end{enumerate}
\end{question}
It is well known that (1) and (2) are true. Note that the statement (3) for $C=\emptyset$ (in this case $n=1$) is equivalent to Conjecture~\ref{conj:LaroseTesson}.
The complexity classes and structures from the open problem are depicted in Figure~\ref{fig:complexityZoo}. 


A powerful tool when working with pp-constructions are \emph{minor conditions}.
A finite structure $\A$ pp-constructs a finite structure $\B$ if and only if every
\emph{minor condition} satisfied by the polymorphism clone of $\A$ is also satisfied by the polymorphism clone of $\B$~\cite{wonderland}; for definitions, see
Section~\ref{sec:h1}. 
In particular, if $\A$ pp-constructs $\B$ and $\A$ has a 4-ary \emph{Siggers} polymorphism, i.e., a homomorphism $f\colon\A^4\to \A$ satisfying
\[f(a,r,e,a)= f(r,a,r,e)\text{ for all }e,a,r\in A,\]
then $\B$ must also have a Siggers polymorphism. Hence, minor conditions can be used as witnesses for proving that $\A$ does not pp-construct $\B$. 

There is a third way to describe $\PPPoset$: a Birkhoff-style approach, extending the concept of a variety in universal algebra to the concept of a variety that is not only closed under homomorphic images, subalgebras, and products, but also closed under taking so-called \emph{reflections} of algebras.
We do not need this perspective in this thesis and refer the reader to~\cite{wonderland}.

Exploring $\PPPoset$ is a massive task. Hence, as often in mathematics we try to solve a simpler problem first. There are some natural ways to restrict the problem:
\begin{enumerate}
    \item restrict the size of the structures;
    \item restrict the signature of the structures.
\end{enumerate}
The first restriction was investigated by Vucaj in his PhD Thesis. He classified all structures with two elements~\cite{PPPoset}. The three-element case is already proving terribly difficult; at the moment Vucaj, Barto, Zhuk, Bodirsky, myself, and others are working on it. Some initial results have been published in~\cite{VucajZhukBodirsky21ThreeElements,vucajzhuk2023submaximal}. 

The second restriction is the one we focus on in this thesis. We will discuss structures with unary signatures in Section~\ref{sec:unarySignature}; they turn out to be quite easy in our setting. Hence, the natural next step is to investigate structures with one binary relation, i.e., \emph{directed graphs} or \emph{digraphs} for short. The poset arising when ordering all finite digraphs by pp-constructability is called $\DGPoset$.
It is known that for every finite relational structure with finite signature one can construct a finite digraph that  satisfies ``almost" the same minor conditions~\cite{FederVardi,BulinDelicJacksonNiven}. For details see Theorem~\ref{thm:StructurToDigraph}. In particular, the CSP over the structure and the CSP over the digraph are log-space equivalent. Using this construction we can show that many questions related to the complexity of  CSPs over finite templates can be answered by studying digraphs. 
For example, the  Larose-Tesson Conjecture for finite structures  is equivalent to the one for finite digraphs. 
In Section~\ref{sec:chaTreesOpenTreesNLConjecture} we present a potential counterexample to this conjecture: a tree whose CSP is neither known to be P-hard nor $\ModL p$-hard, nor is it known to be in NL. 


From the previous paragraph it is clear that studying the poset of all finite digraphs ordered by pp-constructability is almost as difficult as studying the poset of all finite structures ordered by pp-constructability. Hence, in this thesis we focus on specific types of digraphs: smooth digraphs, i.e., digraphs without sources or sinks, tournaments and semicomplete digraphs, orientations of paths and cycles, digraphs with at most four vertices, and orientations of trees. 
%

\section*{Contributions and Structure of the Thesis}

In Chapter~\ref{cha:prelims} we introduce basic concepts in detail; in particular we discuss pp-constructions, minor conditions, and CSPs.

\subsection*{Chapter~\ref{cha:cycles}: Smooth Digraphs}
In this chapter we focus on smooth digraphs. 
Barto, Kozik, and Niven showed that a finite smooth digraph can either pp-construct $\drawKThree$ or is homomorphically equivalent to a disjoint union of directed cycles~\cite{BartoKozikNiven}. Hence, we mainly study disjoint unions of directed cycles in this chapter. We give a description of when one disjoint union of directed cycles pp-constructs another in terms of the prime decomposition of the occurring cycle lengths. In particular, we show the following.

\begin{corollary*}[\ref{cor:PPvsLoopCaBlockerForSa}]
Let $a$ be in $\N^+$ and let $\B$ be a finite structure with finite relational signature. Then 
\begin{align*}
\B\leq \C_a && \text{iff} && \Pol(\B)\not\models\Sigma_{p}\text{ for all $p$ that are prime divisors of $a$},
\end{align*}
where $\C_a$ is the directed cycle of length $a$ and $\Sigma_p$ is the $p$-cyclic loop condition.
\end{corollary*}
Furthermore, we show a version of the theorem of Barto, Kozik, and Niven for pp-constructions.

\begin{corollary*}[\ref{cor:cyclesSquarefree}]
Let $\mathbb{G}$ be a finite smooth digraph. Then either $\mathbb G$ can pp-construct $\drawKThree$ or there is a finite disjoint union of cycles $\C$ whose cycle lengths are square-free such that $\mathbb G$ has the same pp-constructability type as $\C$.
\end{corollary*}
We also describe all digraphs with a quasi Maltsev polymorphism. Any digraph that can be pp-constructed from a digraph with a quasi Maltsev polymorphism also has a quasi Maltsev polymorphism. Hence, this description sheds light on a large part of the upper region of $\DGPoset$.

Most of the content of this chapter has been published in~\cite{StarkeVucajBodirskySmoothDigraphs}.

\subsection*{Chapter~\ref{cha:submax}: Submaximal Digraphs}

A digraph $\G$ lies on the $n$-th level of $\DGPoset$ if a longest chain starting in $\G$ consists of $n$ digraphs.
From the previous chapter we already know all digraphs on the $n$-th level of $\DGPoset$ with a quasi Maltsev polymorphism. In this chapter we show that on the first three levels there is only one digraph that does not have a Maltsev polymorphism, namely $\T_3$, i.e., the transitive tournament on three vertices. We prove that the three topmost levels of $\DGPoset$ are as follows:

\begin{center}
        \begin{tikzpicture}[scale=1.3]
    \node (0) at (2,2)  {$\P_{0} \equiv \C_{1}$};
    \node (1) at (2,1)  {$\P_{1} \equiv \P_2 \equiv \P_3 \equiv \cdots$};
    \node (20) at (0,0)  {$\T_{3}$};
    \node (21) at (1,0)  {$\C_{2}$};
    \node (22) at (2,0)  {$\C_{3}$};
    \node (23) at (3,0)  {$\C_{5}$};
    \node (24) at (4,0)  {$\dots$};
    \path
        (0) edge (1)
        (1) edge (20)
        (1) edge (21)
        (1) edge (22)
        (1) edge (23)
        (1) edge (24)
        ;
    \end{tikzpicture}
\end{center}
where $\P_n$ denotes the directed path with $n$ edges.

The content of this chapter has been published in~\cite{BodirskyStarke2022}.

\subsection*{Chapter~\ref{cha:DatalogIntro}: A Short Overview of Datalog}
In this chapter we recall Datalog, linear Datalog, and symmetric linear Datalog. We give a proof of the well-known fact that containment in one of these classes is closed under primitive positive constructions. Furthermore, we give a quick survey of known sufficient and necessary conditions for containment in one of these classes. 

\subsection*{Chapter~\ref{cha:semicompleteTn}:  Semicomplete Digraphs}
This chapter started out as a continuation of Chapter~\ref{cha:submax}. We wanted to find all digraphs on the fourth level of $\DGPoset$. It is easy to see that the transitive tournaments form an infinite descending chain. We had suspected  that $\T_n$, i.e., the transitive tournament with $n$ vertices, lies on the $n$-th level of $\DGPoset$. This turned out to be wrong. In this chapter we present a digraph that lies strictly between $\T_3$ and $\T_4$ (Example~\ref{exa:T4isNotACoverOfT3}). We show that the infinite descending chain $\T_3\ppstrictlyLess\T_4\ppstrictlyLess\cdots$ has a greatest lower bound which can be represented by the digraph depicted below.
\begin{center}
    \begin{tikzpicture}[scale=0.7,rotate=-90]
    \node[var-b] (3) at (0,3) {};
    \node[var-b] (2) at (0,2) {};
    \node[var-b] (1) at (0,1) {};
    \node[var-b] (0) at (0,0) {};
    \path[->,>=stealth']
        (0) edge (1)
        (1) edge (2)
        (2) edge (3)
        (0) edge[bend left=40] (2)
        (1) edge[bend right=40] (3)
        ;
\end{tikzpicture}
\end{center}
Furthermore, we generalize transitive tournaments to semicomplete digraphs and show that every semicomplete digraph, i.e., a digraph where between any two vertices there is either a directed or an undirected edge, can either pp-construct $\drawKThree$, is a transitive tournament, contains a single two-cycle, or  contains a single three-cycle. These four cases were already studied by Jackson,  Kowalski, and Niven in~\cite{NivenDigraphCSP}. We were not able to fully classify the pp-constructability order for the latter two cases. 

We combine our knowledge about transitive tournaments and about directed cycles to describe the digraphs on the fourth level of $\DGPoset$ that are lower covers of $\C_p$ for some prime $p$.
\begin{theorem*}[\ref{thm:lowerCoversOfSubmaximalCycles}]
Let $p$ be a prime. Then the distinct lower covers of $\C_p$ are  
$\T_3\mathbin{\times}\C_p$, $\C_{p\cdot q_1}$, $\C_{p\cdot q_2}$, $\dots$, where $q_1,q_2,\dots$ are the primes not equal to $p$.
\end{theorem*}

\subsection*{Chapter~\ref{cha:paths}:  Orientations of Paths and Cycles}
In Chapter~\ref{cha:cycles} we studied directed paths and directed cycles. In this chapter we consider orientations of paths and orientations of cycles.  We show that the digraph depicted below
\begin{center}
    \begin{tikzpicture}[scale=0.3]
    \node[var-b] (0) at (0,0) {};
    \node[var-b] (1) at (1,1) {};
    \node[var-b] (2) at (2,2) {};
    \node[var-b] (3) at (3,3) {};
    \node[var-b] (4) at (4,4) {};
    \node[var-b] (5) at (5,3) {};
    \node[var-b] (6) at (6,2) {};
    \node[var-b] (7) at (7,3) {};
    \node[var-b] (8) at (8,2) {};
    \node[var-b] (9) at (9,1) {};
    \node[var-b] (10) at (10,2) {};
    \node[var-b] (11) at (11,3) {};
    \node[var-b] (12) at (12,4) {};
    \node[var-b] (13) at (13,5) {};
    \path[->,>=stealth']
    (0) edge (1)
    (1) edge (2)
    (2) edge (3)
    (3) edge (4)
    (5) edge (4)
    (6) edge (5)
    (6) edge (7)
    (8) edge (7)
    (9) edge (8)
    (9) edge (10)
    (10) edge (11)
    (11) edge (12)
    (12) edge (13)
    ;
    \end{tikzpicture}
\end{center}
is the smallest orientation of a path that does not have a Hagemann Mitschke chain~\cite{HagemannMitschke}.
It is also the smallest orientation of a path such that the complement of its CSP is not in symmetric linear Datalog.  
For orientations of paths with up to 21 vertices that do have a Hagemann Mitschke chain we determine the length of the shortest such chain. Conjecture~\ref{conj:LaroseTesson} states in particular that every CSP over a finite structure in NL is NL-hard or in L.
Combining results of Chapter~\ref{cha:semicompleteTn}, of Kazda~\cite{Kazda-n-permute}, of  Dalmau~\cite{Dalmau_2005LinearDatalog}, Barto, Kozik, and of Willard~\cite{BartoKozikWillard}, we show that this L-NL dichotomy is true for orientations of paths.
\begin{corollary*}[\ref{cor:NLDichotomyHoldForPaths}]
Let $\A$ be a finite relational structure such that the complement of  $\csp(\A)$ is in linear Datalog. This is in particular true if $\A$ is an orientation of a path. Then exactly one of the following is true
\begin{enumerate}
    \item $\A\models\HM n$ for some $n\geq1$ and the complement of $\csp(\A)$ is in symmetric linear Datalog; in this case $\csp(\A)$ is in $\Lclass$, or
    \item $\A\not\models\HM n$ for any $n\geq1$ and $\A\ppleq\Ord$; in this case $\csp(\A)$ is $\NL$-complete. 
\end{enumerate}
\end{corollary*}

\subsection*{Chapter~\ref{cha:4elements}: Digraphs of Order Four}
In previous chapters we have considered digraphs of specific shapes, e.g., smooth digraphs, paths, or cycles. In this chapter we restrict the number of vertices in the graph, similarly to Vucaj and Bodirsky, who described the subposet of $\PPPoset$ consisting of all structures with two elements in the domain~\cite{PPPoset}. We consider  the subposet of $\DGPoset$ which consists of all digraphs with at most four vertices. With the help of a computer program we show that this subposet has 15 different pp-constructability classes. We  give an almost complete description of when a digraph in this subposet can pp-construct another digraph in this subposet. The only remaining question is whether
\begin{align*}
&\begin{tikzpicture}[scale=0.5,baseline=2mm]
\node[bullet] (0) at (0,0){};
\node[bullet] (1) at (0,1){};
\node[bullet] (2) at (1,1){};
\node[bullet] (3) at (1,0){};
\path[->,>=stealth']
(0) edge (1)
(0) edge (2)
(0) edge (3)
(1) edge (0)
(1) edge (2)
(1) edge (3)
(2) edge (3)
;
\end{tikzpicture}
\hspace{10mm}\text{can pp-construct}\hspace{10mm}
\begin{tikzpicture}[scale=0.5,baseline=2mm]
\node[bullet] (0) at (0,0){};
\node[bullet] (1) at (0,1){};
\node[bullet] (2) at (1,1){};
\node[bullet] (3) at (1,0){};
\path[->,>=stealth']
(0) edge (1)
(0) edge (2)
(1) edge (0)
(1) edge (2)
(2) edge (3)
;
\end{tikzpicture}~.
\end{align*}
Whenever a digraph can pp-construct another one we present a pp-construction, and whenever a digraph cannot pp-construct another one we give a separating minor condition.

\subsection*{Chapter~\ref{cha:hardTrees}: The Smallest Hard Trees}
In this last chapter we study orientations of trees. Gutjahr, Welzl, and Woeginger showed in 1992 that there are orientation of trees whose CSP is NP-hard~\cite{GutjahrWW92}. Their example of a tree with an NP-hard CSP has 287 vertices. Since then, smaller and smaller NP-hard trees have been found. We wrote a program to determine the smallest NP-hard trees. It determined that the smallest NP-hard trees have 20 vertices (assuming $\operatorname{P}\neq\operatorname{NP}$). In this chapter we present the key ideas of this program and we present all smallest NP-hard trees. 

The content of this chapter has been published in~\cite{BodirskyBulinStarkeWernthalerSmallestHardTrees}.

\section*{Some Open Questions}
Finally, I want to mention two open problems.
\begin{enumerate}
    \item Are there infinite ascending chains in $\DGPoset$?
\end{enumerate}
This simple sounding question has been haunting us from the very beginning. We show in Chapter~\ref{cha:cycles} that there is no infinite ascending chain of smooth digraphs, i.e., digraphs without sources or sinks. In Chapter~\ref{cha:semicompleteTn} we show that also semicomplete digraphs do not admit such a chain. Whether $\DGPoset$ has such a chain is still open.

\begin{enumerate}[resume]
    \item Is $\DGPoset$ a lattice?
\end{enumerate}
It is well known that $\PPPoset$ is a meet-semilattice. Whether it is also a lattice is still an open question. For $\DGPoset$ it is not even known whether it is a meet semilattice, as it is not clear whether the meet of two digraphs can always be represented by a digraph. In Chapter~\ref{cha:cycles} we show that the subposet of $\DGPoset$ consisting of smooth digraphs is a distributive lattice.  
\chapter{Preliminaries}\label{cha:prelims}

In this chapter we present formal definitions of notions used in this thesis.

\section*{Notation}
\begin{itemize}
    \item For $n\in\N$, we define $[n]\coloneqq\{1,\dots,n\}$. Note that $[0]=\emptyset$.
    \item For $a\geq1$ we denote $\{0,\dots,a-1\}$ by $\Z_a$.
    \item By $\operatorname{Im}(f)$ we denote the image of the function $f$.
    \item By $\dom(f)$ we denote the domain of the function $f$.
    \item By $\lcm$ and $\gcd$ we denote the least common multiple and the greatest common divisor, respectively.
    \item For tuples $\boldsymbol t=(t_1,\dots,t_k)\in A^k$ we usually write the tuple $\boldsymbol t$ in bold font and the entries $t_i$ in normal font.
    \item For a $n$-tuple $\boldsymbol{a} = (a_1,\dots,a_n)$ and a map $\sigma \colon [m]\to [n]$, we
denote the $m$-tuple $(a_{\sigma(1)},\dots,a_{\sigma(m)})$ by $\shiftTuple{\boldsymbol{a}}{\sigma}$.
    \item By $k\equiv_a \ell$ we denote that $a$ divides $k-\ell$. 
\end{itemize}

In the next section we will introduce the poset that this thesis is about. In Sections \ref{sec:h1} and \ref{sec:freestructures} we will present equivalent characterizations.




\section{The pp-constructability Poset}
A \emph{relational signature} is a function $\tau\colon\dom(\tau)\to\N$ that assigns to each \emph{relation symbol} $R$ in the (not necessarily finite) set $\dom(\tau)$ its \emph{arity}. Instead of $R\in\dom(\tau)$ we will usually write $R\in\tau$ or sometimes $R^{(k)}\in\tau$ to indicate that $R\in\dom(\tau)$ and $\tau(R)=k$. A \emph{relational structure $\A$ with signature $\tau$} (or \emph{$\tau$-structure} for short) is a non-empty set $A$ together with a $k$-ary relation $R^{\A}\subseteq A^k$ for each $k$-ary relation symbol in $\tau$. A \emph{digraph} is a structure with one binary relation.

Let $\A=(A;(R^{\A})_{R\in\tau})$ and $\B=(B;(R^{\B})_{R\in\tau})$ be $\tau$-structures with the same relational signature $\tau$. A map $h\colon A \to B$ is a \emph{homomorphism} from $\A$ to $\B$ if it preserves all relations, i.e., for all $R\in \tau$:
\[\text{if } (a_1,\dots,a_n)\in R^{\A} \text{ then } (h(a_1),\dots,h(a_n))\in R^{\B}.\]
We write $\A\to\B$ if there exists a homomorphism from $\A$ to $\B$.  If $\A\to\B$ and $\B\to\A$, then we say that $\A$ and $\B$ are \emph{homomorphically equivalent}.
A \emph{primitive positive} formula (pp-formula) is a first order formula using only existential quantification and conjunctions of \emph{atomic formulas}, i.e., formulas of the form $R(x_1,\dots,x_n)$ or $x\approx y$, and the special formula $\bot$. The formula $\bot(x_1,\dots,x_n)$ always evaluates to the empty set. 
Let $\mathbb A$ be a relational structure and $\phi(x_1,\dots,x_n)$ be a pp-formula. 
Then the relation 
\begin{equation*}
\phi^{\A}\coloneqq\{(a_1,\dots,a_n)\mid \mathbb A \vDash \phi(a_1,\dots,a_n)\}
\end{equation*}
is said to be \emph{pp-definable} in $\mathbb A$.
We say that $\mathbb B$ is a \emph{pp-power} of $\A$ if it is isomorphic to a structure with domain $A^n$, for some $n\in \N^+$\!, whose relations are pp-definable in $\A$ (a $k$-ary relation on $A^n$ is regarded as a $kn$-ary relation on $A$).

Combining the notions of homomorphic equivalence and pp-power we obtain the following definition from~\cite{wonderland}.
\begin{definition}
We say that $\A$ \emph{pp-constructs} $\B$, denoted by $\A\leq\B$, if $\B$ is homomorphically equivalent to a pp-power of $\A$.
\end{definition}

Since pp-constructability is a reflexive and transitive relation on the class of relational structures~\cite{wonderland}, the use of the symbol $\ppleq$ is justified. 
Note that the quasi-order $\ppleq$ naturally induces the equivalence relation 
\[\A \equiv \B \text{ if and only if }\B \ppleq \A \ppleq \B.\]
The equivalence classes of $\equiv$ are called \emph{pp-constructability types}. We denote by $[\A]$ the pp-constructability type of a structure $\A$.
\begin{definition}
We name the poset
\[\PPPoset\coloneqq(\{[\A]\mid \A\ \text{a finite relational structure}\},\leq)\] the \emph{pp-constructability poset}.
\end{definition}

\begin{example}\label{exa:C6LeqC3}
The directed cycle with 6 vertices, $\C_6$, pp-constructs the directed cycle with 3 vertices, $\C_3$. To see this, consider the first pp-power of $\C_6$ given by the pp-formula $\Phi_E(x,y)\coloneqq \exists z.\ E(x,z)\AND E(z,y)$. 
We obtain a structure that consists of two disjoint copies of $\C_3$; which is homomorphically equivalent to $\C_3$; see Figure~\ref{fig:6pp3}. 
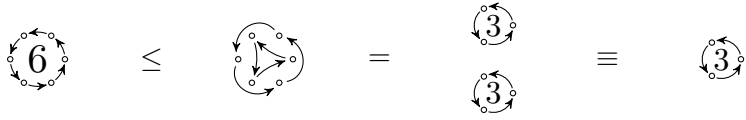
\begin{figure}
    \centering
    \begin{tikzpicture}
    \cycle{6}{(-3,0)}
    
    \def \n {3}
\def \radius {{2*3*0.035+0.15}}

\foreach \s in {1,...,\n}
{
  \draw ({360/\n * (\s - 1)}:\radius) circle (1pt);
  \draw ({360/\n * (\s - 1)+180/\n}:\radius) circle (1pt);
  \draw ({360/\n * (\s - 1)}:\radius) edge[>=stealth',arrows=-{>[bend]},bend left=20,shorten >=1mm,shorten <=1mm] ({360/\n * (\s)}:\radius);
  \draw ({360/\n * (\s - 1)+180/\n}:\radius) edge[>=stealth',arrows=-{>[bend]},bend right = 80,looseness=1.7,shorten >=1mm,shorten <=1mm] ({360/\n * (\s)+180/\n}:\radius);
}

    \cycle{3}{(3,0.45)}
    \cycle{3}{(3,-0.45)}
    
    \cycle{3}{(6,0)}
    
    \node at (-1.5,0) {$\leq$};
    \node at (1.5,0) {$=$};
    \node at (4.5,0) {$\equiv$};
    
    \end{tikzpicture}
    \caption{The structure $\C_6$ pp-constructs $\C_3$.}
    \label{fig:6pp3}
\end{figure}
We will often represent a pp-formula by a graph, very similar to its \emph{canonical database} (see Definition~\ref{def:canonicalDatabase}). For example the graph representing $\Phi_E(x,y)$ is 
\begin{center}
    \begin{tikzpicture}[scale=0.5]
    \node[var-f,label=below:$x$] (x) at (0,0) {};
    \node[var-b] (z) at (1,0) {};
    \node[var-f,label=below:$y$] (y) at (2,0) {};
    \path[->,>=stealth']
    (x) edge (z)
    (z) edge (y)
    ;
    \end{tikzpicture}
\end{center}
where filled vertices represent existentially quantified variables. 
\end{example}

Let us make our first observation about $\PPPoset$ by showing that it has a largest element, i.e., a finite structure that can be pp-constructed by every other finite structure. 
\begin{lemma}\label{lem:loopIsTopElement}
The pp-constructability type  of the loop graph, $[\drawLoop]$, is the top element of $\PPPoset$.
\end{lemma}
\begin{proof}
Let $\A$ be a finite structure. 
Let $\A'$ be the first-pp-power of $\A$ given by the pp-formula $\Phi_E(x,y)\coloneqq (x\approx y)$. Then $\A'$ is the digraph \[(A;\{(a,a)\mid a \in A\}).\] Note that $\A'$ is homomorphically equivalent to a single loop:
\begin{center}
    \begin{tikzpicture}
    \node at (-2.2,0) {$\A$};
       \node at (-1.2,0) {$\ppleq$};
       
\node[circle, fill, scale=0.3] at (0,0) (0) {};
\path (0) edge[out=120, in=60, looseness=20,->,>=stealth'] (0);
\node[circle, fill, scale=0.3] at (1,0.5) (0) {};
\path (0) edge[out=120, in=60, looseness=20,->,>=stealth'] (0);
\node[circle, fill, scale=0.3] at (0.8,-0.3) (0) {};
\path (0) edge[out=120, in=60, looseness=20,->,>=stealth'] (0);
\node[circle, fill, scale=0.3] at (0.4,0.3) (0) {};
\path (0) edge[out=120, in=60, looseness=20,->,>=stealth'] (0);
\node[circle, fill, scale=0.3] at (1.2,0.1) (0) {};
\path (0) edge[out=120, in=60, looseness=20,->,>=stealth'] (0);

    \node at (2.5,0) {$\equiv$};
\node[circle, fill, scale=0.3] at (3.5,0) (0) {};
\path (0) edge[out=120, in=60, looseness=20,->,>=stealth'] (0);
    \end{tikzpicture}
\end{center}
Hence,  $\A\ppleq \drawLoop$ as desired. 
\end{proof}

In graph notation we represent the pp-formula $\Phi_E$ from the previous proof by
\begin{center}
    \begin{tikzpicture}[scale=0.5]
    \node[var-f,label=below:$x$] (x) at (0,0) {};
    \node[var-f,label=below:$y$] (y) at (1,0) {};
    \path[dashed]
    (x) edge (y)
    ;
    \end{tikzpicture}
\end{center}
where the dashed line between $x$ and $y$ indicates that $\Phi_E$ contains the conjunct $x\approx y$.
Note that $\drawLoop\ppleq\drawDot$ using the formula $\bot$. Hence, $[\drawLoop]=[\drawDot]$. The poset $\PPPoset$ also has a bottom element, which is $[\drawKThree]$, i.e., the pp-constructability type of the complete graph on three vertices; see Theorem~\ref{cor:K3IsBottom}. 

\section{Minor Conditions}\label{sec:h1}
As mentioned before, pp-constructability can be characterized algebraically; this characterization will provide the main tool to prove that a structure cannot pp-construct another structure. First, we define the basic notion of a \emph{polymorphism}.
For structures $\A$ and $\B$ define 
\[\Hom(\A,\B)\coloneqq\{f\mid f \text{ a homomorphism from } \A\text{ to } \B\}.\]
For $n\geq 1$, we denote by $\A^n$ the structure
with same signature $\tau$ as $\A$ whose domain is $A^n$ such that for any $k$-ary $R\in \tau$, a tuple $(\boldsymbol{a}_1,\dots,\boldsymbol{a}_k)$ of $n$-tuples is contained in $R^{\A^n}$ if and only if it is contained in $R^{\A}$
componentwise, i.e., $(a_{1j},\dots,a_{kj})\in R^{\A}$ for all $1\leq j\leq n$.

\begin{definition}
For a relational structure $\A$, a \emph{polymorphism of $\A$} is an element of $\bigcup_{n\in\N^+}\Hom(\A^n,\A)$. 
Moreover, we denote by $\Pol(\A)$ the \emph{polymorphism clone} of $\A$, i.e., the set of all polymorphisms of $\A$.  
\end{definition}
A set of operations $\mathcal C$ on a set $A$ is called a \emph{clone} if it contains the projections and is closed under composition, i.e., for every $n$-ary $f\in\mathcal C$ and all $m$-ary $g_1,\dots,g_n\in\mathcal C$ we have that the operation 
\begin{align*}
    f(g_1,\dots,g_n)\colon A^m&\to A\\
    \boldsymbol a&\mapsto f(g_1(\boldsymbol a),\dots,g_n(\boldsymbol a))
\end{align*}
is also contained in $\mathcal C$. 
Note that $\Pol(\A)$ contains all projections $\pi^n_i\colon A^n\to A$ and is closed under composition. Hence, the polymorphism clone of $\A$ is indeed a clone.
We introduce some notation for $f(g_1,\dots,g_n)$ in the case that $g_1,\dots,g_n$ are projections.
\begin{definition}
Let $\sigma \colon [m]\to [n]$ and $f\colon A^m \to A$ be functions. We define the function $f_\sigma\colon A^n\to A$ by the rule
\[f_\sigma(\boldsymbol{a}) \coloneqq f(\shiftTuple{\boldsymbol{a}} {\sigma}),\]
where $\shiftTuple{\boldsymbol{a}} {\sigma}=(a_{\sigma(1)},\dots,a_{\sigma(m)})$.
Any function of the form $f_\sigma$, for some map $\sigma \colon [m]\to [n]$, is called a \emph{minor} of $f$.
\end{definition}
Note that $f_\sigma=f(\pi^n_{\sigma(1)},\dots,\pi^n_{\sigma(m)})$. 
Let $I$ and $J$ are arbitrary sets and $\sigma\colon I\to J$, $f\colon A^I\to A$, and $\boldsymbol{a}\in A^J$. We extend the definitions of $\boldsymbol{a}_\sigma$ and $f_\sigma$ in the obvious way.
Let  $I$ be a finite set and $f\colon \A^I\to\A$ a homomorphism. We say that $f$ is $I$-ary. Then for every bijection $\sigma\colon [|I|] \to I$ the function $f_\sigma$ is a polymorphism of $\A$.
Throughout this thesis we treat all homomorphisms $f\colon \A^I\to\A$, where $I$ is a finite set, as elements of $\Pol(\A)$ by identifying $f$ with $f_\sigma$ for a bijection $\sigma$ (in all such situations the particular choice of $\sigma$ does not matter). 

\begin{definition}
Let $\A$ and $\B$ be structures. An arity-preserving map $\lambda\colon \Pol(\B)\to\Pol(\A)$ is \emph{minor-preserving} if for all $f\colon B^m\to B$ in $\Pol(\B)$ and $\sigma\colon [m]\to [n]$ 
we have
\begin{equation*}
    \lambda(f_\sigma)=(\lambda f)_\sigma.
\end{equation*}
We write $\Pol(\B)\tominor\Pol(\A)$ to denote that there is a minor-preserving map from $\Pol(\B)$ to $\Pol(\A)$.
\end{definition}
Note that every clone homomorphism is a minor-preserving map. Observe the the composition of minor preserving maps is again a minor preserving map. Hence, the existence of minor preserving maps defines a quasi order on the polymorphism clones of finite structures.
The next theorem, restated from~\cite{wonderland}, shows that the concepts presented so far, i.e., pp-constructions and minor-preserving maps, give rise to the same poset.
\begin{theorem}[Theorem 1.3 in~\cite{wonderland}]\label{thm:ppvsminor}
Let $\A$ and $\B$ be finite structures. Then
\begin{align*}
    \B\leq\A&&\text{if and only if}&&\Pol(\B)\tominor\Pol(\A).
\end{align*}
\end{theorem}

We will present the ideas of the proof, in the case that both structures are digraphs, after Theorem~\ref{thm:freestructure}.

\begin{definition}
Let $\sigma\colon [n]\to [r]$ and $\tau\colon [m]\to [r]$ be functions. A \emph{minor identity} is an expression of the form:
\begin{equation*}
    \forall x_1,\dots,x_r( f(x_{\sigma(1)},\dots,x_{\sigma(n)}) \approx g(x_{\tau(1)},\dots,x_{\tau(m)})).
\end{equation*}
\end{definition}
Usually, we write $f(x_{\sigma(1)},\dots,x_{\sigma(n)}) \approx g(x_{\tau(1)},\dots,x_{\tau(m)})$ omitting the universal quantification, or even $f_\sigma\approx g_\tau$ for brevity.
A finite set of minor identities is called \emph{minor condition}.
An important example of such a condition is given in the following definition.
\begin{definition}\label{def:siggers}
The \emph{Siggers condition}, denoted $\Siggers$, consists of the identity
\[f(a,r,e,a)\approx f(r,a,r,e).\]
\end{definition}
A set of functions $F$ \emph{satisfies} a minor condition $\Sigma$, denoted $F\models\Sigma$, if there is a map $\tilde\cdot$ assigning to each function symbol occurring in $\Sigma$ a function in $F$ such that for all $f_\sigma\approx g_\tau\in\Sigma$ we have $\tilde f_\sigma= \tilde g_\tau$. 
Note that we make a distinction between the symbol $\approx$ and $=$ to emphasize the difference between an identity in a formula and an equality of two specific objects.

\begin{definition}\label{def:orderOnConditions}
Let $\Sigma$ and $\Gamma$ be sets of minor identities. We say that $\Sigma$ \emph{implies} $\Gamma$ (and that $\Sigma$ \emph{is stronger then} $\Gamma$), denoted $\Sigma\Rightarrow\Gamma$, if 
\[\mathcal C \models \Sigma \text{ implies } \mathcal C \models \Gamma \text{ for all clones $\mathcal C$.} \] 
If $\Sigma\Rightarrow\Gamma$ and $\Gamma\Rightarrow\Sigma$, then we say that $\Sigma$ and $\Gamma$ are \emph{equivalent}, denoted $\Sigma\Leftrightarrow\Gamma$. We define 
\[[\Sigma]\coloneqq\{\Sigma'\mid \Sigma'\text{ a minor condition, }\Sigma\Leftrightarrow\Sigma'\}.\]
\end{definition}  
We say that a minor condition is \emph{trivial} if it is satisfied by projections on a set $A$ that
contains at least two elements, or, alternatively, if it is implied by any minor condition.
We extend the definition of $\models$ and $\Rightarrow$ to single functions and single minor identities in the obvious way. Hence, if $f$ is a function, $\Sigma$ is a minor identity, and $\Gamma$ is a set of minor identities, then we can write $f\models\Sigma$ instead of $\{f\}\models\{\Sigma\}$ and $\Gamma\Rightarrow\Sigma$ instead of $\Gamma\Rightarrow\{\Sigma\}$.

Observe that if $\lambda\colon\Pol(\B)\to\Pol(\A)$ is minor-preserving and $f_\sigma= g_\tau$, then $\lambda(f)_\sigma=\lambda(g)_\tau$. It follows that minor-preserving maps preserve all minor conditions.
A simple compactness argument shows the following corollary of Theorem~\ref{thm:ppvsminor}.
\begin{corollary}\label{cor:wond}
Let $\A$ and $\B$ be finite structures. Then
\begin{align*}
    \B\ppleq\A&&\text{iff}&&\text{$\Pol(\B)\models\Sigma$ implies $\Pol(\A)\models\Sigma$ for all minor conditions $\Sigma$}.
\end{align*}
\end{corollary}
After Theorem~\ref{thm:freestructure} we give a proof of this corollary for the case that $\A$ and $\B$ are digraphs. The proof does not use compactness but free structures instead.
See Lemma~\ref{lem:pcSatisfyPclc} for an example of how to use minor conditions.

\section{The Polymorphism Preservation Theorem}
The connection of pp-construction and minor-preserving maps presented in the previous section in Theorem~\ref{thm:ppvsminor} might come as a surprise. In this section we look at a simpler version of this connection: that of pp-definable relations and polymorphisms.
First we introduce the following definition to translate between structures and pp-formulas.
\begin{definition}\label{def:canonicalDatabase}
The \emph{canonical database} of a pp-$\tau$-formula $\Phi$ is the $\tau$-structure $\A$ that can be constructed as follows: Let $\Phi'$ be obtained from $\Phi$ by removing all conjuncts of the form $x\approx y$ in $\Phi$ and by identifying variables $x$ and $y$ if there is a conjunct $x\approx y$ in $\Phi$.
Let $V$ be the variables of $\Phi'$. Then $\A$ is the $\tau$-structure with $A=V$ and \[R^{\A}=\{(v_1,\dots,v_k)\mid R(v_1,\dots,v_k)\text{ is a conjunct of }\Phi'\}.\]
The \emph{canonical conjunctive query} of a structure $\A$ with signature $\tau$ is the formula
\[\bigwedge_{R^{(k)}\in\tau}\bigwedge_{\boldsymbol r\in R^{\A\hspace{-4mm}\phantom{(k)}}}R(r_1,\dots,r_k).\]
\end{definition}
Observe that the canonical database of the canonical conjunctive query of a structure $\A$ is again $\A$. 
Note that the graph we write to represent a pp-formula $\Phi$ is similar to the canonical database of $\Phi$ but in general not the same. Formally, instead of identifying variables and ignoring existential quantifiers we introduce a new binary relation for $\approx$ (visualized by dashed edges) and a new unary relation for existentially quantified variables (visualized by filled vertices).
The following preservation theorem is well known. 
\begin{theorem}[\cite{Bodnarchuk1969GaloisTF,Geiger}]\label{thm:ppPolyPreservation}
Let $\A$ be a finite structure and $R$ be some relation on $A$. Then $\Pol(\A)$ preserves $R$ if and only if $R$ is pp-definable in $\A$.
\end{theorem}
\begin{proof}
It is easy to check that pp-definable relation are preserved by polymorphisms. 
Let $R=\{\boldsymbol r_1,\dots,\boldsymbol r_m\}$ be a $k$-ary relation on $A$ that is preserved by all polymorphisms of $\A$. Define $\boldsymbol t_i\coloneqq (r_{1i},\dots,r_{mi})$ for $i\in[k]$ and $j_i\coloneqq \min\{i'\in[k]\mid  \boldsymbol t_i=\boldsymbol t_{i'}\}$.
Note that $(\boldsymbol t_1,\dots,\boldsymbol t_k)\in R^{m}\subseteq A^m$. 
Let $\Phi$ be the canonical conjunctive query of $\A^m$ with $(r_{1i},\dots,r_{mi})$ replaced by $x_{j_i}$ for $i\in[k]$. Let $v_1,\dots,v_n$ be the variables of $\Phi$ that are not equal to any $x_i$. Define 
\[\Phi_R\coloneqq \left(\bigwedge_{\substack{1\leq i<j\leq k\\\boldsymbol t_i=\boldsymbol t_j}}x_i\approx x_j\right) \AND\exists v_1,\dots,v_n\ \Phi(x_1,\dots,x_k,v_1,\dots,v_n).\]
Note that the satisfying assignments of $\Phi$ correspond 1-1 with the $m$-ary polymorphisms of $\A$. Hence, the set $\Phi_R^{\A}$ is equal to
\[\{(f(r_{11},\dots,r_{m1}),\dots,f(r_{1k},\dots,r_{mk}))\mid \text{$f$ is an $m$-ary polymorphism of $\A$}\}.\]
Since $R$ is preserved by all polymorphisms we have that $\Phi_R^{\A}\subseteq R$ and because all projections are polymorphisms we also have $R\subseteq \Phi_R^{\A}$. Therefore, $\Phi_R$ is a pp-definition of $R$.
\end{proof}

The idea of this proof to relate solutions of a formula with polymorphisms of a structure will be reused in Definition~\ref{def:indicatorStructure}.

\section{How I Learned to Love Free Structures }\label{sec:freestructures}
There are two other characterizations of pp-constructions that look utterly incomprehensible at first. However once I finally took the time to understand them I noticed that they are sometimes actually quite useful in practice. For example the pp-constructions in Lemmata~\ref{thm:PPvsLoop} and~\ref{lem:4ElementsT3u1C2EqualsT3u2C2} and the minor condition separating $\OriC{1,3}$ from $\T_3\mathbin{\times}\C_2$ in Section~\ref{sec:CyclesWithTwoChangesInDirection} were found with the help of free structures. 
Nevertheless this characterizations can be skipped without compromising the understanding of the bulk of the rest of this thesis.

\begin{definition}[Definition 4.1 in \cite{Barto_2021FreeStructures}]\label{def:freestructure}
Let $\A$ be a finite relational structure on the set $A=[n]$, and $\mathcal C$ a clone (not necessarily related to $\A$). The \emph{free structure of} $\mathcal C$ \emph{generated by} $\A$ is a relational structure $\F_{\mathcal C}(\A)$ with the same signature as $\A$. Its universe $F_{\mathcal C}(A)$ consists of all $n$-ary operations in $\mathcal C$. For any relation of $\A$, say $R^{\A}=\{\boldsymbol{r}_1,\dots,\boldsymbol{r}_m\}\subseteq A^k$, the relation $R^{\F_{\mathcal C}(\A)}$ is defined as the set of all $k$-tuples $(f_1,\dots,f_k)\in F_{\mathcal C}(A)$ such that there exists an $m$-ary operation $g\in\mathcal C$ that satisfies
\[f_j(x_1,\dots,x_n) = g(x_{{r}_{1j}},\dots,x_{{r}_{mj}}) \text{ for each } j\in[k].\]
\end{definition}
We can also write this identity as $f_j=g_{\pi_j}$, where $\pi_j\colon [m]\to [n],k\mapsto r_{kj}$, i.e., $\pi_j$ maps $k$ to the $j$-th entry of $\boldsymbol r_k$. See Examples~\ref{exa:freeStructureWorks} and~\ref{exa:freeStructureNotWorks} for examples of a free structure generated by a digraph.

\begin{definition}[Section 3.2 in \cite{Barto_2021FreeStructures}]\label{def:freeCondition}
Let $\A$ and $\A'$ be finite relational structures with the same signature $\tau$, such that $A=[n]$. 
For each node $v$ in $A'$ let $f_v$ be an $n$-ary function symbol and for each relation of $\A$, say $R^{\A}=\{\boldsymbol{r}_1,\dots,\boldsymbol{r}_m\}\subseteq A^k$ and each $\boldsymbol v\in R^{\A'}\subseteq (A')^{k}$ let $g^R_{\boldsymbol v}$ be an $m$-ary function symbol. Define the minor condition $\Sigma(\A,\A')$ to consist of all identities
\begin{align*}
f_{v_i}(x_1,\dots,x_n)\approx g^R_{\boldsymbol v}(x_{r_{1i}},\dots,x_{r_{mi}})
\end{align*}
for all $R\in\tau, \boldsymbol v=(v_1,\dots,v_k)\in R^{\A'}$, and $1\leq i\leq k$.
\end{definition}

Note that strictly speaking $\Sigma(\A,\A')$ depends on the order of the enumeration  $R^{\A}=\{\boldsymbol{r}_1,\dots,\boldsymbol{r}_m\}$. However any minor condition we obtain when using a different enumeration is equivalent to $\Sigma(\A,\A')$.
The following theorem links the notion of a free structure to the characterizations of pp-constructability presented in Theorem~\ref{thm:ppvsminor} and Corollary~\ref{cor:wond}.

\begin{theorem}[Theorem 4.12 in~\cite{Barto_2021FreeStructures}]\label{thm:freestructure}
Let $\A$ and $\B$ be finite relational structures. 
Then the following are equivalent
\begin{enumerate}
    \item $\B\ppleq\A$,
    \item $\Pol(\B)\tominor\Pol(\A)$,
    \item $\Pol(\B)\models\Sigma$ implies $\Pol(\A)\models\Sigma$ for all minor conditions $\Sigma$, 
    \item $\Pol(\A)\models\Sigma(\A,\F_{\Pol(\B)}(\A))$, and
    \item $\F_{\Pol(\B)}(\A)\to\A$.
\end{enumerate}
\end{theorem}




We will now present the proof of Theorem~\ref{thm:freestructure} in the case that $\A$ and $\B$ are digraphs and give some examples afterwards. 

\begin{proof}[Proof of Theorem~\ref{thm:freestructure} for digraphs]
$(1)\Rightarrow(2)$ Let $\B'$ be a $k$-th pp-power of $\B$ such that there are homomorphisms $h_a\colon \A\to\B'$ and $h_b\colon \B'\to\A$. Define the map $\lambda\colon \Pol(\B)\to\Pol(\A)$
\begin{align*}
f\mapsto \lambda(f)
\end{align*}
where $\lambda(f)$ is defined as follows: let $a_1,\dots,a_n\in A$ and $h_a(a_i)=(b_{1i},\dots,b_{ki})$ for every $i$, then $\lambda(f)(a_1,\dots,a_n)\coloneqq h_b(f(b_{11},\dots,b_{1n}),\dots,f(b_{k1},\dots,b_{kn}))$. 
With abuse of notation $\lambda(f)$ could be denoted by $h_b\circ f\circ h_a$. 
Using Figure~\ref{fig:lambdaf} it is easy to verify the following two observations: $\lambda$ is minor-preserving and the fact that $h_a$ and $h_b$ are homomorphisms and that $f$ is a polymorphism ensures that $\lambda(f)$ is a polymorphism.

\begin{figure}
    \centering
    \begin{tikzpicture}[scale=0.5]
    \node[var-b,label=above:$a_1$] (a1) at (0,0) {};
    \node[] (a2) at (1.5,0) {$\dots$};
    \node[var-b,label=above:$a_n$] (an) at (3,0) {};
    \node[label=above:{\small $\lambda(f)$}] at (5,0) {$\mapsto$};
    \node[var-b] (an) at (7,0) {};
    
    \node[rotate=-90,label=above:{\small $h_a$}] at (0,-1.5) {$\mapsto$};
    \node[rotate=-90,label=above:{\small $h_a$}] at (3,-1.5) {$\mapsto$};
    \node[rotate=90,label=below:{\small $h_b$}] at (7,-1.5) {$\mapsto$};

    \node[var-b] (a1) at (0,-3) {};
    \node[] (a2) at (1.5,-3) {$\dots$};
    \node[var-b] (an) at (3,-3) {};
    \node[label=above:{\small $f$}] at (5,-3) {$\mapsto$};
    \node[var-b] (an) at (7,-3) {};
    
    \node[rotate=90] at (0,-4) {$\dots$};
    \node[rotate=90] at (3,-4) {$\dots$};
    \node[rotate=90] at (7,-4) {$\dots$};

    \node[var-b] (a1) at (0,-5) {};
    \node[] (a2) at (1.5,-5) {$\dots$};
    \node[var-b] (an) at (3,-5) {};
    \node[label=above:{\small $f$}] at (5,-5) {$\mapsto$};
    \node[var-b] (an) at (7,-5) {};

    \end{tikzpicture}
    \caption{Diagram of $\lambda(f)=h_b\circ f\circ h_a$.}
    \label{fig:lambdaf}
\end{figure}
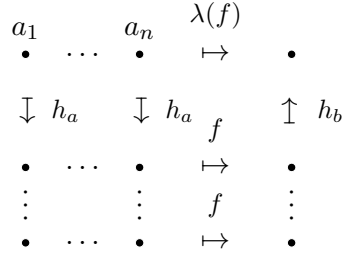

$(2)\Rightarrow(3)$ Let $\lambda\colon \Pol(\B)\to\Pol(\A)$ be a minor-preserving map and $\Sigma$ be some minor condition satisfied by $\Pol(\B)$ using the set of functions symbols $F$. Assume that $\Pol(\B)\models\Sigma$, witnessed be some assignment $\tilde\cdot\colon F\to\Pol(\B)$. Then the assignment $f\mapsto \lambda(\tilde f)$ witnesses that $\Pol(\A)\models\Sigma$.

Before we come to the next implications let us write out the definitions of $\F_{\Pol(\B)}(\A)$ and $\Sigma(\A,\A')$ for digraphs. Let $E(\A)$ denote the set of edges of $\A$.
The set of vertices of $\F_{\Pol(\B)}(\A)$ consists of all homomorphisms from  $\B^A$ to $\B$ and there is an edge from $f_1$ to $f_2$ if there exists a homomorphism $g\colon \B^{E(\A)}\to\B$, 
such that $f_1=g_{\operatorname{source}}$ and $f_2=g_{\operatorname{target}}$, where $\operatorname{source},\operatorname{target}\colon E(\A)\to A$ with  $\operatorname{source}(s,t)= s$ and $\operatorname{target}(s,t)= t$.
Let $\A'$ be a digraph, for each node $\tilde f\in A$ let $f$ be a function symbol,  
and for each edge $(\tilde f_1,\tilde f_2)$ in $\A'$ let $g_{(f_1,f_2)}$ be a function symbol. Define
\begin{align*}
\Sigma(\A,\A')\coloneqq{}&\{f_1\approx (g_{(f_1,f_2)})_{\operatorname{source}}\mid (\tilde f_1,\tilde f_2)\text{ is an edge in }\A'\}\cup{}\\
&\{f_2\approx (g_{(f_1,f_2)})_{\operatorname{target}}\mid (\tilde f_1,\tilde f_2)\text{ is an edge in }\A'\}
\end{align*}

$(3)\Rightarrow(4)$ 
Consider the assignment that maps $f$ to $\tilde f$ and $g_{(f_1,f_2)}$ to some function $g$ witnessing the edge $(\tilde f_1,\tilde f_2)$ in $\F_{\Pol(\B)}(\A)$. This assignment shows that $\Pol(\B)\models\Sigma(\A,\F_{\Pol(\B)}(\A))$. Hence, by (3), $\Pol(\A)\models\Sigma(\A,\F_{\Pol(\B)}(\A))$.

$(4)\Rightarrow(5)$ 
Let $\bar\dot$ be an assignment witnessing $\Pol(\A)\models\Sigma(\A,\F_{\Pol(\B)}(\A))$. Now the map $\tilde f\mapsto \bar f(a\mapsto a)$ is a homomorphism from $\F_{\Pol(\B)}(\A)$ to $\A$, as desired. This map is well defined, since $\tilde f\colon \B^{A}\to\B$ is an element of $\F_{\Pol(\B)}(\A)$, $\bar f\colon \A^{A}\to\A$ is an element of $\F_{\Pol(\A)}(\A)$, and $a\mapsto a\colon\A\to\A$ is an element of $\A^{A}$. 
Although not necessary for the construction it is interesting to note that $\tilde f\mapsto \bar f$ is a homomorphism from $\F_{\Pol(\B)}(\A)$ to $\F_{\Pol(\A)}(\A)$ and that $f\mapsto f(a\mapsto a)$ is a homomorphism from $\F_{\Pol(\A)}(\A)$ to $\A$.

$(5)\Rightarrow(1)$
First, observe that 
$\F_{\Pol(\B)}(\A)$ is up to isolated vertices a $B^A$-th pp-power of $\B$. The elements of this pp-power are from $B^{B^A}$, i.e., functions from $B^A$ to $B$. Hence, it makes sense to denote the variables of the pp-formulas as $f(\boldsymbol b)$ for $\boldsymbol b\in B^A$ (instead of for example $f_{\boldsymbol b}$). Define the pp-formula
\begin{align*}
\Phi_E(f_1,f_2)\coloneqq &\operatorname{isPol}(f_1)\AND \operatorname{isPol}(f_2)\AND \exists g\ \operatorname{isPol}(g)\\
&\AND\bigwedge_{\boldsymbol b\in B^A} f_1(\boldsymbol b)\approx g_{\operatorname{source}}(\boldsymbol b)\AND 
\bigwedge_{\boldsymbol b\in B^A} f_2(\boldsymbol b)\approx g_{\operatorname{target}}(\boldsymbol b)
\intertext{where} 
\operatorname{isPol}(f_i)\coloneqq &\bigwedge_{(\boldsymbol{b_1},\boldsymbol{b_2})\in E(\B^A)} E(f_i(\boldsymbol{b_1}),f_i(\boldsymbol{b_2}))
\intertext{and}
\operatorname{isPol}(g)\coloneqq &\bigwedge_{(\boldsymbol{b_1},\boldsymbol{b_2})\in E(\B^{E(\A)})} E(g(\boldsymbol{b_1}),g(\boldsymbol{b_2})).
\end{align*}
This formula is a lot to take in, especially with the condensed notation. 
But do not worry there is a detailed example below.
The digraph defined by $\Phi_E$ has all maps from $B^A$ to $A$ as vertices. Recall that $\F_{\Pol(\B)}(\A)$ has all homomorphisms from $B^A$ to $A$ as vertices. However both digraphs have exactly the same edges, since $\Phi_E$ ensures in particular that only homomorphisms have edges. Therefore, both graphs are homomorphically equivalent and $\B\ppleq \F_{\Pol(\B)}(\A)$. Furthermore, $\A\hookrightarrow \F_{\Pol(\B)}(\A)$ with $a\mapsto \pi_a$, where $\pi_a(\boldsymbol b)=\boldsymbol b(a)$. By (5), $\F_{\Pol(\B)}(\A)\to\A$. Hence, $\B\ppleq\A$.  Note that we only needed (5) in the very last step, in particular $\B\ppleq \F_{\Pol(\B)}(\A)$ for any finite digraphs $\A$ and $\B$. 

There is one special case we did not yet consider in this construction: the case that $E(\A)$ is empty. Since there are no 0-ary polymorphisms we also have that $\F_{\Pol(\B)}(\A)$ has no edges and $\F_{\Pol(\B)}(\A)\to\A$. In this case we replace $\operatorname{isPol}(g)$ by $\bot$. 
\end{proof}

Observe that whenever $\Pol(\B)$ does not contain a constant polymorphism, then the formula $\bot$ can be replaced by a pp-formula that encodes 'there exists a constant polymorphism'. Hence, $\bot$ is only necessary to have 
$\drawLoop
\ppleq
\drawDot$. 
Let $\C_2\coloneqq (\{a,b\},\{(a,b),(b,a)\})$ and $\P_1\coloneqq(\{0,1\},\{(0,1)\})$. 
\begin{example}\label{exa:freeStructureWorks}
We want to verify that $\C_2\ppleq\P_1$.

(5) To determine $\F_{\Pol(\C_2)}(\P_1)$ we first determine some polymorphisms of $\C_2$.
The set of homomorphisms from $\C_2^{E(\P_1)}$ to $\C_2$ is $\{\operatorname{id}_{\{a,b\}},\overline{\cdot}\}$, where $\overline{a}=b$ and $\overline{b}=a$.
The set of homomorphisms from $\C_2^{\P_1}$ to $\C_2$ is
$\{\pi_0,\pi_1,\overline{\pi_0},\overline{\pi_1}\}$.
Note that 
\begin{align*}
\pi_0(x_0,x_1)&=x_0=\operatorname{id}_{\{a,b\}}(x_0), &\overline{\pi_0}(x_0,x_1)&=\overline{x_0},\\ \pi_1(x_0,x_1)&=x_1=\operatorname{id}_{\{a,b\}}(x_1),&
\overline{\pi_1}(x_0,x_1)&=\overline{x_1}. 
\end{align*}
Hence, $\F_{\Pol(\C_2)}(\P_1)$ consists of the two isolated edges $\pi_0\toEdge\pi_1$ and $\overline{\pi_0}\toEdge\overline{\pi_1}$. This graph has a homomorphism to $\P_1$, therefore $\C_2\ppleq\P_1$. 

(1) The pp-construction from the proof is $\Phi$. We can remove the unnecessary existentially quantified variables and obtain $\Phi'$, simplifying this construction further we obtain $\Phi''$. 
\begin{center}
\begin{tikzpicture}[scale=0.9]

    \node[var-f,label=below:{\small $f_1(aa)$}] (x1) at (0,0) {};
    \node[var-f,label=below:{\small $f_1(ab)$}] (x2) at (1.33,0) {};
    \node[var-f,label=below:{\small $f_1(ba)$}] (x3) at (2.66,0) {};
    \node[var-f,label=below:{\small $f_1(bb)$}] (x4) at (4,0) {};
    \node[var-f,label=above:{\small $f_2(aa)$}] (y1) at (0,2) {};
    \node[var-f,label=above:{\small $f_2(ab)$}] (y2) at (1.33,2) {};
    \node[var-f,label=above:{\small $f_2(ba)$}] (y3) at (2.66,2) {};
    \node[var-f,label=above:{\small $f_2(bb)$}] (y4) at (4,2) {};

    \node[var-b,label=left:{\small $g(a)$}] (g1) at (0,1) {};
    \node[var-b,label=right:{\small $g(b)$}] (g2) at (4,1) {};

    \path
        (x2) edge (x3)
        (x1) edge[bend left] (x4)
        (y2) edge (y3)
        (y1) edge[bend right] (y4)
        (g1) edge (g2)
        (x1) edge[dashed] (g1)
        (x2) edge[dashed] (g1)
        (x3) edge[dashed] (g2)
        (x4) edge[dashed] (g2)
        (y1) edge[dashed] (g1)
        (y2) edge[dashed] (g2)
        (y3) edge[dashed] (g1)
        (y4) edge[dashed] (g2)
        ;
    \node at (2,-1.5) {$\Phi$};
\end{tikzpicture}
\begin{tikzpicture}[scale=0.9]

    \node[var-f,label=below:{\small $f_1(aa)$}] (x1) at (0,0) {};
    \node[var-f,label=below:{\small $f_1(ab)$}] (x2) at (1.33,0) {};
    \node[var-f,label=below:{\small $f_1(ba)$}] (x3) at (2.66,0) {};
    \node[var-f,label=below:{\small $f_1(bb)$}] (x4) at (4,0) {};
    \node[var-f,label=above:{\small $f_2(aa)$}] (y1) at (0,2) {};
    \node[var-f,label=above:{\small $f_2(ab)$}] (y2) at (1.33,2) {};
    \node[var-f,label=above:{\small $f_2(ba)$}] (y3) at (2.66,2) {};
    \node[var-f,label=above:{\small $f_2(bb)$}] (y4) at (4,2) {};

    \path
        (x2) edge (x3)
        (x1) edge[bend left] (x4)
        (y2) edge (y3)
        (y1) edge[bend right] (y4)
        (x1) edge[dashed] (x2)
        (x3) edge[dashed] (x4)
        (y1) edge[dashed,bend right=20] (y3)
        (y2) edge[dashed,bend right=20] (y4)
        (y1) edge[dashed] (x1)
        (x4) edge[dashed] (y4)
        ;
    \node at (2,-1.5) {$\Phi'$};
\end{tikzpicture}
\begin{tikzpicture}[scale=0.9]

    \node[var-f,label=below:$x_1$] (x1) at (0,0) {};
    \node[var-f,label=below:$x_2$] (x2) at (1,0) {};
    \node[var-f,label=above:$y_1$] (y1) at (0,1) {};
    \node[var-f,label=above:$y_2$] (y2) at (1,1) {};
    \path
        (x1) edge (y1)
        (x2) edge[dashed] (x1)
        (x2) edge[dashed] (y2)
        ;
    \node at (0.5,-1.5) {$\Phi''$};
\end{tikzpicture}
\end{center}
The pp-power of $\C_2$ given by $\Phi''$ consists of the two edges $aa\toEdge ba$ and $bb\toEdge ab$ and is isomorphic to $\F_{\Pol(\P_1)}(\C_2)$ and homomorphically equivalent to $\P_1$ with $h_a(0)=aa$, $h_a(1)=ba$, $h_b(xy)=0$ if $x=y$ and $h_b(xy)=1$ if $x=y$.

(2) Let $f\in\Pol(\C_2)$, $\lambda(f)=h_b\circ f\circ h_a$, $x_1,\dots,x_n\in \{0,1\}$, and $a_i=a$ if $x_i=0$ and $a_i=b$ otherwise. Then 
\begin{align*}
    \lambda(f)(x_1,\dots,x_n)=
    \begin{cases}
    0&\text{if }f(a_1,\dots,a_n)=f(a,\dots,a)\\
    1&\text{otherwise}
    \end{cases}
\end{align*}
and $\lambda$ is a minor-preserving map from $\Pol(\C_2)$ to $\Pol(\P_1)$.

(4) We have that $\Sigma(\P_1, \F_{\Pol(\B)}(\A))$ consists of the identities
\begin{align*}
    f_1(x_0,x_1)&\approx g(x_0),&
    f'_1(x_0,x_1)&\approx g'(x_0),\\
    f_2(x_0,x_1)&\approx g(x_1)\text{, and}&
    f'_2(x_0,x_1)&\approx g'(x_1).
\end{align*}
Define an assignment $\tilde\cdot$ to polymorphisms of $\P_1$:
\begin{align*}
\tilde f_1&\coloneqq\tilde f'_1\coloneqq \lambda(\pi_0)=\lambda(\overline{\pi_0})=\pi_0,\\ 
\tilde f_2&\coloneqq\tilde f'_2\coloneqq \lambda(\pi_1)=\lambda(\overline{\pi_1})=\pi_1\text{, and}\\
\tilde g&\coloneqq\tilde g'\coloneqq\lambda(\operatorname{id}_{a,b})=\lambda(\overline{\cdot})=\operatorname{id}_{0,1}.
\end{align*}
The assignment $\tilde\cdot$ shows that $\Pol(\P_1)$ satisfies $\Sigma(\P_1, \F_{\Pol(\B)}(\A))$.
\end{example}

Next we present a minimal example on how to find a simple minor condition separating two structures using free structures.
\begin{example}\label{exa:freeStructureNotWorks}
We want to verify that $\P_1\not\ppleq\C_2$. 

(5) To determine $\F_{\Pol(\P_1)}(\C_2)$ we first determine polymorphisms of $\P_1$.
The set of homomorphisms from $\P_1^{\C_2}$ to $\P_1$ is $\{\pi_{a},\pi_{b},\min,\max\}$.
Note that 
\begin{align*}
\pi_a(x_a,x_b)=x_a&=\pi_{a}(x_a,x_b), &
\pi_b(x_a,x_b)=x_b&=\pi_{b}(x_a,x_b),\\ 
\pi_b(x_a,x_b)=x_b&=\pi_{a}(x_b,x_a),&
\pi_a(x_a,x_b)=x_a&=\pi_{b}(x_b,x_a),\\
\min(x_a,x_b)&=\min(x_b,x_a), & \max(x_a,x_b)&=\max(x_b,x_a).
\end{align*}
Hence, $\F_{\Pol(\P_1)}(\C_2)$ is
\begin{center}
    \begin{tikzpicture}[scale=1]
    \node (0) at (0,0) {$\pi_a$};
    \node (1) at (0,1) {$\pi_b$};
    \node (min) at (1,0) {$\min$};
    \node (max) at (1,1) {$\max$};
    \path (0) edge (1)
    (max) edge[loop above,>=stealth'] (max)
    (min) edge[loop above,>=stealth'] (min)
    ;
    
    \end{tikzpicture}
\end{center}
and does not have a homomorphism to $\C_2$, therefore $\P_1\not\ppleq\C_2$. 

(3) A minimal subgraph of $\F_{\Pol(\P_1)}(\C_2)$ that does not have a homomorphism to $\C_2$ is the loop $\max$. Using this subgraph we get the identity $\Sigma(\C_2,\drawLoop)\equiv\{f(x,y)=f(y,x)\}$ which is satisfied in $\Pol(\P_1)$ and not satisfied in $\Pol(\C_2)$.
\end{example}
In an ongoing research project we successfully used the ideas from Examples~\ref{exa:freeStructureWorks} and~\ref{exa:freeStructureNotWorks} in the program introduced in Chapter~\ref{cha:hardTrees} to verify that $\B\ppleq\A$ or to find a new minor condition that witnesses $\B\not\ppleq\A$ for structures $\A$ and $\B$ with binary signature and three elements.
One more useful observation about free structures is the following lemma. 
\begin{lemma}\label{lem:freeMinorConditionIsSatisfiedIff}
Let $\A$, $\A'$, and $\B$ be finite relational structures. Then 
\begin{align*}
    \Pol(\B)&\models\Sigma(\A,\A')&&\text{if and only if} &\A'&\to\F_{\Pol(\B)}(\A),
    \intertext{in particular}
\Pol(\A)&\models\Sigma(\A,\A')&&\text{if and only if} &\A'&\to\A.
\end{align*}
\end{lemma}
\begin{proof}[Proofsketch]
Let $\bar\dot$ be an assignment witnessing $\Pol(\B)\models\Sigma(\A,\A')$. Now the map $v\mapsto \tilde f_v$ is a homomorphism from $\A'$ to $\A'\to\F_{\Pol(\B)}(\A)$, as desired. 
Let $h\colon \A'\to\F_{\Pol(\B)}(\A)$ be a homomorphism. Then $f_v\mapsto h(v)$ is an assignment witnessing $\Pol(\B)\models\Sigma(\A,\A')$. 
The second part of the lemma follows from the fact that $\F_{\Pol(\A)}(\A)$ is homomorphically equivalent to $\A$.
\end{proof}
Note that we could have used this lemma in the proof of Theorem~\ref{thm:freestructure} for the direction $(4)\Rightarrow(5)$. 
Combining Theorem~\ref{thm:freestructure} with Lemma~\ref{lem:freeMinorConditionIsSatisfiedIff} we obtain the following lemma.
\begin{lemma}
Let $\A$ and $\B$ be finite relational structures and $\mathcal F$ be a set of finite structures such that for every finite structure $\C$ we have $\C\to\A$ if and only if $\F\not\to\C$ for all $\F\in\mathcal{F}$. Then
\begin{align*}
    \B\ppleq\A&&\text{if and only if} &&\Pol(\B)\not\models\Sigma(\A,\F)\text{ for all }\F\in\mathcal{F}.
\end{align*}
\end{lemma}
\begin{proof}
The proof is straight forward:
\begin{align*}
    \B\ppleq\A &\Leftrightarrow \F_{\Pol(\B)}(\A)\to \A\tag{Theorem~\ref{thm:freestructure}}\\
    &\Leftrightarrow \F\not\to \F_{\Pol(\B)}(\A)\text{ for all}\F\in\mathcal{F}\\
    &\Leftrightarrow \Pol(\B)\not\models\Sigma(\A,\F)\text{ for all}\F\in\mathcal{F}.\tag{Lemma~\ref{lem:freeMinorConditionIsSatisfiedIff}}
\end{align*}
As desired.
\end{proof}
Note that if $\mathcal F$ can be chosen to be a one element set $\{\F\}$, then $\B\models\A$ if and only if $\Pol(\B)\not\models\Sigma(\A,\F)$. In this case $\A$ is called a \emph{blocker} for $\Sigma(\A,\F)$. 


\section{Constraint Satisfaction Problems}

This section gives a short introduction to constraint satisfaction problems (CSPs). For more detailed information you can consult \cite{theBodirsky}. All structures in this section are assumed to have a finite relational signature.

\begin{definition}
Let $\A$ be a finite relational structure with finite signature $\tau$. The \emph{Constraint Satisfaction Problem of $\A$}, denoted $\csp(\A)$, is the following decision problem $\{\Instance\mid \Instance\text{ is a finite $\tau$-structure and } \Instance\to\A\}$.
\end{definition}
As we can translate between pp-formulas and structures we can also view pp-formulas as instances of CSPs: $\Phi\in\csp(\A)$ if the canonical database of $\Phi$ has a homomorphism to $\A$ or equivalently if $\Phi$ is satisfiable in $\A$.
Let us make some simple observations.

\begin{observation}
Let $\A$ be a finite structure. Then $\csp(\A)$ is in NP.
\end{observation}

\begin{observation}
Let $\A$ and $\B$ be structures that are homomorphically equivalent. Then $\csp(\A)=\csp(\B)$.
\end{observation}

\begin{example}
Many classical problems from theoretical computer science can be encoded as CSPs:
\begin{enumerate}
    \item $\csp(\K_3)=\{\G\mid\G\text{ is a finite directed 3-colourable graph}\}$
    \item $\csp(\K_n)=\{\G\mid\G\text{ is a finite directed n-colourable graph}\}$ 
    \item 3SAT can be encoded as $\csp(\{0,1\};R_{(a,b,c)}\mid a,b,c\in\{0,1\})$, where $R_{(a,b,c)}\coloneqq \{0, 1\}^3 \setminus \{(a, b, c)\}$.
    \item Horn-3SAT can be encoded as $\csp(\{0,1\};R_{(1,1,0)},\{0\},\{1\})$. 
    \item Solving systems of linear equations over $\Z_p$ can be encoded as \[\csp(\{0,\dots,p-1\};\{(x,y,z)\mid x+y+z\equiv_p 0\},\{0\},\dots,\{p-1\}).\]
    \item $\csp(\{0,1\};\leq,\{0\},\{1\})$ contains a finite graph $\G$ with labels 1 and 0 if and only if $\G$ has no directed path from a node labeled 1 to a node labeled 0. 
    
    \item $\csp(\{0,1\};=,\{0\},\{1\})$ contains a finite graph $\G$ if and only if $\G$ has no path from a node labeled 1 to a node labeled 0. 
    
    \item $\csp(\{0,1,2\};<)$ contains a finite graph $\G$ if and only if $\G$ has no directed path of length at least three. 
\end{enumerate}
Note that the problems in 1.~and 3.~are NP-complete. In 4., we have a P-complete problem. The problem in 5.~is Mod$_p$-complete~\cite{LaroseTesson}. Directed reachability presented in 6.~is NL-complete. Finally, 7.~and 8.~are examples of problems in L. 
\end{example}

It is well known that there is a connection between complexity o CSPs and pp-powers.

\begin{lemma}
Let $\A$ be a finite structure and $\A'$ be a pp-power of $\A$. Then $\csp(\A')\leq_{\operatorname{log}}\csp(\A)$, i.e., there is a log-space reduction from $\csp(\A')$ to $\csp(\A)$. 
\end{lemma}
To explain the idea of the reduction consider the case that $\A'$ is a digraph and the $k$-th pp-power of a structure $\A$ given by the pp-formula $\Phi_E(x_1,\dots,x_k,y_1,\dots,y_k)$. For any instance $\Instance'$ of $\csp(\A')$ we construct an instance $\Instance$ of $\csp(\A)$, such that $\Instance'\to\A'$ if and only if $\Instance\to\A$, as follows: 
\begin{enumerate}
    \item for each element $u$ of $\Instance'$ add $k$ copies of $u$ to $\Instance$, 
    \item for each edge $(u,v)$ in $\Instance'$ add a copy of the canonical database of $\Phi_E$ to $\Instance$ and identify $x_1,\dots,x_k$ and $y_1,\dots,y_k$ with the respective copies of $u$ and $v$ (note that if $\Phi_E$ uses $=$, then this might also identify some copies of $u$ with some copies of $v$).  
\end{enumerate}

As a result we obtain the following well-known connection between pp-constructions and CSPs.

\begin{lemma}
Let $\A$ and $\B$ be finite structures such that $\B\ppleq\A$. Then $\csp(\A)\leq_{\operatorname{log}}\csp(\B)$.
\end{lemma}
The converse is not true, take $\A=
\drawEdge$ and $\B=
\drawLoop$ as counterexample.
In 1999 Feder and Vardi made their famous dichotomy conjecture, which was proven independently by Bulatov and Zhuk in 2017 \cite{ZhukFVConjecture,BulatovFVConjecture}.
\begin{theorem}\label{thm:FederVardiConjSiggerIffNoK3}
Let $\A$ be a finite structure. Then exactly one of the following holds:
\begin{enumerate}
    \item $\A$ can pp-construct {\drawKThree}, and $\csp(\A)$ is $\NP$-hard,
    \item $\A$ has a Siggers polymorphism, and $\csp(\A)$ is in $\Pclass$.
\end{enumerate}
\end{theorem}
Consequently, there are no NP-intermediate CSPs. There is ongoing research by Bodirsky to obtain a similar result for CSPs over special infinite structures~\cite{theBodirsky}. In parallel to work on the dichotomy conjecture there was work on another one that might be easier to proof. Now Conjecture~\ref{conj:LaroseTesson} is the biggest open problem in finite CSPs.
Note that many problems arising around pp-constructions can be phrased as an instance of a CSP:
\begin{itemize}
    \item $\B\ppleq\A$ iff  $\F_{\Pol(\B)}(\A)\in\csp(\A)$ (Theorem~\ref{thm:freestructure}),
    \item $\A$ is homomorphically equivalent to $\B$ iff $\A\in\csp(\B)$ and $\B\in\csp(\A)$,
    \item computing the pp-power of a structure $\A$ given pp-formulas $\Phi_1,\dots,\Phi_n$ requires finding all satisfying assignments of $\Phi_1,\dots, \Phi_n$ under $\A$.
\end{itemize}

Even $\Pol(\A)\models\Sigma$ can be seen as an instance of $\csp(\A)$ in the following way.
\begin{definition}\label{def:indicatorStructure}
Let $\A$ be a finite structure and $\Sigma$ a minor condition. Define the \emph{indicator formula of $\Sigma$ with respect to $\A$} as the pp-formula which contains for each $f_\sigma\approx g_\tau$ in $\Sigma$ the conjunct 
\[\operatorname{isPol}(f)\AND \operatorname{isPol}(g)\AND \bigwedge_{\boldsymbol b}f_\sigma(\boldsymbol b)\approx g_\tau(\boldsymbol b),\]
where $\operatorname{isPol}(f)$ and $\operatorname{isPol}(g)$ are the pp-formulas defined in Section~\ref{sec:freestructures}. The canonical database of the indicator formula is called the \emph{indicator structure of $\Sigma$ with respect to $\A$} and is denoted by $\Ind(\Sigma,\A)$.
\end{definition}
Note that similar to the proof of Theorem~\ref{thm:ppPolyPreservation} the solutions of $\Ind(\Sigma,\A)$ as an instance of $\csp(\A)$ correspond to polymorphisms of $\A$ satisfying $\Sigma$. Hence, the following lemma holds.
\begin{lemma}
We have $\Pol(\A)\models\Sigma$ if and only if $\Ind(\Sigma,\A)\to\A$.
\end{lemma}

\begin{example}
Consider the condition $\Sigma_2=\{f(x,y)\approx f(y,x)\}$,  the directed 3-cycle $\C_3$, and the directed 2-cycle $\C_2$. In Figure~\ref{fig:indicatorGraph} we see why $\Ind(\Sigma_2,\C_3)\to \C_3$ and $\Ind(\Sigma_2,\C_2)\not\to \C_2$. 
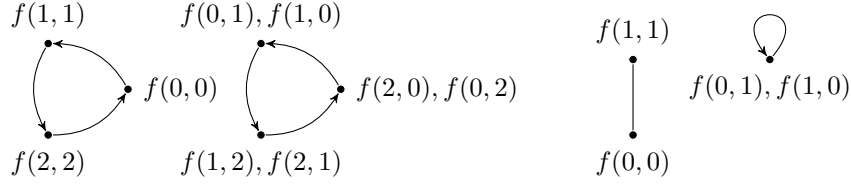
\begin{figure}
    \centering
    \begin{tikzpicture}[scale=0.7]
    \node[var-b,label={right:\small$f(0,0)$}] (0) at (0:1) {};
    \node[var-b,label={above:\small$f(1,1)$}] (1) at (120:1) {};
    \node[var-b,label={below:\small$f(2,2)$}] (2) at (240:1) {};
    
    \path [->,>=stealth']
    (0) edge[bend right] (1)
    (1) edge[bend right] (2)
    (2) edge[bend right] (0)
    ;
    \node[var-b,label={right:\small$f(2,0),f(0,2)$}] (0) at ($(4,0)+(0:1)$) {};
    \node[var-b,label={above:\small$f(0,1),f(1,0)$}] (1) at ($(4,0)+(120:1)$) {};
    \node[var-b,label={below:\small$f(1,2),f(2,1)$}] (2) at ($(4,0)+(240:1)$) {};
    
    \path [->,>=stealth']
    (0) edge[bend right] (1)
    (1) edge[bend right] (2)
    (2) edge[bend right] (0)
    ;
    \end{tikzpicture}
    \hspace{5mm}
    \begin{tikzpicture}
    \node[var-b,label={below:\small$f(0,0)$}] (00) at (0,0) {};
    \node[var-b,label={below:\small$f(0,1),f(1,0)$}] (01) at (1.8,1) {};
    \node[var-b,label={above:\small$f(1,1)$}] (11) at (0,1) {};
    
    \path 
    (01) edge[loop above,>=stealth',in=130,out=50,looseness=30,->] (01)
    (11) edge (00)
    ;
    \end{tikzpicture}
    \caption{The digraphs $\Ind(\Sigma_2,\C_3)$ (left) and $\Ind(\Sigma_2,\C_2)$ (right).}
    \label{fig:indicatorGraph}
\end{figure}
Hence, \[\Pol(\C_3)\models\Sigma_2\text{ and } \Pol(\C_2)\not\models\Sigma_2.\] Note that any homomorphism from $\Ind(\Sigma_2,\C_3)$ to $\C_3$ defines in an obvious way a polymorphism of $\C_3$ satisfying $\Sigma_2$. 
\end{example}



\section{Cores and Idempotent Polymorphisms}\label{sec:coresAndIdempotentcy}
In this section let us look at a very useful observation.
A finite structure $\A$ is called a \emph{core} if $\End(\A)=\Aut(\A)$. Cores are well studied and have very useful properties.
\begin{lemma}
Let $\A$ be a finite structure. Then there exists a core $\A'$ that is homomorphically equivalent to $\A$. Furthermore, this structure $\A'$ is unique up to isomorphism. Hence, we call $\A'$ \emph{the} core of $\A$.
\end{lemma}
An operation $f\colon A^n\to A$ is called \emph{idempotent} if $f(a,\dots,a)=a$ for all $a\in A$.
For an operation $f\colon A^n\to A$ define $\hat f\colon A\to A$ as $\hat f(a)=f(a,\dots,a)$ for all $a$. Note that $f$ is idempotent if and only if $\hat f=\operatorname{id}_A$. Let $\A$ be a structure. Note that if $f\in\Pol(\A)$, then $\hat f\in\End(\A)$. For $h\in\End(\A)$ define $\Pol(\A)|_h\coloneqq\{f\in\Pol(\A)\mid \hat f=h\}$. Note that whenever two operations $f$ and $g$ satisfy $f_\sigma=g_\tau$, then in particular \[\hat f=\hat f_\sigma=\hat g_\tau=\hat g.\] 
In other words: minor identities can only hold between operations that yield the same unary relation. Hence, when searching for minor-preserving maps from $\Pol(\A)$ to $\Pol(\B)$ we can partition $\Pol(\A)$ according to possible unary operations. 

\begin{lemma}\label{lem:minorPresMapsPartitionByEndos}
Let $\A$ and $\B$ be finite structures. Then $\Pol(\A)\tominor\Pol(\B)$ if and only if $\Pol(\A)|_h\tominor\Pol(\B)$ for every $h\in\End(\A)$. 
\end{lemma}

We leave the proof of the following lemma as an exercise to the reader.
\begin{lemma}\label{lem:coresBuildUp}
Let $\A$ be a finite structure and $h,h'\in\End(\A)$ such that there is an $\alpha\in\End(\A)$ with $h=\alpha\circ h'$. Then $\Pol(\A)|_{h'}\to \Pol(\A)|_{h}$, $f\mapsto \alpha\circ f$ is a minor-preserving map. 
If $\alpha$ is an automorphism, then this map is bijective and its inverse is also  minor-preserving.
\end{lemma}



From Theorem~\ref{thm:ppPolyPreservation} we can conclude the following lemma.

\begin{lemma}\label{lem:coresCanPPdefineConstants}
Let $\A$ be a finite structure. Then the following are equivalent:
\begin{enumerate}
    \item $\A$ is a core and for every $a\in A$ the relation $\{a\}$ is pp-definable,
    \item every polymorphism of $\A$ is idempotent, and
    \item $\Pol(\A)=\Pol(\A)|_{\operatorname{id}}$.
\end{enumerate}
\end{lemma}

Combining the previous lemmata we obtain the following well-known result.
\begin{lemma}\label{lem:coresCanPPconstructConstants}
Let $\A$ be finite core structure. Then $\A\ppeq\A_A$, where $\A_A$ is the structure $\A$ together with a unary relation $c_a\coloneqq\{a\}$ for each $a\in A$.
\end{lemma}
\begin{proof}
Clearly $\A_A\ppleq\A$. For the other direction it suffices to show, by Lemma~\ref{lem:coresBuildUp},  that $\Pol(\A)|_{h}\to \Pol(\A_A)$ for all $h\in\End(\A)$. Since $\A$ is a finite core we have that $\End(\A)=\Aut(\A)$. Let $h\in\Aut(\A)$. Then, by Lemma~\ref{lem:minorPresMapsPartitionByEndos}, $\Pol(\A)|_{h}\to \Pol(\A)|_{\operatorname{id}}$. 
By Lemma~\ref{lem:coresCanPPdefineConstants} and Theorem~\ref{thm:ppPolyPreservation} we have $\Pol(\A_A)=\Pol(\A_A)|_{\operatorname{id}}=\Pol(\A)|_{\operatorname{id}}$. Hence, $\Pol(\A)|_{h}\to \Pol(\A_A)$ and therefore $\A\ppleq\A_A$ as desired.
\end{proof}

We will sometimes write ``$\A$ with constants" to denote $\A_A$ and ``$\A$ can, using constants, pp-define" to denote ``$\A_A$ can pp-define". 
Often we will write $(x\approx a)$ instead of $c_a(x)$.
The last two lemmata show that in our pursuit to understand $\PPPoset$ it suffices to only consider structures that are cores and have for every element a constant (i.e., a singleton unary relation containing this element). Although we will mostly talk about digraphs in this thesis, whenever we have a core digraph we will assume that it has a constant for each node.


\section{Structures with Unary Signature}\label{sec:unarySignature}

Before we get to digraphs, let us first look at structures with an even simpler signature, i.e., every relation symbol is unary. 
Let $\Idemp\coloneqq(\{0,1\};c_0,c_1)$. 
With $c_0=\{0\}$ and $c_1=\{1\}$.
Note that $\Pol(\Idemp)$ consists of all idempotent operations on $\{0,1\}$.

\begin{lemma}\label{lem:idempotentCanPPconstructAllUnary}
Let $\A$ be any finite structure with unary signature. Then $\Idemp\ppleq\A$.
\end{lemma}
\begin{proof}
Let $\A$ be a structure with unary signature $\tau$ and $A=\{1,\dots,n\}$. Let $\A'$ be the $(n+1)$-st pp-power of $\Idemp$ given by the formulas $\phi_R$ for $R\in\tau$, where
\begin{align*}
    \phi_R(x_1,\dots,x_{n+1})\coloneqq c_0(x_1)\AND c_1(x_{n+1})\AND \bigwedge_{i\in [n]\setminus  R^{\A}} x_i\approx x_{i+1}.
\end{align*}

For $i\in A$ define the $(n+1)$-tuple $\boldsymbol t_i\coloneqq (0,\dots,0,1,\dots,1)$ containing $i$ zeros. The map $A\to A', i\mapsto \boldsymbol t_i$ is an embedding, since $i\in R^{\A}$ if and only if $\phi_R(\boldsymbol t_i)$ holds in $\A'$.
Next we find a homomorphism from $\A'$ to $\A$. Let $\boldsymbol t$ be an element of $\A'$. If the first entry of $\boldsymbol t$ is not 0 or the last one is not 1 then $\boldsymbol t$ is in no relation and we just map it to 1. Otherwise  we map $\boldsymbol t$ to the number of 0's before the first 1. If $\boldsymbol t$ is mapped to $i$ and is in some relation $R^{\A'}$, then the $i$-th and the $(i+1)$-st entry of $\boldsymbol t$ are different and $\A'\models\phi_R(\boldsymbol t)$ implies $i\in R^{\A}$. Hence, we constructed a homomorphism from $\A'$ to $\A$.
As desired we obtain $\Idemp\ppleq\A'\ppeq\A$. 
\end{proof}
We denote $\drawDot$ by $\P_0$.
Next we show that there is no structure strictly between $\Idemp$ and $\P_0$, no structure with unary signature and even no structure in $\PPPoset$.
\begin{theorem}\label{thm:idempotentIsCoatom}
The structure $\Idemp$ represents the unique coatom of $\PPPoset$. 
\end{theorem}
\begin{proof}
By Lemma~\ref{lem:loopIsTopElement} we have that $\Idemp\ppleq\P_0$.
Note that there is no constant polymorphism in $\Idemp$. Therefore,  $\Pol(\Idemp)\not\models f(x)\approx f(y)$. However $\Pol(\P_0)\models f(x)\approx f(y)$. Hence, $[\Idemp]\neq[\P_0]$ and $\Idemp$ lies strictly below $\P_0$.
Let $\A$ be a finite structure such that $\P_0\not\ppleq\A$. Assume without loss of generality that $\A$ is a core. By Lemma~\ref{lem:coresCanPPconstructConstants}, we can also assume that $\A$ has a constant for each element. Note that $\A$ must have at least two elements $a,b$, as otherwise $\P_0\ppleq\A$. Therefore, $\A$ can pp-construct $\Idemp$ using the formulas $\phi_{c_0}(x)=c_a(x)$ and $\phi_{c_1}(x)=c_b(x)$.
\end{proof}

\begin{corollary}
The finite unary structures ordered by pp-constructability form the following chain with two elements: 
\[[\Idemp]<[\P_0].\]
\end{corollary}

An operation $f$ is called \emph{conservative} if $f(x_1,\dots,x_k)\in\{x_1,\dots,x_k\}$ for all $x_1,\dots,x_k$. Observe that every conservative operation is in particular idempotent. Note that an operation is conservative if and only if it preserves all unary relations.
Lemma~\ref{lem:idempotentCanPPconstructAllUnary} and Theorem~\ref{thm:idempotentIsCoatom} yield the following corollary.
\begin{corollary}
Let $n\geq2$. Then $\Pol(\Idemp)$ is minor equivalent to the clone of all conservative operations on $\{1,\dots,n\}$.
\end{corollary}

\section{Undirected Graphs}
A special case of digraphs are undirected graphs, i.e., structures with one symmetric binary relation that has no loops.
The description of the poset of finite undirected graphs ordered by pp-constructability is an immediate consequence of the following theorem by Barto, Kozik, and Niven.

\begin{theorem}[Theorem 7.1 in \cite{BartoKozikNiven}]\label{thm:BartoKozikNivenSmoothDigraphs}
Let $\G$ be a smooth digraph (i.e., $\G$ has neither sources nor sinks). If $\G$ is a core and not a disjoint union of directed cycles then $\G$ can pp-construct any finite structure.
\end{theorem}

In particular we obtain the following corollary.
\begin{corollary}
Let $\A$ be a finite structure. Then $\K_3\ppleq\A\ppleq\C_1$.
\end{corollary}

\def\scale{0.60}

It is well known that $\K_3\ppleq\G$ for any digraph $\G$. The usual argument goes like this:
Clearly $\drawKThree$ is a core. Hence, we can add constants to $\drawKThree$.
It is easy to show that the polymorphisms of $\drawKThree$ with constants are exactly the projections. Hence, its polymorphism clone has a minor-preserving map into any other polymorphism clone.
However, I often struggled to find a pp-construction. Therefore, in the proof below, I want to share with you a nifty little gadget.
\begin{lemma}\label{lem:KnppdefinesAllGraphs}
Let $n\geq3$ and $\G$ be a digraph without loops and with vertices $0,\dots,n-1$. Then $\K_n$ can, using constants, pp-define $\G$.
\end{lemma}
\begin{proof}
Let $n\geq3$. For pairwise distinct $a,b,c\in\N_{< n}$ let $c_1,\dots,c_{n-3}$ be such that $\{a,b,c,c_1,\dots,c_{n-3}\}=\N_{< n}$ and define the  formula $\Phi_{abc}$ by the following graph representation (all edges are undirected): 
\[
\begin{tikzpicture}[scale=\scale]
    \node[var-f,label=left:$y$] (y) at (0,0) {};
    \node[var-f,label=left:$x$] (x) at (0,-1) {};
    \node[var-b] (x1) at (1,-1) {};
    \node[var-b] (y1) at (1,0) {};
    \node[var-b,label=left:$b$] (x2) at (1,-2) {};
    \node[var-b,label=left:$a$] (y2) at (1,1) {};
    \node[var-b,label=right:$c_1$] (c1) at (2,1) {};
    \node[var-b,label=right:$c_2$] (c2) at (2,0) {};
    \node[rotate=90] at (2,-1) {$\dots$};
    \node[var-b,label=right:$c_{n-3}$] (c3) at (2,-2) {};
    \path[>=stealth'] 
    (y) edge (x)
    (y1) edge (x1)
    (x) edge (x1)
    (x1) edge (x2)
    (y) edge (y1)
    (y1) edge (y2)
    (x1) edge (c1)
    (y1) edge (c1)
    (x1) edge (c2)
    (y1) edge (c2)
    (x1) edge (c3)
    (y1) edge (c3)
    ;
    \node at (1,-3) {$\Phi_{abc}$};
\end{tikzpicture}
\]
Observe that $\Phi_{abc}(x,y)$ is satisfied by all tuples of $\K_n$ that are not of the form $(x,x)$ except $(a,b)$,  
since the neighbours of $a$ and $b$ would both have to take value $c$. Hence, $\Phi_{abc}$ pp-defines the relation 
\[\{(x,y)\in\N_{<n}^2\mid x\neq y\}\setminus\{(a,b)\}.\] 
As the specific choice of $c$ is not relevant we will write $\Phi_{ab}$ instead of $\Phi_{abc}$.
Let $\G=(\N_{<n},E)$ be a digraph without loops and define  \[\overline E\coloneqq \{(u,v)\in\N_{<n}^2 \mid (u,v)\notin E, u\neq v\}.\] 
Then 
\[\Phi(x,y)\coloneqq \bigAND_{(u,v)\in \overline E}\Phi_{uv}(x,y)\]
is a pp-definition of $\G$ in $\K_n$.
\end{proof}

\begin{lemma}\label{lem:K3nisppPowerOfK3}
Let $n\geq1$. Then $\K_{3^n}$ is isomorphic to an $n$-th pp-power of $\K_3$.
\end{lemma}
\begin{proof}
Define 3-ary formula $\Psi$ and the binary formulas $\Psi_0$, $\Psi_1$, and $\Psi_2$ as follows
\[
\begin{tikzpicture}[scale=\scale]
    \node[var-f,label=left:$y$] (y) at (0,0) {};
    \node[var-f,label=left:$x$] (x) at (0,-3) {};
    \node[var-b] (x1) at (0,-1) {};
    \node[var-b] (y1) at (0,-2) {};
    \node[var-f,label=right:$z$] (z) at (1,-1.5) {};
    \path[>=stealth'] 
    (y1) edge (x)
    (x1) edge (y)
    (y1) edge (x1)
    (y1) edge (z)
    (x1) edge (z)
    ;
    \node at (0.5,-4) {$\Psi$};
\end{tikzpicture}
\hspace{10mm}
\begin{tikzpicture}[scale=\scale]
\node[var-f,label=left:$y$] (y1) at (0,0) {};
\node[draw,scale=0.6] (f1) at (0,-1.5) {$0$};
\node[var-f,label=left:$x$] (x1) at (0,-3) {};
\node at (1,-1.5) {$\coloneqq$};
\tikzset{decoration={snake,amplitude=.25mm,segment length=1.3mm,post length=0.6mm,pre length=0.6mm}}
    \draw[decorate] (x1) -- (f1);
    \draw[decorate,->,>=stealth'] (f1) -- (y1);
\node[var-f,label=left:$y$] (y2) at (2,0) {};
\node[var-f,label=left:$x$] (x2) at (2,-3) {};
\path (x2) edge[dashed] (y2);
    \node at (0,-4) {$\Psi_0$};
\end{tikzpicture}\hspace{10mm}
\begin{tikzpicture}[scale=\scale]
\node[var-f,label=left:$y$] (y1) at (0,0) {};
\node[draw,scale=0.6] (f1) at (0,-1.5) {$1$};
\node[var-f,label=left:$x$] (x1) at (0,-3) {};
\node at (1,-1.5) {$\coloneqq$};
\tikzset{decoration={snake,amplitude=.25mm,segment length=1.3mm,post length=0.6mm,pre length=0.6mm}}
    \draw[decorate] (x1) -- (f1);
    \draw[decorate,->,>=stealth'] (f1) -- (y1);
\node at (0,-4) {$\Psi_1$};
\node[var-f,label=left:$y$] (y3) at (2,0) {};
\node[var-f,label=left:$x$] (x3) at (2,-3) {};
\node[var-b,label=right:$2$] (z3) at (3,0) {};
\path (x3) edge (y3)
(z3) edge (y3);
\end{tikzpicture}\hspace{10mm}
\begin{tikzpicture}[scale=\scale]

\node at (0,-4) {$\Psi_2$};
\node[var-f,label=left:$y$] (y1) at (0,0) {};
\node[draw,scale=0.6] (f1) at (0,-1.5) {$2$};
\node[var-f,label=left:$x$] (x1) at (0,-3) {};
\node at (1,-1.5) {$\coloneqq$};
\tikzset{decoration={snake,amplitude=.25mm,segment length=1.3mm,post length=0.6mm,pre length=0.6mm}}
    \draw[decorate] (x1) -- (f1);
    \draw[decorate,->,>=stealth'] (f1) -- (y1);
\node[var-f,label=left:$y$] (y3) at (2,0) {};
\node[var-f,label=left:$x$] (x3) at (2,-3) {};
\node[var-b,label=right:$1$] (z3) at (3,0) {};
\path (x3) edge (y3)
(z3) edge (y3);
\end{tikzpicture}
\]
we use curly arrows as shorthand for $\Psi_0$, $\Psi_1$, and $\Psi_2$. 
Observe that in $\K_3$  for all $a\in\{0,1,2\}$ we have that $\Psi_a(a,x)$ is only satisfied by 0 and that $\Psi$ pp-defines the relation 
\[\{(x,y,z)\in\K_3^3\mid (x=y)\Rightarrow (z=x=y)\}.\]
Furthermore the formula $\exists z'\ \Psi(x,y,z')\AND\Psi_a(z',z)$ has the property that if $x=a=y$, then $z$  has to be 0. And if not $x=a=y$, then are at least two possible satisfying values for $z$.
For $(a_1,\dots,a_n)\in\{0,1,2\}^n$ define the formula
\[
\begin{tikzpicture}[scale=\scale]
    \node[var-f,label=above:$y_1$] (y1) at (0,0) {};
    \node[var-f,label=below:$x_1$] (x1) at (0,-3) {};
    \node[var-b] (a1) at (0,-1) {};
    \node[var-b] (b1) at (0,-2) {};
    \node[var-b] (z1) at (1,-1.5) {};
    \path[>=stealth'] 
    (b1) edge (x1)
    (a1) edge (y1)
    (b1) edge (a1)
    (b1) edge (z1)
    (a1) edge (z1)
    ;
    
    \node[var-f,label=above:$y_2$] (y2) at (2,0) {};
    \node[var-f,label=below:$x_2$] (x2) at (2,-3) {};
    \node[var-b] (a2) at (2,-1) {};
    \node[var-b] (b2) at (2,-2) {};
    \node[var-b] (z2) at (3,-1.5) {};
    \path[>=stealth'] 
    (b2) edge (x2)
    (a2) edge (y2)
    (b2) edge (a2)
    (b2) edge (z2)
    (a2) edge (z2)
    ;
    
    \node[var-f,label=above:$y_3$] (y3) at (4,0) {};
    \node[var-f,label=below:$x_3$] (x3) at (4,-3) {};
    \node[var-b] (a3) at (4,-1) {};
    \node[var-b] (b3) at (4,-2) {};
    \node[var-b] (z3) at (5,-1.5) {};
    \path[>=stealth'] 
    (b3) edge (x3)
    (a3) edge (y3)
    (b3) edge (a3)
    (b3) edge (z3)
    (a3) edge (z3)
    ;
    
    \node[var-f,label=above:$y_4$] (y4) at (6,0) {};
    \node[var-f,label=below:$x_4$] (x4) at (6,-3) {};
    \node[var-b] (a4) at (6,-1) {};
    \node[var-b] (b4) at (6,-2) {};
    \node[var-b] (z4) at (7,-1.5) {};
    \path[>=stealth'] 
    (b4) edge (x4)
    (a4) edge (y4)
    (b4) edge (a4)
    (b4) edge (z4)
    (a4) edge (z4)
    ;

    \node[var-b] (c1) at (1,1.5) {};
    \node[var-b] (d1) at (3,1.5) {};
    \node[var-b] (e1) at (2,2.5) {};
    
    \path
    (c1) edge (d1)
    (d1) edge (e1)
    (e1) edge (c1)
    ;
    
    \node[var-b] (c2) at (4,2.5) {};
    \node[var-b] (d2) at (5,2.5) {};
    \node[var-b] (e2) at (5,1.5) {};
    
    \path
    (c2) edge (d2)
    (d2) edge (e2)
    (e2) edge (c2)
    (e1) edge (c2)
    ;

    \node[var-b] (c4) at (6,2.5) {};
    \node[var-b] (d4) at (7,2.5) {};
    \node[var-b] (e4) at (7,1.5) {};
    
    \path
    (c4) edge (d4)
    (d4) edge (e4)
    (e4) edge (c4)
    (d2) edge (c4)
    ;
    
    \tikzset{decoration={snake,amplitude=.25mm,segment length=1.3mm,post length=0.6mm,pre length=0.6mm}}
    
    \node[draw,scale=0.6] (f44) at (7,0) {$a_{4}$};
    
    \draw[decorate] (z4) -- (f44);
    \draw[decorate,->,>=stealth'] (f44) -- (e4);
    
    \node at (8.5,0) {$\dots$};
    \node at (8,-3) {$\dots$};
    \node at (8.5,-1.5) {$\dots$};
    \node at (9,2.5) {$\dots$};

    \node[var-f,label=above:$y_{n-2}$] (y4) at (10,0) {};
    \node[var-f,label=below:$x_{n-2}$] (x4) at (10,-3) {};
    \node[var-b] (a4) at (10,-1) {};
    \node[var-b] (b4) at (10,-2) {};
    \node[var-b] (z4) at (11,-1.5) {};
    \path[>=stealth'] 
    (b4) edge (x4)
    (a4) edge (y4)
    (b4) edge (a4)
    (b4) edge (z4)
    (a4) edge (z4)
    ;
    
    \node[var-f,label=above:$y_{n-1}$] (y5) at (12,0) {};
    \node[var-f,label=below:$x_{n-1}$] (x5) at (12,-3) {};
    \node[var-b] (a5) at (12,-1) {};
    \node[var-b] (b5) at (12,-2) {};
    \node[var-b] (z5) at (13,-1.5) {};
    \path[>=stealth'] 
    (b5) edge (x5)
    (a5) edge (y5)
    (b5) edge (a5)
    (b5) edge (z5)
    (a5) edge (z5)
    ;

    \node[var-f,label=above:$y_n$] (y6) at (14,0) {};
    \node[var-f,label=below:$x_n$] (x6) at (14,-3) {};
    \node[var-b] (a6) at (14,-1) {};
    \node[var-b] (b6) at (14,-2) {};
    \node[var-b] (z6) at (15,-1.5) {};
    \path[>=stealth'] 
    (b6) edge (x6)
    (a6) edge (y6)
    (b6) edge (a6)
    (b6) edge (z6)
    (a6) edge (z6)
    ;
    
    \node[var-b] (c3) at (13,1.5) {};
    \node[var-b] (d3) at (15,1.5) {};
    \node[var-b] (e3) at (14,2.5) {};
    
    \path
    (c3) edge (d3)
    (d3) edge (e3)
    (e3) edge (c3)
    ;
    
    \node[var-b] (c4) at (10,2.5) {};
    \node[var-b] (d4) at (11,2.5) {};
    \node[var-b] (e4) at (11,1.5) {};
    
    \path
    (c4) edge (d4)
    (d4) edge (e4)
    (e4) edge (c4)
    (e3) edge (c4)
    ;
    
    \node[draw,scale=0.6] (f1) at (1,0) {$a_{1}$};
    \node[draw,scale=0.6] (f2) at (3,0) {$a_{2}$};
    \node[draw,scale=0.6] (f3) at (5,0) {$a_{3}$};

    \node[draw,scale=0.6] (f4) at (11,0) {$a_{n-2}$};
    \node[draw,scale=0.6] (f5) at (13,0) {$a_{n-1}$};
    \node[draw,scale=0.6] (f6) at (15,0) {$a_n$};
    
    \draw[decorate] (z4) -- (f4);
    \draw[decorate,->,>=stealth'] (f4) -- (e4);
    \draw[decorate] (z5) -- (f5);
    \draw[decorate,->,>=stealth'] (f5) -- (c3);
    \draw[decorate] (z6) -- (f6);
    \draw[decorate,->,>=stealth'] (f6) -- (d3);
    \draw[decorate] (z1) -- (f1);
    \draw[decorate,->,>=stealth'] (f1) -- (c1);
    \draw[decorate] (z2) -- (f2);
    \draw[decorate,->,>=stealth'] (f2) -- (d1);
    \draw[decorate] (z3) -- (f3);
    \draw[decorate,->,>=stealth'] (f3) -- (e2);
    
    \node at (7.5,-5) {$\Phi_{(a_1,\dots,a_n)}$};
\end{tikzpicture}
\]
From the previous observations it follows that $\Phi_{(a_1,\dots,a_n)}$ is satisfied by every tuple in $\K_3^{2n}$ except $(a_1,\dots,a_n,a_1,\dots,a_n)$.
Hence, the $n$-th pp-power of $\K_3$ given by the formula 
\[\bigwedge_{(a_1,\dots,a_n)\in\K_3^n} \Phi_{(a_1,\dots,a_n)}(x_1,\dots,x_n,y_1,\dots,y_n)\]
is isomorphic to $\K_{3^n}$.
\end{proof}

The following is an immediate consequence of the previous two lemmata.
\begin{theorem}\label{cor:K3IsBottom}
Let $\G$ be a finite digraph with at most $3^n$ vertices. Then $\G$ is homomorphically equivalent to an $n$-th pp-power of $\K_3$ with constants.  
\end{theorem}
\begin{proof}
Let $\G'$ be the $3^n$-element digraph obtained from $\G$ by adding isolated vertices. Clearly, $\G$ and $\G'$ are homomorphically equivalent. By Lemma~\ref{lem:KnppdefinesAllGraphs} $\G'$ is pp-definable from $\K_{3^n}$, which is, by Lemma~\ref{lem:K3nisppPowerOfK3}, an $n$-th pp-power of $\K_3$. Therefore, $\G'$ is an $n$-th pp-power of $\K_3$ as well.
\end{proof}

Note that the only undirected graph that is a core and not smooth is the graph consisting of a single vertex and no edges. Furthermore, the only undirected graph that is a smooth core and  is a disjoint union of directed cycles is the 2-cycle $\C_2$. Note that $\C_2$ can pp-construct the isolated vertex using the formula $\bot$. Hence, we obtain the following corollary.
\begin{corollary}
The finite undirected graphs ordered by pp-constructability form the following chain with three elements: 
\[[\K_3]<[\C_2]<[\P_0].\]
\end{corollary}

\section{Directed Graphs}
Now we come to the structures that are the focus of this thesis: digraphs, i.e., structures with one binary relation. Let us first make some simple observations.
If $\G$ is a digraph with a loop, then $\G$ is homomorphically equivalent to $\drawLoop$, denoted by $\C_1$. Note that $[\C_1]=[\P_0]$.
Observe that the digraph $\P_1\coloneqq(\{0,1\},\{(0,1)\})$ is in the same pp-constructability class as $\Idemp$. Hence, by Theorem~\ref{thm:idempotentIsCoatom}, $[\P_1]$ is the unique coatom of $\DGPoset$. We also know that $\K_3$ is the bottom element of $\DGPoset$. The goal of this thesis is to shed some light on the interval between $\P_1$ and $\K_3$, see Figure~\ref{fig:DGPosetVersion1}. 
\begin{figure}
    \centering
\begin{tikzpicture}
\node (C1) at (0,0) {$[\C_1]$};
\node (I) at (0,-1) {$[\P_1]$};
\node (K3) at (0,-4) {$[\K_3]$};
\path 
(K3) edge (-1,-3)
(K3) edge (1,-3)
(I) edge (-1,-2)
(I) edge (1,-2)
;
\node at (1,-2.5) {$\vdots$};
\node at (-1,-2.5) {$\vdots$};

\node at (0,-2) {$\dots$};
\node at (0,-3) {$\dots$};

\path (I) edge (C1);
\end{tikzpicture}

\caption{The poset $\DGPoset$ with a huge unknown region in the middle.}
\label{fig:DGPosetVersion1}
\end{figure}

\section{From Structures to Digraphs}
A lot of the results for finite digraphs translate to finite structures. The reason is a theorem from Bul\'in, Deli\'c, Jackson, and Niven. To present it we need some definitions.
Define $\T_3\coloneqq (\{0,1,2\},\{(0,1),(0,2),(1,2)\})$, i.e., the transitive tournament on three vertices. A minor identity \[f(x_1,\dots,x_n)\approx g(y_1,\dots,y_m)\] is \emph{balanced} if $\{x_1,\dots,x_n\}=\{y_1,\dots,y_m\}$.

\begin{theorem}[Theorem 1.4 in~\cite{BulinDelicJacksonNiven}]\label{thm:StructurToDigraph}
Let $\A$ be a finite core structure over a finite signature. Then there exists a finite core digraph $\mathcal{D}(\A)$ such that the following holds:
\begin{enumerate}
    \item $\csp(\A)$ and $\csp(\mathcal{D}(\A))$ are equivalent under log-space reductions,
    \item $\mathcal D(\A)$ can pp-construct $\A$, and
    \item if $\Sigma$ is a minor condition, such that $\Pol(\T_3)\models\Sigma$ and every identity in $\Sigma$ is either balanced or contains at most two variables, then
    \begin{align*}
        &\Pol(\A)\models\Sigma&&\text{if and only if}&&\Pol(\mathcal{D}(\A))\models\Sigma.
    \end{align*}
\end{enumerate}
\end{theorem}

From the minor conditions occurring in this thesis the following ones satisfy the condition in (ii):
\begin{itemize}
    \item $\Sigma_p$ for any prime $p$, see Section~\ref{sec:notationForCycles},
    \item $\Sigma_C$ for any finite set $C\subset\N$,  see Section~\ref{sec:notationForCycles},
    \item $\Majority$, see Definition~\ref{def:nu}, 
    \item $\Siggers$, see Definition~\ref{def:siggers},
    \item $\WNU{3,4}$, see Definition~\ref{def:wnu34},
    \item $\WNU{n}$ for any $n\geq 1$, see Definition~\ref{def:wnu},
    \item $\NU{n}$ for any $n\geq 3$, see Definition~\ref{def:nu},
    \item $\HMcK n$ for any $n\geq 0$, see Definition~\ref{def:hmck},
    \item $\HM n$ for any $n\geq 1$, see Definition~\ref{def:hm},
    \item $\J n$ for any $n\geq0$, see Definition~\ref{def:j},
    \item $\KK n$ for any $n\geq2$, see Definition~\ref{def:kk},
    \item $\TS n$ for any $n\geq 2$, see Definition~\ref{def:ts},
\end{itemize}

and the following ones do not satisfy the condition in (ii):
\begin{itemize}
    \item $\Maltsev$, see Definition~\ref{def:malt},
    \item $\GFS n$ for $n\geq 1$, see Definition~\ref{def:elevatorGFS},
    \item $\NN n$ for $n\geq0$, see Definition~\ref{def:nn}.
\end{itemize}

Recall, that Conjecture~\ref{conj:LaroseTesson} states that if a finite structure $\A$ can neither pp-construct st-Con not $\TLinP$ for any prime $p$, then $\csp(\A)$ is in L. By Corollary~\ref{cor:OrdisLeastLowerBoundForTn} and Theorem~\ref{thm:DLiff34WNUiff3Linp}  this condition can be characterized by the minor conditions $\HM n$ and $\WNU{3,4}$.  Note that these minor conditions satisfy condition (ii) in Theorem~\ref{thm:StructurToDigraph}. Hence, Conjecture~\ref{conj:LaroseTesson} is true for all finite structures if and only if it is true for all finite digraphs.


\chapter[Smooth Digraphs]{Smooth Digraphs\\
\Large{or Cycles Everywhere}}\label{cha:cycles}

In this chapter we will consider a large class of digraphs, i.e.,  \emph{finite smooth digraph}, which are finite directed graphs that have neither sources nor sinks. We will give a comparably simple description of the pp-constructability order on these graphs and of the subposet of $\DGPoset$ consisting of smooth digraphs. In particular, we show that every finite poset embeds into $\DGPoset$. 
Most of the results of this chapter have been published in \cite{StarkeVucajBodirskySmoothDigraphs}. 

A digraph $\G$ is a \emph{disjoint union of cycles} if every node of $\G$ has exactly one incoming and one outgoing edge; in particular, $\G$ is smooth.
Barto, Kozik, and Niven~\cite{BartoKozikNiven} showed  the following dichotomy for smooth digraphs.

\begin{theorem}[Theorem 3.1 in~\cite{BartoKozikNiven}]\label{thm:BartoKozikNiven}
Let $\mathbb G$ be a finite smooth digraph. Then either $\mathbb G$ pp-constructs $\K_3$ or it is homomorphically equivalent to a finite disjoint union of directed cycles.
\end{theorem}

Define $\SDPoset$ to be the subposet of $\PPPoset$ consisting of finite disjoint unions of directed cycles.
We know that $[\K_3]$ is the bottom element of $\PPPoset$. From the previous theorem we know that every non-minimal element in $\PPPoset$, represented by a smooth digraph, can also be represented by a finite disjoint union of directed cycles. Hence, we obtain the following corollary.

\begin{corollary}
The subposet of $\DGPoset$ consisting of smooth digraphs is isomorphic to the poset $\SDPoset$ with a single element $[\K_3]$ added at the bottom.
\end{corollary}

Following this corollary, the big goal of this chapter is to understand the poset $\SDPoset$.
We start by introducing some notation.

\section{Notation for Cycles}\label{sec:notationForCycles}
I like to think of the disjoint union of a cycle of length 3 and one of length 2 as just the set $\{2,3\}$, as this set contains all the important information about this structure. The notation introduced here has the goal of making the transition between the structure and the set as easy as possible.
To any finite set $C\subset\N^+\!$ we associate a finite disjoint union of cycles $\C=(V,E)$ defined by 
\begin{align*}
    V \coloneqq \{(a,k) \mid a\in C, k \in \mathbb Z_a\}  \text{ and }  E \coloneqq \{((a,k),(a,k+_a 1)) \mid a\in C, k \in \mathbb Z_a\}, 
\end{align*} where by $+_a$ we denote the addition modulo $a$.
For the sake of notation we will from now on write $+$ instead of $+_a$; it will be clear from the context to which addition we are referring to.
For any $a_1,\dots,a_n\in\N^+\!$ we write $\C_{a_1,\dots,a_n}$ for the finite disjoint union of cycles associated to the set $\{a_1,\dots,a_n\}$. Note that, for $a\in\N^+$\!, the structure $\C_a$ is a directed cycle of length $a$. For readability we will usually denote the elements of $\C_a$ by $0,\dots,a-1$, as in the original definition of $\C_a$, instead of $(a,0),\dots,(a,a-1)$.
To any finite disjoint union of cycles $\C$ we associate the set $C\coloneqq\{a\mid \C_{a}\hookrightarrow \C\}$.

Be warned, dear reader, that usually $C$ denotes the underlying set of the structure $\C$,
but from now on $\C$ itself will denote the structure as well as the underlying set and $C$
is the set defined above. We hope that this will not lead to any confusion. Also beware
that the finite disjoint union of cycles $\D$ associated to the set associated to a finite disjoint union of cycles $\C$ is
not necessarily isomorphic to $\C$. The structure $\C$ could have multiple cycles of the same length whereas $\D$ may not. However, $\C$ and $\D$ are homomorphically equivalent and thus
represent the same element in $\DGPoset$. 
We define for every $k\in\N^+$ the pp-formula
\begin{align*}
    x\stackrel{k}\toEdge z&\coloneqq \exists y_1,\dots,y_{k-1}.\ E(x,y_1)\wedge E(y_1,y_2)\wedge\dots\wedge E(y_{k-1},z)\\
    &=\tikz[scale=0.6,baseline=-1mm]{
\node[var-f,label=below:{$x$}] (0) at (0,0) {};
\node[var-b,label=below:{$y_1$}] (1) at (1,0) {};
\node[var-b,label=below:{$y_2$}] (2) at (2,0) {};
\node[] (3) at (3,0) {$\dots$};
\node[var-b,label=below:{$y_{k-1}$}] (4) at (4,0) {};
\node[var-f,label=below:{$z$}] (5) at (5,0) {};
\path[>=stealth',->]
        (0) edge (1)
        (1) edge (2)
        (4) edge (5)
        ;
},
\end{align*}
which we will often use in pp-constructions.
For a disjoint union of cycles $\C$; $u,v\in \C$ and $k\in\N^+\!$ we have that $\C\models u\stackrel{k}\toEdge v$ if there is a  directed  path of length $k$ from $u$ to $v$. If $\C$ is clear from the context we abbreviate $\C\models u\stackrel{k}\toEdge v$ by $ u\stackrel{k}\toEdge v$. 
Note that for fixed $u$ and $k$ there is exactly one $v \in \C$ such that $ u\stackrel{k}\toEdge v$; we denote this element $v$ by $u+k$.
The \emph{$k$-th relational power of $\C$} is  the digraph $(\C,\{(u,u+k)\mid u\in\C\})$. 
We denote the map $\C\to\C$, $u\mapsto u+1$ by $\sigma_{\C}$. Note that $\sigma_{\C} \in \operatorname{Aut}(\C)$.

Observe that the cyclic loop condition $\Sigma_{\C}$, introduced in Definition~\ref{def:SigmaG}, is $f\approx f_{\sigma_{\C}}$ and that the edge relation of $\C$ is $\{(u,\sigma_{\C}(u))\mid u\in\C\}$. If $\C$ is the finite disjoint union of cycles associated to the set $C$, then we write $\sigma_{\! C}$ and $\Sigma_{C}$ instead of $\sigma_{\C}$ and $\Sigma_{\C}$, respectively. Note that $\sigma_{\! C}(a,k) = (a,k+1)$ for any $a\in C$ and $k\in\{0,\dots,a-1\}$.

We want to remark that any finite power $\mathbb C^n$ of a disjoint union of cycles $\mathbb C$ is again a disjoint union of cycles. Hence, for an element $\boldsymbol{t}\in\mathbb C^n$, $k\in\N^+\!$, we have that $\boldsymbol{t}+k$ is already defined, furthermore \[\boldsymbol{t}+k=(t_1+k,\dots,t_n+k).\]

For $a,c\in\N^+\!$ define $a\dotdiv c\coloneqq\frac{a}{\gcd(a,c)}$. Note that $a\dotdiv c$ is always a natural number.
The choice of the symbol $\dotdiv$ is meant to emphasize that $a\dotdiv c$ is the numerator of the fraction $a\div c$
in reduced form. Roughly speaking, the operation $\dotdiv$ should be understood as ``divide as much as you can". 
The operation $\dotdiv$ has the following useful properties.
\begin{lemma}\label{lem:dotdiv}
For all $a,b,c\in\N^+\!$ we have
\begin{enumerate}
    \item $a\dotdiv (a\dotdiv c)=\gcd(a,c)$,
    \item $(a\dotdiv b)\dotdiv c=a\dotdiv (b\cdot c)$,
    \item $\gcd(a\dotdiv c,c\dotdiv a)=1$, and
    \item $a \dotdiv c=1$ if and only if $a$ divides $c$. 
\end{enumerate}
\end{lemma}
\begin{proof}
Let $a,b,c\in\N^+\!$. Simply applying the definitions we obtain:
\[a\dotdiv (a\dotdiv c)=\frac{a}{\gcd(a,a\dotdiv c)}=\frac{a}{\gcd\left(a,\frac{a}{\gcd(a,c)}\right)}
    =\frac{a}{\frac{a}{\gcd(a,c)}}
    =\gcd(a,c)\]
The reader may verify the other statements.
\end{proof}
For a finite disjoint union of cycles $\C$ and $c\in\N^+\!$, we let $\C\dotdiv c$ denote the finite disjoint union of cycles associated to the set $C\dotdiv c\coloneqq\{a\dotdiv c\mid a\in C\}$.

\begin{lemma}\label{lem:relPowEqualsCdotdivc}
Let $\C$ be a finite disjoint union of cycles and $c\in\N^+\!$. The $c$-th relational power of $\C$, i.e., the digraph with the same elements as $\C$ that has an edge from $u$ to $v$ whenever $\C$ has a directed path of length $c$ from $u$ to $v$,  is homomorphically equivalent to $\C\dotdiv c$.
\end{lemma}
\begin{proof}
Note that it suffices to verify the claim for cycles.
Let $\C=\C_a$. Then $\C\dotdiv c$ consists of one cycle of length $a\dotdiv c$ and the  $c$-th relational power of $\C$ consists of $\frac{a}{a\dotdiv c}=\gcd(a,c)$ many cycles of length $a\dotdiv c$. Hence, they are homomorphically equivalent.
\end{proof}

\begin{example}\label{ex.division}
Recall that Example~\ref{exa:C6LeqC3} showed $\C_6\ppleq\C_3$ using the formula $x\stackrel{2}{\toEdge}y$. 
Note that $\C_6$ pp-constructs $\C_2$ with the formula $\Phi_E(x,y)\coloneqq x\stackrel{3}{\toEdge}y$.
In general, the first pp-power of a disjoint union of cycles $\C$ given by the formula $x\stackrel{c}{\toEdge}y$ is the $c$-th relational power of $\C$, which is, by Lemma~\ref{lem:relPowEqualsCdotdivc}, homomorphically equivalent to $\C\dotdiv c$. Therefore, $\C\leq \C\dotdiv c$.
\end{example}

\section{Loop Conditions}

In general, showing that there is no minor-preserving map from $\Pol(\B)$ to $\Pol(\A)$ is a rather complicated task.
However, a result by Barto, Bul\'in, Krokhin, and Opr\v{s}al provides a concrete minor condition to check~\cite{Barto_2021FreeStructures}.
We show that, for disjoint unions of cycles $\mathbb G$, $\mathbb H$, whenever $\mathbb G\nleq\mathbb H$ there is a single minor identity with only one function symbol witnessing this. minor identities of this form have been studied in the literature and are called \emph{loop conditions}~\cite{olsak-loop, olsak-strong}. 
 
\begin{definition}
Let $\sigma,\tau\colon [m] \to [n]$ be maps. A \emph{loop condition} is a minor identity of the form
\begin{equation*}\label{Sigma}
    f_\sigma \approx f_\tau.
\end{equation*}
\end{definition}
To any loop condition $\Sigma$ we can assign a digraph in a natural way:
\begin{definition}
Let $\sigma,\tau\colon [m] \to [n]$ be maps and let $\Sigma$ be the loop condition, given by the identity $f_\sigma \approx f_\tau$. We define the digraph 
\[\G_\Sigma\coloneqq
([n],\{(\sigma(i),\tau(i))\mid i\in [m]\}).\] 
\end{definition}

\begin{examplex}\label{ex:G_Sigma} Some loop conditions and the corresponding digraphs.
\begin{itemize}
    \item Let $\Sigma_S$ be the loop condition $f(x,y,x,z,y,z)\approx f(y,x,z,x,z,y)$. Then $\mathbb G_{\Sigma_S}$ is isomorphic to $\K_3$.
    \item Let $\Sigma_3$ be the loop condition $f(x,y,z) \approx f(y,z,x)$. Then $\mathbb G_{\Sigma_3}$ is isomorphic to a directed cycle of length 3. \eoe
\end{itemize}
\end{examplex}
Observe that, for every digraph $\mathbb G$, all loop conditions $\Sigma$ such that $\mathbb G_{\Sigma}$ is isomorphic to $\mathbb G$ are equivalent. For convenience, we will from now on allow any finite set in place of $[m]$ and $[n]$. This allows us to construct a concrete loop condition from any graph.
\begin{definition}\label{def:SigmaG}
Let $\mathbb G=(V,E)$ be a digraph. We define the loop condition $\Sigma_{\mathbb G}\coloneqq (f_\sigma\approx f_\tau)$, where $\sigma,\tau\colon E\to V$ with $\sigma(u,v)=u$ and $\tau(u,v)=v$.
\end{definition} 
Observe that $\mathbb G_{\Sigma_\mathbb G}=\mathbb G$.
The name loop condition is justified by the following observation. If $\mathbb G$ is a finite graph such that $\Pol(\mathbb G)$ satisfies
$\Sigma$ and $\mathbb G_\Sigma \to \mathbb G$, then $\mathbb G$ has a loop. 
Consequently, if $\mathbb G$ does not have a loop, then $\Pol(\mathbb G)\not\models\Sigma_{\mathbb G}$. 
If $\mathbb G_\Sigma$ itself has a loop, then there is an $i$ with $\sigma(i)=\tau(i)$ and a structure $\A$ satisfies $\Sigma$ with the projection $\pi_i\colon\boldsymbol{a}\mapsto a_i$ and therefore $\Sigma$ is trivial. 
If $\mathbb G_\Sigma$ is a disjoint union of directed cycles, then we say that $\Sigma$ is a \emph{cyclic loop condition}. For instance, the identity $\Sigma_3$ in Example~\ref{ex:G_Sigma} is a cyclic loop condition.

\section[Unions of Prime Cycles]{Unions of Prime Cycles and Prime Cyclic Loop Conditions}\label{sec:primefun}

This section has two goals.
Firstly, we want to understand $\mPCL$ and $\PCPoset$, defined below, as these posets will be used to describe $\SDPoset$, i.e., the poset of \emph{finite disjoint union of directed cycles ordered by pp-constructability}.
Secondly, this section gently introduces the reader to techniques used in Sections~ \ref{sec:conditions} and~\ref{sec:structures}  without many of the difficulties that come with the more general case.

\begin{definition}
A cyclic loop condition $\Sigma_P$ is called a \emph{prime cyclic loop condition} if $P$ is a set of primes. 
A disjoint unions of cycles $\mathbb P$ is called a \emph{disjoint unions of prime cycles} if $P$ is a set of primes.
We define the posets $\mPCL$ and $\PCPoset$ as 
\begin{align*}
    \mPCL&\coloneqq (\{[\Sigma] \mid \Sigma \text{ a prime cyclic loop condition}\},\Leftarrow),
    \\
    \PCPoset&\coloneqq (\{[\mathbb P]\mid \mathbb P\text{ a disjoint union of prime cycles}\},\leq).
\end{align*}

\end{definition} 

To gain a better understanding of the objects we are working with, we reproduce a well-known fact about cycles and cyclic loop conditions, presented in Lemma~\ref{lem:pcSatisfyPclc}. We start with a simple example.
\begin{example}\label{exa:C3notmodelsS3}
Consider the digraph $\C_3$. Observe that $\Pol(\C_3) \models \Sigma_2$ as witnessed by the polymorphism
\[f(x,y) = 2\cdot(x + y) \pmod 3.\]
On the other hand, $\Pol(\C_3) \not\models \Sigma_3$. 
Assume that $f$ is a polymorphism of $\C_3$ satisfying $\Sigma_3$, then
\[
\begin{tikzpicture}
\node at (0,0.9) {$f((3,0),(3,1),(3,2))$};
\node at (0,0) {$f((3,1),(3,2),(3,0))$};
\node at (2.0,0.9) {$=a$};
\node at (2.0,0) {$=a$};
\path[<-,>=stealth'] 
    (-0.9,0.25) edge (-0.9,0.6)
    (0.15,0.25) edge (0.15,0.6)
    (1.2,0.25) edge (1.2,0.6)
    (2.17,0.25) edge (2.17,0.6)
    ;
\end{tikzpicture}\]
and  $(a,a)$ is a loop, a contradiction.
\end{example}

\begin{lemma}\label{lem:pcSatisfyPclc}
Let $p,q$ be primes. Then $\Pol(\C_{q})\models\Sigma_p$ if and only if $p\neq q$.
\end{lemma}
\begin{proof}
If $p\neq q$, then there is an $n\in \N^+\!$ such that $p\cdot n\equiv_q1$. The map
\[f(x_1,\dots,x_p) = n\cdot(x_1 +\ldots + x_p) \pmod q\]
is a polymorphism of $\C_q$ satisfying $\Sigma_p$.

Assume that $f$ is a polymorphism of $\C_q$ satisfying $\Sigma_p$, then
\[f(0,\dots,q-2,q-1) = a = f(1,\dots,q-1,0)\]
and $(a,a)$ is a loop, a contradiction.
\end{proof}
From Corollary~\ref{cor:wond} it is easy to see that the digraphs $\C_2, \C_3, \C_5,\dots$ represent an infinite antichain in $\PCPoset$ and that the conditions $\Sigma_2,\Sigma_3,\Sigma_5,\dots$ represent an infinite antichain in $\mPCL$. 
In order to describe these two posets we generalize Lemma~\ref{lem:pcSatisfyPclc} to disjoint unions of prime cycles in Lemma~\ref{lem:duofpcSatisfyPclc}. 
First, we present an example that (hopefully) helps to better understand the polymorphisms of disjoint unions of cycles. 

\begin{example}
Let us examine the binary polymorphisms of $\C_{2,3}$ (see Figure~\ref{fig:23poly}). 
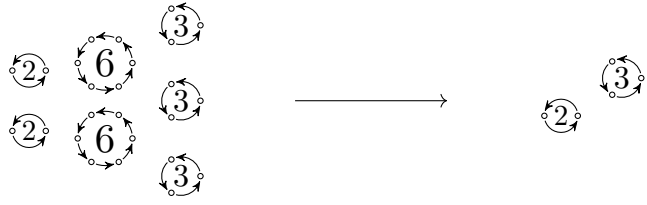
\begin{figure}
    \centering
    \begin{tikzpicture}
    \cycle{2}{(0,-0.4)}
    \cycle{2}{(0,0.4)}
    
    \cycle{3}{(2,0)}
    \cycle{3}{(2,1)}
    \cycle{3}{(2,-1)}
    
    \cycle{6}{(1,-0.5)}
    \cycle{6}{(1,0.5)}
    
    \path (3.5,0) edge[->] (5.5,0);
    
    \cycle{2}{(7,-0.2)}
    \cycle{3}{(7.8,0.3)}
    \end{tikzpicture}
    \caption{Just a nice image to keep you motivated to understand the set $\Hom(\C_{2,3}^2,\C_{2,3})$, i.e., the binary polymorphisms of $\C_{2,3}$.} 
    \label{fig:23poly}
\end{figure} 
Every element in $\Hom(\C_{2,3}^2,\C_{2,3})$ is build from homomorphisms from the connected components of $\C_{2,3}^2$  to $\C_{2,3}$. 
Hence, \[\left|\Hom(\C_{2,3}^2,\C_{2,3})\right| = 2^2\cdot  3^3\cdot 5^2.\]
More generally, let $\C$ be a finite disjoint union of cycles and $n\in \N^+$\!. Let $\textbf G$ be the subgroup of the symmetric group (even automorphism group) on $\C^n$ generated by $+1\colon\boldsymbol t\mapsto(\boldsymbol t+1)$.   
For every  orbit  of $\C^n$ under $\textbf G$ pick  a representative and denote the set of representatives by $T$. Note that every connected component of $\C^n$ contains exactly one element of $T$. Let $f\colon T\to\C$ be such that for every $\boldsymbol t=((a_1,k_1),\dots,(a_n,k_n))\in T$ we have that $f(\boldsymbol t)$ lies in a cycle whose length divides $\ell_{\boldsymbol t}\coloneqq\lcm(a_1,\dots,a_n)$. Since every $\boldsymbol t\in T$ lies in a cycle of length $\ell_{\boldsymbol t}$ we have that $f$ can be uniquely extended to a polymorphism $\tilde f$ of $\C$. Furthermore, 
\[\tilde f((+1)^d(\boldsymbol t))=f(\boldsymbol t)+d\]
for all $\boldsymbol t\in T$ and $d\in\N$. Note that all $n$-ary polymorphisms of $\C$ can be constructed this way.
\end{example}


\begin{lemma}\label{lem:duofpcSatisfyPclc}
Let $\mathbb Q$ be a disjoint union of prime cycles and let $p$ be a prime. Then $\Pol(\mathbb Q)\models \Sigma_p$ if and only if $p\notin Q$.
\end{lemma}
\begin{proof}
We assume without loss of generality that $\mathbb Q$ is the disjoint union of cycles associated to $Q$.

($\Rightarrow$) 
Let $f$ be a polymorphism of $\mathbb Q$ satisfying $\Sigma_p$. Assume $p\in Q$. Then 
\[f((p,0),\dots,(p,p-2),(p,p-1)) = a = f((p,1),\dots,(p,p-1),(p,0))\]
and $(a,a)$ is a loop, a contradiction.

($\Leftarrow$)
For this direction we construct a polymorphism $f\colon \mathbb Q^{\C_p} \to \mathbb Q$ of $\mathbb Q$ satisfying $\Sigma_p$. Here $\C_p$ stands for the set of elements of $\C_p$.
Let $\textbf G$ be the subgroup of the symmetric group on $\mathbb Q^{\C_p}$ generated by $\sigma_{\! p}\colon \boldsymbol t \mapsto \shiftTuple{\boldsymbol t }{\sigma_{\! p}}$ and $+1\colon \boldsymbol t\mapsto (\boldsymbol t+1)$. Recall that 
\begin{align*}
    (\shiftTuple{\boldsymbol t} {\sigma_{\! p}})_{(p,k)}=\boldsymbol t_{(p,k+1)},&& +1=\sigma_{ (\mathbb Q^{\C_p})}, && ((+1)(\boldsymbol t))_{(p,k)}=\boldsymbol t_{(p,k)}+1.
\end{align*}
Hence, $\sigma_{\! p}\circ (+1)=(+1)\circ \sigma_{\! p}$, $\textbf G$ is commutative, and every element of $\textbf G$ is of the form $\sigma_{\! p}^c\circ (+1)^d$ for some $c, d \in\N$.
For every orbit of $\mathbb Q^{\C_p}$ under $\textbf G$ pick a representative and denote the set of representatives by $T$. Let $\boldsymbol t\in T$ and $(q,k)=\boldsymbol t_{(p,0)}$.
Note that the orbit of $\boldsymbol t$ is a disjoint union of cycles of a fixed length $n$ and that $q$ divides $n$. Define $f$ on the orbit of $\boldsymbol t$ as 
\[f\left((\sigma_{\! p}^c\circ (+1)^d)(\boldsymbol t)\right)
=f\left(\shiftTuple {(\boldsymbol t+d)}{\sigma_{\! p}^c}\right)
\coloneqq (q,d)\text{ for all $c,d$.}\] 
To show that $f$ is well defined on the orbit of $\boldsymbol t$ it suffices to prove that $\shiftTuple{(\boldsymbol t+d)}{\sigma_{\! p}^c}=\shiftTuple{(\boldsymbol t+\ell)}{\sigma_{\! p}^{m}}$ implies $d\equiv_p\ell $ for all $d,c,\ell,m\in\N$. Without loss of generality we can assume that $m=\ell=0$. 
Observe that $\boldsymbol t=\shiftTuple{(\boldsymbol t+d)}{\sigma_{\! p}^c}$ implies
$\boldsymbol t_{(p,k\cdot c)}=\boldsymbol t_{(p,(k+1)\cdot c)}+d$ for all $k$.
Considering that $p\cdot c\equiv_{p}0$ we have 
\[\boldsymbol t_{(p,0)}=\boldsymbol t_{(p,c)}+d=
\dots=\boldsymbol t_{(p,(p-1)\cdot c)}+(p-1)\cdot d=\boldsymbol t_{(p,0)}+p\cdot d.\] 
Hence, $p\cdot d\equiv_{q}0$. Since $p\notin Q$ we have that $p$ and $q$ are coprime. Therefore, $d\equiv_{q}0$ as desired.

Repeating this for every $\boldsymbol t\in T$ defines $f$ on $\mathbb Q^{\C_p}$.
If $\boldsymbol r\stackrel{1}\toEdge \boldsymbol s$, then 
\[\boldsymbol s=(\boldsymbol r+1)=(+1)(\boldsymbol r)\] 
and $f(\boldsymbol r)\stackrel{1}\toEdge f(\boldsymbol s)$. Hence, $f$ is a polymorphism of $\mathbb Q$. Furthermore, $f=f_{\sigma_{\! p}}$ by definition.
\end{proof}

Understanding whether a loop condition implies another one is helpful for describing $\mPCL$ and $\PCPoset$.
This problem has already been studied before.
A sufficient condition is given in the following result from Ol\v{s}ak~\cite{olsak-loop}. 
\begin{theorem}[Corollary 1 in~\cite{olsak-loop}]\label{thm:loopolsak}
Let $\Sigma$ and  $\Gamma$ be loop conditions. If  $\mathbb{G}_{\Sigma}\to \mathbb{G}_{\Gamma}$, then $\Sigma\Rightarrow\Gamma$.
\end{theorem}

The idea of the proof is to show that if $h\colon \mathbb{G}_{\Sigma}\to \mathbb{G}_{\Gamma}$ and $f\models\Sigma$, then $f_{h'}\models\Gamma$, where $h'\colon E(\G)\to E(\H), (u,v)\mapsto (h(u),h(v))$. 
We can use Theorem~\ref{thm:loopolsak} to easily deduce implications between cyclic loop conditions. However, for convenience, we state some results explicitly.
\begin{lemma}\label{lem:orderOnLoopcond} Let $C\subset\N^+\!$ be  finite and $c,d\in\N^+$\!.
\begin{enumerate}
    \item We have $\Sigma_C\Rightarrow\Sigma_{C\dotdiv c}$. In particular, $\Sigma_{C\dotdiv c}\Rightarrow \Sigma_{C\dotdiv (c\cdot d)}$.
    \item We have $\Sigma_C\Rightarrow\Sigma_{C\cupdot \{d\}}$.
    \item If $d$ is a multiple of an element of $C$, then $\Sigma_{C\cupdot \{d\}}\Leftrightarrow\Sigma_C$.
\end{enumerate}
\end{lemma}

Using these results we can characterize the orders of $\mPCL$ and $\PCPoset$.
\begin{lemma}\label{cor:implicationSinglePClc}
Let $\Sigma_P, \Sigma_Q$ be prime cyclic loop conditions. Then the following are equivalent:
\begin{multicols}{3}
\begin{enumerate}[label=(\arabic*)]
    \item $\Sigma_P\Rightarrow\Sigma_Q$,
    \item $\mathbb Q\leq \mathbb P$,
    \item $\Pol(\mathbb Q)\not\models\Sigma_P$,
    \item $P\subseteq Q$,
    \item $\mathbb P\to\mathbb Q$.
    \item[]
\end{enumerate}
\end{multicols}
\end{lemma}
\begin{proof}
We show (1) $\Rightarrow$ (3) $\Rightarrow$ (4)  $\Rightarrow$ (5)  $\Rightarrow$ (1) and (4)  $\Rightarrow$ (2) $\Rightarrow$ (3).

(1) $\Rightarrow$ (3)
Since $\Pol(\mathbb Q)\not\models\Sigma_Q$ we have, by assumption, $\Pol(\mathbb Q)\not\models\Sigma_P$.

(3) $\Rightarrow$ (4)
We show the contraposition. 
Since $P\nsubseteq Q$ there is a $p\in P\setminus Q$. By Lemma~\ref{lem:duofpcSatisfyPclc} we have $\Pol(\mathbb Q)\models\Sigma_p$. By Lemma~\ref{lem:orderOnLoopcond} we have $\Sigma_p\Rightarrow\Sigma_P$. Hence, $\Pol(\mathbb Q)\models\Sigma_P$.

(4) $\Rightarrow$ (5)
The identity-function is even an embedding from $\mathbb P$ to $\mathbb Q$.

(5) $\Rightarrow$ (1)
Follows from Theorem~\ref{thm:loopolsak}.

(4) $\Rightarrow$ (2)
Let $\C$ be the first pp-power of $\mathbb Q$ given by the formula \[\Phi_E(x,y)\coloneqq x\stackrel{c}{\toEdge}x\wedge x\stackrel{1}\toEdge y\] 
where $c\coloneqq \prod_{p\in P}p$. 
The set of edges of $\C$ is the subset of edges of $\mathbb Q$ containing all edges that lie in a cycle whose length divides $c$.
Note that $q\in Q$ divides $c$ if and only if $q\in P$. Hence, $\C_q\hookrightarrow\C$ if and only if $q\in P$.
Therefore, $\C$ consists of  $\mathbb P$ and possibly some isolated points and is homomorphically equivalent to $\mathbb P$.

(2) $\Rightarrow$ (3)
By Corollary \ref{cor:wond}, every minor identity satisfied by $\mathbb Q$ is also satisfied by $\mathbb P$. Since $\Pol(\mathbb P)\not\models\Sigma_P$ we have that $\Pol(\mathbb Q)\not\models\Sigma_P$.
\end{proof}

As a consequence of the previous corollary 
any element $[\Sigma_P]$ of $\mPCL$ can be represented by exactly one prime cyclic loop condition. Hence, we will from now on identify $[\Sigma_P]\in\mPCL$ with $\Sigma_P$. Analogously, we identify $[\mathbb P]\in\PCPoset$ with $\mathbb P$.
Furthermore, we obtain a simple description of $\mPCL$ and $\PCPoset$.

\begin{corollary}\label{cor:mPclEqPPC}
The following holds:
\[\mPCL\simeq\PCPoset\simeq (\{P\mid P\text{ a finite nonempty set of primes}\},\supseteq).\]
\end{corollary}





\section{Cyclic Loop Conditions}\label{sec:conditions}
The goal of this section is to give a comprehensible description of the poset of cyclic loop conditions and of the poset of sets of cyclic loop conditions ordered by strength. The results from this section will help us to describe the poset of
pp-constructability types of finite disjoint unions of cycles, $\SDPoset$.
We start by giving a description of the implication order on cyclic loop conditions in Theorem~\ref{thm:implicationSingleClc}. Using a compactness argument we extend this description to  sets of cyclic loop conditions  in Corollary~\ref{thm:implicationClc}. 
We show that every cyclic loop condition is equivalent to a set of prime cyclic loop conditions. Finally, in Corollary~\ref{cor:characterizationMCLAndmCL}, we give a comprehensible description of the poset of sets of prime cyclic loop conditions ordered by their strength for clones; in particular, this gives a description for the poset of cyclic loop conditions ordered by their strength for clones.

\subsection{Single Cyclic Loop Conditions}\label{ssc:clc}

\begin{definition}
We introduce the poset $\mCL$ as follows  
\begin{align*}
    \mCL&\coloneqq (\{[\Sigma] \mid \Sigma \text{ a cyclic loop condition}\},\Leftarrow).
\end{align*}
\end{definition}

One way to show that one cyclic loop condition is stronger than another is presented in the following example.

\begin{example}
Let $\mathcal{C}$ be a clone with domain $D$ such that $\mathcal C\models \Sigma_4$. Then, by definition, there is some $f\in\mathcal C$ satisfying 
\[f(a,b,c,d)=f(b,c,d,a)\text{ for all }a,b,c,d\in D.\]
The function $g\colon (a,b)\mapsto f(a,b,a,b)$ satisfies
\[g(a,b)=g(b,a)\text{ for all }a,b\in D.\]
Hence, $\mathcal C\models \Sigma_2$. Therefore, $\Sigma_4\Rightarrow\Sigma_2$.
\end{example}

Alternatively, $\Sigma_4\Rightarrow\Sigma_2$ follows from Theorem \ref{thm:loopolsak} and the fact that $\Cyc4\to\Cyc2$.
However, unlike for prime cyclic loop conditions, not every implication between cyclic loop conditions can be shown using Theorem~\ref{thm:loopolsak} as seen in the following example. 

\begin{example}\label{exa:2implies4}
Let $\mathcal{C}$ be a clone with domain $D$ such that $\mathcal C\models \Sigma_2$. Then there is some $g\in\mathcal C$ satisfying 
\[g(a,b)=g(b,a)\text{ for all }a,b\in D.\]
The function $f\colon (a,b,c,d)\mapsto g(g(a,b),g(c,d))$ satisfies the identity
\begin{align*}
    f(a,b,c,d)&=g(g(a,b),g(c,d))\\&=g(g(a,b),g(d,c))
    \\&= g(g(d,c),g(a,b))=f(d,c,a,b)
\end{align*}

 for all $a,b,c,d\in D$.
 Note that the digraph corresponding to this identity is isomorphic to $\Cyc4$.
 Hence, $\mathcal C\models \Sigma_4$. Therefore, $\Sigma_2\Rightarrow\Sigma_4$. However, $\mathbb{G}_{\Sigma_2}=\Cyc 2 \not\to \Cyc 4 = \mathbb{G}_{\Sigma_4}$.
\end{example}

The next goal is to weaken the condition in Theorem~\ref{thm:loopolsak} such that the converse implication also holds. First we generalize Example~\ref{exa:2implies4} by proving that $\{\Sigma_C,\Sigma_D\} \Rightarrow\Sigma_{C\cdot D}$ in Corollary~\ref{cor:bulletLoopConditions}. This corollary generalizes Proposition 2.2 in~\cite{Barto_modularity} 
from cyclic loop conditions with only one cycle to cyclic loop conditions in general. To prove this corollary we introduce some notation.  
The function $f$ constructed in Example~\ref{exa:2implies4} is the so-called \emph{star product} of $g$ with itself \cite{absorption}.

\begin{definition}
Let $A$ be a set, $f\colon A^n\to A$ and $g\colon A^m\to A$ be maps. The \emph{star product} of $f$ and $g$ is the function $(f\star g)\colon A^{n\cdot m}\to A$ defined by
\[(x_1,\dots,x_{n\cdot m})\mapsto f(g(x_1,\dots,x_m),\dots,g(x_{(n-1)\cdot m+1},\dots,x_{n\cdot m})).\]
For functions $f\colon A^I\to A$ and $g\colon A^J\to A$, where $I$ and $J$ are finite sets, we define the star product $(f\star g)\colon A^{I\mathbin\times J}\to A$ by
\[\boldsymbol t\mapsto f(i\mapsto g(j\mapsto  t_{(i,j)})).\]
\end{definition}
Note that the second definition extends the first one in the following sense: Let $f\colon A^I\to A$ and $g\colon A^J\to A$ be functions, where $I=[n]$ and $J=[m]$. 
Define $\tilde f\colon A^n\to A$ and $\tilde g\colon A^m\to A$ as $\tilde f(t_1,\dots,t_n)\coloneqq f(\boldsymbol t)$ and $\tilde g(s_1,\dots,s_m)\coloneqq g(\boldsymbol s)$ for all $\boldsymbol t\in A^I,\boldsymbol s\in A^J$. For $\sigma\colon I\mathbin\times J\to[n\cdot m]$ with $\sigma(i,j)=(i-1)\cdot m+j$ we have 
\[(f\star g)_\sigma(\boldsymbol t)=(\tilde f\star \tilde g)(t_1,\dots,t_{n\cdot m}).\]

Using the star product we can easily show the following lemma. 
\begin{lemma}\label{lem:CDimpliesCtimesD}
Let $C,D\subseteq\N^+$\! be finite. Then $\{\Sigma_C,\Sigma_D\}\Rightarrow \Sigma_{\C\times \D}$.
\end{lemma}
\begin{proof}
Let $\mathcal{C}$ be a clone such that $\mathcal{C}\models \Sigma_C$ and $\mathcal{C}\models \Sigma_D$. Then there are $f,g\in\mathcal C$ with $f\models \Sigma_C$ and $g\models\Sigma_D$.
Furthermore 
\[(f\mathbin\star g)(\boldsymbol t)= f((a,k)\mapsto g((b,\ell)\mapsto t_{((a,k),(b,\ell))})).\]

Observe that $f\mathbin\star g$ satisfies
\[f\mathbin\star g=f_{\sigma_{\!C}}\mathbin\star g=f_{\sigma_{\!C}}\mathbin\star g_{\sigma_{\!D}}=(f\mathbin\star g)_{\sigma_{\C\mathbin\times\D}}.\]

Hence, $f\mathbin\star g\models \Sigma_{\C\mathbin\times\D}$.  Therefore,  $\{\Sigma_C,\Sigma_D\}\Rightarrow \Sigma_{\C\times \D}$. 
\end{proof}
Note that this lemma in particular implies that $\{\Sigma_n,\Sigma_m\}\Rightarrow \Sigma_{\lcm(n,m)}$ for any $n,m\in\N^+$\!.
The following example points out the difference between Example~\ref{exa:2implies4} and the special case of $C=D=\{2\}$ in Lemma~\ref{lem:CDimpliesCtimesD}.
\begin{example}\label{exa:TimesvsBullet}
Let $g\models\Sigma_2$. Then the digraph corresponding to the identity 
\[g\mathbin\star g=g(g,g)=g_{\sigma_2}(g_{\sigma_2},g_{\sigma_2})=g_{\sigma_2} \mathbin\star g_{\sigma_2}\] 
is $\Cyc2\mathbin\times\Cyc2$, which consists of two cycles of length 2. Hence, $g\mathbin\star g$ satisfies $\Sigma_{\Cyc2\mathbin\times\Cyc2}$, which is equivalent to $\Sigma_2$.
Observe that in order to show $\Sigma_2\Rightarrow\Sigma_4$, in Example~\ref{exa:2implies4}, we used the identity 
\[g(g,g)=g_{\sigma_2}(g,g_{\sigma_2}).\] 
This identity corresponds to a digraph $\mathbb G$ isomorphic to $\Cyc4$. Hence, $g\mathbin\star g$ also satisfies $\Sigma_{\mathbb G}$, which is equivalent to $\Sigma_4$. Note that $\mathbb G$ and $\Cyc2\mathbin\times\Cyc2$ have the same vertices.
\end{example}

We now introduce a new edge relation on $\C\mathbin\times\D$, which in the case of $\C=\D=\Cyc2$ yields the graph $\mathbb G$ from Example~\ref{exa:TimesvsBullet}.

\begin{definition}
Let $C,D\subset\N^+\!$ be finite.
Define $\C\mathbin\bullet\D$ as the finite disjoint union of cycles $(V,E)$ with  $V=\C\times\D$ and  $E=\{(\boldsymbol t,\sigma_{\C\mathbin\bullet\D}(\boldsymbol t))\mid \boldsymbol t\in \C\times\D\}$, where $\sigma_{\C\mathbin\bullet\D}\colon \C\mathbin\bullet\D\to \C\mathbin\bullet\D $ is defined as
\[\sigma_{\C\mathbin\bullet\D}((a,k),(b,\ell))\coloneqq    
\begin{cases}
    ((a,0),(b,\ell+1))&\text{if }k = a-1\\
    ((a,k+1),(b,\ell)) & \text{otherwise.}
\end{cases}\]
Define 
$\C^{\bullet 0}\coloneqq\Cyc1$ and $\C^{\bullet (k+1)}\coloneqq\C \mathbin\bullet \C^{\bullet k}$ for $k\in\N$.
\end{definition}
For example, $\Cyc{2}\mathbin\bullet\Cyc{2}$ is isomorphic to $\Cyc{4}$ and $\C^{\bullet 1}$ is isomorphic to $\C$. 
Observe that the set associated to  $\C\mathbin\bullet\D$ is 
\[C\cdot D\coloneqq \{a\cdot b\mid a\in C, b\in D\},\]
whereas, the set associated to $\C\times\D$ is $\{\lcm(a,b)\mid a\in C,b\in D\}$. Therefore,  $\C\mathbin\bullet\D\to \C\times\D$. Hence, by Theorem~\ref{thm:loopolsak}, the following lemma is stronger than Lemma~\ref{lem:CDimpliesCtimesD}.

\begin{lemma}\label{lem:CbulletDequalsCD}
Let $C,D\subset \N^+\!$ finite. Then $\{\Sigma_C,\Sigma_D\}\Rightarrow \Sigma_{C\cdot D}$.
\end{lemma}
\begin{proof}
Let $\mathcal C$ be a clone such that $\mathcal C\models\Sigma_{C}$ and $\mathcal C\models\Sigma_{D}$. Then there are $f,g\in\mathcal C$ with $f\models \Sigma_{C}$ and $g\models \Sigma_{D}$.
Let $\{a_1,\dots,a_n\}=C$.
For any $i\in[n]$ define the following permutations on $\C\mathbin\bullet\D$
\begin{align*}
\tau\colon ((a,k),(b,\ell))&\mapsto 
((a,k+1),(b,\ell))
\\
\tau_{i}\colon ((a,k),(b,\ell))&\mapsto
\begin{cases}
((a,k),(b,\ell+1)) & \text{if }a= a_i, k=a-1 \\
((a,k),(b,\ell)) & \text{otherwise.}
\end{cases}
\end{align*}

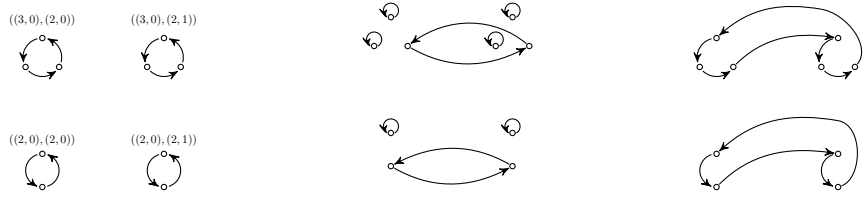
\begin{figure}
    \centering
    \begin{tikzpicture}
        \node at (0+90:0.5) [scale=0.4] {$((3,0),(2,0))$};
        \node at (0,-1.1) [scale=0.4] {$((2,0),(2,0))$};
        \cycleEmpty{3}{(0,0)}{{1,2,3}}
        \cycleEmpty{2}{(0,-1.5)}{{1,2}}
        
        \node at ($(1.6,0)+(0+90:0.5)$) [scale=0.4] {$((3,0),(2,1))$};
        \node at (1.6,-1.1) [scale=0.4] {$((2,0),(2,1))$};
        \cycleEmpty{3}{(1.6,0)}{{1,2,3}}
        \cycleEmpty{2}{(1.6,-1.5)}{{1,2}}
    
    \end{tikzpicture}
    \hspace{1.9cm}
    \begin{tikzpicture}
        \cycleEmpty{3}{(0,0)}{{}}
        \cycleEmpty{2}{(0,-1.5)}{{}}
        
        \cycleEmpty{3}{(1.6,0)}{{}}
        \cycleEmpty{2}{(1.6,-1.5)}{{}}
    \def \radiusA {{3*0.035+0.15}}
    \def \radiusB {{2*0.035+0.15}}
        
        \path[->,>=stealth']
            (-30:\radiusA) edge[bend right,shorten >=0.5mm,shorten <=0.5mm] ($(1.6,0)+(-30:\radiusA)$)
            ($(1.6,0)+(-30:\radiusA)$) edge[bend right=40,shorten >=0.5mm,shorten <=0.5mm] (-30:\radiusA)
            ($(0,-1.5)+(90+180:\radiusB)$) edge[bend right,shorten >=0.5mm,shorten <=0.5mm] ($(1.6,-1.5)+(90+180:\radiusB)$)
            ($(1.6,-1.5)+(90+180:\radiusB)$) edge[bend right,shorten >=0.5mm,shorten <=0.5mm] ($(0,-1.5)+(90+180:\radiusB)$)
            ;
            \draw[>=stealth',arrows=-{>[bend]}] ($(1.6+0.05,0)+(90+120:\radiusA)$)
    arc ({-60}:{240}:1mm);
            \draw[>=stealth',arrows=-{>[bend]}] ($(1.6+0.05,0)+(90:\radiusA)$)
    arc ({-60}:{240}:1mm);
            \draw[>=stealth',arrows=-{>[bend]}] ($(0.05,0)+(90+120:\radiusA)$)
    arc ({-60}:{240}:1mm);
            \draw[>=stealth',arrows=-{>[bend]}] ($(0.05,0)+(90:\radiusA)$)
    arc ({-60}:{240}:1mm);
            \draw[>=stealth',arrows=-{>[bend]}] ($(0.05,-1.5)+(90:\radiusB)$)
    arc ({-60}:{240}:1mm);
            \draw[>=stealth',arrows=-{>[bend]}] ($(1.6+0.05,-1.5)+(90:\radiusB)$)
    arc ({-60}:{240}:1mm);
    \end{tikzpicture}
    \hspace{1.9cm}
    \begin{tikzpicture}

        \cycleEmpty{3}{(0,0)}{{1,2}}
        \cycleEmpty{2}{(0,-1.5)}{{1}}
        
        \cycleEmpty{3}{(1.6,0)}{{1,2}}
        \cycleEmpty{2}{(1.6,-1.5)}{{1}}
        
        \def \radiusA {{3*0.035+0.15}}
        \def \radiusB {{2*0.035+0.15}}

        \path[->,>=stealth']

            (90+240:\radiusA) edge[out=45,in=170,shorten >=0.5mm,shorten <=0.5mm] ($(1.6,0)+(90:\radiusA)$)

            ($(0,-1.5)+(90+180:\radiusB)$) edge[out=45,in=170,shorten >=0.5mm,shorten <=0.5mm] ($(1.6,-1.5)+(90:\radiusB)$)
        ;
        \path 
        (90:\radiusA) edge[<-,>=stealth',out=45,in=170,shorten <=0.5mm] ($(1.6,0)+(90:\radiusA)+(90:\radiusA)-(90+240:\radiusA)$)
        
        ($(1.6,0)+(90:\radiusA)+(90:\radiusA)-(90+240:\radiusA)$) edge[out=-10,in=50,shorten >=0.5mm] ($(1.6,0)+(90+240:\radiusA)$)

        ($(0,-1.5)+(90:\radiusB)$) edge[<-,>=stealth',out=45,in=170,shorten <=0.5mm] ($(1.6,-1.5)+(90:\radiusB)+(90:\radiusB)-(90+180:\radiusB)$)
        
        ($(1.6,-1.5)+(90:\radiusB)+(90:\radiusB)-(90+180:\radiusB)$) edge[out=-10,in=20,shorten >=0.5mm] ($(1.6,-1.5)+(90+180:\radiusB)$)
        
        ;
    
    \end{tikzpicture}
    \caption{Graphs of the permutations $\tau$ (left), $\tau_{1}\circ\tau_{2}$ (middle), and $\tau\circ\tau_{1}\circ\tau_{2}$ (right) on $\Cyc{2,3}\mathbin\bullet\Cyc{2}$.
}
    \label{fig:CBulletD}
\end{figure}
Observe that, since $f=f_{\sigma_{\!C}}$ and $g=g_{\sigma_{\!D}}$, we have \[f\star g=(f\star g)_\tau=(f\star g)_{\tau_i}\text{ for all $i\in [n]$.}\] 
Now we show that
\[\tau\circ\tau_{1}\circ\cdots\circ\tau_{n}
=\sigma_{\C\mathbin\bullet \D}.\]
First, the reader can verify that this equality holds for the example in Figure~\ref{fig:CBulletD}.
To prove that the equality holds in general let $(a,k)\in\C$ and $(b,\ell)\in\D$.
We have that
\begin{align*}
    (\tau\circ\tau_{1}\circ\cdots\circ\tau_{n})((a,k),(b,\ell))
    &=
    \begin{cases}
    \tau((a,k),(b,\ell+1))&\text{if }k=a-1\\
    \tau((a,k),(b,\ell)) & \text{otherwise}
    \end{cases}\\
    &=\begin{cases}
    ((a,0),(b,\ell+1))&\text{if }k=a-1\\
    ((a,k+1),(b,\ell)) & \text{otherwise}
    \end{cases}\\
    &=\sigma_{\C\mathbin{\bullet}\D}((a,k),(b,\ell)).
\end{align*}
Hence, $(f\star g)=(f\star g)_{\sigma_{\C\mathbin\bullet\D}}$ and since $C\cdot D$ is associated to $\C\mathbin\bullet\D$ there is, by Theorem~\ref{thm:loopolsak}, also an element in $\mathcal C$ that satisfies $\Sigma_{C\cdot D}$.
\end{proof}

Observe that $\C\mathbin\bullet\D\to\C$ and $\C\mathbin\bullet\D\to\D$. Therefore, by Theorem~\ref{thm:loopolsak}, we have $\Sigma_{C\cdot D}\Rightarrow\Sigma_C$ and $\Sigma_{C\cdot D}\Rightarrow\Sigma_D$. Hence, the implication in Lemma~\ref{lem:CbulletDequalsCD} is actually an equivalence. As a consequence we obtain the following corollary. 

\begin{corollary}\label{cor:bulletLoopConditions}
Let $C,C_1,\dots,C_n\subset \N^+\!$ be finite. Then 
\[\{\Sigma_{C_1},\dots,\Sigma_{C_n}\}\Leftrightarrow\Sigma_{C_1\cdot\ldots\cdot C_n}\] 
holds. In particular, $\Sigma_C\Leftrightarrow\Sigma_{C^k}$ for every $k\in\N^+\!$.
\end{corollary}

We now weaken the condition in Theorem~\ref{thm:loopolsak} such that it gives a full characterization of the implication order on $\mCL$. 
\begin{theorem}\label{thm:implicationSingleClc}

Let $C,D\subset\N^+\!$ finite. Then the following are equivalent:
\begin{enumerate}[label=(\arabic*)]
    \item We have $\Sigma_C\Rightarrow\Sigma_D$.
    \item For all $c\in \N^+\!$ we have $\Pol(\D\dotdiv c)\models\Sigma_C$ implies $\Pol(\D\dotdiv c)\models\Sigma_D$.
    \item For every $a\in C$ there exist $b\in D$ and $k\in\N$ such that $b$ divides $a^k$.
    \item There exists $k\in\N^+\!$ such that $\C^{\bullet k}\to \D$.
\end{enumerate}
\end{theorem}

In order to show (2) $\Rightarrow$ (3) we need to know more about the connection between disjoint unions of cycles and cyclic loop conditions. This connection is investigated later in Section~\ref{sec:structures}. As this section is devoted exclusively to cyclic loop conditions we will provide the proof of this implication later. However, it should be noted that the result required from section~\ref{sec:structures} is not based on the current section, so there is no circular reasoning.

\begin{proof}
The direction (1) $\Rightarrow$ (2) is clear from the definition of the implication order. 

(2) $\Rightarrow$ (3) Will follow from Lemma~\ref{lem:missingImplicationSingleClc}.

(3) $\Rightarrow$ (4) For $a\in C$ let $k_a\in\N$ be such that there is a $b\in D$ with $b$ divides  $a^{k_a}$. Let $n=|C|$ and $k=\max\{k_a\mid a\in C\}$. Note that any number in $C^{n\cdot k}$ must be a multiple of $a^k$ for some $a\in C$. Hence, $\C^{\bullet (n\cdot k)}\to \D$.

(4) $\Rightarrow$ (1) 
Let $k\in\N^+\!$ be such that $\C^{\bullet k}\to \D$. Then
\[\Sigma_C\stackrel{\ast}\Leftrightarrow\Sigma_{C^k}\stackrel{\ast\ast}\Rightarrow\Sigma_{D},\]
where $\ast$ and $\ast\ast$ hold by Corollary~\ref{cor:bulletLoopConditions} and Theorem~\ref{thm:loopolsak}, respectively.
\end{proof}

\begin{remark} Note that it follows from (3) that the problem
\begin{align*}
    &\text{Input: two finite sets $C,D\subset\N^+$}\\
    &\text{Output: Does $\Sigma_C\Rightarrow\Sigma_D$ hold?}
\end{align*}
is decidable; in fact, this can be decided in polynomial time even if the integers in C and D are given in binary.
\end{remark}
Define the  function $\Flat$ from $\N^+$ to $\N^+$ as
\begin{align*}
    p_1^{\alpha_1}\cdot\ldots\cdot p_n^{\alpha_n}&\mapsto p_1\cdot\ldots\cdot p_n.
\end{align*}
We extend $\Flat$ to subsets of $\N^+$ in the obvious way. Note that (3) in Theorem~\ref{thm:implicationSingleClc} is equivalent to
\[ \text{(3') For every $a\in C$ there exist $b\in D$ such that $\Flat(b)$ divides $\Flat(a)$.}\]
Hence, every cyclic loop condition is equivalent to a cyclic loop condition where only square-free numbers occur. 
\begin{corollary}\label{cor:flat}
For every finite $C\subset\N^+\!$ we have that $\Sigma_C\Leftrightarrow\Sigma_{\Flat(C)}$. 
In particular, we have $\Sigma_{\C\times\D}\Leftrightarrow\Sigma_{\C\mathbin\bullet\D}$.
\end{corollary}

\subsection{Sets of Cyclic Loop Conditions}\label{ssc:setsOfClc}
The next goal is to understand the implication order on sets of cyclic loop conditions. 
\begin{definition}
We introduce the poset $\MCL$ as follows  

\begin{align*}
    \MCL&\coloneqq (\{[\Sigma]\mid \Sigma \text{ a set of cyclic loop conditions}\},\Leftarrow).
\end{align*}
\end{definition}

We already understand this order on finite sets of cyclic loop conditions: By Corollary \ref{cor:bulletLoopConditions}, every finite set of cyclic loop conditions is equivalent to a single cyclic loop condition and we know how to compare  single cyclic loop conditions by Theorem \ref{thm:implicationSingleClc}.
Using the compactness theorem for first-order logic we will show that in order to determine the order on infinite sets it suffices to consider their finite subsets. 
\begin{theorem}\label{thm:compactness}
Let $\Gamma$, $\Sigma$ be sets of minor identities, where $\Sigma$ is finite. Then $\Gamma\Rightarrow\Sigma$ if and only if there is a finite $\Gamma'\subseteq \Gamma$ such that $\Gamma'\Rightarrow\Sigma$.
\end{theorem}
\begin{proof}
Let $\Sigma=\{(f_1)_{\sigma_1}\approx (g_1)_{\tau_1},\dots,(f_k)_{\sigma_k}\approx (g_k)_{\tau_k}\}$ and let $\kappa$ be the set of function symbols occurring in $\Gamma$. We construct a first-order theory $T_{\Gamma,\Sigma}$ which is satisfiable if and only if $\Gamma\not\Rightarrow\Sigma$. 
The goal of the construction is that every model of $T_{\Gamma,\Sigma}$ encodes a clone which witnesses $\Gamma\not\Rightarrow\Sigma$ and 
conversely every clone which witnesses $\Gamma\not\Rightarrow\Sigma$ is a model of $T_{\Gamma,\Sigma}$.

The signature of $T_{\Gamma,\Sigma}$ is  $\kappa\cup\{D,F_n,E_n\mid n\in\N^+\}$, where
\begin{itemize}
    \item $\kappa$ is a set of constant symbols (intended to denote the operations satisfying $\Gamma$),
    \item $D$ is a unary relation symbol (intended to denote the domain),
    \item $F_n$ is a unary relation symbol (intended to denote the $n$-ary operations) for every $n\in\N^+\!$,
    \item $E_n$ is an $(n+2)$-ary relation symbol (intended to denote the graphs of all $n$-ary operations) for every $n\in\N^+\!$.
\end{itemize}

For the sake of readability we set the following abbreviations:
\begin{itemize}
    \item $\forall x\in R: \Phi(x)$ abbreviates $\forall x( R(x)\Rightarrow \Phi(x))$, for every unary relation symbol $R$
    \item $\exists x\in R: \Phi(x)$ abbreviates $\exists x( R(x)\wedge \Phi(x))$, for every unary relation symbol $R$
    \item $f(\boldsymbol x)\approx x_i$ abbreviates $E_n(f,x_1,\dots,x_n,x_i)$.
\end{itemize}
   
The theory $T_{\Gamma,\Sigma}$ consists of the following sentences
\begin{itemize}
    \item $\exists x (D(x))$,
    \item $\forall f\in F_n\ \forall x_1,\dots,x_n\in D\ \exists! y\in D :E_n(f,x_1,\dots,x_n,y)$ for every
    \item $\exists f\in F_n\ \forall \boldsymbol x\in D^n:f(\boldsymbol x)\approx x_i$ for every $n\in\N^+$ and $i\in[n]$,
    \item $\forall f\in F_n\ \forall g_1,\dots,g_n\in F_m\ \exists h\in F_m$
    \[\forall \boldsymbol x\in D^m :f(g_1(\boldsymbol x),\dots,g_n(\boldsymbol x))\approx h(\boldsymbol x)\] for every $n,m\in\N^+\!$,
    \item $F_n(f)$ for every $n$-ary $f\in \kappa$, $n\in\N^+\!$,
    \item $\forall x_1,\dots, x_r \in D: f(x_{\sigma(1)},\dots,x_{\sigma(n)})\approx g(x_{\tau(1)},\dots,x_{\tau(m)})$ for all identities $f_\sigma\approx g_\tau \in \Gamma$,
    \item $\forall f_1\in F_{n_1}\ \forall g_1\in F_{m_1}\dots \forall f_k\in F_{n_k}\ \forall g_k\in F_{m_k}:$
    \[\bigvee_{i=1}^k \exists x_1,\dots, x_{r_i} \in D:f_i(x_{\sigma_i(1)},\dots,x_{\sigma_i(n_i)})\not\approx g_i(x_{\tau_i(1)},\dots,x_{\tau_i(m_i)})\] 
    assume w.l.o.g. that $f_1,g_1,\dots,f_k,g_k\notin\kappa$.
\end{itemize}
Now we show that $T_{\Gamma,\Sigma}$ is unsatisfiable if and only if $\Gamma\Rightarrow\Sigma$.
If $\Gamma\not\Rightarrow\Sigma$, then there is a clone $\mathcal C$ over some domain $C$ such that $\mathcal C\models\Gamma$, witnessed by a function $\tilde\cdot\colon \tau\to\mathcal C$, and $\mathcal C\not\models \Sigma$. Let $\A$ be the structure with domain $C\cup \mathcal{C}$ and 
\begin{align*}
    D^{\A}&\coloneqq C,\ 
F_n^{\A}\coloneqq\{f\in\mathcal C\mid f\text{ is $n$-ary}\},\ f^{\A}\coloneqq\tilde f,\text{ and }\\
E_n^{\A}&\coloneqq\{(f,x_1,\dots,x_n,y)\mid f\in F_n^{\A}, 
x_1,\dots,x_n,y\in D, f(x_1,\dots,x_n)=y\}
\end{align*}for every $n\in\N^+\!$ and every $f\in\kappa$. 
By construction $\A\models T_{\Gamma,\Sigma}$. Hence, $T_{\Gamma,\Sigma}$ is satisfiable.

If $T_{\Gamma,\Sigma}$ is satisfiable, then it has some model $\A$. Let $C=D^{\A}$. For every $n\in\N^+\!$ and $f\in F_n^{\A}$ let $\bar f\colon C^n\to C$ be the operation with graph \[\{(c_1,\dots,c_n,d)\mid (f,c_1,\dots,c_n,d)\in E_n^{\A}\}.\]
Define $\mathcal{C}\coloneqq\{\bar f\mid n\in\N^+, f\in F_n^{\A}\}$. By construction $\mathcal{C}$ is a clone over domain $C$ that satisfies $\Gamma$, witnessed by the assignment $f\mapsto \bar f^{\A}$, and does not satisfy $\Sigma$. Hence, $\Gamma\not\Rightarrow\Sigma$.

If $\Gamma\Rightarrow\Sigma$ we have that theory $T_{\Gamma,\Sigma}$ is unsatisfiable. Hence, by compactness, there is a finite $T'\subseteq T$ that is also unsatisfiable. Since $T'$ is finite there must be a finite $\Gamma'\subseteq\Gamma$  such that $T'\subseteq T_{\Gamma',\Sigma}$. Hence, $T_{\Gamma',\Sigma}$ is unsatisfiable and $\Gamma'\Rightarrow\Sigma$.
\end{proof}

Note that the restriction to minor identities is not essential. The proof of Theorem~\ref{thm:compactness} can easily be adapted to the case where $\Gamma$ and $\Sigma$ are sets of identities (not necessarily minor). 
From Theorem~\ref{thm:compactness},  Corollary~\ref{cor:bulletLoopConditions}, and Theorem~\ref{thm:implicationSingleClc} we obtain the following corollary.

\begin{corollary}\label{thm:implicationClc}
Let $\Gamma$ be a set of cyclic loop conditions and $D\subset\N^+\!$ be finite. Then the following are equivalent:
\begin{enumerate}[label=(\arabic*)]
    \item $\Gamma\Rightarrow\Sigma_D$.
    
    \item There is a finite $\Gamma'\subseteq\Gamma$ such that $\Gamma'\Rightarrow \Sigma_D$.
    \item There is a finite $\Gamma'\subseteq\Gamma$ such that for all $c\in \N^+\!$ we have $\Pol(\D\dotdiv c)\models\Gamma'$ implies $\Pol(\D\dotdiv c)\models\Sigma_D$.
    \item There are cyclic loop conditions $\Sigma_{C_1},\dots,\Sigma_{C_n}\in\Gamma$ such that for every $a_1\in C_1,\dots,a_n\in C_n$ there exist $b\in D$ and $k_1,\dots,k_n\in\N^+\!$ such that $b$ divides $a_1^{k_1}\cdot\ldots\cdot a_n^{k_n}$.    
    \item There are cyclic loop conditions $\Sigma_{C_1},\dots,\Sigma_{C_n}\in\Gamma$ and $k_1,\dots,k_n\in\N^+\!$ with $\C_1^{\bullet k_1}\mathbin\bullet\cdots\mathbin\bullet\C_n^{\bullet k_n} \to \D$.
\end{enumerate}
\end{corollary}

Theorem~\ref{thm:implicationSingleClc} and~\ref{thm:implicationClc} provide a simple characterization of the implication order on $\mCL$ and $\MCL$, respectively. However, we did not yet achieve the goal of obtaining a comprehensible description of these posets.

\subsection{Irreducible Cyclic Loop Conditions}\label{ssc:irredClc}

We take a second look at Corollary~\ref{cor:bulletLoopConditions}. It states in particular that the cyclic loop condition $\Sigma_{C_1\cdot\ldots\cdot C_n}$ is equivalent to the set $\{\Sigma_{C_1},\dots,\Sigma_{C_n}\}$. This suggests that one can replace a single cyclic loop condition by a set of possibly simpler cyclic loop conditions.

\begin{definition}
Let $\Sigma$ be a cyclic loop condition. A \emph{decomposition} of $\Sigma$ is a set of pairwise incomparable cyclic loop conditions $\Gamma$ such that $\Sigma$ and $\Gamma$ are equivalent. The condition $\Sigma$ is called \emph{irreducible} if every decomposition of $\Sigma$ contains exactly one element.
\end{definition}

\begin{example}\label{exa:decomp}
Some examples of decompositions.
\begin{multicols}{2}
\begin{itemize}
    \item $\Sigma_6 \Leftrightarrow\{\Sigma_2, \Sigma_3\}$
    \item $\Sigma_{30} \Leftrightarrow \{\Sigma_2, \Sigma_{15}\} \Leftrightarrow {}$\hbox to 1cm {$\{\Sigma_2, \Sigma_3, \Sigma_5\}$}
\item $\Sigma_{6,20} \Leftrightarrow \{ \Sigma_2, \Sigma_{3,5}\}$
\item     $\Sigma_{2,15} \Leftrightarrow \{\Sigma_{2,3},  \Sigma_{2,5}\}$
\end{itemize}
\end{multicols}
\begin{itemize}
    \item $\Sigma_{14,15}\Leftrightarrow\{\Sigma_{2,3},\Sigma_{2,5},\Sigma_{3,7},\Sigma_{5,7}\}$
\end{itemize}
Note that all decompositions except for $\{\Sigma_2, \Sigma_{15}\}$ contain only irreducible cyclic loop conditions.
\end{example}

Observe that all irreducible  cyclic loop conditions from the example are prime cyclic loop conditions.

\begin{lemma}\label{lem:pclcAreImpliedBySingleCond}
Let $\Sigma_P$ be a prime cyclic loop condition and let $\Gamma$ be a set of cyclic loop conditions such that $\Gamma\Rightarrow\Sigma_P$. Then there is a $\Sigma\in\Gamma$ with $\Sigma\Rightarrow\Sigma_P$.
\end{lemma}
\begin{proof}
By Corollary~\ref{cor:flat}, we may  assume all numbers occurring in conditions in $\Gamma$ are square-free. 

Since $\Gamma\Rightarrow\Sigma_P$ we have, by Corollary~\ref{thm:implicationClc}, that there are there are $\Sigma_{C_1},\dots,\Sigma_{C_n}\in\Gamma$  with $\C_1\mathbin\bullet\cdots\mathbin\bullet\C_n \to \mathbb P$.
Assume that $\C_i\not\to\mathbb P$ for any $i$. Then for every $i$ there is an $a_i\in\C_i$ such that no $p\in P$ divides $a_i$. Hence, there is a cycle of length $a_1\cdot\ldots\cdot a_n$ in $\C_1\mathbin\bullet\cdots\mathbin\bullet\C_n$ and no $p\in P$ divides $a_1\cdot\ldots\cdot a_n$, a contradiction.
Therefore, there must be an $i$ such that $\C_i\to\mathbb P$. Hence, by Theorem~\ref{thm:loopolsak}, we have that $\Sigma_{C_i}\Rightarrow\Sigma_P$.
\end{proof}

\begin{corollary}\label{lem:pclcAreIrred}
Every prime cyclic loop condition is irreducible.
\end{corollary}

The next step is to show that every cyclic loop condition can be decomposed into a set of prime cyclic loop conditions. The next example gives the idea how to find this decomposition.

\begin{examplex}\label{exa:620}
Let us consider $\Sigma_{6,20}$ from Example \ref{exa:decomp} and let $\Gamma$ be a decomposition of $\Sigma_{6,20}$. Then every condition in $\Gamma$ is implied by $\Sigma_{6,20}$. 
By Lemma~\ref{lem:orderOnLoopcond}, we know that $\Sigma_{6,20}$ implies $\Sigma_{\{6,20\}\dotdiv c}$ for every $c\in\N^+\!$. All nontrivial conditions of this form are presented in Figure~\ref{fig:exampleDotdivPrimeConditions}. 
\begin{figure}
\begin{minipage}{0.48\textwidth}
     \centering
    \begin{tikzpicture}
    \node (01) at (0,-2) {$10$};
    
    \node (10) at (-1.5,-1) {$3$};
    \node (11) at (0,-1) {$5$};
    \node (12) at (1.5,-1) {$2$};
    
    \node (20) at (0,0) {$1$};
    \node (02) at (1.5,-2) {$4$};
    \node (00) at (-1.5,-2) {$15$};
    \path 
        (20) edge (10)
        (20) edge (11)
        (20) edge (12)
        
        (10) edge (00) 
        (11) edge (00) 
        (11) edge (01)
        (12) edge (01)
        (12) edge (02)
        ;
    \end{tikzpicture}
   \end{minipage}\hfill
   \begin{minipage}{0.48\textwidth}
     \centering
    \begin{tikzpicture}
    \node (01) at (0,-2) {$\boldsymbol{\Sigma_{3,2}}$};
    
    \node (10) at (-1.5,-1) {$\Sigma_{2,20}$};
    \node (11) at (0,-1) {$\Sigma_{6,4}$};
    \node (12) at (1.5,-1) {$\Sigma_{3,10}$};
    
    \node (20) at (0,0) {$\Sigma_{6,20}$};
    \node (02) at (1.5,-2) {$\boldsymbol{\Sigma_{3,5}}$};
    \node (00) at (-1.5,-2) {$\boldsymbol{\Sigma_{2,4}}$};
    
    \node at (-0.8,-1) {$\Leftrightarrow$};
    \node[rotate=90] at (-1.5,-1.5) {$\Leftrightarrow$};
    \node[rotate=45] at (-0.8,-1.5) {$\Leftrightarrow$};
    \path 
        (20) edge (10)
        (20) edge (11)
        (20) edge (12)
        (11) edge (01)
        (12) edge (01)
        (12) edge (02)
        ;
    \end{tikzpicture}
   \end{minipage}
   \caption{Some numbers ordered by divisibility (left) and corresponding cyclic loop conditions of the form $\Sigma_{\{6,20\}\dotdiv c}$ ordered by strength
   (right).}\label{fig:exampleDotdivPrimeConditions}
\end{figure}
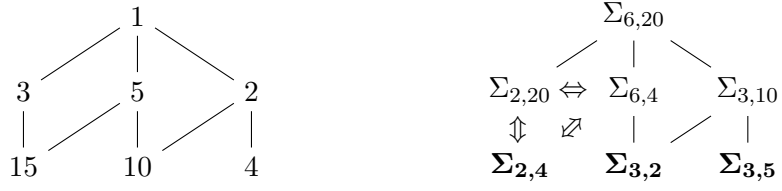
We would like to emphasize the following two observations: 
\begin{enumerate}
    \item the conditions  $\Sigma_{2,4},\Sigma_{3,2}, \Sigma_{3,5}$ 
    printed in bold at the bottom of the figure
    are equivalent to prime cyclic loop conditions and
    \item $\Sigma_{6,20}\Leftrightarrow\{\Sigma_{2,4},\Sigma_{3,2}, \Sigma_{3,5}\}\Leftrightarrow \{ \Sigma_2,\Sigma_{3,2}, \Sigma_{3,5}\}\Leftrightarrow \{ \Sigma_2, \Sigma_{3,5}\}$. \hfill$\triangle$
\end{enumerate} 
\end{examplex}

\begin{definition}\label{def:minimal}

Let $C\subseteq\N^+\!$ be a finite set. A number $c\in\N^+\!$ is \emph{maximal for} $C$ if $1\notin(C\dotdiv c)$ and $1\in{C\dotdiv (c\cdot d)}$ for all $d>1$ dividing $\lcm(C\dotdiv c)$.
\end{definition}
The numbers 15, 10, and 4 are maximal for $\{6,20\}$.
We now show that the two observations in Example~\ref{exa:620} hold in general.

\begin{lemma}\label{lem:minimalIsPrime}
Let $C\subseteq\N^+\!$ be a finite set and let  $c\in\N^+\!$ be maximal for $C$. 
Then $\Sigma_{C\dotdiv c}$ is equivalent to a prime cyclic loop condition.

\end{lemma}
\begin{proof}

We show that $\Sigma_{C\dotdiv c}$ is equivalent to a prime cyclic loop condition by showing that every $a\in (C\dotdiv c)$ is a multiple of some prime in $C\dotdiv c$.
Let $a$ be in $ (C\dotdiv c)$ and $p$ be a prime divisor of $a$ (which exists since $a\neq1$). 
Since $c$ is maximal for $C$ we have $1\notin(C\dotdiv c)$ and $1\in (C\dotdiv (c\cdot p))$. Hence, $p\in (C\dotdiv c)$ and  $a$ is a multiple of a prime in $ (C\dotdiv c)$, as desired.
\end{proof}

\begin{theorem}\label{thm:clcAreSetsOfPclc}
Every cyclic loop condition is equivalent to a finite set of prime cyclic loop conditions.
\end{theorem}
\begin{proof}
Let $C\subset\N^+\!$ be finite and consider the cyclic loop condition $\Sigma_C$. Without loss of generality assume that $C$ is square-free. 
Let $c_1,\dots, c_n\in\N^+\!$ be such that 
\[\{\Sigma_{C\dotdiv c_1},\dots,\Sigma_{C\dotdiv c_n}\}=\{\Sigma_{C\dotdiv c}\mid c\text{ is maximal for }C\}.\] 
By Lemma~\ref{lem:orderOnLoopcond}, we have that  $\Sigma_C\Rightarrow\{\Sigma_{C\dotdiv c_1},\dots,\Sigma_{C\dotdiv c_n}\}$.
From Lemma~\ref{lem:minimalIsPrime} we know that for every $i$ the condition  $\Sigma_{C\dotdiv c_i}$ is equivalent to a prime cyclic loop condition. 
Assume  for contradiction that 
\[(\C\dotdiv c_1)\mathbin\bullet\dots\mathbin\bullet(\C\dotdiv c_n)\not\to\C.\]
Hence, there are $a_1,\dots,a_n$ with $a_i\in C\dotdiv c_i$ such that no $a\in C$ divides $c := a_1\cdot\ldots\cdot a_n$.
Therefore, $1\notin C\dotdiv c$ and the cyclic loop condition $\Sigma_{C\dotdiv c}$ is not trivial.
Let $d\in \N^+\!$ be such that $c\cdot d$ is maximal for $C$. 
Then there is an $i$ such that $\Sigma_{C\dotdiv c\cdot d}= \Sigma_{C\dotdiv c_i}$. 
Note that the number $a_i$ cannot be in $C\dotdiv c\cdot d$,  because $C$ is square-free. Therefore, $C\dotdiv c\cdot d\neq C\dotdiv c_i$ for all $i$. A contradiction.
By Corollary~\ref{thm:implicationClc} it follows that 
$\{\Sigma_{C\dotdiv c_1},\dots,\Sigma_{C\dotdiv c_n}\}\Rightarrow\Sigma_C$, 
as desired.
\end{proof}

Note that, as seen in Example~\ref{exa:620}, the set \[\Gamma\coloneqq\{\Sigma_{C\dotdiv c}\mid c\text{ is maximal for }C\}\] is not necessarily a decomposition of $\Sigma_C$. However, the set consisting of the strongest conditions in $\Gamma$ is a decomposition of $\Sigma_C$. This observation together with Theorem~\ref{thm:clcAreSetsOfPclc} and Corollary~\ref{lem:pclcAreIrred} yields the following corollary.

\begin{corollary}
Let $\Sigma$ be a cyclic loop condition. Then
\begin{enumerate}
    \item $\Sigma$ is irreducible iff $\Sigma$ is equivalent to a prime cyclic loop condition and
    \item there is a  finite decomposition of $\Sigma$ which consists of prime cyclic loop conditions.
\end{enumerate}
\end{corollary}

The next step towards finding a comprehensible description of $\mCL$ and $\MCL$ is to understand 
the poset of sets of prime cyclic loop conditions ordered by strength.

\begin{definition}
We introduce the poset $\MPCL$ as follows
\begin{align*}
\MPCL&\coloneqq (\{[\Sigma]\mid \Sigma \text{ a set of prime cyclic loop conditions}\},\Leftarrow).
\end{align*}
\end{definition} 

From Theorem~\ref{thm:clcAreSetsOfPclc} we obtain the following corollary.
\begin{corollary}
We have that $\MCL=\MPCL$.
\end{corollary}
 
Combining Lemma~\ref{lem:pclcAreImpliedBySingleCond} and Corollary~\ref{cor:implicationSinglePClc} we obtain the following result.
\begin{corollary}\label{cor:implicationPclc}
Let $\Gamma$ be a set of prime cyclic loop conditions and $Q$ a finite nonempty set of primes. Then the following are equivalent:
\begin{enumerate}[label=(\arabic*)]
    \item We have $\Gamma\Rightarrow\Sigma_Q$.
    
    \item There is a $\Sigma_P\in\Gamma$ such that $\Sigma_P\Rightarrow \Sigma_Q$.
    \item There is a $\Sigma_P\in\Gamma$ such that  $\Pol(\mathbb Q)\not\models\Sigma_P$.
    \item There is a $\Sigma_{P}\in\Gamma$  with $P\subseteq Q$.
    \item There is a $\Sigma_{P}\in\Gamma$  with $\mathbb P \to \mathbb Q$.
    
\end{enumerate}
\end{corollary}

Let $\Gamma$ be a set of cyclic loop conditions and $\Gamma'\coloneqq\{\Sigma_P\in\mPCL\mid \Gamma\Rightarrow\Sigma_P \}$. Then $\Gamma\Leftrightarrow\Gamma'$ and $\Gamma'$ is a downset of $\mPCL$. Corollary~\ref{cor:implicationPclc} implies that if $\Gamma$ is a downset  of $\mPCL$, then $\Gamma=\Gamma'$. Hence, the elements of $\MPCL$ correspond to downsets of $\mPCL$. 

\begin{corollary}\label{cor:characterizationMCLAndmCL}
We have that 
$\MCL=\MPCL\simeq (\operatorname{Downsets Of}(\mPCL),\subseteq)$ and $\mCL\simeq(\operatorname{Finitely Generated Downsets Of}(\mPCL),\subseteq)$.
\end{corollary}

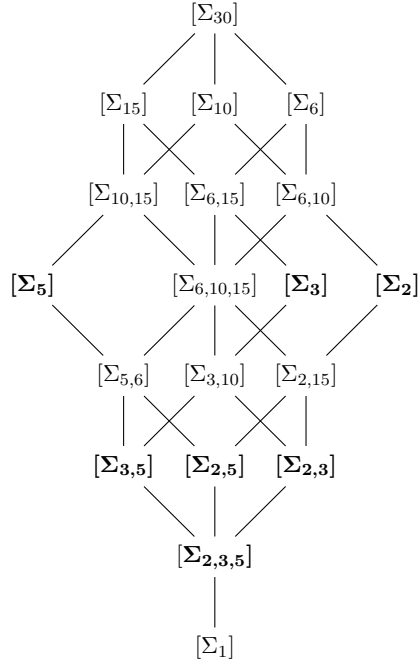
\begin{figure}
    \centering
        \begin{tikzpicture}[scale=0.6]
\def\myScale{0.8}

\node[scale=\myScale] (0) at (2,-2)  {$[\Sigma_{1}]$};
\node[scale=\myScale] (00) at (2,0)  {$[\boldsymbol{\Sigma_{2,3,5}]}$};
\node[scale=\myScale] (11) at (0,2)  {$\boldsymbol{[\Sigma_{3,5}]}$};
\node[scale=\myScale] (12) at (2,2) {$\boldsymbol{[\Sigma_{2,5}]}$};
\node[scale=\myScale] (13) at (4,2) {$\boldsymbol{[\Sigma_{2,3}]}$};
\node[scale=\myScale] (21) at (0,4) {$[\Sigma_{5,6}]$};
\node[scale=\myScale] (22) at (2,4) {$[\Sigma_{3,10}]$};
\node[scale=\myScale] (23) at (4,4) {$[\Sigma_{2,15}]$};
\node[scale=\myScale] (31) at (-2,6) {$\boldsymbol{[\Sigma_{5}]}$};
\node[scale=\myScale] (33) at (4,6) {$\boldsymbol{[\Sigma_{3}]}$};
\node[scale=\myScale] (34) at (6,6) {$\boldsymbol{[\Sigma_{2}]}$};
\node[scale=\myScale] (32) at (2,6) {$[\Sigma_{6,10,15}]$};
\node[scale=\myScale] (41) at (0,8) {$[\Sigma_{10,15}]$};
\node[scale=\myScale] (42) at (2,8) {$[\Sigma_{6,15}]$};
\node[scale=\myScale] (43) at (4,8) {$[\Sigma_{6,10}]$};
\node[scale=\myScale] (51) at (0,10) {$[\Sigma_{15}]$};
\node[scale=\myScale] (52) at (2,10) {$[\Sigma_{10}]$};
\node[scale=\myScale] (53) at (4,10) {$[\Sigma_{6}]$};
\node[scale=\myScale] (60) at (2,12) {$[\Sigma_{30}]$};
\path 
    (0)  edge (00)
    (00) edge (11)
    (00) edge (12)
    (00) edge (13)
    (11) edge (21)
    (11) edge (22)
    (12) edge (21)
    (12) edge (23)
    (13) edge (22)
    (13) edge (23)
    (21) edge (31)
    (21) edge (32)
    (22) edge (32)
    (22) edge (33)
    (23) edge (32)
    (23) edge (34)
    (31) edge (41)
    (32) edge (41)
    (32) edge (42)
    (32) edge (43)
    (33) edge (42)
    (34) edge (43)
    (41) edge (51)
    (41) edge (52)
    (42) edge (51)
    (42) edge (53)
    (43) edge (52)
    (43) edge (53)
    (51) edge (60)
    (52) edge (60)
    (53) edge (60)
    ;
\end{tikzpicture}
    \caption{The poset $\mCL$ restricted to cyclic loop conditions using 2, 3, and 5.}
    \label{fig:conditionPoset}
\end{figure}
A subposet of $\mCL$ is represented in Figure~\ref{fig:conditionPoset}.
We will discuss more properties of the posets $\mCL$ and $\MPCL$ in Section~\ref{sec:lattice}; prime cyclic loop conditions in bold.

\section{Unions of Cycles and Cyclic Loop Conditions}\label{sec:structures}

In this section we describe $\SDPoset$, namely the subposet of $\PPPoset$ containing the pp-constructability types of disjoint unions of cycles.

Firstly, in Lemma~\ref{lem:duofcSatisfyClc} we characterize $\models$ on disjoint unions of cycles and cyclic loop conditions. Using this characterization we prove the missing implication in Theorem~\ref{thm:implicationSingleClc}.
Secondly, in Lemma~\ref{thm:PPvsLoop}, we show that cyclic loop conditions suffice to
determine the pp-constructability order on finite disjoint unions of cycles.
This allows us to view a finite disjoint union of cycles as the set of prime cyclic loop conditions that it satisfies. These sets are always downsets of $\mPCL$. In Lemma~\ref{lem:surjectivety} we show when a downset is realized by a finite disjoint union of cycles. Combining these results we obtain a complete description of $\SDPoset$ as downsets of $\mPCL$.

\subsection{Characterization of the Satisfaction Relation on Unions of Cycles and Cyclic Loop Conditions}\label{ssc:models}
In Lemma~\ref{lem:duofcSatisfyClc} we generalize Lemma~\ref{lem:duofpcSatisfyPclc} to disjoint unions of cycles and cyclic loop conditions. 
First, recall that in Example~\ref{exa:C3notmodelsS3} we showed that $\Cyc3\not\models\Sigma_3$ by constructing a suitable tuple. We adopt the same idea in the following example.
\begin{example}\label{exa:C56notmodelsS25}
We show that the structure $\Cyc{5, 6}$ does not satisfy the cyclic loop condition $\Sigma_{2,5}$.
Define
\begin{align*}
\boldsymbol t\coloneqq{}& ((6,0),(6,3),(5,0),(5,3),(5,1),(5,4),(5,2)).
    \intertext{Note that $\boldsymbol t$ is contained in a cycle of length 30 in $\Cyc{5,6}^7$ and that}
\boldsymbol t+3={}& ((6,3),(6,0),(5,3),(5,1),(5,4),(5,2),(5,0)).
\end{align*} 
Since $\Cyc{5,6}^7$ is isomorphic to $\Cyc{5,6}^{\Cyc{2,5}}$ we view $\boldsymbol t$ as an element of $\Cyc{5,6}^{\Cyc{2,5}}$. Then  $\boldsymbol t+3=\boldsymbol t_{\sigma_{2,5}}$ and any $f\in\Pol(\Cyc{5,6})$ with $f\models\Sigma_{2,5}$ satisfies
\begin{align*}
f(\boldsymbol t)
=f_{\sigma_{2,5}}(\boldsymbol t)
=f(\shiftTuple{\boldsymbol t}{\sigma_{2,5}})
=f(\boldsymbol t+3).
\end{align*}
Therefore, $f$ would map $\boldsymbol t$ to an element of $\Cyc{5,6}$ that lies in a cycle whose length divides 3. Neither 5 nor 6 divide 3. Hence, such an $f$ does not exist and $\Cyc{5,6}\not\models\Sigma_{2,5}$. 
\end{example}

The following lemma shows that the line of reasoning used in Example~\ref{exa:C56notmodelsS25} is the only obstacle for a cyclic loop condition to be satisfied by a disjoint union of cycles.
It characterises the relation $\models$ on disjoint unions of cycles and cyclic loop conditions. Therefore, it will be a main tool in many of the subsequent proofs in this article.

\begin{lemma}\label{lem:duofcSatisfyClc}
Let $\C$ be a finite disjoint union of cycles and $D\subset\N^+\!$ finite. We have $\Pol(\C)\models \Sigma_D$ if and only if for all maps $h\colon D\to C$ there is an $a\in C$ such that
\[a\text{ divides }\lcm(\{h(b)\dotdiv b \mid {b\in D}\}).\]
\end{lemma}
The proof of this lemma is very similar to that of Lemma \ref{lem:duofpcSatisfyPclc}, it just has more technical difficulties.

\begin{proof}
Assume without loss of generality that $\C$ is the finite disjoint union of cycles associated to the set $C\subset \N^+$\!.

$(\Rightarrow)$
Consider a map $h\colon D\to C$ and a polymorphism $f\colon \C^{\D}\to\C$ of $\C$ satisfying $\Sigma_D$. Define 
\[c \coloneqq \lcm(\{h(b)\dotdiv b\mid{b\in D}\})\]
and the tuple $\boldsymbol t\in\C^{\D}$ as $\boldsymbol t_{(b,k)}\coloneqq (h(b),c\cdot k)$ for $(b,k)\in\D$.
We will show $f(\boldsymbol t)=f(\boldsymbol t+c)$.
Observe that $(\shiftTuple{\boldsymbol t}{\sigma_{\! D}})_{(b,k)}=\boldsymbol t_{\sigma_{\! D}(b,k)}=\boldsymbol t_{(b,k+1)}$ for all $(b,k)\in\D$.
Let $b\in D$. We have, by the definition of $c$, that $h(b)$ divides $c\cdot b$. Hence, $c\cdot b\equiv_{h(b)}0$ and 
\[(\boldsymbol t_{\sigma_{\! D}})_{(b,b-1)}=\boldsymbol t_{(b,0)}=(h(b),0)=(h(b),c\cdot b)=(h(b),c\cdot (b-1)+c)= \boldsymbol t_{(b,b-1)}+c.\] 
Therefore, $\shiftTuple{\boldsymbol t}{\sigma_{\! D}}=\boldsymbol t+c$, which just says that  $\boldsymbol t$ maps neighbouring points in $\D$ to points in $\C$ that are connected with a path of length $c$.
Furthermore, we have
\begin{align*}
f(\boldsymbol t)
=f_{\sigma_{\! D}}(\boldsymbol t)
=f(\shiftTuple{\boldsymbol t}{\sigma_{\! D}})
=f(\boldsymbol t+c).
\end{align*}

Since $f$ is a polymorphism we have that $f(\boldsymbol t) \stackrel{c}{\toEdge}f(\boldsymbol t+c)=f(\boldsymbol t)$. Hence, $f(\boldsymbol t)$ is in a cycle whose length divides $c$.

$(\Leftarrow)$
For this direction we construct a polymorphism $f\colon \C^{\D} \to \C$ of $\C$ satisfying $\Sigma_D$.
Let $\textbf G$ be the subgroup of the symmetric group on $\C^{\D}$ generated by $\sigma_{\! D}\colon \boldsymbol t \mapsto \shiftTuple{\boldsymbol t }{\sigma_{\! D}}$ and $+1\colon \boldsymbol t\mapsto (\boldsymbol t+1)$. Recall that
\begin{align*}
    (\shiftTuple{\boldsymbol t} {\sigma_{\! D}})_{(b,k)}=\boldsymbol t_{(b,k+1)},&& +1=\sigma_{\C^{\D}}, && ((+1)(\boldsymbol t))_{(b,k)}=\boldsymbol t_{(b,k)}+1.
\end{align*}
Hence, $\sigma_{\! D}\circ (+1)=(+1)\circ \sigma_{\! D}$, $\textbf G$ is commutative, and every element of $\textbf G$ is of the form $\sigma_{\! D}^c\circ (+1)^d$ for some $c, d \in\N$.
For every orbit of $\C^{\D}$ under $\textbf G$ pick a representative and denote the set of representatives by $T$. Let $\boldsymbol t\in T$ and define the map $h\colon D\to C$ such that for every $b\in D$ there is a $k\in \mathbb Z_{h(b)}$ with $\boldsymbol t_{(b,0)}=(h(b),k)$. 
By assumption there is an $a_{\boldsymbol t}\in C$ such that
\begin{equation}
a_{\boldsymbol t}\text{ divides } \lcm(\{h(b)\dotdiv b\mid {b\in D}\}).\label{equ:aDivides}
\end{equation}
Note that the orbit of $\boldsymbol t$ is a disjoint union of cycles of a fixed length $n$ and that $a_{\boldsymbol t}$ divides $n$. Define $f$ on the orbit of $\boldsymbol t$ as 
\[f\left((\sigma_{\! D}^c\circ (+1)^d)(\boldsymbol t)\right)
=f\left(\shiftTuple {(\boldsymbol t+d)}{\sigma_{\! D}^c}\right)
\coloneqq (a_{\boldsymbol t},d)\text{ for all $c,d$.}\] 
To show that $f$ is well defined on the orbit of $\boldsymbol t$ it suffices to prove that $\shiftTuple{(\boldsymbol t+d)}{\sigma_{\! D}^c}=\shiftTuple{(\boldsymbol t+\ell)}{\sigma_{\! D}^{m}}$ implies $d\equiv_{a_{\boldsymbol t}}\ell $ for all $d,c,\ell,m\in\N$. Without loss of generality we can assume that $m=\ell=0$. 
Fix some $b\in D$. We want to show that $(h(b)\dotdiv b)$ divides $d$.
Observe that $\boldsymbol t=\shiftTuple{(\boldsymbol t+d)}{\sigma_{\! D}^c}$ implies
\[\boldsymbol t_{(b,k\cdot c)}=\boldsymbol t_{(b,(k+1)\cdot c)}+d \text{ for all $k$.}\]
Considering that $(b\dotdiv c)\cdot c\equiv_{b}0$ we have 
\[\boldsymbol t_{(b,0)}=\boldsymbol t_{(b,c)}+d=\boldsymbol t_{(b,2\cdot c)}+2\cdot d=\dots=\boldsymbol t_{(b,0)}+(b\dotdiv c)\cdot d.\] 
Hence, $(b\dotdiv c)\cdot d\equiv_{h(b)}0$ and there is some $k\in\N$ such that $d=\frac{k\cdot h(b)}{b\dotdiv c}$. Therefore, $h(b)\dotdiv (b\dotdiv c)$ divides $d$ and also $h(b)\dotdiv b$ divides $d$. Since this holds for all $b\in D$ we have, by \eqref{equ:aDivides}, that $a_{\boldsymbol t}$ divides $d$ as desired.

Repeating this for every $\boldsymbol t\in T$ defines $f$ on $\C^{\D}$. The function $f$ is well defined since the orbits partition  $\C^{\D}$.
If $\boldsymbol r\stackrel{1}\toEdge \boldsymbol s$, then $\boldsymbol s=(\boldsymbol r+1)=(+1)(\boldsymbol r)$ and $f(\boldsymbol r)\stackrel{1}\toEdge f(\boldsymbol s)$. Hence, $f$ is a polymorphism of $\C$. Furthermore, $f=f_{\sigma_{\! D}}$ by definition.
\end{proof}

\begin{examplex}
Some applications of Lemma~\ref{lem:duofcSatisfyClc}:
\begin{itemize}
    \item $\Cyc{10} \models\Sigma_{2,5}$, since $h(2)=h(5)=10$ is the only map from $\{2,5\}$ to $\{10\}$ and 10 divides $10=\lcm(h(2)\dotdiv2,h(5)\dotdiv 5)$,
    \item $\Cyc{n}\not\models\Sigma_n$ for $n>1$, witnessed by $h(n)=n$, since $n$ does not divide $1=h(n)\dotdiv n$,
    \item $\Cyc{5, 6} \not\models\Sigma_{2,5}$, witnessed by $h(2)=6$, $h(5)=5$, since neither 5 nor 6 divide $3=\lcm(h(2)\dotdiv2,h(5)\dotdiv 5)$.
\hfill$\triangle$
\end{itemize}
\end{examplex}

As a consequence of Lemma~\ref{lem:duofcSatisfyClc} we obtain the following lemma.
\begin{lemma}\label{lem:AdoesntSatisfyOwnLoops}
Let $\C$ be a finite disjoint union of cycles and $c\in\N^+$\!. Then
\begin{align*}
    \Pol(\C)\models\Sigma_{C\dotdiv c} &&\text{if and only if}&&1\in(C\dotdiv c).
\end{align*}
\end{lemma}
\begin{proof}
The $\Leftarrow$ direction is clear. If $1\in(C\dotdiv c)$, then $\Sigma_{C\dotdiv c}$ is trivial since it is satisfied by the projection $\pi\colon \C^{\C\dotdiv c}\to \C$, $\boldsymbol{t}\mapsto \boldsymbol{t}_{(1,0)}$.

For the $\Rightarrow$ direction  
let $h\colon (C\dotdiv c)\to C$ be some map with $h(b)\dotdiv c=b$ for every $b\in C\dotdiv c$. For instance, $h(b)=\min\{a\in C\mid a\dotdiv c=b \}$.
Note that $h(d\dotdiv c)=d$ for all $d\in\operatorname{Im}(h)$. We apply Lemma~\ref{lem:duofcSatisfyClc} to $h$ and obtain some $\Cyc{a}\hookrightarrow \C$ such that $a$ divides 
\begin{align*}
\lcm(\{h(b)\dotdiv b \mid {b\in (C\dotdiv c)}\})&=
\lcm(\{d\dotdiv(d\dotdiv c)\mid{d\in \operatorname{Im}(h)}\})\\
&=\lcm(\{\gcd(d,c)\mid {d \in  \operatorname{Im}(h)}\})    
\end{align*}
which divides $c$. Therefore, $(a\dotdiv c)=1\in (C\dotdiv c)$.
\end{proof}

Note that we did not use any of the results from Section~\ref{sec:conditions} to prove Lemmata~\ref{lem:duofcSatisfyClc} and~\ref{lem:AdoesntSatisfyOwnLoops}.
Using these lemmata we can prove the missing implication from Theorem~\ref{thm:implicationSingleClc}. 

\begin{lemma}\label{lem:missingImplicationSingleClc}
Let $C,D\subset\N^+\!$ finite. We have that (2) implies (3).
\begin{enumerate}[label=(\arabic*)]
\setcounter{enumi}{1}
\item For all $c\in \N^+\!$ we have $\Pol(\D\dotdiv c)\models\Sigma_C$ implies $\Pol(\D\dotdiv c)\models\Sigma_D$.
\item For every $a\in C$ there exist $b\in D$ and $k\in\N$ such that $b$ divides $a^k$.
\end{enumerate}
\end{lemma}

\begin{proof}
We show the contraposition. Let $a\in C$ be such that no $b\in D$ divides $a^k$ for any $k\in\N$. Define $c\coloneqq a^k$ where $k$ is the highest prime power of any number in $D$, i.e., $k$ is the largest $\ell\in\N$ for which there is a prime $p$ and an $b\in D$ such that $p^{\ell}$ divides $b$. Consider the structure $\D\dotdiv c$. By construction $1\notin (D\dotdiv c)$. Hence, by Lemma~\ref{lem:AdoesntSatisfyOwnLoops}, $\Pol(\D\dotdiv c)\not\models\Sigma_D$.

Using Lemma~\ref{lem:duofcSatisfyClc} we prove  $\Pol(\D\dotdiv c)\models\Sigma_C$. Let $h\colon C\to D\dotdiv c$ be a map. Choose any $b\in D$ such that $h(a)=b\dotdiv c$. Then, by construction of $c$, 
\[h(a)=b\dotdiv c=b\dotdiv a^k=b\dotdiv a^k\dotdiv a=h(a)\dotdiv a\]
and $h(a)$ divides $\lcm(\{h(\tilde a)\dotdiv \tilde a\mid{\tilde a\in C}\})$. Applying Lemma~\ref{lem:duofcSatisfyClc} yields $\Pol(\D\dotdiv c)\models\Sigma_C$ as desired. 
\end{proof}

\subsection{On the pp-constructability of  Unions of Cycles}\label{ssc:ppOrder}

Now we have all the necessary ingredients to prove the connection between cyclic loop conditions and pp-constructions in $\SDPoset$ stated in Lemma~\ref{thm:PPvsLoop}.
In particular,  we
show that cyclic loop conditions suffice to separate disjoint unions of cycles.
We suggest to look at the following concrete pp-constructions first.
Recall that for every $k\in\N^+$ we abbreviate the pp-formula
\[\exists y_1,\dots,y_{k-1}.\ E(x,y_1)\wedge E(y_1,y_2)\wedge\dots\wedge E(y_{k-1},z)\]
by $x\stackrel{k}\toEdge y$.
\begin{example}\label{ex.components}
The digraph $\Cyc{2,3}$ pp-constructs $\Cyc6$. Consider the second pp-power of $\Cyc{2,3}$ given by the pp-formula:
\[\Phi_E(x_1,x_2,y_1,y_2)\coloneqq (x_1\stackrel{1}{\toEdge}y_1) \wedge (x_2\stackrel{1}{\toEdge}y_2) \wedge (x_1\stackrel{2}{\toEdge}x_1) \wedge (x_2\stackrel{3}{\toEdge}x_2).\]
The resulting structure, which consists of one copy of $\Cyc6$ and 19 isolated points, is homomorphically equivalent to $\Cyc6$, and therefore $\Cyc{2,3}\leq\Cyc6$.
\end{example}

\begin{example}\label{ex.jakubconstruction}
The digraph $\Cyc3$ pp-constructs $\Cyc9$. Consider the third pp-power of $\Cyc3$ given by the formula:
\begin{align*}
    \Phi_E(x_1,x_2,x_3,y_1,y_2,y_3)\coloneqq (x_2\approx y_1) \wedge 
    (x_3\approx y_2)\wedge (x_1\stackrel{1}{\toEdge}y_3).
    \tag{$\ast$}
\end{align*}
Let us denote the resulting structure by $\C$. There is an edge $\boldsymbol s\toEdge\boldsymbol t$ in $\C$ if the tuple $\boldsymbol t$ is obtained from $\boldsymbol s$ by first increasing the first entry and then shifting all entries cyclically, see Figure~\ref{fig:incShift}. 
\begin{figure}
    \centering
    \begin{tikzpicture}
    \node at (-0.7,0) {$\boldsymbol s=$};
    \node at (0,0) {$\begin{pmatrix}
    a\\b\\c
    \end{pmatrix}$};
    \node at (3,0) {$\begin{pmatrix}
    a+1\\b\\c
    \end{pmatrix}$};
    \node at (6,0) {$\begin{pmatrix}
    b\\c\\a+1
    \end{pmatrix}$};
    \node at (7,0) {$=\boldsymbol t$};
    
    \path[->] 
    (0.55,0) edge node[above] {increase} (2.2,0)
    (3.8,0) edge node[above] {shift} (5.2,0);
    \end{tikzpicture}
    \caption{The shape of tuples $\boldsymbol s$ and $\boldsymbol t$ in the edge-relation defined by the pp-formula $(\ast)$.
    }
    \label{fig:incShift}
\end{figure}
With this it is clear that for every element $\boldsymbol t$ in $\C$ we have $\boldsymbol t\stackrel{9}{\toEdge}\boldsymbol t$. 
It turns out that $\C$ consists of three copies of $\Cyc9$. Hence, $\Cyc3\leq\Cyc9$ and even $\Cyc3\equiv\Cyc9$.

Note that the third pp-power of $\Cyc2$ given by the formula $(\ast)$ is not homomorphically equivalent to $\Cyc6$; instead, it is isomorphic to $\Cyc{2,6}$, which is homomorphically equivalent to $\Cyc2$.
\end{example}

Although it is neither clear nor necessary we would like to mention that the pp-construction in the proof of the following lemma is essentially just a combination of the three constructions we saw in the Examples~\ref{ex.division},~\ref{ex.components} and~\ref{ex.jakubconstruction}. 

\begin{lemma}\label{thm:PPvsLoop}
Let $\C$ be a finite disjoint union of cycles and let $\B$ be a finite structure with finite relational signature $\tau$. Then 
\begin{align*}
    \B\leq\mathbb C && \text{if and only if} &&
    \tikz[baseline=-2pt]{
    \node[text width=7cm] {$\Pol(\B)\models\Sigma_{C\dotdiv c}$ implies $\Pol(\C)\models\Sigma_{C\dotdiv c}$ for all $c$ that divide $\lcm(C)$.};}
\end{align*}
\end{lemma}
We remark that the first part of the following proof is a specific instance of the proof of 
$ \F_{\Pol(\B)}(\A)\to\A$ implies  $\B\leq\A$
in Theorem~\ref{thm:freestructure}.
As the reader might not be familiar with free structures, we present a self-contained proof.
\begin{proof}
We show both directions separately.

$(\Rightarrow)$ 
Since $\Sigma_{C\dotdiv c}$ is a minor condition, this direction follows 
from Corollary~\ref{cor:wond}.
$(\Leftarrow)$ Assume without loss of generality that $\C$ is the structure associated to $C$.  
Let $\F$ be the $\left|\B^{\C}\right|$-th pp-power of $\B$ defined by the formula
\begin{align*}
\Phi_E(x,y)\coloneqq\bigwedge\left\{
x_{\scalebox{0.77}{$\shiftTuple{\boldsymbol t}{{\sigma_{\! C}}}$}}\approx 
y_{\boldsymbol t} \wedge \Phi_R(x)\ \middle|\ \boldsymbol t\in \B^{\C}, R\in\tau \right\}\!,
\end{align*}
where for every $k$-ary   $R\in\tau$ we have that $\Phi_R(x)$ is defined as
\begin{align*}
\bigwedge\left\{R(x_{\boldsymbol t_1},\dots,x_{\boldsymbol t_k})\ \middle|\ \boldsymbol t_1,\dots,\boldsymbol t_k\in\B^{\C},\text{ $( t_{1u},\dots, t_{ku})\in R^{\B}$ for all $u\in\C$}\right\}\!.
\end{align*}
We can think of the elements of $\F$ as maps from $\B^{\C}$ to $\B$. 
The formula $\Phi_R(f)$ holds if and only if $f$ preserves $R^{\B}$. Note that $f$ preserves $R^{\B}$ if and only if $f_{\sigma_{\!C}}$ preserves $R^{\B}$. Hence, $\Phi_E$ ensures that all elements of $\F$ that are not polymorphisms of $\B$ are isolated points. On the other hand polymorphisms $f$ of $\B$ that are in $\F$ have exactly one in-neighbour, namely $f_{\sigma^{-1}_{\! C}}$, and one out-neighbour, namely $f_{\sigma_{\! C}}$. 
Hence, $\F$ is homomorphically equivalent to a disjoint union of cycles, i.e., the structure $\F$ without isolated points. 
Furthermore, all cycles in $\F$ are of the form 
\[ f\toEdge f_{\sigma_{\! C}}\toEdge f_{\sigma^2_{\! C}}\toEdge \dots \toEdge f_{\sigma^k_{\! C}}=f \]
for some $k\in\N$.

We show that $\F$ and $\C$ are homomorphically equivalent by proving the following two statements:
\begin{enumerate}
\item For all $a\in\N^+\!$ we have $\Cyc{a}\hookrightarrow \C$ implies $\Cyc{a} \hookrightarrow \F$ and
\item for all $c\in\N^+\!$ we have $\Cyc{c} \hookrightarrow \F$ implies $\Cyc{c}\to \C$.
\end{enumerate}
First statement: Suppose that $\Cyc{a}\hookrightarrow \C$. Then the polymorphism 
\[\pi_{(a,0)}\colon\B^{\C}\to\B,  \boldsymbol t\mapsto \boldsymbol t_{(a,0)}\] 
lies in the following cycle of length $a$ in $\F$:
\[\pi_{(a,0)}\toEdge\pi_{(a,1)}\toEdge\dots\toEdge\pi_{(a,a-1)}\toEdge\pi_{(a,0)}.\]
Second statement: Suppose that $\Cyc{c}\hookrightarrow \F$ and let $f$ be a polymorphism in a cycle of length $c$ in $\F$. Then $f=f_{\sigma^c_{\! C}}$. Let $\D$ be the $c$-th relational power of $\C$. Note that $\sigma_{\D}=\sigma^c_{\! C}$. Hence,  $f\models \Sigma_{\D}$. By Lemma~\ref{lem:relPowEqualsCdotdivc}, the digraphs $\D$ and $\C\dotdiv c$ are homomorphically equivalent. Therefore, $\Pol(\B)\models\Sigma_{C\dotdiv c}$ and, by assumption, $\Pol(\C)\models\Sigma_{C\dotdiv c}$ as well. Applying Lemma~\ref{lem:AdoesntSatisfyOwnLoops} we conclude that $1\in(C\dotdiv c)$. Hence, there is some $a\in C$ such that $a$ divides $c$ and $\Cyc{c}\to \Cyc{a}\hookrightarrow \C$.
It follows that $\F$ and $\C$ are homomorphically equivalent. Hence, $\C$ is pp-constructable from $\B$.
\end{proof}

The construction in the proof was discovered by Jakub Opr\v{s}al (see~\cite{Barto_2021FreeStructures} or~\cite{oprsal2018taylors} for more details). We thank him for explaining it to us.
Note that, following the notation introduced in Definition~\ref{def:freestructure}, the structure $\F$, after removing all isolated points, is $\F_{\Pol(\B)}(\C)$.
In particular this lemma shows the following. 
\begin{corollary}\label{cor:PPvsLoopCaBlockerForSa}
Let $a$ be in $\N^+$ and let $\B$ be a finite structure with finite relational signature. Then 
\begin{align*}
\B\leq \C_a && \text{iff} && \Pol(\B)\not\models\Sigma_{p}\text{ for all $p$ that are prime divisors of $a$}.
\end{align*}
\end{corollary}

\begin{example}\label{ex.6,20,15}
Suppose we want to test whether a structure $\B$ can  pp-construct $\Cyc{6,20,15}$.
By Lemma~\ref{lem:AdoesntSatisfyOwnLoops}, $\Pol(\Cyc{6,20,15})$ does not satisfy any non-trivial loop conditions of the form $\Sigma_{\{6,20,15\}\dotdiv c}$. Hence, by Lemma~\ref{thm:PPvsLoop}, to verify that $\B\leq \Cyc{6,20,15}$ we only have to check whether $\B$ satisfies none of the cyclic loop conditions $\Sigma_{6,20,15}$, $\Sigma_{2,20,5}$, $\Sigma_{3,10,15}$, $\Sigma_{6,4,3}$, $\Sigma_{3,5,15}$, $\Sigma_{3,2,3}$, which are the non-trivial ones of the form $\Sigma_{\{6,20,15\}\dotdiv c}$. By Theorem \ref{thm:implicationSingleClc}, these conditions are equivalent to  $\Sigma_{6,10,15}$, $\Sigma_{2,5}$, $\Sigma_{3,10}$, $\Sigma_{2,3}$, $\Sigma_{3,5}$, and $\Sigma_{2,3}$, respectively.

We show that the disjoint union of cycles $\Cyc{2,3,5}$ can pp-construct $\Cyc{6,20,15}$. First we check that $\Cyc{2,3,5}\not\models \Sigma_{6,10,15}$. Consider the map $h(6)=2$, $h(10)=2$, $h(15)=3$. We have
\[\lcm(2\dotdiv 6, 2\dotdiv 10, 3\dotdiv 15)=\lcm(1,1,1)=1.\]
Clearly, neither $2$ nor $3$ nor $5$ divide 1. Hence, by Lemma~\ref{lem:duofcSatisfyClc}, we have $\Cyc{2,3,5}\not\models \Sigma_{6,10,15}$. Similarly, $\Cyc{2,3,5}$ does not satisfy the other four loop conditions. Therefore, $\Cyc{2,3,5}\leq\Cyc{6,20,15}$. 
On the other hand, the structure $\Cyc{2,15}$ cannot pp-construct $\Cyc{6,20,15}$ since it satisfies  $\Sigma_{3,5}$. 
\end{example}

By  Theorem~\ref{thm:clcAreSetsOfPclc}  every cyclic loop condition is equivalent to a set of prime cyclic loop conditions; we thus obtain that even prime cyclic loop conditions suffice to separate disjoint unions of cycles.

\begin{corollary}\label{thm:classificationPPCon}
Let $\mathbb C$ be a finite disjoint union of cycles and $\B$ be a finite structure with finite relational signature. Then
\begin{align*}
    \B\leq\mathbb C && \text{if and only if} &&
    \tikz[baseline=-2pt]{
    \node[text width=6cm] {$\Pol(\B)\models\Sigma$ implies $\Pol(\C)\models\Sigma$ for all prime cyclic loop conditions $\Sigma$.};}
\end{align*}


\end{corollary}
Note that we can also prove this corollary using easier results. The cyclic loop conditions that are considered in Lemma~\ref{thm:PPvsLoop} to test whether $\B$ pp-constructs $\C$ have the form $\Sigma_{C\dotdiv c}$. Note that if $1\notin C$ then $\C$ does not satisfy any of the conditions $\Sigma_{C\dotdiv c}$. Hence, to verify the condition in the lemma it suffices to check whether $\B$ does not satisfy any condition that is minimal in $\{\Sigma_{C\dotdiv c}\mid c\in\N^+\!, 1\notin(C\dotdiv c)\}$. 
Any such minimal condition $\Sigma_{C\dotdiv c}$ has a $c$ that is maximal for $C$ and is, by Lemma~\ref{lem:minimalIsPrime}, equivalent to a prime cyclic loop condition. Combining this observations with Lemma~\ref{lem:duofcSatisfyClc} we obtain a practical way to verify whether two disjoint unions of cycles can pp-construct each other.
\begin{corollary}
Let $\C$ and $\D$ be finite disjoint unions of cycles. Then $\D\ppleq\C$ if and only if for all $c$ such that every number in $C\dotdiv c$ is a multiple of some prime in $C\dotdiv c$ there is a map $h\colon  (C\dotdiv c)\to D$ such that for all $a\in D$ 
\[a\text{ does not divide }\lcm(\{h(b)\dotdiv b \mid {b\in (C\dotdiv c)}\}).\]
\end{corollary}

\subsection{Characterizing Unions of Cycles in Terms of Prime Cyclic Loop Conditions}\label{ssc:SDPoset}
As a consequence of Corollary~\ref{thm:classificationPPCon}, every element $[\C]$ of $\SDPoset$ is uniquely determined by the set of prime cyclic loop conditions that $\C$ satisfies. Hence, the map  
\[\PL \colon \C \mapsto \{\Sigma\mid \text{$\Sigma$ a prime cyclic loop condition and } \Pol(\C)\models\Sigma \}\] 
is injective. 
The next goal is to determine the image of $\PL$. 
First we simplify the characterization from Lemma~\ref{lem:duofcSatisfyClc} for prime cyclic loop conditions.

\begin{lemma}\label{lem:duofcSatisfyPclc}
Let $\C$ be a finite disjoint union of cycles and let $\Sigma_P$ be a prime cyclic loop condition. Then the following are equivalent:
\begin{enumerate}[label=(\arabic*)]
    \item $\Pol(\C)\not\models\Sigma_P$.
    \item There is a $c\in\N^+\!$ such that $1\notin (C\dotdiv c)$ and $P\subseteq (C\dotdiv c)$.
    \item There is a $c\in\N^+\!$ such that $\Sigma_{C\dotdiv c}$ is non-trivial and $\Sigma_P\Rightarrow\Sigma_{C\dotdiv c}$.
\end{enumerate} 
\end{lemma}
\begin{proof}
(1) $\Rightarrow$ (2) 
By Lemma~\ref{lem:duofcSatisfyClc} $\Pol(\C)\not\models\Sigma_P$ implies that there is a map $h\colon P\to C$ such that no $a\in C$ divides $c$, where $c\coloneqq \lcm(\{h(p)\dotdiv p\mid {p\in P}\})$.
In particular, every $p$ divides $h(p)$.
Note that if $a\dotdiv c=1$, then $\gcd(a,c)=a$ and  $a$ divides $c$.
Hence, we have $1\notin(C\dotdiv c)$. 
Let $p\in P$. Then $h(p)$ does not divide $c$. Hence, $h(p)\dotdiv c\neq1$. Furthermore, $h(p)\dotdiv c$ divides 
\[h(p)\dotdiv (h(p)\dotdiv p)=\gcd(h(p),p)=p.\]
Therefore,  $h(p)\dotdiv c=p$ for all $p$ and $P\subseteq (C\dotdiv c)$.

The directions
(2) $\Rightarrow$ (3) and (3) $\Rightarrow$ (1) follow from Lemma~\ref{lem:orderOnLoopcond} and  Lemma~\ref{lem:AdoesntSatisfyOwnLoops}, respectively.
\end{proof}

\begin{lemma}\label{lem:surjectivetyUpperBound}
For every disjoint union of cycles $\C$ we have that $\PL(\C)$ is a cofinite downset of $\mPCL$.
\end{lemma}
\begin{proof}
Let $P$ be the set of all prime divisors of $\lcm(C)$. Using Lemma~\ref{lem:duofcSatisfyPclc} we conclude that every prime cyclic loop condition $\Sigma_S$ with $\Pol(\C)\not\models\Sigma_S$ satisfies $S\subseteq P$. Since $P$ is finite there are only finitely many prime cyclic loop conditions that are not satisfied by $\C$. 
\end{proof}

Next we show that every cofinite downset of $\mPCL$ is also realized by some disjoint union of cycles.
\begin{lemma}\label{lem:surjectivety}
Let $\Gamma$ be a cofinite downset of $\mPCL$ and $\Downset_{\min}$ the set of minimal prime cyclic loop conditions of $\mPCL\setminus\Downset$. Then 
\[\C\coloneqq\bigtimes_{\Sigma\in\Downset_{\min}} \mathbb G_\Sigma\]
is a finite disjoint union of cycles and $\PL(\C)=\Gamma$.
\end{lemma}
\begin{proof}
Let $P$ denote the set $\bigcup\{T\mid \Sigma_T\in\Downset_{\min}\}$, which contains the primes occurring in $\Downset_{\min}$. 
Since $\Downset$ is cofinite we have that $\Downset_{\min}$ is finite. Hence, $\C$ is a finite disjoint union of cycles and
\[ C= \left\{
\lcm(\{p_T\mid { \Sigma_T\in\Downset_{\min}}\})
~\middle|~ 
\text{$p_T\in T$ for every $\Sigma_T\in\Downset_{\min}$}\right\}\!.\] 
We prove that $\PL(\C)=\Downset$. 

$(\subseteq)$  
Let $\Sigma_S\in\mPCL\setminus\Downset$. Since $\PL(\C)$ is closed under implication, we can assume $\Sigma_S$ to be minimal in $\mPCL\setminus\Downset$. Define $c_S\coloneqq \prod(P\setminus S)$. 
We want to apply Lemma~\ref{lem:duofcSatisfyPclc} to show $\Pol(\C)\not\models\Sigma_S$. 
Firstly, note that, since $\Sigma_S\in\Downset_{\min}$, any $a\in C$ is a multiple of some prime $p\in S$. Furthermore, this $p$ does not divide $c_S$. Hence, $a\dotdiv c_S\not=1$ and $1\notin (C\dotdiv c_S)$. 
Secondly, let $p\in S$. Since $\Sigma_S$ is minimal we have that for every other $\Sigma_T\in\Downset_{\min}$ there exists a $p_T\in T\setminus S$. Define $p_S\coloneqq p$ and $a\coloneqq \lcm(\{p_T\mid { \Sigma_T\in\Downset_{\min}}\})$. Then $a\in C$ and $a\dotdiv c_S=p$. Therefore, $p\in (C\dotdiv c_S)$. Hence, $S\subseteq (C\dotdiv c_S)$ and, by Lemma~\ref{lem:duofcSatisfyPclc}, $\Pol(\C)\not\models \Sigma_S$ as desired.

$(\supseteq)$ Let $\Sigma_S\in\mPCL\setminus\PL(\C)$. Since $\Pol(\C)\not\models\Sigma_S$, by Lemma~\ref{lem:duofcSatisfyPclc}, we have that $S$ is contained in a set of the form $(C\dotdiv c)$.
Next we show that there is some $\Sigma_{S_c}\in\Downset_{\min}$ such that $S_c$ is contained in $(C\dotdiv c)$. 
Assume for contradiction that for every $\Sigma_T\in\Downset_{\min}$ the set $T$ is not contained in $(C\dotdiv c)$. Let $p_T$ be a witness of this fact.  
Note that $p_T\notin (C\dotdiv c)$ implies $p_T$ divides $c$. Then $a\coloneqq \lcm(\{p_T\mid { \Sigma_T\in\Downset_{\min}}\})\in C$ but $a\dotdiv c=1$, a contradiction. Hence, there is a $\Sigma_{S_c}\in\Downset_{\min}$ such that $S_c\subseteq (C\dotdiv c)$.

We show that $S\subseteq S_c$. Let $p\in S\subseteq (C\dotdiv c)$. Then there is some $a\in C$ such that $a\dotdiv c=p$. Again $a$ is of the form $\lcm(\{p_T\mid{ T\in\Downset_{\min}}\})$. Note that, since all numbers in $C$ are square-free, no element from $S_c$ can divide $c$. Hence, $p=p_{S_c}\in S_c$.
Therefore, $S\subseteq S_c$ and $\Sigma_S$ implies $\Sigma_{S_c}$. Since $\Downset$ is implication-closed and $\Sigma_{S_c}\notin\Gamma$ we conclude that $\Sigma_S\notin \Downset$. This yields $\PL(\C)=\Downset$, as desired.
\end{proof}

The following corollaries are immediate consequences of Corollary~\ref{thm:classificationPPCon}, Lemma~\ref{lem:surjectivetyUpperBound}, and Lemma~\ref{lem:surjectivety}. 

\begin{corollary}
Let $\C$ be a finite disjoint union of cycles. 
Then there are unions of prime cycles $\C_{P_1},\dots,\C_{P_n}$ such that $\C\ppleq \C_{P_1},\dots,\C\ppleq \C_{P_n}$ and $\C\ppeq \C_{P_1}\mathbin{\times}\dots\mathbin{\times}\C_{P_n}$.
In particular $\C$ has the same pp-constructability type as a disjoint union of cycles whose  cycle  lengths  are  square-free.
\end{corollary}

Combining this with Theorem~\ref{thm:BartoKozikNiven} and Corollary~\ref{thm:classificationPoset} we obtain the following dichotomy for smooth digraphs.

\begin{corollary}\label{cor:cyclesSquarefree}
Let $\mathbb{G}$ be a finite smooth digraph. Then either $[\mathbb G]=[\K_3]$ or there is a finite disjoint union of cycles $\C$ whose cycle lengths are square-free such that $[\mathbb G]= [\C]$.
\end{corollary}

We finally obtain a classification of $\SDPoset$.
\begin{corollary}
\label{thm:classificationPoset}
The map
\begin{align*}
    \SDPoset&\to\MPCL\\
    [\C]&\mapsto[\PL(\C)]
\end{align*}
is well defined and an embedding of posets. Its image consists of the elements that can be represented by a cofinite downset of $\mPCL$. Put differently, \[\SDPoset\simeq(\operatorname{CofiniteDownsetsOf}(\mPCL),\subseteq).\]
\end{corollary}

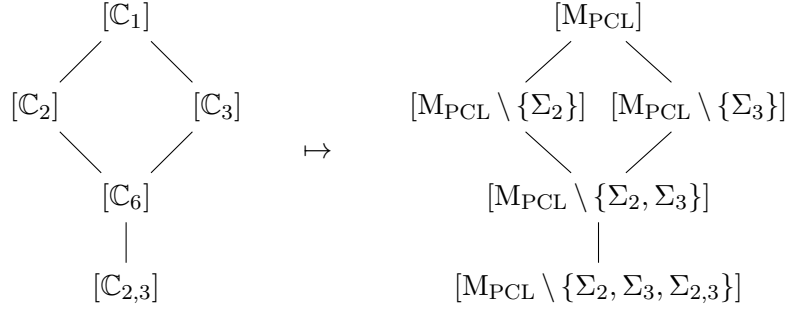
\begin{figure}
    \centering
    \begin{tikzpicture}[scale=0.6]

\node (0) at (2,-2)  {$[\Cyc{2,3}]$};
\node (00) at (2,0)  {$[\Cyc6]$};
\node (11) at (0,2)  {$[\Cyc2]$};
\node (13) at (4,2) {$[\Cyc3]$};
\node (22) at (2,4) {$[\Cyc1]$};

\node at (6.2,1) {$\mapsto$};
\path 
    (0)  edge (00)
    (00) edge (11)
    (00) edge (13)
    (11) edge (22)
    (13) edge (22)
    ;
\end{tikzpicture}
\hspace{5mm}
\begin{tikzpicture}[scale=0.6]

\node (0) at (2,-2)  {$[\mPCL\setminus\left\{\Sigma_{2},\Sigma_{3},\Sigma_{2,3}\right\}]$};
\node (00) at (2,0)  {$[\mPCL\setminus\left\{\Sigma_{2},\Sigma_{3}\right\}]$};
\node (11) at (-0.2,2)  {$[\mPCL\setminus\left\{\Sigma_{2}\right\}]$};
\node (13) at (4.2,2) {$[\mPCL\setminus\left\{\Sigma_{3}\right\}]$};
\node (22) at (2,4) {$[\mPCL]$};

\path 
    (0)  edge (00)
    (00) edge (11)
    (00) edge (13)
    (11) edge (22)
    (13) edge (22)
    ;
\end{tikzpicture}
\caption{ The embedding from Corollary~\ref{thm:classificationPoset} restricted to disjoint unions of cycles of length 2 and 3.}
    \label{fig:example23}
\end{figure}

To better understand what the map from Corollary~\ref{thm:classificationPoset} does, have a look at the illustration in Figure~\ref{fig:example23}.
We can give an explicit description of $\PL(\C)$.

\begin{lemma}\label{lem:pclCharacterization}
Let $\C$ be a finite disjoint union of cycles. Then
\[\PL(\C)=\mPCL\setminus\operatorname{UpsetOf}(\{\Sigma_{C\dotdiv c}\mid c\text{ is maximal for $C$}\}).\]
\end{lemma}
\begin{proof}
($\subseteq$) 
By Lemma~\ref{lem:AdoesntSatisfyOwnLoops}, we have that $\C$ does not satisfy $\Sigma_{C\dotdiv c}$ for any $c$ that is maximal for $C$. Hence, no condition in  \[\operatorname{UpsetOf}(\{\Sigma_{C\dotdiv c}\mid c\text{ is maximal for $C$}\})\] is satisfied by $\C$.

($\supseteq$) Let $\Sigma_P$ be a prime cyclic loop condition such that $\Pol(\C)\not\models\Sigma_P$. Then, by Lemma~\ref{lem:duofcSatisfyPclc}, there is a $c\in\N^+\!$ such that $\Sigma_{C\dotdiv c}$ is non-trivial and $\Sigma_P\Rightarrow\Sigma_{C\dotdiv c}$. Note that we can choose $c$ to be maximal for $C$. Hence, the condition $\Sigma_P$ is in $\operatorname{UpsetOf}(\{\Sigma_{C\dotdiv c}\mid c\text{ is maximal for $C$}\})$.
\end{proof}

For a given finite disjoint union of cycles, Lemmata~\ref{lem:surjectivety} and~\ref{lem:pclCharacterization}
provide a method to construct a finite disjoint union of cycles of square-free length with the same pp-constructability type which can be carried out by hand on small examples, as illustrated by the following example.
\begin{example}\label{exa:620Structure}
Consider the structure $\Cyc{6,20}$.  The set from Lemma~\ref{lem:pclCharacterization} containing all conditions of the form $\Sigma_{\{6,20\}\dotdiv c}$ with $c$ maximal for $\{6,20\}$ is $\{\Sigma_{2,4},\Sigma_{3,2},\Sigma_{3,5}\}$. This set is equivalent to   $\{\Sigma_{2},\Sigma_{3,2},\Sigma_{3,5}\}$.
Note that $\Sigma_2\Rightarrow\Sigma_{3,2}$. Hence,
\begin{align*}
\PL(\Cyc{6,20})&\makebox[7mm]{${}\stackrel{\ref{lem:pclCharacterization}}{=}{}$}\operatorname{UpsetOf}(\{\Sigma_{2},\Sigma_{3,2},\Sigma_{3,5}\})\\
&\makebox[7mm]{${}={}$}\operatorname{UpsetOf}(\{\Sigma_{3,2},\Sigma_{3,5}\})\\
&\makebox[7mm]{${}\stackrel{\ref{lem:surjectivety}}{=}{}$}\PL(\Cyc{3,10}).
\end{align*}
Therefore, $\Cyc{6,20}\equiv\Cyc{3,10}$. 
\end{example}

We have a similar situation for cyclic loop conditions: Corollary \ref{cor:flat} in particular states that every cyclic loop condition is equivalent to one that only uses square-free numbers.
Recall from Example~\ref{exa:620} that \[\Sigma_{6,20}\Leftrightarrow\{\Sigma_{2},\Sigma_{3,2},\Sigma_{3,5}\}\Leftrightarrow\{\Sigma_{2},\Sigma_{3,5}\}\Leftrightarrow\Sigma_{6,10}.\]

The different behaviors of disjoint unions of cycles and cyclic loop conditions when constructing square-free representatives is explained by the following observations. 
Let $C\subset\N$ finite and 
\[S=\{\{ p\in (C\dotdiv c)\mid p\text{ is prime}\}\mid \text{$c$ is maximal for $C$}\}.\]  
Let $S_{\min}$ ($S_{\max}$) be the set of minimal (maximal) elements of $S$ with respect to inclusion. Then  
\begin{itemize}
    \item $\C$ has the same pp-constructability type as $\bigtimes_{P\in S_{\min}}\mathbb P$,
    \item $\Sigma_C\Leftrightarrow\{\Sigma_P\mid P\in S_{\max}\}\Leftrightarrow\Sigma_{\Flat(C)}$, and
    \item if $C$ contains only square-free numbers, then $S=S_{\min}=S_{\max}$ and $S$ is an antichain.
\end{itemize}
Compare this to  Example~\ref{exa:620Structure}, where $C=\{6,20\}$. In this case we have that $S=\{\{2\},\{3,2\},\{3,5\}\}$, $S_{\min}=\{\{3,2\},\{3,5\}\}$, and $S_{\max}=\{\{2\},\{3,5\}\}$.

\section[The Lattices of Unions of Cycles]{The Lattices of Unions of Cycles and Cyclic Loop Conditions}\label{sec:lattice}

The characterizations from Corollary~\ref{cor:characterizationMCLAndmCL} and  Theorem~\ref{thm:classificationPoset} suggest that the posets $\MPCL$ and $\SDPoset$ can be described lattice-theoretically. 
Let $X$ be a finite set. The \emph{free distributive lattice generated by $X$}, denoted by $\mathcal{F}_D(X)$, is the poset of all downsets of $2^X$. The fact that $\mathcal{F}_D(X)$ is a free distributive lattice can be found in Theorem~13 in~\cite{birkhoff1940lattice}.   
\begin{figure}
    \centering
    \begin{tikzpicture}[scale=0.54]
\def\myScale{0.8}

\node[scale=\myScale] (0) at (2,-2)  {$[\Cyc{2,3,5}]$};
\node[scale=\myScale] (00) at (2,0)  {$[\Cyc{6,10,15}]$};
\node[scale=\myScale] (11) at (0,2)  {$[\Cyc{2,15}]$};
\node[scale=\myScale] (12) at (2,2) {$[\Cyc{3,10}]$};
\node[scale=\myScale] (13) at (4,2) {$[\Cyc{5,6}]$};
\node[scale=\myScale] (21) at (0,4) {$[\Cyc{10,15}]$};
\node[scale=\myScale] (22) at (2,4) {$[\Cyc{6,15}]$};
\node[scale=\myScale] (23) at (4,4) {$[\Cyc{6,10}]$};
\node[scale=\myScale] (31) at (-2,6) {$[\Cyc{2,3}]$};
\node[scale=\myScale] (33) at (4,6) {$[\Cyc{2,5}]$};
\node[scale=\myScale] (34) at (6,6) {$[\Cyc{3,5}]$};
\node[scale=\myScale] (32) at (2,6) {$[\Cyc{30}]$};
\node[scale=\myScale] (41) at (0,8) {$[\Cyc{6}]$};
\node[scale=\myScale] (42) at (2,8) {$[\Cyc{10}]$};
\node[scale=\myScale] (43) at (4,8) {$[\Cyc{15}]$};
\node[scale=\myScale] (51) at (0,10) {$[\Cyc{2}]$};
\node[scale=\myScale] (52) at (2,10) {$[\Cyc{3}]$};
\node[scale=\myScale] (53) at (4,10) {$[\Cyc{5}]$};
\node[scale=\myScale] (60) at (2,12) {$[\Cyc{1}]$};
\path 
    (0)  edge (00)
    (00) edge (11)
    (00) edge (12)
    (00) edge (13)
    (11) edge (21)
    (11) edge (22)
    (12) edge (21)
    (12) edge (23)
    (13) edge (22)
    (13) edge (23)
    (21) edge (31)
    (21) edge (32)
    (22) edge (32)
    (22) edge (33)
    (23) edge (32)
    (23) edge (34)
    (31) edge (41)
    (32) edge (41)
    (32) edge (42)
    (32) edge (43)
    (33) edge (42)
    (34) edge (43)
    (41) edge (51)
    (41) edge (52)
    (42) edge (51)
    (42) edge (53)
    (43) edge (52)
    (43) edge (53)
    (51) edge (60)
    (52) edge (60)
    (53) edge (60)
    ;
\end{tikzpicture}
\hspace{1cm}
\begin{tikzpicture}[scale=0.54]
\def\myScale{0.4}

\node at (2,-2.2)  {};

\node[circle, draw,scale=\myScale,label={right:$2^{\{2,3,5\}}$}] (0) at (2,-2)  {};
\node[circle, draw,scale=\myScale,label={right:$2^{\{2,3,5\}}\setminus\{2,3,5\}$}] (00) at (2,0)  {};
\node[circle, draw,scale=\myScale] (11) at (0,2)  {};
\node[circle, draw,scale=\myScale] (12) at (2,2) {};
\node[circle, draw,scale=\myScale] (13) at (4,2) {};
\node[circle, draw,scale=\myScale] (21) at (0,4) {};
\node[circle, draw,scale=\myScale] (22) at (2,4) {};
\node[circle, draw,scale=\myScale] (23) at (4,4) {};
\node[circle, fill,scale=\myScale] (31) at (-2,6) {};
\node[circle, fill,scale=\myScale] (33) at (4,6) {};
\node[circle, fill,scale=\myScale,label={right:$2^{\{3,5\}}$}] (34) at (6,6) {};
\node[circle, draw,scale=\myScale] (32) at (2,6) {};
\node[circle, draw,scale=\myScale] (41) at (0,8) {};
\node[circle, draw,scale=\myScale] (42) at (2,8) {};
\node[circle, draw,scale=\myScale,label={right:$\{\{3\},\{5\},\emptyset\}$}] (43) at (4,8) {};
\node[circle, draw,scale=\myScale] (51) at (0,10) {};
\node[circle, draw,scale=\myScale] (52) at (2,10) {};
\node[circle, draw,scale=\myScale,label={right:$\{\{5\},\emptyset\}$}] (53) at (4,10) {};
\node[circle, draw,scale=\myScale,label={right:$\{\emptyset\}$}] (60) at (2,12) {};
\node[circle, draw,scale=\myScale,label={right:$\emptyset$}] (70) at (2,14)  {};

\path 
    (0)  edge (00)
    (00) edge (11)
    (00) edge (12)
    (00) edge (13)
    (11) edge (21)
    (11) edge (22)
    (12) edge (21)
    (12) edge (23)
    (13) edge (22)
    (13) edge (23)
    (21) edge (31)
    (21) edge (32)
    (22) edge (32)
    (22) edge (33)
    (23) edge (32)
    (23) edge (34)
    (31) edge (41)
    (32) edge (41)
    (32) edge (42)
    (32) edge (43)
    (33) edge (42)
    (34) edge (43)
    (41) edge (51)
    (41) edge (52)
    (42) edge (51)
    (42) edge (53)
    (43) edge (52)
    (43) edge (53)
    (51) edge (60)
    (52) edge (60)
    (53) edge (60)
    (60) edge (70)
    ;
\end{tikzpicture}
\caption{The poset $\SDPoset$ restricted to disjoint unions of cycles using 2, 3, and 5 (left). The poset $\mathcal{F}_D(\{2,3,5\})$ \textbf{upside down}  (right).} 
    \label{fig:freedist3}
\end{figure}
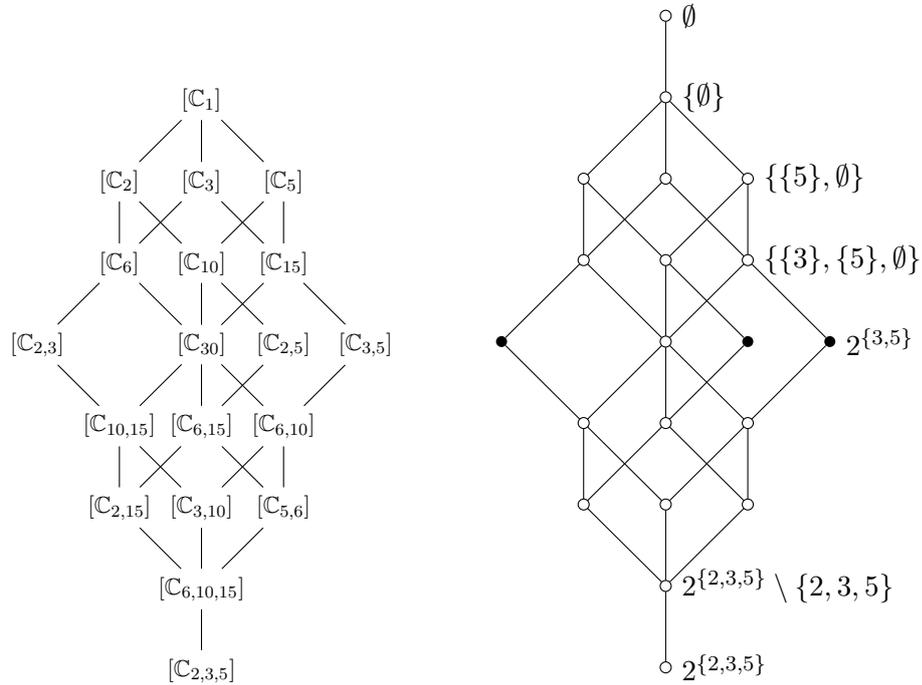
Observe that the poset in Figure~\ref{fig:freedist3} (left) is isomorphic to the free distributive lattice generated by $\{2,3,5\}$ after removing the top element. 
More generally, whenever we restrict $\SDPoset$ to disjoint unions of cycles using only a fixed finite set of $n$ primes, then the resulting poset is isomorphic to the free distributive lattice generated by an $n$-element set, $\mathcal F_D(n)$, (again, after removing the top element).
Note that if we add the $(n+1)$-st prime to each of the unions of cycles then we have $\C_1$ as top element. Hence, we obtain the following corollary. 
\begin{corollary}\label{cor:DGPosetEmbedsAllFreeDistLattices}
The free distributive lattice on $n$ generators (as a poset) embeds into $\DGPoset$.
\end{corollary}
Using this corollary we 
get the following result.
\begin{theorem}\label{thm:DGPosetEmbedsAllFinitePosets}
Every finite poset embeds into $\DGPoset$.
\end{theorem}
\begin{proof}
Let $(P,\preceq)$ be a finite poset. Then the map
\begin{align*}
    P&\to \mathcal F_D(P)\\
    p&\mapsto \{Q\subseteq P \mid Q\subseteq \{q\in P\mid q\preceq p\}\}
\end{align*}
is an embedding. Hence, by Corollary~\ref{cor:DGPosetEmbedsAllFreeDistLattices}, $(P,\preceq)$ also embeds into the poset $\DGPoset$. 
\end{proof}

Consider the power set $2^X$ of a countably infinite set $X$ ordered by inclusion. We denote by $\FD$ the poset of all downsets of $2^X$ ordered by inclusion.
Markowsky proved that $\FD$ is the free completely distributive lattice, i.e., complete and distributive over infinite meets and joins, on countably many generators~\cite{Markowsky}. 
The generating set consists of the principal downsets generated by $X\setminus\{x\}$ for $x\in X$. If we choose $X$ to be the set of all primes, then the following corollary is an immediate consequence of Corollaries~\ref{cor:characterizationMCLAndmCL} and~\ref{thm:classificationPoset}. First, for every finite set $C\subset\N^+$ define $\operatorname{Drop}(\Sigma_C)\coloneqq C$.
\begin{corollary}\label{thm:gorgeous}
The following holds
\begin{align*}
    \SDPoset&\hookrightarrow\MPCL\hookrightarrow\FD.
\end{align*}
One choice for the second embedding $\iota$ is to map an element of $\MPCL$, represented by a set of prime cyclic loop conditions $\Gamma$, to the downset of $\operatorname{Drop}(\Gamma)$.

\end{corollary}
Note that, since $\operatorname{Drop}(\Gamma)$ only contains finite sets it is necessary to take the downset to obtain an element of $\FD$ whose nonempty elements also contain infinite sets.
Furthermore, the image of $\iota$ is closed under finite meets and joins. Hence, the subposet $\iota(\MPCL)$ is also a sublattice of $\FD$.  Therefore, we obtain the following corollary.
\begin{corollary}\label{cor:cyclesFormALattice}
The posets $\SDPoset$ and $\MPCL$ are distributive lattices.
\end{corollary}

One could wonder whether $\SDPoset$ also distributes over infinite meets and joins. This is not the case as shown by this counterexample, provided (personal communication) by Friedrich Martin Schneider,
\[[\Cyc2]\vee\bigwedge_{\substack{p\text{ an odd}\\\text{prime}}}[\Cyc p]=[\Cyc2]\neq[\Cyc1]=\bigwedge_{\substack{p\text{ an odd}\\\text{prime}}} (\Cyc2\vee\Cyc p).\]

Here we summarize the results from
Corollaries~\ref{cor:characterizationMCLAndmCL} and~\ref{cor:mPclEqPPC} and Theorem~\ref{thm:classificationPoset} about the posets that have been studied in this article.
\begin{itemize}
    \item $\mPCL\simeq \PCPoset\simeq (\{P\mid P\text{ a finite nonempty set of primes}\},\supseteq)$
    \item $\mCL\simeq(\operatorname{Finitely Generated Downsets Of}(\mPCL),\subseteq)$.
    \item $\MCL=\MPCL\simeq(\operatorname{ Downsets Of}(\mPCL),\subseteq)$.
    \item $\SDPoset\simeq(\operatorname{Cofinite Downsets Of}(\mPCL),\subseteq)$.
\end{itemize}

Note that, $\mCL$ is isomorphic to a sublattice of $\MPCL$. Hence, $\mCL$ is a distributive lattice, as well. The meet and join of $\mCL$ can be described in the following way.

\begin{lemma}\label{lem:meetAndJoinInMCL}
Let $[\Sigma_C], [\Sigma_D]\in\mCL$. Then
\begin{enumerate}
    \item $[\Sigma_C]\wedge[\Sigma_D]=[\Sigma_{C\cup D}]$ and
    \item $[\Sigma_C]\vee[\Sigma_D]=[\Sigma_{C\cdot D}]=[\Sigma_{\C\mathbin\times \D}]$.
\end{enumerate} 
\end{lemma}
\begin{proof}
Let $\Sigma_E$ and $\Sigma_{E'}$ be cyclic loop conditions.
Assume without loss of generality that the sets $C$, $D$, $E$, and $E'$ contain only square-free numbers.

Clearly, $\C\to\C\cup\D$, $\D\to \C\cup\D$, $\C\mathbin\times\D\to\C$, and $\C\mathbin\times\D\to\D$. Suppose that $\Sigma_C\Rightarrow\Sigma_E$, $\Sigma_D\Rightarrow\Sigma_E$, 
$\Sigma_{E'}\Rightarrow\Sigma_{C}$, and $\Sigma_{E'}\Rightarrow\Sigma_{D}$. Then, by Theorem~\ref{thm:implicationSingleClc}, there are homomorphisms $f,g,f',$ and $g'$ as in Figure~\ref{fig:meetAndJoin}. 
\begin{figure}
    \centering
    \begin{tikzpicture}
        \node (C) at (-1,0) {$\C$};
        \node (D) at (1,0) {$\D$};
        \node (CuD) at (0,-1) {$\C\cup\D$};
        \node (E) at (0,-2.414) {$\mathbb E$};
        \node (CxD) at (0,1) {$\C\mathbin\times\D$};
        \node (Ep) at (0,2.414) {$\mathbb E'$};
        
        \path[->,>=stealth']
            (C) edge (CuD)
            (D) edge (CuD)
            (CxD) edge (C)
            (CxD) edge (D)
            (Ep) edge node[right] {$h'$} (CxD)
            (Ep) edge[bend right] node[left] {$f'$} (C)
            (Ep) edge[bend left] node[right] {$g'$} (D)
            (CuD) edge node[right] {$h$}(E)
            (C) edge[bend right] node[left] {$f$}(E)
            (D) edge[bend left] node[right] {$g$}(E)
            ;
        \node at (3,0) {implies};
        \node (C) at (-1+6,0) {$\Sigma_C$};
        \node (D) at (1+6,0) {$\Sigma_D$};
        \node (CuD) at (0+6,-1) {$\Sigma_{C\cup D}$};
        \node (E) at (0+6,-2.414) {$\Sigma_E$};
        \node (CxD) at (0+6,1) {$\Sigma_{C\cdot D}$};
        \node (Ep) at (0+6,2.414) {$\Sigma_{E'}$};
        
        \node[rotate=-90] at (0+6,-1.707) {$\Rightarrow$};
        \node[rotate=-90] at (0+6,1.707) {$\Rightarrow$};
        \node[rotate=-70] at (-0.7+6,-1.5) {$\Rightarrow$};
        \node[rotate=-110] at (0.7+6,-1.5) {$\Rightarrow$};
        \node[rotate=-45] at (-0.5+6,-0.5) {$\Rightarrow$};
        \node[rotate=-135] at (0.5+6,-0.5) {$\Rightarrow$};
        \node[rotate=-110] at (-0.7+6,1.5) {$\Rightarrow$};
        \node[rotate=-70] at (0.7+6,1.5) {$\Rightarrow$};
        \node[rotate=-135] at (-0.5+6,0.5) {$\Rightarrow$};
        \node[rotate=-45] at (0.5+6,0.5) {$\Rightarrow$};
        
    \end{tikzpicture}
    \caption{Determining the meet and join of $[\Sigma_C]$ and $[\Sigma_D]$.}
    \label{fig:meetAndJoin}
\end{figure}
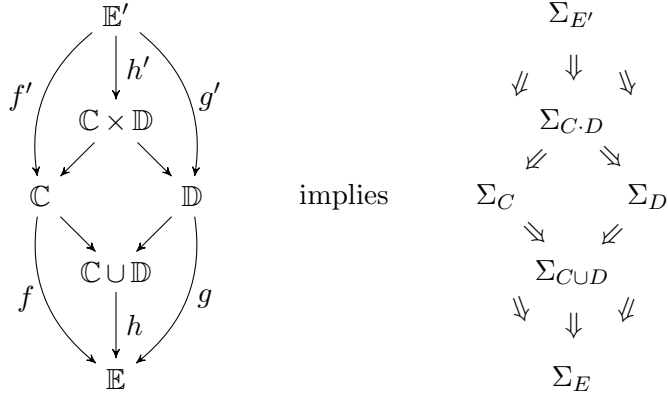
The maps 
\begin{align*}
    h\colon\C\cup\D&\to \mathbb E & h'\colon\mathbb E'&\to\C\mathbin\times\D\\
u&\mapsto
\begin{cases}
f(u)&\text{if }u\in\C\\
g(u)&\text{otherwise}
\end{cases} &
u&\mapsto (f(u),g(u))
\end{align*}
are homomorphism as well. Hence, $[\Sigma_{C\cup D}]$ is the meet and $[\Sigma_{C\mathbin\times D}]$ is the join of $[\Sigma_C]$ and $[\Sigma_D]$. 
\end{proof}

Using Lemmata~\ref{lem:surjectivety} and~\ref{lem:pclCharacterization} we can describe the meet of $\SDPoset$ of disjoint unions of prime cycles.
\begin{corollary}\label{cor:meetToTimesCycles}
For any finite antichain $Q$ of $\PCPoset$ we have that
\[\bigwedge_{\mathbb P\in Q} [\mathbb P]=\left[\bigtimes_{\mathbb P\in Q} \mathbb P\right].\]
\end{corollary}

Note that this connection between $\times$ and $\wedge$ does not hold in general as seen in the following example.
\begin{examplex}
Use Figure~\ref{fig:freedist3} to verify that we have
\begin{itemize}
    \item $[\Cyc{30}\mathbin\times\Cyc{2,3}]=[\Cyc{30}]\neq[\Cyc{10,15}]=[\Cyc{30}]\wedge[\Cyc{2,3}]$ and
    \item $[\Cyc{10,15}\mathbin\times\Cyc{6,10}]=[\Cyc{10}]\neq[\Cyc{3,10}]=[\Cyc{10,15}]\wedge[\Cyc{6,10}]$.\eoe
\end{itemize}
\end{examplex}

Note that Corollary~\ref{cor:meetToTimesCycles} and Lemma~\ref{lem:surjectivety} imply that the following conjecture is true in the case that $\A$ and $\B$ are disjoint unions of cycles.
\begin{conjecture}
Let $\A$ and $\B$ be finite relational structures. Then there are $\A_1,\dots,\A_n\ppleq\A$ and $\B_1,\dots,\B_m\ppleq\B$ such that 
\[[\A]\wedge[\B]=[\A_1\mathbin{\times}\dots\mathbin{\times}\A_n\mathbin{\times}\B_1\mathbin{\times}\dots\mathbin{\times}\B_m].\]
\end{conjecture}

Observe that a prime cyclic loop condition $[\Sigma_P]\in\mCL$ corresponds to a downset of $\mPCL$, which is generated by a single element, namely $\Sigma_P$. Analogously, a disjoint union of prime cycles $[\mathbb P]\in\SDPoset$ corresponds to the complement of an upset of $\mPCL$, which is generated by a single element, namely $\Sigma_P$. 

\begin{corollary}\label{cor:joinAndMeetIrred}
The join irreducible elements of $\mCL\setminus\{[\Sigma_1]\}$ are exactly the elements that can be represented by a prime cyclic loop condition.
The meet irreducible elements of $\SDPoset\setminus\{[\Cyc 1]\}$ are exactly the elements that can be represented by a finite disjoint union of prime cycles.
\end{corollary}

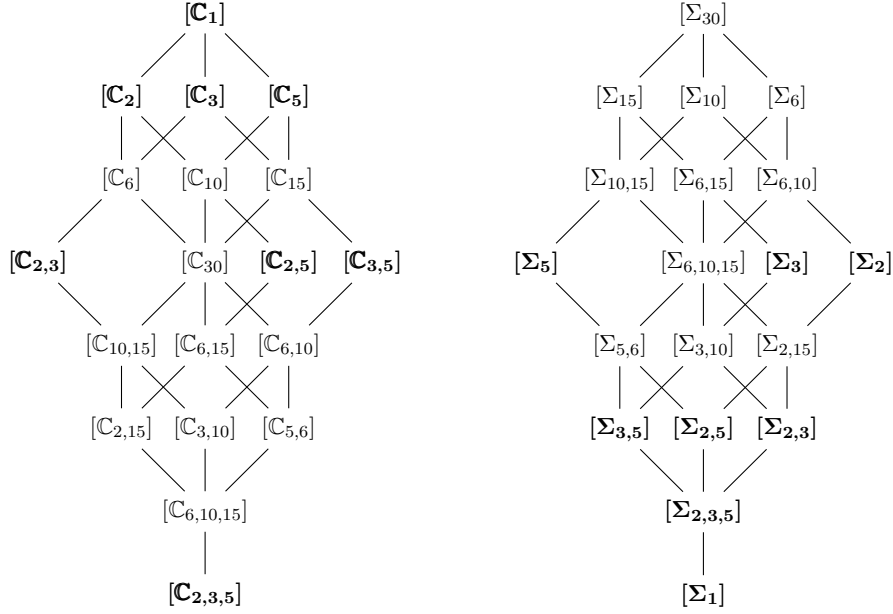
\begin{figure}
    \centering
    \begin{tikzpicture}[scale=0.55]
\def\myScale{0.8}

\node[scale=\myScale] (0) at (2,-2)  {$\boldsymbol{[\pmb{\C}_{2,3,5}]}$};
\node[scale=\myScale] (00) at (2,0)  {$[\Cyc{6,10,15}]$};
\node[scale=\myScale] (11) at (0,2)  {$[\Cyc{2,15}]$};
\node[scale=\myScale] (12) at (2,2) {$[\Cyc{3,10}]$};
\node[scale=\myScale] (13) at (4,2) {$[\Cyc{5,6}]$};
\node[scale=\myScale] (21) at (0,4) {$[\Cyc{10,15}]$};
\node[scale=\myScale] (22) at (2,4) {$[\Cyc{6,15}]$};
\node[scale=\myScale] (23) at (4,4) {$[\Cyc{6,10}]$};
\node[scale=\myScale] (31) at (-2,6) {$\boldsymbol{[\pmb{\C}_{2,3}]}$};
\node[scale=\myScale] (33) at (4,6) {$\boldsymbol{[\pmb{\C}_{2,5}]}$};
\node[scale=\myScale] (34) at (6,6) {$\boldsymbol{[\pmb{\C}_{3,5}]}$};
\node[scale=\myScale] (32) at (2,6) {$[\Cyc{30}]$};
\node[scale=\myScale] (41) at (0,8) {$[\Cyc{6}]$};
\node[scale=\myScale] (42) at (2,8) {$[\Cyc{10}]$};
\node[scale=\myScale] (43) at (4,8) {$[\Cyc{15}]$};
\node[scale=\myScale] (51) at (0,10) {$\boldsymbol{[\pmb{\C}_{2}]}$};
\node[scale=\myScale] (52) at (2,10) {$\boldsymbol{[\pmb{\C}_{3}]}$};
\node[scale=\myScale] (53) at (4,10) {$\boldsymbol{[\pmb{\C}_{5}]}$};
\node[scale=\myScale] (60) at (2,12) {$\boldsymbol{[\pmb{\C}_{1}]}$};
\path 
    (0)  edge (00)
    (00) edge (11)
    (00) edge (12)
    (00) edge (13)
    (11) edge (21)
    (11) edge (22)
    (12) edge (21)
    (12) edge (23)
    (13) edge (22)
    (13) edge (23)
    (21) edge (31)
    (21) edge (32)
    (22) edge (32)
    (22) edge (33)
    (23) edge (32)
    (23) edge (34)
    (31) edge (41)
    (32) edge (41)
    (32) edge (42)
    (32) edge (43)
    (33) edge (42)
    (34) edge (43)
    (41) edge (51)
    (41) edge (52)
    (42) edge (51)
    (42) edge (53)
    (43) edge (52)
    (43) edge (53)
    (51) edge (60)
    (52) edge (60)
    (53) edge (60)
    ;
\end{tikzpicture}
\hspace{1cm}
    \begin{tikzpicture}[scale=0.55]
\def\myScale{0.8}

\node[scale=\myScale] (0) at (2,-2)  {$\boldsymbol{[\Sigma_{1}]}$};
\node[scale=\myScale] (00) at (2,0)  {$\boldsymbol{[\Sigma_{2,3,5}]}$};
\node[scale=\myScale] (11) at (0,2)  {$\boldsymbol{[\Sigma_{3,5}]}$};
\node[scale=\myScale] (12) at (2,2) {$\boldsymbol{[\Sigma_{2,5}]}$};
\node[scale=\myScale] (13) at (4,2) {$\boldsymbol{[\Sigma_{2,3}]}$};
\node[scale=\myScale] (21) at (0,4) {$[\Sigma_{5,6}]$};
\node[scale=\myScale] (22) at (2,4) {$[\Sigma_{3,10}]$};
\node[scale=\myScale] (23) at (4,4) {$[\Sigma_{2,15}]$};
\node[scale=\myScale] (31) at (-2,6) {$\boldsymbol{[\Sigma_{5}]}$};
\node[scale=\myScale] (33) at (4,6) {$\boldsymbol{[\Sigma_{3}]}$};
\node[scale=\myScale] (34) at (6,6) {$\boldsymbol{[\Sigma_{2}]}$};
\node[scale=\myScale] (32) at (2,6) {$[\Sigma_{6,10,15}]$};
\node[scale=\myScale] (41) at (0,8) {$[\Sigma_{10,15}]$};
\node[scale=\myScale] (42) at (2,8) {$[\Sigma_{6,15}]$};
\node[scale=\myScale] (43) at (4,8) {$[\Sigma_{6,10}]$};
\node[scale=\myScale] (51) at (0,10) {$[\Sigma_{15}]$};
\node[scale=\myScale] (52) at (2,10) {$[\Sigma_{10}]$};
\node[scale=\myScale] (53) at (4,10) {$[\Sigma_{6}]$};
\node[scale=\myScale] (60) at (2,12) {$[\Sigma_{30}]$};
\path 
    (0)  edge (00)
    (00) edge (11)
    (00) edge (12)
    (00) edge (13)
    (11) edge (21)
    (11) edge (22)
    (12) edge (21)
    (12) edge (23)
    (13) edge (22)
    (13) edge (23)
    (21) edge (31)
    (21) edge (32)
    (22) edge (32)
    (22) edge (33)
    (23) edge (32)
    (23) edge (34)
    (31) edge (41)
    (32) edge (41)
    (32) edge (42)
    (32) edge (43)
    (33) edge (42)
    (34) edge (43)
    (41) edge (51)
    (41) edge (52)
    (42) edge (51)
    (42) edge (53)
    (43) edge (52)
    (43) edge (53)
    (51) edge (60)
    (52) edge (60)
    (53) edge (60)
    ;
\end{tikzpicture}
\caption{Poset $\SDPoset$ restricted to disjoint unions of cycles using 2, 3, and 5; meet irreducible elements in bold (left). Poset $\mCL$ restricted to cyclic loop conditions using 2, 3, and 5; join irreducible elements in bold (right).}
    \label{fig:SDPvsmCL}
\end{figure}

Let $P$ be a finite set of primes and $\N^+_P$ be the set of all positive natural numbers whose prime decomposition only uses numbers from $P$. 
Observe that
the subposet of $\SDPoset$ consisting of $\{[\C]\mid C\subset \N_P^+ \text{ finite}\}$ and the subposet of $\mCL$ consisting of $\{[\Sigma_C]\mid C\subset \N_P^+  \text{ finite}\}$ are isomorphic. 
The map 
\[[\C]\mapsto[\{\Sigma_D\mid D\subset\N_P^+\text{ finite, }\Pol(\C)\models\Sigma_D\}]\]
is an isomorphism.
See for example Figure~\ref{fig:SDPvsmCL}. 
However, the posets $\mCL$ and $\SDPoset$ are not isomorphic since $\mCL$ has infinite ascending chains, e.g. $([\Sigma_{p_1\cdot\ldots\cdot p_n}])_{n\in\N}$, and $\SDPoset$ does not.

So far in this chapter we have described the poset $\SDPoset$, i.e., the subposet of $\DGPoset$ where every element is a pp-constructability type of some finite disjoint union of cycles. From the provided description it follows that $\SDPoset$ is a distributive lattice and that it contains infinite descending chains and infinite antichains but no infinite ascending chains. Some of these properties are inherited by $\PPPoset$.
For instance, it follows that $\PPPoset$ contains  infinite antichains and infinite descending chains; the latter was already known from the description of $\mathfrak{P}_{\operatorname{Boole}}$~\cite{PPPoset}.

\section{Digraphs with a Maltsev Polymorphism}
In this section we consider all digraphs that satisfy a certain minor condition. 
\begin{definition}\label{def:malt}
The \emph{quasi Maltsev} condition $\Maltsev$ is the minor condition
\[f(y,y,x)\approx f(x,x,x)\approx f(x,y,y).\]
\end{definition}


 
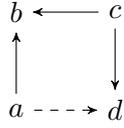
\begin{figure}
    \centering
    \begin{tikzpicture}[scale=1.3]
    \node (a) at (0,0) {$a$};
    \node (b) at (0,1) {$b$};
    \node (c) at (1,1) {$c$};
    \node (d) at (1,0) {$d$};
    \path[->,>=stealth']
        (a) edge (b)
        (c) edge (b)
        (c) edge (d)
        (a) edge[dashed] (d)
        ;
    \end{tikzpicture}
    \caption{Rectangularity in digraphs.}
    \label{fig:rect}
\end{figure}

Define $\mathfrak{P}_{\Maltsev}$ to be the subposet of $\DGPoset$ that consists of all finite digraphs that do have a quasi Maltsev polymorphism.
A finite digraph $H$ is called 
\emph{$k$-rectangular} if whenever $H$ contains directed paths of length $k$ from $a$ to $b$, from $c$ to $b$, 
and from $c$ to $d$, then also from $a$ to $d$; see Figure~\ref{fig:rect}. 
A digraph $H$ is called \emph{totally rectangular} if it is $k$-rectangular for all $k \geq 1$. 
For $k\geq0$ define $\P_k$, the \emph{directed path of length $k$}, as 
\begin{align*}
    &(\{0,\dots,k\},\{(0,1),\dots,(k-1,k)\}).
\end{align*}
\begin{example}\label{exa:cyclesHaveMaltsev}
Let $a\in\mathbb N^+$\!. 
Then the ternary operation 
\begin{align*}
\Z_a&\to\Z_a\\    
(x_1,x_2,x_3) &\mapsto x_1 - x_2 + x_3 \pmod a
\end{align*}
is an idempotent Maltsev polymorphism of $\C_a$. Hence, $\Pol(\C_a) \models \Maltsev$. Note that $\C_a$ is a core and that it is totally rectangular.
\end{example}
Similarly, we can show that any disjoint union of cycles has an idempotent Maltsev polymorphism.
Hence, $\SDPoset$ is a subposet of $\mathfrak{P}_{\Maltsev}$.
To determine the digraphs in $\mathfrak{P}_{\Maltsev}$ that are not in $\SDPoset$ we use 
the following well-known 
statements that hold for any finite digraph $\G$.


\begin{enumerate}
    \item $\G$ has a Maltsev polymorphism if and only if the core of $\G$ has an idempotent Maltsev polymorphism.
    \item $\G$ is totally rectangular if and only if
it has an idempotent Maltsev polymorphism. 
    \item If $\G$ is totally rectangular then the core of $\G$ is either a directed path or a disjoint union of directed cycles (Lemma 3.10 in~\cite{CarvalhoEgriJacksonNiven}).
    \item If $\G$ is a directed path or a disjoint union of directed cycles then $\G$ is totally rectangular.
\end{enumerate}
The proofs of these statements are left as an exercise to the reader. Combining 1.,2.,3., and 4.~we obtain the following lemma:
\begin{lemma}\label{lem:maltIffRect}
Let $\G$ be a finite digraph. Then $\Pol(\G)\models\Maltsev$ if and only if the core of $\G$ is a directed path or a disjoint union of directed cycles.
\end{lemma}
Our goal is to understand the pp-constructability order on digraphs $\G$ that do have a Maltsev polymorphism. 
We already understand disjoint unions of directed cycles. Recall that $\P_1$ is the unique lower cover of $\C_1$ and that $[\C_1]=[\P_0]$. The following lemma shows that all $\P_k$ collapse for $k\geq 1$.

\begin{lemma}\label{lem:IdempConstructPaths}
For all $k\geq 1$ we have that $\P_1\ppeq\P_k$.
\end{lemma}
\begin{proof}
Let $k\geq 1$. Then $\P_{k+1}$ can pp-construct $\P_k$ with a 1-dimensional pp-construction using the formula $\exists z.\ E(x,y)\AND E(y,z)$ which is 
\begin{center}
    
\begin{tikzpicture}[scale =0.5]
\node[var-b] (1) at (0,2) {};
\node[var-f,label=left:{$x$}] (x) at (0,0) {};
\node[var-f,label=left:{$y$}] (y) at (0,1) {};
\path[>=stealth',->] 
        (x) edge (y)
        (y) edge (1)
        ;
\end{tikzpicture}
\end{center}
in graph notation. Let $\G$ be the k-dimensional pp-power of $\P_1$ given by the $2k$-ary pp-formula $E(x_1,y_k)\AND (x_2\approx y_1)\AND\dots\AND (x_k\approx y_{k-1})$ which is
\begin{center}
\begin{tikzpicture}[scale=0.5]
\node[var-f,label=below:{$x_1$}] (x1) at (0,0) {};
\node[var-f,label=below:{$x_2$}] (x2) at (1,0) {};
\node[var-f,label=below:{$x_3$}] (x3) at (2,0) {};
\node (x4) at (3,0) {$\dots$};
\node[var-f,label=below:{$x_k$}] (x5) at (4,0) {};

\node[var-f,label=above:{$y_1$}] (y1) at (0,1) {};
\node[var-f,label=above:{$y_2$}] (y2) at (1,1) {};
\node[var-f,label=above:{$y_3$}] (y3) at (2,1) {};
\node (y4) at (3,1) {$\dots$};
\node[var-f,label=above:{$y_k$}] (y5) at (4,1) {};
\path[>=stealth'] 
        (x1) edge[->] (y5)
        (x2) edge[dashed] (y1)
        (x3) edge[dashed] (y2)
        (x4) edge[dashed] (y3)
        (x5) edge[dashed] (y4)
        ;
\end{tikzpicture}
\end{center}
in graph notation. 
Then $\G$ contains the following path of length $k$:
\begin{center}
    \begin{tikzpicture}
    \node (0) at (0,0) {$(0,0,\dots,0,0)$};
    \node (1) at (0,1) {$(0,0,\dots,0,1)$};
    \node (2) at (0,2) {$\vdots$};
    \node (3) at (0,3) {$(0,1,\dots,1,1)$};
    \node (4) at (0,4) {$(1,1,\dots,1,1)$};
    \path[->,>=stealth']
        (0) edge (1)
        (1) edge (2)
        (2) edge (3)
        (3) edge (4)
        ;
    \end{tikzpicture}
\end{center}
and $\P_{k}\rightarrow \G$. Note that for every edge $(u,v)$ in $\G$ the tuple $v$ contains exactly one 1 more than $u$. Hence, the map $(x_1,\dots,x_k)\mapsto x_1+\dots+x_k$ is a homomorphism from $\G$ to $\P_{k}$. Therefore, $\G\homeq\P_{k}$ and $\P_1\ppleq \P_{k}$.
\end{proof}

We conclude that in $\mathfrak{P}_{\Maltsev}$ the element $[\P_1]$ is the only one that cannot be represented by a disjoint union of cycles.
Note that the $n$-ary $\max$ function is a totally symmetric polymorphism of $\P_1$. Hence, $\P_1$ satisfies all cyclic loop conditions. Additionally,  $[\P_1]\neq[\C_1]$. Therefore, we have that $[\P_1]$ cannot be represented by a disjoint union of cycles. Note that for any digraph $\G\in\DGPoset$ such that $\H\ppleq\G$ for some digraph $\H\in\mathfrak{P}_{\Maltsev}$ we have that $\G$ must have a Maltsev polymorphism and therefore $\G\in\mathfrak{P}_{\Maltsev}$. From these observations we obtain the following corollary.

\begin{corollary}\label{cor:classificationMaltsevPoset}
The poset $\mathfrak{P}_{\Maltsev}$ is isomorphic to the poset $\SDPoset$ with a single element $[\P_1]$ inserted as unique lower cover of $[\C_1]$. The covering relation of $\mathfrak{P}_{\Maltsev}$ is the covering relation of $\DGPoset$ restricted to $\mathfrak{P}_{\Maltsev}$.
\end{corollary}

We conclude that $\C_2,\C_3,\C_5,\dots$ are lower covers of $[\P_1]$ in $\DGPoset$. A natural next question to ask is: are there other lower covers of $[\P_1]$? The answer is yes and will be discussed in Chapter~\ref{cha:submax}.



\section[Concluding Remarks]{Concluding Remarks and Directions for Further Research}








This thesis does not focus on \emph{conservative} CSPs, i.e., CSPs where the polymorphisms of the template structure preserve all possible unary relations. However, we will show in the following that a directed cycle $\C$ has the same pp-constructibility type as $\C$ expanded by all unary relations. For the proof we need the following notion. Let $f\colon A^I\to A$ be a function and let $i\in I$. Then the index $i$ is \emph{essential} for $f$ if there exist tuples $a,a'\in A^I$ that only differ in the $i$-th entry such that $f(a)\neq f(a')$. Otherwise, $i$ is \emph{inessential} for $f$. Note that $f(a)$ only depends on the entries of $a$ at essential indices.

\begin{lemma}
    Let $\C$ be a directed cycle. Then $\C$ can pp-construct the structure $\C$ expanded by all unary relations.
\end{lemma}
The proof of this lemma in the original version of this thesis is incorrect. 
\begin{proof}
    Without loss of generality $\C=\C_n$. 
    Then $\sigma=(0\,\dots\,n-1)$ is the permutation of $\{0,\dots,n-1\}$ with $\C_n=(\{0,\dots,n-1\},\{(a,\sigma(a)\mid a\in C\})$. 
    Let $\tilde\C_n$ be the structure $\C_n$ expanded by all unary relations. Clearly, $\tilde\C_n\ppleq\C_n$. By Theorem~\ref{thm:freestructure} it suffices to show $\F_{\Pol(\C_n)}(\tilde\C_n)\to\tilde\C_n$ in order to proof $\C_n\ppleq\tilde\C_n$. 
    We construct the homomorphism from $\F_{\Pol(\C_n)}(\tilde\C_n)$ to $\tilde\C_n$ componentwise. Let $f$ be in the connected component $\D$ of $\F_{\Pol(\C_n)}(\tilde\C_n)$. Note that all elements of $\D$ are of the form $f_{\sigma^k}$ for some $k$. Since $\C_n\not\models\Sigma_k$ for any divisor $k>1$ of $n$, we have that the smallest $k\geq 1$ with $f=f_{\sigma^k}$ is $n$. Hence, $\D$ is a directed cycle of length $n$ with some unary relations. 
    Let $U$ be a unary relation of $\tilde\C_n$ and let $U_F$ be the corresponding unary relation of $\F_{\Pol(\C_n)}(\tilde\C_n)$.
    Note that an $n$-ary polymorphism $f$ of $\C_n$ is in $U_F$ if there is a $g\in\Pol(\C_n)$ such that $f(x_0,\dots,x_{n-1})=g(x_{i_1},\dots,x_{i_{\ell}})$ with $\{x_{i_1},\dots,x_{i_{\ell}}\}= U$, i.e., if all essential indices of $f$ are contained in $U$. In particular we have for any essential index $a$ of $f$ that if $f\in U_F$ then $a\in U$.
    
    Let $a\in C$ such that the index $a$ is essential for $f$. Define the map $\alpha\colon \D\to \tilde\C, f_{\sigma^k}\mapsto a+k$. Since the edge relation of $\D$ is a cycle of length $n$ we have that $\alpha$ is well defined and preserves the edge relation. Note that for all $k$ we have that $f_{\sigma^k}$ is mapped to an element of $\tilde\C$ that is an essential index for $f_{\sigma^k}$. Hence, $\alpha$ preserves all unary relations. 
\end{proof}
This lemma can not be generalized to disjoint unions of cycles as seen in the following example.
\begin{example}
    Let $\tilde\C_{6,14,15}$ be $\C_{6,14,15}$ expanded by all unary relations. Use Lemma~\ref{lem:duofcSatisfyClc} to verify that $\C_{6,14,15}\models\Sigma_{5,7}$ and $\C_{14,15}\not\models\Sigma_{5,7}$. Note that any polymorphism of $\tilde\C$ restricted to $\C_{14,15}$ is a polymorphism of $\C_{14,15}$ (satisfying the same minor identities). Hence $\tilde\C_{6,14,15}\not\models\Sigma_{5,7}$. 
\end{example}

There are also interesting questions about random unions of cycles. For each $n$ we consider a uniform distribution over all unions of cycles with vertices $\{1,\dots,n\}$ and we denote the \emph{random union of cycles with $n$ vertices} by $\mathbb R\!\operatorname{andomUC}_n$.
Note that $\mathbb R\!\operatorname{andomUC}_n$ can also be seen as the random permutation of $[n]$.

\begin{question}
Let $\C$ be a union of cycles. What is the probability that $\mathbb R\!\operatorname{andomUC}_n$ can pp-construct $\C$ as $n$ tends to infinity?
\end{question}

Note  that there are $(k-1)!$ ways to make a cycle of length $k$ using a fixed set of $k$ vertices.
A back of an envelop calculation suggests that 
\[P(\C_k \not\hookrightarrow \R\!\operatorname{andomUC}_n)\approx\left(1-\frac{(k-1)!}{n^k}\right)^{\binom nk}\to e^{-\frac{1}{k}}.\]

For $k=1$ I verified the formula computationally using OEIS A002467. Hence, if the formula is indeed true, then $P(  \R\!\operatorname{andomUC}_n\ppeq\C_1)\to 1-\frac{1}{e}$.
How about the probability that $\R\!\operatorname{andomUC}_n$  pp-constructs $\C_2$?
Note that $P( \R\!\operatorname{andomUC}_n\ppleq \C_2)$ is equal to
\[\sum_{k\in2\cdot\N^+ } P(\C_k \hookrightarrow \R\!\operatorname{andomUC}_n, \C_1,\C_2,\C_4,\dots,\C_{k-2} \not\hookrightarrow \R\!\operatorname{andomUC}_n).\]
If the all the events $\C_k\not\hookrightarrow\R\!\operatorname{andomUC}_n$ are independent, then as $n$ goes to infinity this tends to
\[\sum_{k\in2\cdot \N^+ }(1-e^{-\frac{1}{k}})\cdot e^{-1}\cdot e^{-\frac{1}{2}}\cdot \dots\cdot e^{-\frac{1}{k-2}}.\]
Whether this is actually true and can be simplified further is left for future research.


\chapter{Submaximal Digraphs}\label{cha:submax}

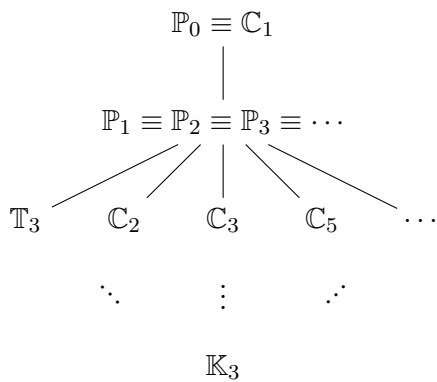
\begin{figure}
    \centering
    \begin{tikzpicture}[scale=1.3]
    \node (0) at (2,2)  {$\P_{0} \equiv \C_{1}$};
    \node (1) at (2,1)  {$\P_{1} \equiv \P_2 \equiv \P_3 \equiv \cdots$};
    \node (20) at (0,0)  {$\T_{3}$};
    \node (21) at (1,0)  {$\C_{2}$};
    \node (22) at (2,0)  {$\C_{3}$};
    \node (23) at (3,0)  {$\C_{5}$};
    \node (24) at (4,0)  {$\dots$};

    \node (3) at (2,-0.7)  {$\vdots$};
    \node[rotate = 45] (3) at (0.8,-0.7)  {$\vdots$};
    \node[rotate = -45] (3) at (3.2,-0.7)  {$\vdots$};

    \node (4) at (2,-1.5)  {$\K_3$};
    
    \path
        (0) edge (1)
        (1) edge (20)
        (1) edge (21)
        (1) edge (22)
        (1) edge (23)
        (1) edge (24)
        ;
    
    \end{tikzpicture}
    \caption{The pp constructibility poset on finite digraphs.}
    \label{fig:main}
\end{figure}

Our goal in this chapter is to classify the maximal elements below $\P_1$ in $\DGPoset$, we call such lower covers \emph{submaximal elements}. The results presented in this chapter have been published in \cite{BodirskyStarke2022}.
For $n\in\N$ let $\C_n$ denote the directed cycle of length $n$, $\P_n$ denote the directed path of length $n$, and $\T_n$ denote the strict order on $\{0,\dots,n-1\}$. 
The following theorem is our main result and will be shown in the remainder of this chapter; see Figure~\ref{fig:main}.

\begin{theorem}\label{thm:submaximalGraphs}
The submaximal elements of $\DGPoset$ are precisely $\T_3$, $\C_2$, $\C_3$, $\C_5$, $\dots$.
If $\G$ is a finite digraph that cannot be pp-constructed by $\P_1$, then $\G\leq \T_3$ or $\G\leq\C_p$ for some prime $p$.
\end{theorem}

\section{Submaximal Digraphs and Minor Conditions}
We first discuss which of 
the minor conditions that we have encountered 
are satisfied by the polymorphisms of 
the digraphs that appear in Theorem~\ref{thm:submaximalGraphs}. 
Recall from Lemma~\ref{lem:pcSatisfyPclc}, Example~\ref{exa:cyclesHaveMaltsev},  Lemmata~\ref{lem:maltIffRect} and~\ref{lem:IdempConstructPaths}, and Corollary~\ref{cor:PPvsLoopCaBlockerForSa} from Chapter~\ref{cha:cycles}
\begin{itemize}
    \item that for any primes
$p$ and $q$ we have that
$\Pol(\C_p) \models \Sigma_q$ if and only if $p \neq q$,
\item that $\C_p\models\Maltsev$ for all primes $p$, 
\item that for any finite digraph $\G$ we have that $\G\models\Maltsev$ if and only if the core of $\G$ is a directed path or a disjoint union of directed cycles,
\item that $\P_1\ppeq\P_2\ppeq\P_3\ppeq\dots$, and
\item that for any finite structure $\B$ we have that $\B\leq \C_a$ if and only if  $\Pol(\B)\not\models\Sigma_{a}$.
\end{itemize}



\begin{lemma}\label{lem:T3Conditions}
We have $\Pol(\T_3) \not \models \Maltsev$ and $\Pol(\T_3) \models \Sigma_n$ for every $n \in {\mathbb N}^+$. 
\end{lemma}
\begin{proof}
The operation $(x_1,\dots,x_n) \mapsto \max(x_1,\dots,x_n)$ is a polymorphism of $\T_3$ that satisfies $\Sigma_n$. 
Secondly, $\T_3 = (\{0,1,2\},E)$ is not 1-rectangular, witnessed by $(1,2),(0,2),(0,1) \in E$ and $(1,1) \notin E$. Therefore, $\T_3\not\models\Maltsev$. 
\end{proof}


The following theorem states that 
the digraph $\P_1$ is the unique smallest element of 
$\DGPoset$ that satisfies $\Maltsev$ and $\Sigma_p$ for all primes $p$. 

\begin{theorem}\label{thm:maltAndCyclImplyIdempotent}
Let $\G$ be a finite digraph
that satisfies $\Maltsev$
and $\Sigma_p$ for all primes $p$. Then 
$\P_1 \leq \G$. 
\end{theorem}


\begin{proof}

By Lemma~\ref{lem:maltIffRect} there are two cases to consider:
the first is that $\G$ is homomorphically equivalent to $\P_k$ for some $k$. Then $\P_1\leq \G$ by Lemma~\ref{lem:IdempConstructPaths}. 

The second case is that $\G$ is homomorphically equivalent to  a disjoint union of directed cycles.
Without loss of generality we may assume that $\G$ itself is a disjoint union of directed cycles.   
Let $(a_0,\dots,a_{\ell-1})$ be a shortest cycle in $\G$. Let $p$ be a prime and $k\in\N^+\!$ such that $p\cdot k=\ell$, and
let $f\in\Pol(\G)$ be a function that witnesses that  $\Pol(\G)\models\Sigma_p$. 
Then
\begin{align*}
    f(a_0,a_k,\dots,a_{(p-1)\cdot k})=a=f(a_k,a_{2 k},\dots,a_{0}).
\end{align*}
Since $f$ is a polymorphism of $\G$ there is a directed path of length $k$ from $a$ to $a$. Thus, $\G$ contains a directed cycle whose length divides $k$, which  contradicts the assumption that $\ell$ is the length of a shortest directed cycle in $\G$. Therefore, $\ell$ has no prime divisors, and $\ell=1$. So $\G$ contains a loop and is thus homomorphically equivalent to $\C_1$; it follows that $\P_1 \leq \G$. 
\end{proof}

\section{Proof of the Main Result}



Recall from Corollary~\ref{cor:PPvsLoopCaBlockerForSa} that $\C_p$ is a blocker for $\Sigma_p$.
We show that the same is true for $\T_3$ and $\Maltsev$.

\begin{figure}
    \centering
    \begin{tikzpicture}
    \node (0) at (0,0) {$(c,d)$};
    \node (1) at (0,1.5) {$(a,d)$};
    \node (2) at (0,3) {$(a,b)$};
    
    \path[->,>=stealth']
        (0) edge (1)
        (1) edge (2)
        (0) edge[bend left=42] (2)
        ;
    \end{tikzpicture}
    \caption{The pp-construction of $\T_3$ in the proof of Lemma~\ref{lem:notSatMalt}.}
    \label{fig:T3constr}
\end{figure}
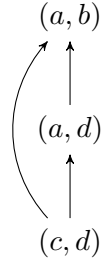

\begin{lemma}\label{lem:notSatMalt}
Let $\G$ be a finite digraph. 
If $\Pol(\G) \not \models \Maltsev$, then $\G \leq \T_3$. 
\end{lemma}
\begin{proof}
Assume without loss of generality that $\G = (V,E)$ is a core. Since $\G\not\models\Maltsev$ we have that $\G$ is not totally rectangular. Hence, there are vertices $a,b,c,d \in V$ such that in $\G$ there are directed paths of length $k$ from $a$ to $b$, from $c$ to $b$, from $c$ to $d$, and there is no directed path of length $k$ from $a$ to $d$. Note that by Lemma~\ref{lem:coresCanPPconstructConstants} we are allowed to use constants in pp-constructions. 
Recall that $x\stackrel{k}\toEdge y$
is a shortcut for the pp-formula 
$\exists u_1,\dots,u_{k-1}.\ E(x,u_1) \wedge E(u_1,u_2) \wedge \cdots \wedge E(u_{k-1},y)$. 
Consider the second pp-power of $\G$ given by the formula
\begin{align*}
    \phi(x_1,x_2,y_1,y_2)&\coloneqq x_1\stackrel{k}\toEdge y_2 \AND (x_2\approx d) \AND (y_1\approx a).
\end{align*}
Let $\Hb$ be the resulting digraph. Consider the vertices 
\[\text{$v_0\coloneqq(c,d)$, $v_1\coloneqq(a,d)$, and $v_2\coloneqq(a,b)$ of $\Hb$.}\] 
Note that the only vertex of $\Hb$ that can have incoming and outgoing edges is $v_1$. Since there is no path of length $k$ from $a$ to $d$ the vertex $v_1$ does not have a loop. Furthermore, $\Hb$ has the edges $(v_0,v_1), (v_1,v_2),$ and $(v_0,v_2)$ (see Figure~\ref{fig:T3constr}). 
Hence, $i\mapsto v_i$ is an embedding of $\T_3$ into $\Hb$. 
Let $V_0$ be the set of all vertices in $\Hb$ that have outgoing edges and $V_2$ be the set of all vertices in $\Hb$ that have incoming edges. Let $V_1$ denote the set $(V_0\cap V_2)\cup(H\setminus(V_0\cup V_2))$. Note that $V_1$ consists of $v_1$ and all isolated vertices. Clearly, $V_0\setminus V_2$, $V_1$, and $V_2\setminus V_0$ form a partition of $H$ and
the map
\begin{align*}
    v\mapsto
    \begin{cases}
    2&\text{if }v\in V_2\setminus V_0\\
    1&\text{if }v\in V_1\\
    0&\text{if }v\in V_0\setminus V_2
    \end{cases}
\end{align*}
is a homomorphism from $\Hb$ to $\T_3$. Hence, $\G \leq \T_3$.
\end{proof}

Now we can classify the lower covers of $\P_1$ in $\DGPoset$.
\begin{proof}[Proof of Theorem~\ref{thm:submaximalGraphs}]
Let $\G$ be a digraph such that $\P_1\not\ppleq \G$. 
Theorem~\ref{thm:maltAndCyclImplyIdempotent} implies that either $\Pol(\G)$ does not satisfy $\Maltsev$ or that it does not satisfy $\Sigma_p$ for some prime $p$. In the first case $\G \leq \T_3$, by Lemma~\ref{lem:notSatMalt}. 
In the second case $\G\leq\C_p$, by Corollary~\ref{cor:PPvsLoopCaBlockerForSa}.
Hence, all submaximal elements of $\DGPoset$ are contained in 
$\{\T_3, \C_2 ,\C_3, \C_5, \dots\}$. 
Lemma~\ref{lem:pcSatisfyPclc}, Example~\ref{exa:cyclesHaveMaltsev},  and Lemma~\ref{lem:T3Conditions} imply that these digraphs form an antichain in $\DGPoset$. Hence, each of these digraphs represents a different submaximal element of $\DGPoset$. 
\end{proof}


Note that our result implies the following. 

\begin{corollary}\label{cor:maltAndCyclImplyIdempotent}
If a finite digraph $\G$ satisfies $\Maltsev$, $\Sigma_2$, $\Sigma_3$, $\Sigma_5$, $\dots$, then any minor condition satisfied by $\Pol(\P_1)$ is also satisfied by $\Pol(\G)$.
\end{corollary}

The statement of Corollary~\ref{cor:maltAndCyclImplyIdempotent}  may also be phrased as
\[\{\Maltsev, \Sigma_2, \Sigma_3, \Sigma_5, \dots\}\subseteq \Sigma(\G) \quad \Rightarrow \quad  \Sigma(\P_1)\subseteq \Sigma(\G),\]
where $\Sigma(\G)$ and $\Sigma(\P_1)$ denote the sets of minor conditions that hold in $\G$ and in $\P_1$ respectively. 

\begin{remark}
We do not know whether Corollary~\ref{cor:maltAndCyclImplyIdempotent} holds for arbitrary clones of operations on a finite set, instead of just clones of the form $\Pol(\G)$ for a finite digraph $\G$. However, the statement is false for clones of operations on an infinite set, as illustrated by the clone of operations on ${\mathbb Q}$ of the form
$(x_1,\dots,x_n) \mapsto a_1x_1+\cdots+a_nx_n$ for 
$a_1,\dots,a_n \in {\mathbb Q}$ such that $a_1 + \cdots + a_n = 1$. 
This clone satisfies $\Sigma_n$ for every $n \in \mathbb N$,
and contains the function $(x_1,x_2,x_3) \mapsto x_1 - x_2 + x_3$, so it also satisfies $\Maltsev$. 
However, it is easy to see that it does not contain an operation $f$ that satisfies
\[f(x,x,y) = f(y,y,x) = f(x,y,y) = f(y,x,x)\]
for all $x,y \in {\mathbb Q}$;
however, this minor condition is satisfied by $\Pol(\P_2)$ (for example by $f = \max$). 
\end{remark}

\begin{remark}
Many, but not all the statements that we have shown also apply to \emph{infinite} digraphs. In the statement
\begin{align*}
    \Hb \ppleq \G &&\text{if and only if} && \Sigma(\Hb)\subseteq \Sigma(\G),
\end{align*}
only the forward direction holds if $\G$ and $\Hb$ are infinite digraphs; however, in this text we only used the forward direction. 

Every digraph with a Maltsev polymorphism is totally rectangular even if the digraph
is infinite. And every rectangular digraph
is homomorphically equivalent to 
\begin{itemize}
    \item some $\P_n$,
    \item some (possible infinite) disjoint union of directed cycles,
    \item 
    $\P^\infty\coloneqq(\N,\{(u,u+1)\mid u\in\N\})$,
    \item 
$\P_\infty\coloneqq(\N,\{(u+1,u)\mid u\in\N\})$, 
    \item 
the disjoint union $\P_\infty + \P^\infty$ of $\P_\infty$ and $\P^\infty$, or
\item 
  $\P^\infty_\infty\coloneqq(\Z,\{(u,u+1)\mid u\in\Z\})$.
\end{itemize}  
All of these graphs have a Maltsev polymorphism. 
Infinite disjoint unions of cycles are clearly not submaximal. 
Observe that, $\P_1$ cannot pp-construct the core digraphs $\P_{\infty}$, $\P^{\infty}$, $\P_\infty + \P^\infty$, and $\P^\infty_\infty$, 
and these graphs can pp-construct $\P_1$. Clearly $\P_{\infty}$ and $\P^{\infty}$ pp-construct each other. 
We do not know whether these graphs are submaximal in the class of all digraphs.  
\end{remark}

\section[Concluding Remarks]{Concluding Remarks and Directions for Further Research}


Recall that we have described the pp-constructibility poset on disjoint unions of cycles  in Chapter~\ref{cha:cycles}. 
Note that this description combined with the result of the current chapter shows that when exploring $\DGPoset$ it only remains to explore the interval between $\K_3$ and $\T_3$: 
if a digraph $\Hb$ does not have a Maltsev polymorphism, then we proved that it is below $\T_3$ (and above $\K_3$);
otherwise, it is homomorphically equivalent to a directed path or a disjoint union of cycles and thus falls into the region that has already been fully described. However it is still not clear what the lower covers of $\C_2,\C_3,\C_5,\dots$ are that do not have a Maltsev polymorphism. This question will be answered in Lemma~\ref{thm:lowerCoversOfSubmaximalCycles}. The following problem remains open.

\begin{question}\label{pro:lowerCoversOfTthree}
What are the lower covers of $\T_3$?
\end{question}
We will present a candidate for a lower cover in Problem~\ref{que:candidateLowerCoverOfT3}.



\chapter{A Short Overview of Datalog}\label{cha:DatalogIntro} 
In this chapter we will consider some of the most commonly used Datalog fragments and their well-known connections to pp-constructions and minor conditions.

\begin{definition}
    A \emph{Datalog program} is a tuple $\prog P=(\tau,\sigma,\mathcal R,G)$, where $\tau$ (\emph{EDBs}) and $\sigma$ (\emph{IDBs}) are disjoint signatures, $G\in\sigma$ is the \emph{goal relation}, and $\mathcal R$ is a set of \emph{rules} of the form
\[P(x_1,\dots,x_n)\dashv\Phi(x_1,\dots,x_n,y_1,\dots,y_m),\]
where $P\in\sigma$, called the \emph{head} of the rule, and $\Phi$, called the \emph{body} of the rule, is a $(\tau\cup\sigma)$-formula of the form 
\[R_1(z_{1,1},\dots,z_{1,n_1})\AND\dots\AND R_k(z_{k,1},\dots,z_{k,n_k}).\]
It is common to separate the conjuncts of $\Phi$  by ``," instead of ``$\AND$"~\cite{Dalmau_2005LinearDatalog}.
Let $\B,\B'$ be $(\tau\cup\sigma)$-structures. 
The program $\prog P$ can \emph{derive} $\B'$ from $\B$ in one step, denoted $\B\vdash_{\prog P}\B'$, if $B=B'$ and there is a rule $P(x_1,\dots,x_n)\dashv\Phi$ in $\mathcal  R$ and a satisfying assignment $h$ of $\Phi$  such that 
\[P^{\B'}=P^{\B}\cupdot\{(h(x_1),\dots,h(x_n))\}\]
and $Q^{\B'}=Q^{\B}$ for all $Q\in (\sigma\cup\tau)\setminus\{P\}$. We denote the transitive reflexive closure of $\vdash_{\prog P}$ by $\vdash_{\prog P}^\ast$.
The \emph{language} recognized by $\prog P$, denoted $L(\prog P)$, is the class of all finite $\tau$-structures $\B$ on which $\prog P$ can derive the goal relation, i.e., there exists a $(\tau\cup\sigma)$-structure $\B'$ such that $G^{\B'}\neq\emptyset$ and $\B\vdash_{\prog P}^\ast\B'$, where $\B$ is seen as the $(\tau\cup\sigma)$-expansion of $\B$ with $P^{\B}=\emptyset$ for all $P\in\sigma$.

A Datalog program $\prog P$ is called \emph{linear} if for all rules $P(x_1,\dots,x_n)\dashv\Phi$ in $\mathcal R$ the formula $\Phi$ has at most one conjunct with an IDB. 
A linear Datalog program $\prog P$ is called \emph{symmetric} if for every rule 
\begin{align*}
P(x_1,\dots,x_n)&\dashv\Phi\AND Q(y_1,\dots,y_m) 
\intertext{in $\mathcal R$, where $Q$ is an IDB and $\Phi$ contains no IDBs, we have that $\mathcal{ R}$ also contains the rule} 
Q(y_1,\dots,y_m)&\dashv\Phi\AND P(x_1,\dots,x_n).
\end{align*}

\end{definition}

We denote the class of all languages recognized by ((symmetric) linear) Datalog programs by  \emph{((symmetric) linear) Datalog}. 
For a finite $\tau$-structure $\A$ let $\cocsp(\A)$ denote the complement of $\csp(\A)$, i.e., the class of all finite $\tau$-structures that do not have a homomorphism into $\A$. There are close connections between the introduced Datalog fragments, CSPs, and pp-constructions.
It is well known that ((symmetric) linear) Datalog is closed under pp-constructions, i.e., if $\cocsp(\A)$ is in ((symmetric) linear) Datalog and $\A\ppleq\B$ then $\cocsp(\B)$ is in ((symmetric) linear) Datalog. However, I was not able to find a reference for this fact. Hence, we present a proof here. We roughly follow the proof that Datalog is closed under pp-construction given by Bodirsky in Section 8.3 in~\cite{theBodirsky}.

\begin{lemma}\label{lam:DLoneEDBinBody}
For every ((symmetric) linear) Datalog program $\prog P$ there is a ((symmetric) linear) Datalog program $\prog P'$ recognizing the same language such that each rule of $\prog P'$ has at most one EDB in its body.
\end{lemma}
\begin{proof}
From $\prog P$ we obtain a new Datalog program $\prog P'$ by replacing each rule 
\begin{align*}
P(a_1,\dots,a_n)&\dashv\Phi(b_1,\dots,b_k)\AND R(c_1,\dots,c_m)
\intertext{in $\mathcal R$, where $R\in\tau$ and $\Phi$ contains at least one EDB, by the two rules}
P_{\Phi}(b_1,\dots,b_k)&\dashv \Phi(b_1,\dots,b_k)\text{ and}\\
P(a_1,\dots,a_n)&\dashv P_{\Phi}(b_1,\dots,b_k)\AND R(c_1,\dots,c_m).   
\end{align*}
Moreover, if $\prog P$ is symmetric linear then we also add the symmetric counterparts of these new rules.
Clearly, $L(\prog P)=L(\prog P')$ and if $\prog P$ is (symmetric) linear, then $\prog P'$ is also (symmetric) linear. Repeating this construction multiple times yields the desired ((symmetric) linear) Datalog program. 
\end{proof}

\begin{lemma}\label{lem:DLclosedUnderPPDefinitions}
Let $\cocsp(\A)$ be in ((symmetric) linear) Datalog and $R$ be a pp-definable relation of $\A$. Then $\cocsp(\A,R)$ is in ((symmetric) linear) Datalog.
\end{lemma}
\begin{proof}
Let $\prog P=(\tau,\sigma,\mathcal R,G)$ be a ((symmetric) linear) Datalog program recognizing $\cocsp(\A)$. It suffices to show that we can add $=$, conjunctions, and projections of relations in $\A$ and stay in ((symmetric) linear) Datalog.

\underline{Adding $=$} 
Let $E$ be the new binary relation symbol used for $=$.
By Lemma~\ref{lam:DLoneEDBinBody} we can assume that each rule in $\mathcal R$ has at most one EDB in its body. Define $\prog P'=(\tau\cup\{E\},\sigma,\mathcal{R}\cup \mathcal R',G)$, where $\mathcal R'$ contains for each $P^{(n)}\in\sigma$ and $i\in\{1,\dots,n\}$ the rules
\begin{align*}
    P(a_1,\dots,a_n)&\dashv P(a_1,\dots,a_i',\dots,,a_n)\AND E(a_i', a_i)\text{ and}\\
    P(a_1,\dots,a_n)&\dashv P(a_1,\dots,a_i',\dots,,a_n)\AND E(a_i, a_i').
\end{align*}
It is easy to verify that $\prog P'$ recognizes $\cocsp(\A,=)$, note that for this to be true it is necessary that we assumed there to be at most one EDB in each body of a rule. Moreover, if $\prog P$ is (symmetric) linear, then $\prog P'$ is (symmetric) linear.

\underline{Adding conjunctions} 
Let $S^{\A}$ be the conjunction of two relations of $\A$ given by the pp-formula $\psi(x_1\dots,x_n)=R_1(y_1,\dots,y_{m_1})\AND R_2(z_1,\dots,z_{m_2})$, where the variables $y_1,\dots,y_{m_1},z_1,\dots,z_{m_2}$ are not necessarily distinct.
Define $\prog P'=(\tau\cup\{S\},\sigma,\mathcal{R}\cup \mathcal R',G)$, where $\mathcal R'$ is constructed as follows: 
For every rule \[P(a_1,\dots,a_n)\dashv\Phi(b_1,\dots,b_k)\AND R_1(c_1,\dots,c_{m_1})\] in $\mathcal R$, where $R\in\tau$, let $d_1,\dots,d_n,e_1,\dots,e_{m_2}$ be variables such that 
\begin{itemize}
    \item for all $1\leq i\leq n,1\leq j\leq m_1$ we have $x_i=y_j$ implies
$d_i=c_j$ and
\item for all $1\leq i\leq n,1\leq j\leq m_2$ we have $x_i=z_j$ implies $d_i=e_j$. 
\end{itemize} Note that
\[\psi(d_1\dots,d_n)=R_1(c_1,\dots,c_{m_1})\AND R_2(e_1,\dots,e_{m_2}).\]
Add to $\mathcal R'$ a new rule
\[P(a_1,\dots,a_n)\dashv\Phi(b_1,\dots,b_k)\AND S(d_1,\dots,d_n).\]
Analogously add new rules to $\mathcal R'$ for all rules in $\mathcal R$ with $R_2$ in their body. Clearly the resulting program $\prog P'$ recognizes $\cocsp(\A,S)$. Moreover, if $\prog P$ is (symmetric) linear, then $\prog P'$ is (symmetric) linear.

\underline{Adding projections} Let $\S^{\A}$ be a projection of a relation in $\A$ given by the pp-formula $\psi(x_1,\dots,x_m)=\exists x.\ R(x,x_1,\dots,x_m)$.
By Lemma~\ref{lam:DLoneEDBinBody} we can assume that each rule in $\mathcal R$ has at most one EDB in its body. 
Define $\prog P'=(\tau\cup\{S\},\sigma\cup\sigma',\mathcal{R}\cup \mathcal R',G)$, where $\sigma'$ and $\mathcal R'$ are constructed as follows: For every $P\in\sigma^{(k)}$ and $I\subseteq\{1,\dots,k\}$ add a new $(k+(n-1)\cdot|I|)$-ary relation symbol $P_I$ to $\sigma'$. 
For every rule 
\[P(a_1,\dots,a_n)\dashv\Phi(b_1,\dots,b_k)\AND R(c,c_1,\dots,c_{m})\] 
in $\mathcal{R}$ and $V\subseteq\{b_1,\dots,b_k\}\setminus\{c_1,\dots,c_m\}$ let $\boldsymbol t_x$ be an $n$-tuple of new variables for every $x\in V\cup\{c\}$. Let $I\coloneqq\{i\mid 1\leq i\leq n,a_i\in V\}$, let $h$ be a map such that  $h(x)=\boldsymbol t_x$ if $x\in V$ and $h(x)=x$ otherwise, and let $\Phi_V$ be the pp-formula obtained from $\Phi$ by replacing every IDB $P'(b_1',\dots,b'_{k'})$ with $P'_J(h(b_1'),\dots,h(b'_{k'}))$, where $J=\{i\mid 1\leq i\leq n,b'_i\in V\}$. Now we add new rules to $\mathcal R'$:
\begin{align*}
\intertext{If $c\notin V$ and $c\notin \{a_1,\dots,a_n\}$, then add the rule}
P_I(h(a_1),\dots,h(a_n))&\dashv \Phi_V(h(b_1),\dots,h(b_k))\AND S(c_1,\dots,c_m).\\
\intertext{If $c\notin V$ and $c=a_i$ for some $i$, then add the rule}
P_{I\cup\{i\}}(h(a_1),\dots,\boldsymbol t_c,\dots,h(a_n))&\dashv \Phi_V(h(b_1),\dots,h(b_k))\AND S(c_1,\dots,c_m).\\
\intertext{If $c\in V$, then add the rule}
P_I(h(a_1),\dots,h(a_n))&\dashv \Phi_V(h(b_1),\dots,h(b_k)).
\end{align*}
Similarly, for every rule $P(a_1,\dots,a_n)\dashv\Phi(b_1,\dots,b_k)\AND R'(c_1,\dots,c_{m'})$ in $\mathcal{R}$ (where $R'$ could also be $R$) and $V\subseteq\{b_1,\dots,b_k\}\setminus\{c_1,\dots,c_{m'}\}$ let $\boldsymbol t_x$, $h$, $I$, and $\Phi_V$ be defined as before and add the following rule to $\mathcal R'$
\[P_I(h(a_1),\dots,h(a_n))\dashv \Phi_V(h(b_1),\dots,h(b_k))\AND R'(c_1,\dots,c_{m'}).\]
Let $\Instance'$ be an instance of $\csp(\A,S)$. From $\Instance'$ we obtain an instance $\Instance$ of $\csp(\A)$ by dropping the relation $S$ and adding a new element $a_{(c_1,\dots,c_m)}$ for each $(c_1,\dots,c_m)\in S^{\Instance'}$ with $(a_{(c_1,\dots,c_m)},c_1,\dots,c_m)\in R^{\Instance}$.
Let $(a_1,\dots,a_n)$ be a tuple in $\Instance$ and let $I\subseteq \{1,\dots,n\}$ contain the indices of the newly added elements in this tuple and let $h$ be a map such that $h(a_i)=a_i$ if $i\notin I$ and $h(a_i)=(c_1,\dots,c_m)$ if $a_i$ is the newly added element $a_{(c_1,\dots,c_m)}$.
It is left as an exercise for the reader to verify that $\prog P$ can derive an IDB $P$ on $(a_1,\dots,a_n)$ if and only if  $\prog P'$ can derive the IDB $P_I$ on $(h(a_1),\dots,h(a_n))$.
Hence, $\prog P$ derives $G$ on $\Instance$ if and only if $\prog P'$ derives $G$ on $\Instance'$. Therefore $\prog P'$ recognizes $\cocsp(\A,S)$.  Moreover, if $\prog P$ is (symmetric) linear, then $\prog P'$ is (symmetric) linear.
\end{proof}

\begin{theorem}\label{thm:DLClosedUnderPPPowers}
Let $\cocsp(\A)$ be in ((symmetric) linear) Datalog and $\B$ be a pp-power of $\A$. Then $\cocsp(\B)$ is in ((symmetric) linear) Datalog.
\end{theorem}
\begin{proof}
Let $\prog P=(\tau,\sigma,\mathcal R,G)$ be a ((symmetric) linear) Datalog program recognizing $\cocsp(\A)$ and
let $\B$ be a $\tilde\tau$-structure and a $p$-th pp-power of $\A$.
By Lemma~\ref{lam:DLoneEDBinBody} we can assume that each rule in $\mathcal R$ has at most one EDB in its body.  
For every $k$-ary relationsymbol $\tilde R$ in $\tau$ we have that $\tilde R^{\B}$ as an $(n\cdot k)$-ary relation of $\A$ is pp-definable in $\A$. Hence, we can assume by $(1)$, $(2)$, and $(3)$ that this relation is already in $\A$ and has the relationsymbol $R$. Note that we can remove relations from $\A$ and all rules from $\mathcal R$ that use the removed relationsymbols of the removed relations to obtain a Datalog program for the reduct of $\A$. Therefore we can assume without loss of generality that $\tau=\{R\mid \tilde R\in\tilde \tau\}$ and for all $R\in\tau$ 
\begin{align*}
    &(c_{1,1},\dots,c_{1,p},\dots,c_{m,1},\dots,c_{m, p})\in R^{\A} \text{ if and only if }\\ &((c_{1,1},\dots,c_{1,p}),\dots,(c_{m,1},\dots,c_{m, p}))\in\tilde R^{\B}.
\end{align*}
As an intermediate step we construct the program $\prog P_1=(\tau,\sigma_1,\mathcal R_1,G)$, 
where $\mathcal R_1$ is obtained from $\mathcal R$ by replacing each rule  
\begin{align*}
P(a_1,\dots,a_n)&\dashv\Phi(b_1,\dots,b_k)\AND R(c_1,\dots,c_{m})
    \intertext{in $\mathcal R_1$ where the $c_1,\dots,c_m$ are not pairwise distinct by the two rules}
P_{\Phi,c_1,\dots,c_m}(a_1,\dots,a_n,d_1,\dots,d_{m})&\dashv\Phi(b_1,\dots,b_k)\AND R(d_1,\dots,d_{m})\\
P(a_1,\dots,a_n)&\dashv P_{\Phi,c_1,\dots,c_m}(b_1,\dots,b_k,c_1,\dots,c_m)
\end{align*}
where the $d_i$'s are new pairwise distinct variables and $P_{\Phi,c_1,\dots,c_m}$ is a new IDB.
If $\prog P$ is symmetric we also add the symmetric counterparts of the new rules to $\mathcal R_1$. 
Note that for every EDB $R(c_1,\dots,c_{m})$ that occurs in a body of a rule of $\prog P_1$ we have that the $c_i$'s are pairwise distinct. Furthermore, $L(\prog P)=L(\prog P_1)$ and that if $\prog P$ is (symmetric) linear, then $\prog P_1$ is (symmetric) linear. 

Now construct a program $\prog P_2=(\tilde \tau,\sigma_2,\mathcal R_2,G)$, where $\mathcal R_2$ is constructed as follows.
For every rule 
\begin{align*}
P(a_1,\dots,a_n)&\dashv \Phi(b_1,\dots,b_k)\AND R(c_{1,1},\dots,c_{1,p},\dots,c_{m,1},\dots,c_{m, p})
\intertext{ in $\mathcal R_2$    
 and $i_1,\dots,i_n,j_1,\dots,j_k\in\{1,\dots,p\}$ let $t_1,\dots,t_m$ be variables such that 
if $c_{i,j}=b_{\ell}$ then $j_{\ell}=j$ and $s_{\ell}=t_i$ and
if $c_{i,j}=a_{\ell}$ then $i_{\ell}=j$ and $r_{\ell}=t_i$.
Now add to $\mathcal{R}_2$ the rule}
P_{i_1,\dots,i_n}(r_1,\dots,r_n)&\dashv \Phi_{j_1,\dots,j_k}(s_1,\dots,s_k)\AND \tilde R(t_1,\dots,t_m),
\end{align*}
where $\Phi_{j_1,\dots,j_k}$ is obtained from $\Phi$ by replacing every conjunct $Q(b_{\ell_1},\dots,b_{\ell_{k'}})$ by $Q_{j_{\ell_1},\dots,j_{\ell_{k'}}}(s_{\ell_1},\dots,s_{\ell_{k'}})$.
For example let $q=2$. Then for the rule
\begin{align*}
    P(a_1,a_2)&\dashv Q(a_1,b_3)\AND Q'(a_1,b_2)\AND R(a_1,b_2,a_2,c_4)
    \intertext{the following two rules would be added to $\mathcal R_2$ }
    P_{1,1}(t_1,t_2)&\dashv Q_{1,1}(t_1,s)\AND Q'_{1,2}(t_1,t_1) \AND \tilde R(t_1,t_2)\text{ and}\\
    P_{1,1}(t_1,t_2)&\dashv Q_{1,2}(t_1,s)\AND Q'_{1,2}(t_1,t_1) \AND \tilde R(t_1,t_2).
\end{align*}
We are still not done. We also need to add to $\mathcal R_2$ the rule
\[P_{i_1,\dots,i_j,\dots,i_n}(r_1,\dots,t_k,\dots,r_n)\dashv
\begin{tabular}{l}
$\phantom{\AND{}}P_{i_1,\dots,i'_j,\dots,i_n}(r_1,\dots,t_\ell,\dots,r_n)$\\ 
$\AND\tilde R(t_1,\dots,t_m),$
\end{tabular}\]
whenever
\[\A\models \forall c_{1,1},\dots,c_{m, p}\ R(c_{1,1},\dots,c_{1,p},\dots,c_{m,1},\dots,c_{m, p})\Rightarrow (c_{k,i_j}\approx c_{\ell,i_j'}).\]
Note that, if $\prog P_1$ is (symmetric) linear, then $\prog P_2$ is (symmetric)
linear as well.
It is left as a cumbersome exercise to the reader to verify that $\prog P_2$ indeed recognizes $\cocsp(\B)$. 
\end{proof}

Combining the previous theorem with the fact that if $\A$ and $\B$ are homomorphically equivalent, then $\cocsp(\A)=\cocsp(\B)$, we obtain the desired corollary.

\begin{corollary}\label{cor:DLclosedUnderPPConstructions}
Datalog, linear Datalog, and symmetric linear Datalog are closed under pp-constructions, i.e., if $\A\ppleq\B$ and $\cocsp(\A)$ in ((symmetric) linear) Datalog, then $\cocsp(\B)$ in ((symmetric) linear) Datalog.
\end{corollary}

\section{Containment in Datalog}
The structure $\TLinP$ has the domain
$D = \{0,\dots,p-1\}$ where $p$ is some prime, the relation $\{(x,y,z) \mid x+y+z\equiv 0 \pmod p\}$, 
and the relation $\{x\}$ for every $x \in D$. 
It is well known that 
$\cocsp(\TLinP)$ is not in Datalog~\cite{FederVardi}. 

\begin{definition}\label{def:wnu34}
The \emph{3-4 weak near-unanimity condition} consists of the identities
\begin{align*}
    &f(y,x,x)\approx f(x,y,x)\approx f(x,x,y)\\
    \approx{} &g(x,x,x,y)\approx g(x,x,y,x)\approx g(x,y,x,x)\approx g(y,x,x,x).
\end{align*}
We abbreviate the condition by  $\WNU{3,4}$.
\end{definition}

We have the following well-known characterization of CSPs in Datalog.
\begin{theorem}[\cite{BoundedWidthJournal,Maltsev-Cond}, Theorem 47 in \cite{Pol}]\label{thm:DLiff34WNUiff3Linp}
Let $\A$ be a finite relational structure. Then the following are equivalent. 
\begin{enumerate}
    \item $\cocsp(\A)$ is in Datalog, in particular $\csp(\A)$ is in P,
\item $\A\not\ppleq\TLinP$ for all primes $p$, and 
\item $\A\models\WNU{3,4}$.
\end{enumerate}
\end{theorem}


\section{Containment in Linear Datalog}
\label{sect:NL}
In this section we present some sufficient and some necessary conditions for containment of $\cocsp(\A)$ in linear Datalog. 

\begin{definition}\label{def:nu}
For $n\geq 3$, the \emph{$n$-ary quasi near-unanimity condition}, denoted $\NU n$, is
\[f(x,x,\dots,x)\approx f(y,x,\dots,x) \approx f(x,y,\dots,x) \approx \cdots \approx f(x,x,\dots,y).\]
The \emph{quasi majority condition}, denoted $\Majority$, is $\NU3$.
\end{definition}

The existence of a near-unanimity polymorphism characterizes \emph{bounded strict width} \cite{FederVardi}. More importantly for us, a near-unanimity polymorphism is sufficient to put $\cocsp(\A)$ in linear Datalog, using the following two results. Barto, Kozik, and Willard proved that finite structures with finite relational signature and a near-unanimity polymorphism have bounded pathwidth duality~\cite{BartoKozikWillard}. Dalmau proved that bounded pathwidth duality implies containment in NL~\cite{Dalmau_2005LinearDatalog}. It is well known that linear Datalog is contained in NL. The following theorem gives a sufficient condition for containment in linear Datalog.

\begin{theorem}[\cite{Barto-cd,BartoKozikWillard,Dalmau_2005LinearDatalog}]\label{thm:NUimpliesInLinDL}
Let $\A$ be a finite relational structure. If there is an $n\geq3$ with $\A\models\NU n$, 
then $\cocsp(\A)$ is in linear Datalog, and in particular $\cocsp(\A)$ is in NL. 
\end{theorem}


The structure $\HornSAT$ has the domain $\{0,1\}$ and a ternary relation $\{0,1\}^3 \setminus \{(1,1,0)\}$,
and the two unary relations $\{0\}$ and $\{1\})$.
It is well known that $\csp(\HornSAT)$ is P-complete, i.e., complete for the complexity class P under deterministic log-space reductions.

\begin{definition}\label{def:hmck}
Let $n\geq0$. The \emph{quasi Hobby-McKenzie condition of length $n$}, denoted $\HMcK n$, is equal to $\Maltsev$ if $n=0$ and otherwise consists of the identities 
\begin{align*}
    d_1(x,y,y) & \approx d_1(x,x,x) \\
    d_i(x,y,y) & \approx d_{i+1}(x,y,y) & \text{for even } i < n \\
    d_i(x,x,y) & \approx d_{i+1}(x,x,y) & \text{for odd } i < n \\
    d_i(x,y,x) & \approx d_{i+1}(x,y,x) & \text{for odd } i < n \\
    d_n(x,y,y) & \approx p(x,y,y) \\
    p(x,x,y) & \approx e_1(x,x,y) \\
    e_i(x,y,y) & \approx e_{i+1}(x,y,y) & \text{for odd } i < n \\
    e_i(x,x,y) & \approx e_{i+1}(x,x,y) & \text{for even } i < n \\
    e_i(x,y,x) & \approx e_{i+1}(x,y,x) & \text{for odd } i < n \\
    e_{n}(x,y,y) & \approx e_{n}(y,y,y) & \text{if $n$ is odd}\\
    e_{n}(x,y,x) & \approx e_{n}(x,x,x) & \text{if $n$ is odd}\\
    e_{n}(x,x,y) & \approx e_n(y,y,y) & \text{if $n$ is even.}
\end{align*}
\end{definition}

Observe that our definition varies slightly from the one Hobby and McKenzie give in Theorem 9.8 in~\cite{HobbyMcKenzie}. We have replaced the identities $x\approx d_0(x,y,z)$ and $d_0(x,y,y)\approx d_1(x,y,y)$ by the identity 
\[d_1(x,y,y)\approx d_1(x,x,x).\] Similarly, we have removed the last function from the chain of $e_i$'s. Clearly, both conditions are satisfied by the same finite core structures. Observe that our minor condition satisfies (ii) in Theorem~\ref{thm:StructurToDigraph}. In Figure~\ref{fig:HMcKChain} we present a visual representation of the Hobby-McKenzie condition of length 3.

\begin{figure}
    \centering
    \begin{tikzpicture}[node distance=13mm]
    \node (0) {$\pi_1$};
    \node[right of =0] (1) {$d_1$};
    \node[right of =1] (2) {$d_2$};
    \node[right of =2] (3) {$d_3$};
    \node[right of =3] (4) {$p$};
    \node[right of =4] (5) {$e_1$};
    \node[right of =5] (6) {$e_2$};
    \node[right of =6] (7) {$e_3$};
    \node[right of =7] (8) {$\pi_3$};
    \path
    (0) edge node[above] {$xyy$} (1)
(1) edge node[above,align=left] {$xxy$ \\$xyx$} (2)
(2) edge node[above] {$xyy$} (3)
(3) edge node[above] {$xyy$} (4)
(4) edge node[above] {$xxy$} (5)
(5) edge node[above,align=left] {$xyy$ \\$xyx$} (6)
(6) edge node[above] {$xxy$} (7)
(7) edge node[above,align=left] {$xyy$ \\$xyx$} (8)
    ;
    \end{tikzpicture}
    \caption{Visualization of $\HMcK 3$.}
    \label{fig:HMcKChain}
\end{figure}


\begin{theorem}[Consequence of Theorem 9.8 in~\cite{HobbyMcKenzie}]\label{thm:HMcKiffnoHorn3SAT}
A finite relational structure does not satisfy $\HMcK n$ for any $n \geq 1$ if and only if it can pp-construct $\HornSAT$.
\end{theorem}
The following theorem gives a necessary condition for containment in linear Datalog.
\begin{theorem}[Theorem~4.2 in \cite{LaroseTesson}]
Let $\A$ be a finite relational core structure. If $\mathcal{V}(\Pol(\A))$ admits the semilattice type, then $\cocsp(\A)$ is not in linear Datalog. 
\end{theorem}
Let $\A$ be a finite core structure; By Definition~4.10 in~\cite{HobbyMcKenzie} the algebra $\Pol(\A)$ is of \emph{semilattice type} if it is polynomially equivalent to a two-element semilattice. The clone $\Pol(\HornSAT)$ is generated by the binary semilattice operation $\min$. Hence, $\Pol(\HornSAT)$ is of semilattice type and we conclude the following well-known Theorem.
\begin{theorem}[Theorem 3.6 in \cite{AfratiCosmadakis}]
We have that $\cocsp(\HornSAT)$ is not in linear Datalog.
\end{theorem}
Since linear Datalog is closed under pp-constructions we obtain the following corollary.
\begin{corollary}\label{cor:linDLimpliesNotHornSat}
Let $\A$ be a finite relational structure. If $\cocsp(\A)$ is in linear Datalog, then $\A\not\ppleq\HornSAT$ and $\A\models\HMcK n$ for some $n\geq 1$. 
\end{corollary}

\section{Containment in Symmetric Linear Datalog}
\label{sect:L}

Define  $\Ord\coloneqq(\{0,1\};\leq,\{0\},\{1\})$. This structure is also known as st-Con, since its CSP encodes the complement of the directed connectivity (or reachability) problem, a typical NL-hard problem. 

\begin{definition}\label{def:hm}
For $n \geq 1$, the \emph{quasi Hagemann-Mitschke condition of length $n$}, called $\HM n$, consists of the identities 
\begin{align*}
p_1(x,x,x) & \approx p_1(x,y,y) \\
p_{i}(x,x,y) & \approx p_{i+1}(x,y,y) && \text{ for all } i \in \{1,\dots,n-1\} \\
p_n(x,x,y) & \approx p_n(y,y,y).
\end{align*}
\end{definition}

Note that $\HM1=\Maltsev$ and that $\HM n$ implies $\HM{n+1}$ for every $n \geq 1$. 
Combining known results we obtain the following characterisation of $\cocsp$'s in symmetric linear Datalog.


\begin{theorem}\label{thm:symLinDLiffHM}
Let $\A$ be a finite relational structure with $\cocsp(\A)$ in linear Datalog. Then the following are equivalent
\begin{enumerate}
    \item $\cocsp(\A)$ is in symmetric linear  Datalog, in particular $\csp(\A)$ is in $\Lclass$,
    \item $\A\models\HM n$ for some $n \geq 1$, and
    \item $\A\not\ppleq\Ord$.
\end{enumerate}
\end{theorem}
\begin{proof}
$1\Rightarrow3$ Assume that $\cocsp(\A)$ is in symmetric linear Datalog and that  $\A\leq\Ord$. Then, by Corollary~\ref{cor:DLclosedUnderPPConstructions}, $\cocsp(\Ord)$ in symmetric linear Datalog. This contradicts Theorem 12 in \cite{EgristConNotinSymLinDL}, which states that $\cocsp(\Ord)$ is not in symmetric linear Datalog.

$3\Rightarrow2$ 
This is well known, for a proof see for example Corollary~\ref{cor:OrdisLeastLowerBoundForTn} and Theorem~\ref{thm:TnIsHMblocker}.

$2\Rightarrow1$ 
This is exactly the statement of Theorem 2.7 in \cite{Kazda-n-permute}.
\end{proof}

Note that $2\Rightarrow1$ is the only implication which needs the assumption that $\cocsp(\A)$ is in linear Datalog.

\chapter{Semicomplete Digraphs}\label{cha:semicompleteTn}

A digraph $\G$ is \emph{semicomplete} if any two nodes in $\G$ are connected by either a directed or an undirected edge. We call $\G$ a \emph{tournament} if any two nodes in $\G$ are connected by a directed edge.
In this chapter we present a generalization of $\T_3$ and rectangularity, i.e, transitive tournaments and braidedness. We show that the transitive tournaments form an infinite descending chain and find a digraph which is a lower cover of this chain. Furthermore, we attempt to classify all semicomplete digraphs.  The main motivation for this chapter was to answer Problem~\ref{pro:lowerCoversOfTthree} and classify the lower covers of the submaximal elements in $\DGPoset$. We managed to determine the lower covers of directed cycles of prime length (Theorem~\ref{thm:lowerCoversOfSubmaximalCycles}). However, the lower covers of $\T_3$ remain elusive.

\section{Transitive Tournaments}
Before we come to general semicomplete graphs we consider a special case.
For any $n\in\N^+$ define the \emph{transitive tournament on $n$ nodes} 
\[\T_n\coloneqq(\{0,\dots,n-1\},\{(i,j)\mid 0\leq i<j<n\}).\]
The CSP of $\T_n$ is very simple, i.e., for every finite digraph $\G$ we have that $\G\to\T_n$ if and only if $\P_n\not\to\G$. 

\begin{lemma}
Let $n\in\N^+$\!. Then $\T_{n+1}\ppleq\T_n$.
\end{lemma}
\begin{proof}
The first pp-power of $\T_{n+1}$ given by 
\begin{center}
    \begin{tikzpicture}[scale=0.5]
    \node[var-f,label=left:$x$] (x) at (0,-1) {};
    \node[var-b] (a) at (0,1) {};
    \node[var-f,label=left:$y$] (y) at (0,0) {};
    \path[->,>=stealth']
        (x) edge (y)
        (y) edge (a)
        ;
\end{tikzpicture}
\end{center}
is isomorphic to the disjoint union of $\T_n$ with an isolated vertex.
\end{proof}

We can also separate the $\T_n$ from one another. To do this we use Hagemann-Mitschke conditions (see Definition~\ref{def:hm}) and the notion of braids which was introduced by Kazda in Definition 4.2 in~\cite{Kazda-n-permute}.
Let $\A$ be a structure. An \emph{$n$-braid of $\A$} consists of pp-formulas $\phi_1(x_1,y_1),\dots,\phi_n(x_n,y_n)$ together with a tuple $(a_0,\dots,a_n,b_0,\dots,b_n)$ from $A$ such that 
\[\A\models \phi_{i}(a_{i-1},a_i)\AND \phi_{i}(b_{i-1},b_i)\AND \phi_{i}(b_{i-1},a_i)\text{ for all $1\leq i\leq n$.}\]
This $n$-braid can be seen as the structure $\B$ that has a homomorphism into $\A_A$  expanded by the relations $\phi_1^{\A},\dots,\phi_n^{\A}$, where $\B$ is
\begin{center}
    \begin{tikzpicture}[scale = 1]
    \node[var-b,label=above:$a_0$] (10) at (0,1) {};
    \node[var-b,label=below:$b_0$] (00) at (0,0) {};
    \node[var-b,label=above:$a_1$] (11) at (1,1) {};
    \node[var-b,label=below:$b_1$] (01) at (1,0) {};
    \node[] (12) at (2,1) {$\dots$};
    \node[] (02) at (2,0) {$\dots$};
    \node[var-b,label=above:$a_{n-1}$] (13) at (3,1) {};
    \node[var-b,label=below:$b_{n-1}$] (03) at (3,0) {};
    \node[var-b,label=above:$a_n$] (14) at (4,1) {};
    \node[var-b,label=below:$b_n$] (04) at (4,0) {};
\tikzset{decoration={snake,amplitude=.25mm,segment length=1.3mm,post length=0.7mm,pre length=0.7mm}}
    \draw[decorate,->,>=stealth'] (00) -- node[above] {\small$1$} (01);
    \draw[decorate,->,>=stealth'] (00) -- node[above] {\small$1$} (11);
    \draw[decorate,->,>=stealth'] (10) -- node[above] {\small$1$} (11);
    \draw[decorate,->,>=stealth'] (03) -- node[above] {\small$n$} (04);
    \draw[decorate,->,>=stealth'] (03) -- node[above] {\small$n$} (14);
    \draw[decorate,->,>=stealth'] (13) -- node[above] {\small$n$} (14);
\end{tikzpicture}
\end{center}
and \begin{tikzpicture}[scale=0.7]
    \node[var-f,label=above:$x$] (x1) at (0,0) {};
    \node[var-f,label=above:$y$] (x2) at (1,0) {};
\tikzset{decoration={snake,amplitude=.25mm,segment length=1.3mm,post length=0.7mm,pre length=0.7mm}}
    \draw[decorate,->,>=stealth'] (x1) -- node[above] {\small$i$} (x2);
\end{tikzpicture} 
is an abbreviation for 
the relation $\{(a,b)\mid \A\models\phi_i(a,b)\}$. 
An $n$-braid $\phi_1(x_1,y_1),\dots,\phi_n(x_n,y_n)$ with $(a_0,\dots,a_n;b_0,\dots,b_n)$ of $\A$ is \emph{satisfied} by a tuple $(c_0,\dots,c_n)$ if $c_0=a_0$, $c_n=b_n$, and 
\[\A\models\phi_1(c_0,c_1)\AND\phi_2(c_1,c_2)\AND\dots\AND \phi_n(c_{n-1},c_n), \]
i.e., that
\begin{center}
    \begin{tikzpicture}[scale = 1]
    \node[var-b,label=above:$a_0$] (10) at (0,1) {};
    \node[var-b,label=above:$c_1$] (11) at (1,1) {};
    \node[] (12) at (2,1) {$\dots$};
    \node[var-b,label=above:$c_{n-1}$] (13) at (3,1) {};
    \node[var-b,label=above:$b_n$] (14) at (4,1) {};
\tikzset{decoration={snake,amplitude=.25mm,segment length=1.3mm,post length=0.7mm,pre length=0.7mm}}
\draw[decorate,->,>=stealth'] (10) -- node[above] {\small$1$} (11);
\draw[decorate,->,>=stealth'] (13) -- node[above] {\small$n$} (14);
\end{tikzpicture}
\end{center}
is a path in $\A$ expanded by $\phi_1^{\A},\dots,\phi_n^{\A}$.
A structure $\A$ is called \emph{$n$-braided} if any $n$-braid of $\A$ can be satisfied. 
Note that if a digraph is $1$-braided then it is totally rectangular (see Chapter~\ref{cha:cycles}). Observe that from any $n$-braid
\begin{itemize}
    \item $\phi_1(x_1,y_1),\phi_2(x_2,y_2),\dots,\phi_n(x_n,y_n)$ with 
    \item 
$(a_0,a_1,a_2\dots,a_n,b_0;b_1,b_2\dots,b_n)$
\end{itemize} 
we can construct an $(n-1)$-braid 
\begin{itemize}
    \item $\exists z.\ \phi_1(x_1,z)\AND\phi_2(z,y_1),\dots,\phi_n(x_n,y_n)$ with
    \item $(a_0,a_2,\dots,a_n;b_0,b_2,\dots,b_n)$.
\end{itemize}Hence, $n$-braidedness implies $(n+1)$-braidedness.

The $n$-braid of $\A$ \emph{generated} by $( a_0,\dots, a_n; b_0,\dots, b_n)$ consists of the pp-formulas \[\phi_i( x_i, y_i)\coloneqq \bigAND\{ \phi( x_i, y_i)\mid \A\models\phi_{i}( a_{i-1}, a_i)\AND \phi_{i}( b_{i-1}, b_i)\AND \phi_{i}( b_{i-1}, a_i)\}\]
together with the generating tuple.
Note that, by Theorem~\ref{thm:ppPolyPreservation}, we have 
\[\phi_i^{\A}=\{(f( a_{i-1},b_{i-1},b_{i-1}),f(a_{i},a_{i},b_{i}))\mid f\in \Hom(\A^3,\A)\}.\]
Hence, if an $n$-braid of $\A$ generated by $( a_0,\dots, a_n; b_0,\dots, b_n)$ is satisfied by $(a_0,c_1,\dots,c_{n-1},b_n)$, then there are ternary polymorphisms $p_1,\dots,p_n$ of $\A$ such that
\begin{align*}
a_0 & = p_1(a_0,b_0,b_0) \\
p_{i}(a_i,a_i,b_i) & =c_i= p_{i+1}(a_i,b_i,b_i) && \text{ for all } i \in \{1,\dots,n-1\} \\
p_n(a_n,a_n,b_n) & = b_n.
\end{align*}
It is important to note that $p_1,\dots,p_n$ do not have to satisfy $\HM{n}$. 
To show that $n$-braidedness implies the existence of a Hagemann-Mitschke chain of length $n$ we need to generalize the definition of $n$-braids to allow tuples instead of single elements.
A \emph{$(k_0,\dots,k_n)$-multi-braid of $\A$} consists of 
\begin{itemize}
    \item pp-formulas $\phi_1(\boldsymbol x_1,\boldsymbol y_1),\dots,\phi_n(\boldsymbol x_n,\boldsymbol y_n)$ together with 
    \item a tuple $(\boldsymbol a_0,\dots,\boldsymbol a_n;\boldsymbol b_0,\dots,\boldsymbol b_n)$ from $A$,
\end{itemize} where the $\boldsymbol x_i$ are $k_{i-1}$-tuples, the $\boldsymbol y_i$ are $k_i$-tuples, and the $\boldsymbol a_i,\boldsymbol b_i$ are $k_i$-tuples, such that 
\[\A\models \phi_{i}(\boldsymbol a_{i-1},\boldsymbol a_i)\AND \phi_{i}(\boldsymbol b_{i-1},\boldsymbol b_i)\AND \phi_{i}(\boldsymbol b_{i-1},\boldsymbol a_i)\text{ for all $1\leq i\leq n$.}\]
The notions \emph{braided}, \emph{satisfied}, and \emph{generated} are generalized to multi-braids in the obvious way.
We show that $n$-braidedness and multi-braidedness for $n$-tuples coincide. The proof of the following theorem was provided (personal communication) by  Alexandr Kazda.  Kazda and Valeriote proved a version of this theorem that generalizes braidedness to other testing patterns~\cite{kazdaValerioteMultibraidedness}.
\begin{theorem}\label{thm:braidedEqualsMultibraided}
Let $\A$ be a finite core structure and $n,k_0,\dots,k_n\in\N^+$\!. Then $\A$ is $n$-braided if and only if $\A$ is $(k_0,\dots,k_n)$-multi-braided.
\end{theorem}
\begin{proof}
Note that $\A$ is $n$-braided if and only if $\A$ is $(1,\dots,1)$-multi-braided, where $(1,\dots,1)$ contains $n+1$ many 1's.
We proof by induction that $(1,\dots,1)$-multi-braidedness implies $(k_0,\dots,k_n)$-multi-braidedness. 
Assume that $\A$ is $(k_0,\dots,k_n)$-multi-braided and let $i\in\{0,\dots,n\}$. We show that $\A$ is $(k_0,\dots,k_i+1,\dots,k_n)$-multi-braided.
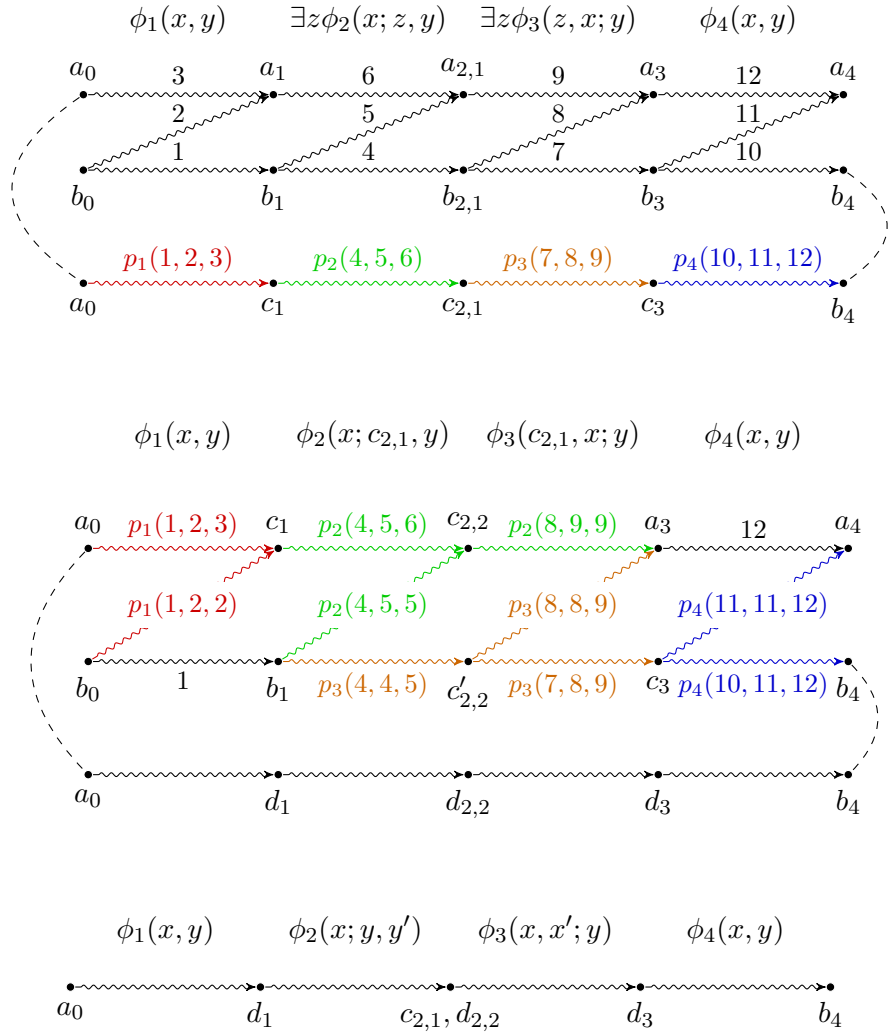
\begin{figure}
    \centering
    \begin{tikzpicture}[xscale = 2.5]
    \node[var-b,label=above:$a_0$] (10) at (0,1) {};
    \node[var-b,label=below:$b_0$] (00) at (0,0) {};

    \node[var-b,label=above:$a_{1}$] (11) at (1,1) {};
    \node[var-b,label=below:$b_{1}$] (01) at (1,0) {};

    \node[var-b,label=above:$a_{2,1}$] (12) at (2,1) {};
    \node[var-b,label=below:$b_{2,1}$] (02) at (2,0) {};

    \node[var-b,label=above:$a_{3}$] (13) at (3,1) {};
    \node[var-b,label=below:$b_{3}$] (03) at (3,0) {};
    
    \node[var-b,label=above:$a_4$] (14) at (4,1) {};
    \node[var-b,label=below:$b_4$] (04) at (4,0) {};
\tikzset{decoration={snake,amplitude=.25mm,segment length=1.3mm,post length=0.7mm,pre length=0.7mm}}
    \draw[decorate,->,>=stealth'] (00) -- node[above] {\small$1$} (01);
    \draw[decorate,->,>=stealth'] (00) -- node[above] {\small$2$} (11);
    \draw[decorate,->,>=stealth'] (10) -- node[above] {\small$3$} (11);
    
    \draw[decorate,->,>=stealth'] (01) -- node[above] {\small$4$} (02);
    \draw[decorate,->,>=stealth'] (01) -- node[above] {\small$5$} (12);
    \draw[decorate,->,>=stealth'] (11) -- node[above] {\small$6$} (12);
    
    \draw[decorate,->,>=stealth'] (02) -- node[above] {\small$7$} (03);
    \draw[decorate,->,>=stealth'] (02) -- node[above] {\small$8$} (13);
    \draw[decorate,->,>=stealth'] (12) -- node[above] {\small$9$} (13);
    
    \draw[decorate,->,>=stealth'] (03) -- node[above] {\small$10$} (04);
    \draw[decorate,->,>=stealth'] (03) -- node[above] {\small$11$} (14);
    \draw[decorate,->,>=stealth'] (13) -- node[above] {\small$12$} (14);

    \node[var-b,label=below:$a_0$] (20) at (0,-1.5) {};
    \node[var-b,label=below:$c_{1}$] (21) at (1,-1.5) {};
    \node[var-b,label=below:$c_{2,1}$] (22) at (2,-1.5) {};
    \node[var-b,label=below:$c_3$] (23) at (3,-1.5) {};
    \node[var-b,label=below:$b_4$] (24) at (4,-1.5) {};
    
    \draw[decorate,->,>=stealth',red!80!black] (20) -- node[above] {\small$p_1(1,2,3)$} (21);
    \draw[decorate,->,>=stealth',green!80!black] (21) -- node[above] {\small$p_2(4,5,6)$} (22);
    \draw[decorate,->,>=stealth',orange!80!black] (22) -- node[above] {\small$p_3(7,8,9)$} (23);
    \draw[decorate,->,>=stealth',blue!80!black] (23) -- node[above] {\small$p_4(10,11,12)$} (24);
    
    \node at (0.5,2) {$\phi_1(x,y)$};
    \node at (1.5,2.) {$\exists z\phi_2(x;z,y)$};
    \node at (2.5,2.) {$\exists z\phi_3(z,x;y)$};
    \node at (3.5,2.) {$\phi_4(x,y)$};

    \path 
    (10) edge[dashed,bend right] (20)
    (04) edge[dashed,bend left] (24)
    ;
    
\end{tikzpicture}
\vspace{10mm}

\begin{tikzpicture}[xscale = 2.5, yscale=1.5]
    \node[var-b,label=above:$a_0$] (10) at (0,1) {};
    \node[var-b,label=below:$b_0$] (00) at (0,0) {};

    \node[var-b,label=above:$c_1$] (11) at (1,1) {};
    \node[var-b,label=below:$b_1$] (01) at (1,0) {};

    \node[var-b,label=above:$c_{2,2}$] (12) at (2,1) {};
    \node[var-b,label=below:$c'_{2,2}$] (02) at (2,0) {};

    \node[var-b,label=above:$a_{3}$] (13) at (3,1) {};
    \node[var-b,label=below:$c_{3}$] (03) at (3,0) {};
    
    \node[var-b,label=above:$a_4$] (14) at (4,1) {};
    \node[var-b,label=below:$b_4$] (04) at (4,0) {};
\tikzset{decoration={snake,amplitude=.25mm,segment length=1.3mm,post length=0.7mm,pre length=0.7mm}}
    \draw[decorate,->,>=stealth'] (00) -- node[below] {\small$1$} (01);
    \draw[decorate,->,>=stealth',red!80!black] (00) -- node[fill=white] {\small$p_1(1,2,2)$} (11);
    \draw[decorate,->,>=stealth',red!80!black] (10) -- node[above] {\small$p_1(1,2,3)$} (11);
    
    \draw[decorate,->,>=stealth',orange!80!black] (01) -- node[below] {\small$p_3(4,4,5)$} (02);
    \draw[decorate,->,>=stealth',green!80!black] (01) -- node[fill=white] {\small$p_2(4,5,5)$} (12);
    \draw[decorate,->,>=stealth',green!80!black] (11) -- node[above] {\small$p_2(4,5,6)$} (12);
    
    \draw[decorate,->,>=stealth',orange!80!black] (02) -- node[below] {\small$p_3(7,8,9)$} (03);
    \draw[decorate,->,>=stealth',orange!80!black] (02) -- node[fill=white] {\small$p_3(8,8,9)$} (13);
    \draw[decorate,->,>=stealth',green!80!black] (12) -- node[above] {\small$p_2(8,9,9)$} (13);
    
    \draw[decorate,->,>=stealth',blue!80!black] (03) -- node[below] {\small$p_4(10,11,12)$} (04);
    \draw[decorate,->,>=stealth',blue!80!black] (03) -- node[fill=white] {\small$p_4(11,11,12)$} (14);
    \draw[decorate,->,>=stealth'] (13) -- node[above] {\small$12$} (14);

    \node[var-b,label=below:$a_0$] (20) at (0,-1.) {};
    \node[var-b,label=below:$d_1$] (21) at (1,-1.) {};
    \node[var-b,label=below:$d_{2,2}$] (22) at (2,-1.) {};
    \node[var-b,label=below:$d_3$] (23) at (3,-1.) {};
    \node[var-b,label=below:$b_4$] (24) at (4,-1.) {};
    
    \draw[decorate,->,>=stealth'] (20) --  (21);
    \draw[decorate,->,>=stealth'] (21) --  (22);
    \draw[decorate,->,>=stealth'] (22) --  (23);
    \draw[decorate,->,>=stealth'] (23) --  (24);
    
    
    \node at (0.5,2) {$\phi_1(x,y)$};
    \node at (1.5,2.) {$\phi_2(x;c_{2,1},y)$};
    \node at (2.5,2.) {$\phi_3(c_{2,1},x;y)$};
    \node at (3.5,2.) {$\phi_4(x,y)$};
    
    \path 
    (10) edge[dashed,bend right] (20)
    (04) edge[dashed,bend left] (24)
    ;
    
\end{tikzpicture}

\vspace{10mm}

\begin{tikzpicture}[xscale = 2.5, yscale=1.5]
    \tikzset{decoration={snake,amplitude=.25mm,segment length=1.3mm,post length=0.7mm,pre length=0.7mm}}

    \node[var-b,label=below:$a_0$] (20) at (0,1) {};
    \node[var-b,label=below:$d_1$] (21) at (1,1) {};
    \node[var-b,label={below:$c_{2,1},d_{2,2}$}] (22) at (2,1) {};
    \node[var-b,label=below:$d_3$] (23) at (3,1) {};
    \node[var-b,label=below:$b_4$] (24) at (4,1) {};
    
    \draw[decorate,->,>=stealth'] (20) --  (21);
    \draw[decorate,->,>=stealth'] (21) --  (22);
    \draw[decorate,->,>=stealth'] (22) --  (23);
    \draw[decorate,->,>=stealth'] (23) --  (24);
    
    
    \node at (0.5,1.5) {$\phi_1(x,y)$};
    \node at (1.5,1.5) {$\phi_2(x;y,y')$};
    \node at (2.5,1.5) {$\phi_3(x,x';y)$};
    \node at (3.5,1.5) {$\phi_4(x,y)$};

\end{tikzpicture}

    \caption{The multi-braids $B^{\exists}$ with edges named $1$ to $12$ (top) and $B^{c_{2,1}}$ (middle). The satisfying path of $B$ (bottom), where $c_{2,2}=p_2(a_{2,2},a_{2,2},b_{2,2})$ and  $c'_{2,2}=p_3(a_{2,2},b_{2,2},b_{2,2})$.}
    \label{fig:braidednessVSmultibraidedness}
\end{figure}
In Figure~\ref{fig:braidednessVSmultibraidedness} there is a sketch of the proof for the case $n=4$, $i=2$, $k_0=k_1=k_2=k_3=k_4=1$.
Consider a $(k_0,\dots,k_i+1,\dots,k_n)$-multi-braid  $B$ consisting of
\begin{itemize}
    \item $\phi_1(\boldsymbol x_1,\boldsymbol y_1),\dots,\phi_i(\boldsymbol x_i;\boldsymbol y'_i),\phi_{i+1}(\boldsymbol x'_{i+1};\boldsymbol y_{i+1}),\dots,\phi_n(\boldsymbol x_n,\boldsymbol y_n)$, where \\ $\boldsymbol y'_i=\boldsymbol y_i,y_{i,k_i+1}$ and $\boldsymbol x'_{i+1}=\boldsymbol x_{i+1},x_{i+1,k_i+1}$, together   with 
\item $\boldsymbol t\coloneqq(\boldsymbol a_0,\dots,(\boldsymbol a_i,a_{i,k_i+1}),\dots,\boldsymbol a_n;\boldsymbol b_0,\dots,(\boldsymbol b_i,b_{i,k_i+1}),\dots,\boldsymbol b_n)$.
\end{itemize} 
We can assume without loss of generality that the multi-braid $B$ is generated by $\boldsymbol t$. 
Define the $(k_0,\dots,k_i,\dots,k_n)$-multi-braid $B^{\exists}$ to consist of
\begin{itemize}
    \item $\phi_1(\boldsymbol x_1,\boldsymbol y_1),\dots,\exists z\ \phi_i(\boldsymbol x_i;\boldsymbol y_i,z),\exists z\ \phi_{i+1}(\boldsymbol x_{i+1},z;\boldsymbol y_{i+1}),\dots,\phi_n(\boldsymbol x_n,\boldsymbol y_n)$ \\ together with 
\item $\boldsymbol t^{\exists}\coloneqq(\boldsymbol a_0,\dots,\boldsymbol a_i,\dots,\boldsymbol a_n;\boldsymbol b_0,\dots,\boldsymbol b_i,\dots,\boldsymbol b_n)$.
\end{itemize} 
By induction hypothesis there exists a satisfying tuple \[(\boldsymbol c_0,\boldsymbol c_1,\dots,\boldsymbol c_i,\dots,\boldsymbol c_{n-1},\boldsymbol c_n),\] where $\boldsymbol c_0=\boldsymbol a_0$ and $\boldsymbol c_n=\boldsymbol a_n$. 
Observe that by definition of $B^{\exists}$ and since $B$ is generated by $\boldsymbol t$ we also have that $B^{\exists}$ is generated by $\boldsymbol t^{\exists}$. Hence, by Theorem~\ref{thm:ppPolyPreservation}, there are ternary polymorphisms $p_1,\dots,p_n$ of $\A$ such that 
$p_j(\boldsymbol a_{j-1},\boldsymbol b_{j-1},\boldsymbol b_{j-1})=\boldsymbol c_{j-1}$ and $p_j(\boldsymbol a_{j},\boldsymbol a_{j},\boldsymbol b_{j})=\boldsymbol c_{j}$ for all $j\in\{1,\dots,n\}$.
Define 
\begin{align*}
    c_{i,k_i+1}&\coloneqq p_i(a_{i,k_i+1},a_{i,k_i+1},b_{i,k_i+1}) \text{ and}\\  c'_{i,k_i+1}&\coloneqq p_{i+1}(a_{i,k_i+1},b_{i,k_i+1},b_{i,k_i+1}).
\end{align*}
To cover the cases $i=0$ and $i=n$ we say that $p_0=\pi^3_1$ and $p_{n+1}=\pi^3_3$.
Observe that
$\A\models\phi_i(\boldsymbol c_{i-1},\boldsymbol c_i,c_{i,k_i+1})$ and $\A\models\phi_{i+1}(\boldsymbol c_{i},c'_{i,k_i+1},\boldsymbol c_{i+1})$. 
Since $\A$ is a core we can assume without loss of generality that we can use constants in pp-formulas. Observe that the following $B^{\boldsymbol c_i}$ is a $(k_0,\dots,1,\dots,k_n)$-multi-braid of $\A$:
\begin{itemize}
    \item $\phi_1(\boldsymbol x_1,\boldsymbol y_1),\dots,\phi_i(\boldsymbol x_i;\boldsymbol c_i,z), \phi_{i+1}(\boldsymbol c_{i},x;\boldsymbol y_{i+1}),\dots,\phi_n(\boldsymbol x_n,\boldsymbol y_n)$ with
\item the tuple
\begin{align*}
 \boldsymbol t^{\boldsymbol c_i}\coloneqq(&\boldsymbol a_0,\boldsymbol c_1,\dots,\boldsymbol c_{i-1},c_{i,k_i+1},\boldsymbol a_{i+1}\dots,\boldsymbol a_{n-1},\boldsymbol a_n;\\
&\boldsymbol b_0,\boldsymbol b_1,\dots,\boldsymbol b_{i-1},c'_{i,k_i+1},\boldsymbol c_{i+1}, \dots,\boldsymbol c_{n-1},\boldsymbol b_n).   
\end{align*}
\end{itemize}
Note that this is not trivial, for example for $j<i$: $\A\models\phi_j(\boldsymbol b_{j-1},\boldsymbol c_j)$ since $p_j(\boldsymbol b_{j-1},\boldsymbol b_{j-1},\boldsymbol b_{j-1})=\boldsymbol b_{j-1}$ and $p_j(\boldsymbol a_{j},\boldsymbol a_{j},\boldsymbol b_{j})=\boldsymbol c_{j}$. The existence of the other edges can be argued similarly, for a full picture see Figure~\ref{fig:braidednessVSmultibraidedness}.
Hence, we can apply the induction hypothesis again, this time on $B^{\boldsymbol c_i}$ to obtain a satisfying tuple 
$(\boldsymbol a_0,\boldsymbol d_1,\dots, d_{i,k_i+1},\dots,\boldsymbol d_{n-1},\boldsymbol b_n)$. Observe that 
\[\text{$\A\models\phi_i(\boldsymbol d_{i-1};\boldsymbol c_i,c_{i,k_i+1})$ and 
$\A\models\phi_{i+1}(\boldsymbol c_i,c_{i,k_i+1};\boldsymbol d_i)$.}\] Hence, $(\boldsymbol a_0,\boldsymbol d_1,\dots, (\boldsymbol c_i,d_{i,k_i+1}),\dots,\boldsymbol d_{n-1},\boldsymbol b_n)$ is a path satisfying the original multi-braid $B$.
\end{proof}
Note that if a structure is $n$-braided, then its core is also $n$-braided.
Hence, the above theorem shows in particular that $n$-braidedness is preserved by taking pp-constructions of cores.  
\begin{corollary}\label{cor:braidednessIsPreservedByPPConstruction}
Let $\A$ and $\B$ be finite core structures such that $\A\ppleq \B$. If $\A$ is $n$-braided, then $\B$ is $n$-braided. 
\end{corollary}


The following theorem is well known, see for example Theorem 8.4 in~\cite{FreeseV09}
and Observation 6 in~\cite{kazdaValerioteMultibraidedness}, but has not been presented in the language of pp-constructions. 

\begin{theorem}\label{thm:TnIsHMblocker}
Let $\A$ be a finite core structure and $n\in\N^+$\!. Then the following are equivalent:
\begin{enumerate}
    \item $\A\ppleq \T_{n+2}$,
    \item $\A\not\models\HM n$, and
    \item $\A$ is not $n$-braided.
\end{enumerate}
\end{theorem}
\begin{proof}
In this proof we show the directions
$(1)\Rightarrow(3)\Rightarrow(2)\Rightarrow(3)\Rightarrow(1)$.

$(1) \Rightarrow (3)$ 
From Corollary~\ref{cor:braidednessIsPreservedByPPConstruction} we know that $n$-braidedness is preserved by pp-constructions. Hence, to show that $\A$ is not $n$-braided it suffices to prove that $\T_{n+2}$ is not $n$-braided. Consider the $n$-braid
\begin{center}
    \begin{tikzpicture}[scale = 0.7]
    \node[var-b,label=above:1] (10) at (0,1) {};
    \node[var-b,label=below:\strut0] (00) at (0,0) {};
    \node[var-b,label=above:2] (11) at (1,1) {};
    \node[var-b,label=below:\strut1] (01) at (1,0) {};
    \node[] (12) at (2,1) {$\dots$};
    \node[] (02) at (2,0) {$\dots$};
    \node[var-b,label=above:\small$n$] (13) at (3,1) {};
    \node[var-b,label=below:\strut\small$n-1$] (03) at (3,0) {};
    \node[var-b,label=above:\small$n+1$] (14) at (4,1) {};
    \node[var-b,label=below:\strut\small$n$] (04) at (4,0) {};
    \path[->,>=stealth']
        (00) edge (01)
        (00) edge (11)
        (10) edge (11)
        (03) edge (04)
        (03) edge (14)
        (13) edge (14)
        ;
\end{tikzpicture}
\end{center}
and note that there is no path of length $n$ from 1 to $n$ in $\T_{n+2}$. Hence, $\T_{n+2}$ is not $n$-braided.

$(3) \Rightarrow (2)$ We show the contraposition. Let $\A$ be a structure satisfying $\HM n$ witnessed by polymorphisms $p_1,\dots,p_n$ and let
\begin{itemize}
    \item $\phi_1(x_1,y_1),\dots,\phi_n(x_n,y_n)$ with
    \item $(a_0,\dots,a_n;b_0,\dots,b_n)$
\end{itemize} be an $n$-braid of $\A$. For $1\leq i<n$ define $c_i\in A$ as $p_i(a_i,a_i,b_i)$. Note that $c_i=p_{i+1}(a_i,b_i,b_i)$. It is easy to see that
\[\A\models \phi_1(a_0,c_1)\AND\phi_2(c_1,c_2)\AND\dots\AND \phi_n(c_{n-1},b_n).\]

$(2)\Rightarrow(3)$ We show the contraposition. Let $\A$ be an $n$-braided finite core structure with elements $a_1,\dots,a_m$. Then, by Theorem~\ref{thm:braidedEqualsMultibraided}, the structure $\A$ is $(m^2,\dots,m^2)$-multi-braided. 
Let 
\begin{align*}
    \boldsymbol t&\coloneqq (a_1,a_1,\dots,a_1,a_2,a_2,\dots,a_2,\dots,a_m,a_m,\dots,a_m)\text{ and}\\
    \boldsymbol s&\coloneqq (a_1,a_2,\dots,a_m,a_1,a_2,\dots,a_m,\dots,a_1,a_2,\dots,a_m).
\end{align*}
Consider the $(m^2,\dots,m^2)$-multi-braid (with $n+1$ many $m^2$'s) generated by 
$(\b t,\b s,\b t,\dots;\b s,\b t,\b s,\dots )$. This multi-braid is satisfied by some path $(\b c_0,\dots,\b c_n)$ and
there are polymorphisms $p_1,\dots,p_n$ of $\A$ such that 
$p_i(t_k,s_k,s_k)=c_{i-1,k}$ and $p_i(s_k,s_k,t_k)=c_{i,k}$ for all $i\in\{1,\dots,n\}$ and $k\in\{1,\dots,m\}$.
Hence, $p_1(t_i,s_i,s_i)=t_i$, $p_1(s_i,s_i,t_i)=p_2(s_i,t_i,t_i)$, $\dots$ for all $i\in\{1,\dots,m^2\}$.
By the choice of $\b s$ and $\b t$ we have that for every $a,b\in A$ there is an $i$ with $s_i=a$ and $t_i=b$. Therefore, $\{p_1,\dots,p_n\}\models\HM n$ as desired.


$(3) \Rightarrow (1)$ Let $\A$ be a structure and let $\phi_1(x_1,y_1),\dots,\phi_n(x_n,y_n)$ with $(a_0,\dots,a_n,b_0,\dots,b_n)$ be an $n$-braid of $\A$ that witnesses that $\A$ is not $n$-braided. Let $\Hb$ be the $(n+1)$-st pp-power of $\A$ given by
\begin{center}
    \begin{tikzpicture}[scale = 1]
    \node[var-b,label=above:$a_0$] (10) at (0,1) {};
    \node[var-b,label=below:\strut$x_0$] (00) at (0,0) {};
    \node[var-b,label=above:$y_1$] (11) at (1,1) {};
    \node[var-b,label=below:\strut$x_1$] (01) at (1,0) {};
    \node[] (12) at (2,1) {$\dots$};
    \node[] (02) at (2,0) {$\dots$};
    \node[var-b,label=above:$y_{n-1}$] (13) at (3,1) {};
    \node[var-b,label=below:\strut$x_{n-1}$] (03) at (3,0) {};
    \node[var-b,label=above:$y_n$] (14) at (4,1) {};
    \node[var-b,label=below:\strut$b_n$] (04) at (4,0) {};
\tikzset{decoration={snake,amplitude=.25mm,segment length=1.3mm,post length=0.7mm,pre length=0.7mm}}
    \draw[decorate,->,>=stealth'] (00) -- node[above] {\small$1$} (01);
    \draw[decorate,->,>=stealth'] (00) -- node[above] {\small$1$} (11);
    \draw[decorate,->,>=stealth'] (10) -- node[above] {\small$1$} (11);
    \draw[decorate,->,>=stealth'] (03) -- node[above] {\small$n$} (04);
    \draw[decorate,->,>=stealth'] (03) -- node[above] {\small$n$} (14);
    \draw[decorate,->,>=stealth'] (13) -- node[above] {\small$n$} (14);
\end{tikzpicture}
\end{center}
where \begin{tikzpicture}[scale=0.7]
    \node[var-f,label=above:$x$] (x1) at (0,0) {};
    \node[var-f,label=above:$y$] (x2) at (1,0) {};
\tikzset{decoration={snake,amplitude=.25mm,segment length=1.3mm,post length=0.7mm,pre length=0.7mm}}
    \draw[decorate,->,>=stealth'] (x1) -- node[above] {\small$i$} (x2);
\end{tikzpicture} 
is an abbreviation for 
the relation $\{(a,b)\mid \A\models\phi_i(a,b)\}$. Note that this relation is pp-definable.
For $-1\leq i\leq n$ define \[v_i\coloneqq (a_0,\dots,a_i,b_{i+1},\dots,b_n).\]
Note that the substructure of $\H$ induced by $v_{-1},\dots,v_n$ is isomorphic to $\T_{n+2}$. If $\H$ would have a path $u_0,\dots,u_{n+2}$ then \[a_0=u_{1,0},u_{2,1},\dots,u_{n+1,n}=b_n\] would contradict our witness for $\A$ not being $n$-braided. Hence, $\Hb$ has no path of length $n+2$ and therefore has a homomorphism to $\T_{n+2}$. We conclude that $\Hb$ is homomorphically equivalent to $\T_{n+2}$ and that $\A\ppleq\T_{n+2}$.
\end{proof}
As a consequence of this theorem we obtain that $\T_{n+2}$ is a blocker for $\HM {n}$.
To separate the $\T_n$'s we still need the following lemma.

\begin{lemma}
Let $n\in\N^+$\!. Then $\T_{n+2}\models\HM{n+1}$.
\end{lemma}
\begin{proof}
Example 5.4 in~\cite{NivenDigraphCSP} shows that the tournament $\T_{n+2}$ is congruence $(n+2)$-permutable. Theorem~2 in~\cite{HagemannMitschke} states that a structure is  congruence $(n+2)$-permutable if and only if it satisfies $\HM{n+1}$.
\end{proof}

Consequently, $\T_3\ppstrictlyGreater\T_4\ppstrictlyGreater\cdots$ is an infinite descending chain. Next we want to find the lower bound of this chain. Recall from Chapter~\ref{cha:DatalogIntro} the structure $\Ord=(\{0,1\};\leq,\{0\},\{1\})$. 
\begin{lemma}\label{lem:OrdCanPPconstructTn}
For any $n\in\N^+$ we have that $\Ord\ppleq\T_n$.
\end{lemma}
\begin{proof}
Clearly, $\Ord\ppleq \T_1,\T_2$. 
Let $n\in\N^+$. Consider the following $n$-braid 
\begin{center}
    \begin{tikzpicture}[scale = 0.5]
    \node[var-b,label=above:1] (10) at (0,1) {};
    \node[var-b,label=below:0] (00) at (0,0) {};
    \node[var-b,label=above:1] (11) at (1,1) {};
    \node[var-b,label=below:0] (01) at (1,0) {};
    \node[] (12) at (2,1) {$\dots$};
    \node[] (02) at (2,0) {$\dots$};
    \node[var-b,label=above:1] (13) at (3,1) {};
    \node[var-b,label=below:0] (03) at (3,0) {};
    \node[var-b,label=above:1] (14) at (4,1) {};
    \node[var-b,label=below:0] (04) at (4,0) {};
    \path[->,>=stealth']
        (00) edge (01)
        (00) edge (11)
        (10) edge (11)
        (03) edge (04)
        (03) edge (14)
        (13) edge (14)
        ;
\end{tikzpicture}
\end{center}
Since there is no directed path from $1$ to $0$ in $\Ord$ we have that $\Ord$ is not $n$-braided. Hence, by Theorem \ref{thm:TnIsHMblocker}, $\Ord\ppleq\T_{n+2}$.
\end{proof}

What makes $\Ord$ interesting is the following lemma, which shows that $\Ord$ is the greatest lower bound of the chain $\T_3\ppstrictlyGreater\T_4\ppstrictlyGreater\cdots$. This fact was originally discovered by Hobby, McKenzie, Larose, and Tesson~\cite{LaroseTesson,HobbyMcKenzie}.

\begin{lemma}\label{lem:noHMimpliesACanppDefineOrder}
Let $\A$ be a finite core structure with $\A\not\models\HM {|A|\cdot(|A|-1)^2}$. Then $\A$ can, using constants, pp-define a structure that is homomorphically equivalent to $\Ord$; in particular, $\A\ppleq\Ord$.
\end{lemma}
\begin{proof}
Let $k\coloneqq |A|\cdot(|A|-1)^2$.
Since $\A$ is a core we know from Theorem~\ref{thm:TnIsHMblocker} that $\A$ is not $k$-braided. Note that for any $k$-braid that witnesses that $\A$ is not $k$-braided all the tuples $(a_i,b_i)$ must satisfy $a_i\neq b_i$. Hence, such a braid can have at most $|A|\cdot(|A|-1)$ many different tuples. Since the $k$-braid has $k+1$ many tuples there must be a tuple $(a,b)\in A^2$ such that at least $|A|$ many $(a_i,b_i)$ are equal to $(a,b)$. Let $n\coloneqq|A|-1$. Then there is an $(a,b)\in A^2$ and an $n$-braid $\phi_1(x_1,y_1),\dots,\phi_n(x_n,y_n)$ with $(a_0,\dots,a_n,b_0,\dots,b_n)$ that witnesses that $\A$ is not $n$-braided such that 
\[(a_0,\dots,a_n,b_0,\dots,b_n)=(a,\dots,a,b,\dots,b).\]
Let $\Hb$ be the digraph with constants $a$ and $b$ and pp-definable edge relation 
\begin{align*}
    \{(c,d)\mid \A\models\phi_1(c,d)\AND\dots\AND\phi_n(c,d) \}.
\end{align*}

Observe that $\Hb$ has the edges $(b,b),(b,a)$, and $(a,a)$, but not $(a,b)$. Hence, the substructure of $\Hb$ induced by $\{a,b\}$ is isomorphic to $\Ord$. Assume there is a directed path $a,c_1,\dots,c_m,b$ in $\Hb$. By taking shortcuts and adding $b$'s at the end we can assume that this path is of length $|A|-1$ and therefore $m=|A|-2=n-1$. 
Hence, $a,c_1,\dots,c_m,b$ 
contradicts the witness that $\A$ is not $n$-braided.
Therefore, there is no directed path from $a$ to $b$ in $\Hb$ and $\Hb$ is homomorphically equivalent to $\Ord$.
\end{proof}

As an immediate consequence of Theorem~\ref{thm:TnIsHMblocker} and Lemmata~\ref{lem:noHMimpliesACanppDefineOrder} and~\ref{lem:OrdCanPPconstructTn} we obtain the following well-known corollary.
\begin{corollary}\label{cor:OrdisLeastLowerBoundForTn}
Let $\A$ be a finite structure. Then 
\begin{align*}
    &\text{$\A\ppleq\T_n$ for all $n$}&& \text{if and only if} &&\A\leq\Ord.
\end{align*}
\end{corollary}


Since we are only working with digraphs the next question is whether the structure $\Ord$ can be represented by  a digraph.
It turns out there is such a digraph that even has just 4 vertices.
Define $\T_4^{\setminus(03)}\coloneqq(\{0,1,2,3\},<\setminus\{(0,3)\})$.

\begin{lemma}\label{lem:OrdDGppeqOrd}
We have that $\T_4^{\setminus(03)}\ppeq\Ord$.
\end{lemma}
\begin{proof}
First we show $\T_4^{\setminus(03)}\ppleq\Ord$.
Note that the constants $1$ and $2$ are pp-definable in $\T_4^{\setminus(03)}$ and consider the digraph $\Hb$ with constants 1 and 2 and pp-definable edge relation
\[\{(a,b)\mid \exists z_1,z_2\in\{0,1,2,3\}
\text{ such that }
a\toEdge z_1\fromEdge z_2\toEdge b
\text{ in }\T_4^{\setminus(03)}\}.\]
It is easy to see that $\Hb$ and $\Ord$ are homomorphically equivalent.

Next we show $\Ord\ppleq \T_4^{\setminus(03)}$.
The third pp-power of $\Ord$ given by
\begin{center}
\begin{tikzpicture}[scale=0.5,baseline=1mm]
    \node[var-f,label=below:\strut$x_1$] (0) at (0,0) {};
    \node[var-f,label=below:\strut$1$] (1) at (1,0) {};
    \node[var-f,label=below:\strut$x_3$] (2) at (2,0) {};
    \node[var-f,label=above:$0$] (4) at  (0,1) {};
    \node[var-f,label=above:$y_2$] (5) at (1,1) {};
    \node[var-f,label=above:$y_3$] (6) at (2,1) {};
    \path[->,>=stealth']
        (0) edge (5)
        (6) edge (0)
        (5) edge (2)
        ;
\end{tikzpicture}
\hspace{0.2cm}has the subgraph\hspace{0.2cm}
\begin{tikzpicture}[scale=0.7,baseline=10mm]
    \node[var-b,label=right:$000$] (3) at (0,3) {};
    \node[var-b,label=right:$010$] (2) at (0,2) {};
    \node[var-b,label=right:$011$] (1) at (0,1) {};
    \node[var-b,label=right:$111$] (0) at (0,0) {};
    \path[->,>=stealth']
        (0) edge (1)
        (1) edge (2)
        (2) edge (3)
        (0) edge[bend left] (2)
        (1) edge[bend left] (3)
        ;
\end{tikzpicture}
\end{center}
induced by the non-isolated vertices.  
Hence, $\Ord\ppeq \T_4^{\setminus(03)}$.
\end{proof}



\section{Semicomplete Digraphs}
We have classified all transitive tournaments. In this section we show that every semicomplete digraph, i.e., digraphs where between any two vertices there is either a directed or an undirected edge, can either pp-construct $\K_3$, is a transitive tournament, contains a single two-cycle, or contains a single three-cycle. These families of digraphs where already studied by Jackson,  Kowalski, and Niven in in~\cite{NivenDigraphCSP}. We were not able to fully classify the pp-constructability order for the latter two families. 

Let $\G$ be a digraph. Then we call $\G$ \emph{strongly connected} if any two distinct vertices of $\G$ are connected via a directed path. 
Next we show a well-known result.
\begin{lemma}
Let $\G$ be a semicomplete digraph. Then there is a unique map $\height_{\G}\colon V(\G)\to\N$ such that 
\begin{enumerate}
    \item the image of $\height_{\G}$ is equal to $\{0,\dots,n\}$ for some $n\in\N$,
    \item $\height_{\G}(u)<\height_{\G}(v)$ implies $u\toEdge v$, and
    \item $\height_{\G}^{-1}(n)$ is empty or induces a subgraph of $\G$ which is strongly connected.
\end{enumerate}
\end{lemma}
\begin{proof}
We show the claim by induction on the number of vertices in $\G$. 
Let $V_1,\dots,V_n\subseteq V(\G)$ pairwise distinct such that each $V_i$ induces a maximal strongly connected subgraphs  of $\G$. Clearly, $V_1,\dots,V_n$ form a partition of $V(\G)$. 
Define the digraph $\H$ with 
\begin{align*}
    V(\H)&=\{V_1,\dots,V_n\}\text{ and}\\
    E(\H)&=\{(V_i,V_j)\mid i\neq j\text{ and $\exists\ v_i\in V_i$, $v_j\in V_j$ with $(v_i,v_j)\in E(\G)$}\}.
\end{align*}
Since $V_1,\dots,V_n$ are maximal we have that $\H$ is acyclic. That $\G$ is semicomplete implies that $\H$ is also semicomplete. Every acyclic semicomplete digraph with $n$ vertices is a transitive tournament and therefore isomorphic to $\T_n$.
Let $h\colon\H\to\T_n$ be an isomorphism. For any $v\in V(\G)$ let $[v]$ denote the unique $V_i$ that contains $v$. 
Define
\begin{align*}
    \height_{\G}\colon V(\G)&\to\N\\
    v&\mapsto h([v]) 
\end{align*} 
Clearly, $\height_{\G}$ satisfies 1., 2., and 3. 
It is easy to see that the map $\height_{\G}$ is also unique.
\end{proof}

If $\G$ is clear from the context we just write $\height$ instead of $\height_{\G}$. For example $\height_{\T_n}(k)=k$ for all $0\leq k< n$.
Define $f\colon \Z\cup\{a,b,c\}\to\Z$, $a,b,c\mapsto 0$, and $n\mapsto n$ for $n\in\Z$.
Let $n,k\geq0$. Define the semicomplete digraphs
\begin{itemize}
    \item $\nCk nk2$ that has vertices $-n,\dots,-1,a,b,1,\dots,k$, an edge from $u$ to $v$ whenever $f(u)<f(v)$, and additionally the edges $a\toEdge b\toEdge a$ 
     and 
    \item $\nCk nk3$ that has vertices $-n,\dots,-1,a,b,c,1,\dots,k$, an edge from $u$ to $v$ whenever $f(u)<f(v)$, and additionally the edges $a\toEdge b\toEdge c\toEdge a$. 
\end{itemize}
See Figure~\ref{fig:C2nkAndC3nk} for examples. 
\begin{figure}
    \centering
    \begin{tikzpicture}[scale=0.7]
    \node[var-b,label=right:$2$] (3) at (0,3) {};
    \node[var-b,label=right:$1$] (2) at (0,2) {};
    \node[var-b,label=left:$a$] (1a) at (-0.5,1) {};
    \node[var-b,label=right:$b$] (1b) at (0.5,1) {};
    \node[var-b,label=right:$-1$] (0) at (0,0) {};
    \path[->,>=stealth']
        (0) edge (1a)
        (0) edge (1b)
        (1a) edge[<->] (1b)
        (1a) edge (2)
        (1b) edge (2)
        (2) edge (3)
        (0) edge (2)
        (1a) edge[bend left] (3)
        (1b) edge[bend right] (3)
        (0) edge[bend left=60] (3)
        (0) edge[bend right=60] (3)
        ;
    \end{tikzpicture}
    \hspace{2cm}
    \begin{tikzpicture}[scale=0.7]
    \node[var-b,label=right:$1$] (3) at (0,3) {};
    \node[var-b,label=right:$c$] (1c) at (0,2) {};
    \node[var-b,label=left:$a$] (1a) at (-0.5,1) {};
    \node[var-b,label=right:$b$] (1b) at (0.5,1) {};
    \path[->,>=stealth']
        (1a) edge (1b)
        (1b) edge (1c)
        (1c) edge (1a)
        (1a) edge[bend left] (3)
        (1b) edge[bend right] (3)
        (1c) edge (3)
        ;
    \end{tikzpicture}
    \caption{The digraphs $\nCk 122$ (left) and $\nCk 013$ (right).}
    \label{fig:C2nkAndC3nk}
\end{figure}
Note that $f(v)+n=\height(v)$ for every node $v$ in $\nCk nki$.
Recall from Theorem~\ref{thm:BartoKozikNiven} that any smooth digraph that is a core and is not a disjoint union of cycles is in the same pp-constructability class as $\K_3$. 
Using this we partition semicomplete digraphs into four classes.

\begin{theorem}
Let $\G$ be a semicomplete digraph. Then $\G$ is
\begin{enumerate}
    \item  isomorphic to $\T_n$ for some $n\geq 1$,
    \item  isomorphic to $\nCk nk2$ for some $n,k\geq 0$,
    \item  isomorphic to $\nCk nk3$ for some $n,k\geq 0$, or
    \item  in the same pp-constructability class as $\K_3$.
\end{enumerate}
\end{theorem}
\begin{proof}
Let $\G$ be a semicomplete digraph. 
If for every $n\in\N$ the set $\height^{-1}(n)$ contains at most one element, then $\G$ is isomorphic to $\T_{|V(\G)|}$ and we are in case (1). Hence, assume that there is an $n\in\N$ such that $\height^{-1}(n)$ contains more than one element. Let $n_{\min}\in\N$ be minimal such that $|\height^{-1}(n)|\geq 2$ and $n_{\max}\in\N$ be maximal such that $|\height^{-1}(n)|\geq 2$. Using the formula
\begin{center}
\begin{tikzpicture}[scale=0.5]
    \node[var-b] (a0) at (-4,0) {};
    \node[var-b] (a1) at (-3,0) {};
    \node at (-2,0) {$\dots$};
    \node[var-b] (a2) at (-1,0) {};
    \node[var-f,label=above:$x$] (x) at (0,0) {};
    \node[var-f,label=above:$y$] (y) at (1,0) {};
    \node[var-b] (b0) at (2,0) {};
    \node at (3,0) {$\dots$};
    \node[var-b] (b1) at (4,0) {};
    \node[var-b] (b2) at (5,0) {};
    \path[->,>=stealth']
        (a0) edge (a1)
        (a2) edge (x)
        (x) edge (y)
        (y) edge (b0)
        (b1) edge (b2)
        ;
\end{tikzpicture}
\end{center}
we see that the smooth subgraph $\H$ of $\G$ induced by $\height^{-1}(\{n_{\min},\dots,n_{\max}\})$ is pp-constructable from $\G$. Note that $\H$ is semicomplete. Hence, $\H$ is a core. There are several cases to consider. If $\H$ is not a disjoint union of cycles then, by Theorem~\ref{thm:BartoKozikNiven}, $\H$ can pp-construct $\K_3$ and we are in case (4). If $\H$ is a disjoint union of cycles, then, since $\H$ is semicomplete, $\H$ is isomorphic to $\C_2$ or $\C_3$. In both cases $n_{\min}=n_{\max}$ and $\G$ is isomorphic to either $\nCk nk2$ or $\nCk nk3$ for some $n,k\geq 0$.
\end{proof}

\begin{lemma}\label{lem:nCkleqT}
Let $n,k\geq 0$ and $i\in\{2,3\}$. Then the following statements are true.
\begin{enumerate}
    \item $\nCk nki=\nCk kni$.
    \item $\nCk {n+1}ki\ppleq \nCk nki$ and $\nCk n{k+1}i\ppleq \nCk nki$; in particular $\nCk nki\ppleq \C_i$.
    \item $\nCk nki\ppleq\T_{2\cdot\max(n,k)+2}$.
    \item $\nCk nk2\models\Sigma_3$ and $\nCk nk3\models\Sigma_2$.
\end{enumerate}
\end{lemma}
\begin{proof}
Let $n,k\geq0$ and $i\in\{2,3\}$. It is easy to see that (1) and (2) hold. For (3) assume without loss of generality that $k\geq n$. Consider the following $2k$-braid
\begin{center}
    \begin{tikzpicture}[scale = 1.0]
    \node[var-b,label=above:$a$] (10) at (0,1) {};
    \node[var-b,label=below:\strut$k$] (00) at (0,0) {};
    \node[var-b,label=above:$b$] (11) at (1,1) {};
    \node[var-b,label=below:\strut$k-1$] (01) at (1,0) {};
    \node[]  at (2,1) {$\dots$};
    \node[]  at (2,0) {$\dots$};
    \node[var-b,label=above:$c$] (13) at (3,1) {};
    \node[var-b,label=below:\strut1] (03) at (3,0) {};
    \node[var-b,label=above:$d$] (14) at (4,1) {};
    \node[var-b,label=below:\strut$c$] (04) at (4,0) {};
    
    \node[var-b] (m1) at (4,0.7) {};
    \node[var-b] (m0) at (4,0.3) {};
    
    \node[var-b,label=above:$1$] (15) at (5,1) {};
    \node[var-b,label=below:\strut$d$] (05) at (5,0) {};
    \node[]  at (6,1) {$\dots$};
    \node[]  at (6,0) {$\dots$};
    \node[var-b,label=above:$k-1$] (16) at (7,1) {};
    \node[var-b,label=below:\strut$a$] (06) at (7,0) {};
    \node[var-b,label=above:$k$] (17) at (8,1) {};
    \node[var-b,label=below:\strut$b$] (07) at (8,0) {};
    \path[<-,>=stealth']
        (00) edge (01)
        (00) edge (11)
        (10) edge (11)
        (03) edge (04)
        (03) edge (14)
        (13) edge (14)
        ;
\path[->,>=stealth']
        (04) edge (05)
        (04) edge (15)
        (14) edge (15)
        (06) edge (07)
        (06) edge (17)
        (16) edge (17)
        ;
\path   (04) edge (m0)
        (14) edge (m1)
        ;

\end{tikzpicture}
\end{center}
where $c=a$ and $d=b$ if $k$ is odd and $c=b$ and $d=a$ if $k$ is even. Note that there is no vertex $v$ such that $v$ has an undirected edge to some neighbour and  $a\stackrel{k}{\fromEdge}v\stackrel{k}{\toEdge}b$. Hence, $\nCk nki$ is not $2k$-braided and by Theorem~\ref{thm:TnIsHMblocker} we have $\nCk nki\ppleq \T_{2\cdot k +2}$.
To show (4) recall from Lemma~\ref{lem:pcSatisfyPclc} that $\C_2\models\Sigma_3$. Let $f$ be a polymorphism of $\C_2$ satisfying $\Sigma_3$. Observe that the map
\begin{align*}
    (x,y,z)\mapsto
    \begin{cases}
    f(x,y,z)&\text{if }\{x,y,z\}=\{a,b\}\\
    \max(\{x,y,z\}\setminus\{a,b\})&\text{if } \{x,y,z\}\neq\{a,b\}
    \end{cases}
\end{align*}
is a polymorphism of $\nCk nk2$ satisfying $\Sigma_3$. We can show $\nCk nk3\models\Sigma_2$ analogously.
\end{proof}
From the previous lemma we know that $\nCk nki\not\models\HM{2\cdot\max(n,k)}$.
That $\nCk nki\models\HM{2\cdot\max(n,k)+1}$ was shown by Jackson, Kowalski, and Niven in Theorem~8.1 in \cite{NivenDigraphCSP}. Hence, if $n\leq k$ then $\nCk n{k+1}i$ lies strictly below $\nCk nki$. Using (4) in Lemma~\ref{lem:nCkleqT} we can separate $\nCk nk2$ and $\nCk nk3$ for any $n,k\geq0$. 
We can show the following perhaps surprising fact.
\begin{lemma}\label{lem:0C1iPPeq1C1i}
We have that $\nCk 01i\ppeq\nCk 11i$ for $i\in\{2,3\}$.
\end{lemma}
\begin{proof}
Let $\H_i$ be the pp-power of $\nCk 01i$ given by $\Phi_i$:
\begin{center}
\begin{tikzpicture}[scale=0.5]
\node[var-f,label=below:$x_1$] (x1) at (0,0) {};
    \node[var-f,label=below:$x_2$] (x2) at (0+1,0) {};
    \node[var-b] (a) at (0+0.5,1) {};
    \node[var-f,label=above:$y_1$] (y1) at (0,2) {};
    \node[var-f,label=above:$y_2$] (y2) at (0+1,2) {};
    \path[<-,>=stealth']
        (x1) edge (x2)
        (x1) edge (a)
        (y2) edge (y1)
        (y2) edge (a)
        ;
    \node at (0.5,-2) {$\Phi_{2}$};
\end{tikzpicture}
\hspace{20mm}
\begin{tikzpicture}[scale=0.5]
    \node[var-f,label=below:$x_1$] (x1) at (0,0) {};
    \node[var-f,label=below:$x_2$] (x2) at (0+1,0) {};
    \node[var-f,label=below:$x_3$] (x3) at (0+2,0) {};
    \node[var-b] (a) at (0+1,1) {};
    \node[var-f,label=above:$y_1$] (y1) at (0,2) {};
    \node[var-f,label=above:$y_2$] (y2) at (0+1,2) {};
    \node[var-f,label=above:$y_3$] (y3) at (0+2,2) {};
    \path[->,>=stealth']
        (x2) edge (x3)
        (x3) edge[bend right] (x1)
        (x1) edge (a)
        (y1) edge (y2)
        (y2) edge (y3)
        (y3) edge (a)
        ;
    \node at (1,-2) {$\Phi_{3}$};
\end{tikzpicture}
\end{center}
Observe that $\{1a,ab,ba,a1\}$ induces a subgraph of $\H_2$ that is the core of $\H_2$ and is isomorphic to $\nCk 112$. Similarly, $\{1bc,abc,bca,cab,ab1\}$ induces a subgraph of $\H_3$ that is the core of $\H_3$ and is isomorphic to $\nCk 113$. 
\end{proof}
Note that this construction can be generalized to show that $\nCk 01n\ppeq \nCk 11n$ for any $n\geq2$, where $\nCk 01n$ and $\nCk 11n$ are defined in the obvious way. However, we cannot simply generalize it to show $\nCk 0ki\ppeq\nCk kki$. Whether we can separate $\nCk 0ki$, $\nCk 1ki, \dots, \nCk kki$ is not clear. So for now the poset of the digraphs $\nCk nki$'s looks as in Figure~\ref{fig:nCkPoset} and we have the following open problem.

\begin{question}
Let $k\geq2$ and $i\in\{2,3\}$. Is $\nCk 0ki\ppeq\nCk 1ki\ppeq \dots\ppeq \nCk kki$?
\end{question}

\begin{figure}
    \centering
    \begin{tikzpicture}[scale=1.3]
    
    \node (00) at (0,1) {$\C_i$};
    \node (01) at (0,0) {$\nCk 01i\ppeq\nCk 11i$};
    \node (02) at (0,-1) {$\nCk 02i$};
    \node (12) at (1,-2) {$\nCk 12i$};
    \node (22) at (2,-3) {$\nCk 22i$};
    
    \node (03) at (0,-2) {$\nCk 03i$};
    \node (13) at (1,-3) {$\nCk 13i$};
    \node (23) at (2,-4) {$\nCk 23i$};
    \node (33) at (3,-5) {$\nCk 33i$};
    
    \node (04) at (0,-3) {$\nCk 04i$};
    \node (14) at (1,-4) {$\nCk 14i$};
    \node[rotate=135] at (2,-5) {$\dots$};
    
    \node (05) at (0,-4) {$\nCk 05i$};
    \node[rotate=135] at (1,-5) {$\dots$};
    
    \node[rotate=90] at (0,-5) {$\dots$};
    \path[->]
        (01) edge (00)
        (02) edge (01)
        (03) edge (02)
        (04) edge (03)
        (05) edge (04)
        
        (14) edge (04)
        (13) edge (03)
        (12) edge (02)
        
        (23) edge (13)
        (22) edge (12)
        
        (33) edge (23)
        
        (13) edge (12)
        (14) edge (13)
        
        (23) edge (22)
    ;
    
    \end{tikzpicture}
    \caption{Poset of $\nCk nki$'s, where $\nCk nki$ with $n>k$ are omitted.}
    \label{fig:nCkPoset}
\end{figure}



\section{Complexity of Semicomplete 
Digraphs}
In Theorem 8.1 in~\cite{NivenDigraphCSP} it was shown that $\csp(\nCk nki)$ is in L for all $n,k\geq0$, $i\in\{2,3\}$.
In this section we show that $\cocsp(\nCk nki)$ is also in symmetric linear Datalog for all $n,k\geq0,i\in\{2,3\}$. There is an introduction to Datalog in Chapter~\ref{cha:DatalogIntro}.
We will use  Theorem~\ref{thm:symLinDLiffHM} to show that $\cocsp(\nCk nki)$ is in symmetric linear Datalog. To do so we need to show that $\nCk nki\models\HM n$ for some $n$ and that $\cocsp(\nCk nki)$ is in linear Datalog.
We already know from Theorem 8.1 in~\cite{NivenDigraphCSP} that  $\nCk nki\models\HM{2\cdot\max(n,k)+2}$. Hence, it only remains to show that $\cocsp(\nCk nki)$ is in linear Datalog.
One of the prominent methods for proving containment in linear Datalog is to show the existence of a near unanimity polymorphism (Theorem~\ref{thm:NUimpliesInLinDL}). However, this method does not work in the case of $\cocsp(\nCk nki)$ if $\max(n,k)\geq2$. There are not many examples of $\cocsp$'s in linear Datalog that do not have a near unanimity polymorphism. Hence, we will prove that $\nCk nki$ has no near unanimity polymorphism if $\max(n,k)\geq2$, although it does not help us prove that $\cocsp(\nCk nki)$ is in linear Datalog. 

\begin{theorem}[Theorem 25 in \cite{FederVardi},
Theorem 8.5.11 in \cite{theBodirsky}]\label{thm:NUiffNaryRelationIsConjuct}
Let $\A$ be a finite structure. Then $\A\models\NU n$ if and only if for every pp-definable $n$-ary relation $R$ in $\A$  we have
\begin{align*}
    R=\{(a_1,\dots,a_n)\mid \forall 1\leq i\leq n\ \exists a\in A.\  (a_1,\dots,a_{i-1},a,a_{i+1}\dots,a_n)\in R\}.
\end{align*}
\end{theorem}

The following lemma was proven (personal communication) by Sebastian Meyer.
\begin{lemma}\label{lem:nCkNoNUs}
Let $n,k\geq0$ and $m\geq 3$. Then $\nCk nki\models\NU m$ if and only if $\max(n,k)\leq 1$.
\end{lemma}
\begin{proof}
The fully symmetric majority operation of $\nCk 012$ that maps $(a,b,1)$ to 1 is a 3-ary near unanimity polymorphism of $\nCk 012$. Hence, $\nCk 012\models\NU m$ for all $m\geq 3$.
Let $f$ be the majority operation of $\nCk 013$ that 
maps tuples $(x,y,z)$ with three different entries to $x$ if $x,y,z\in\{a,b,c\}$ and to 1 otherwise. Then $f$ is a polymorphism of $\nCk 013$ that satisfies $\NU3$. 

Let $m\geq 3$. Any $\nCk nki$ with $\max(n,k)\geq 2$ can pp-construct $\nCk 02i$. Hence, it suffices to show $\nCk02i\not\models\NU m$. We can use the same proof for $i=2$ and $i=3$.
Let $R$ be the $m$-ary relation defined by the pp-formula
\begin{center}
\begin{tikzpicture}[scale=0.5]
\node[var-f,label=below:$x_1$] (00) at (0,0) {};
\node[var-b] (01) at (1,0) {};
\node[var-b] (02) at (2,0) {};
\node at (3,0) {$\dots$};
\node[var-b] (04) at (4,0) {};
\node[var-f,label=below:$x_n$] (05) at (5,0) {};

\node[var-b] (10) at (0.5,1) {};
\node[var-b] (11) at (1.5,1) {};
\node at (3,1) {$\dots$};
\node[var-b] (14) at (4.5,1) {};

\node[var-f,label=above:$x_2$] (20) at (0.5,2) {};
\node[var-f,label=above:$x_3$] (21) at (1.5,2) {};
\node at (3,2) {$\dots$};
\node[var-f,label=above:$x_{n-1}$] (24) at (4.5,2) {};
\path[->,>=stealth]
(00) edge (10)
(01) edge (10)
(10) edge (20)

(01) edge (11)
(02) edge (11)
(11) edge (21)

(04) edge (14)
(05) edge (14)
(14) edge (24)
;

\end{tikzpicture}
\end{center}
Observe that 
\begin{align*}
    (b,1,1,\dots,1,a)&\in R\\
    (a,2,1,\dots,1,a)&\in R\\
    (a,1,2,\dots,1,a)&\in R\\
    \vdots\\
    (a,1,1,\dots,2,a)&\in R\\
    (a,1,1,\dots,1,b)&\in R\\
\end{align*}
and that $(a,1,1,\dots,1,b)\notin R$. Therefore $R$ does not satisfy the condition in Theorem~\ref{thm:NUiffNaryRelationIsConjuct}. Hence, $\nCk 02i\not\models\NU m$.
\end{proof}

Since we cannot use $\operatorname{NU}$'s to show that $\cocsp(\nCk nki)$ is in linear Datalog  we will use a more direct approach. 
We will introduce a new construction that preserves containment in linear Datalog and reduce the problem to showing that $\cocsp(\C_i)$ is in linear Datalog.
Let $\G=(V,E)$ be a digraph. Define
\begin{align*}
    \G^{\top}&\coloneqq(V\cupdot\{\top\},E\cup\{(v,\top)\mid v \in V\})\text{ and}\\ \G^{\bot}&\coloneqq(V\cupdot\{\bot\},E\cup\{(\bot,v)\mid v \in V\}).
\end{align*}
Given a Datalog program $\prog P=(\{E^{(2)}\},\sigma,\mathcal R,G)$ define the Datalog program $\prog P^{\top}\coloneqq (\{E^{(2)}\},\sigma,\mathcal R^{\top},G)$, where $\mathcal R^{\top}$ is obtained from $\mathcal R$ by replacing each rule 
\begin{align*}
    P(x_1,\dots,x_n)&\vdash\Phi(x_1,\dots,x_n,y_1,\dots,y_m)
    \intertext{in $\mathcal R$ by the rule}
    P(x_1,\dots,x_n)&\vdash\Phi(x_1,\dots,x_n,y_1,\dots,y_m)\AND E(x_1,u_1)\AND\dots\AND E(x_n,u_n)\\&\phantom{{}\vdash{}}\AND E(y_1,v_1)\AND\dots\AND E(y_m,v_m),
\end{align*}
where $u_1,\dots,u_n,v_1,\dots,v_m$ are new variables.
\begin{lemma}\label{lem:CSPGTopvsCSPG}
Let $\G$ and $\H$ be digraphs such that $\H$ has at least one edge. 
Then $\H\to\G^{\top}$ if and only if $\H'\to\G$, where $\H'$ is obtained from $\H$ by removing all sinks. 
\end{lemma}
\begin{proof}
Let $h\colon \H\to\G^{\top}$ be a homomorphism. Since $\top$ is a sink in $\G^{\top}$ we have that $h^{-1}(\top)$ contains only sinks of $\H$. Hence, $h|_{V(\H')}$ is a homomorphism from $\H\to\G$. 
Conversely, if $h'\colon\H'\to\G$ is a homomorphism then the extension of $h'$ to $\H$ that maps all sinks of $\H$ to $\top$ is a homomorphism from $\H\to\G^{\top}$.
\end{proof}

\begin{lemma}\label{lem:DatalogTopvsDatalog}
Let $\G$ be a digraph. If a Datalog program $\prog P$ recognizes $\cocsp(\G)$, then the Datalog program $\prog P^{\top}$ recognizes $\cocsp(\G^{\top})$.
\end{lemma}
\begin{proof}
Let $\prog P=(\{E\},\sigma,R,G)$ be a Datalog program recognizing $\cocsp(\G)$, let $\S$, $\T$ be $(\{E\}\cup\sigma)$-structures with at least one $E$-edge, and let 
\[P(x_1,\dots,x_n)\dashv\Phi\] be a rule of $\prog P^{\top}$. Let $\A$ be the canonical database of $\Phi$ and let $\S', \T', \A'$ be obtained by removing all sinks from $\S, \T, \A$, respectively. Then there is a rule $P(x_1,\dots,x_n)\dashv\Phi'$ in $\prog P$ such that the canonical database of $\Phi'$ is $\A'$. Similar to the proof of Lemma~\ref{lem:CSPGTopvsCSPG} we can show $(\ast)$:
\begin{align*}
    \S&{}\vdash_{\prog P^{\top}}\T\text{ with rule }P(x_1,\dots,x_n)\dashv\Phi\\
    \text{iff }\S'&{}\vdash_{\prog P}\T'\text{ with rule }P(x_1,\dots,x_n)\dashv\Phi'.
\end{align*}
Let $\H$ be a digraph with at least one edge. Then
\begin{align*}
    &\H \in\cocsp(\G^{\top})\\
    \Leftrightarrow{} & \H' \in\cocsp(\G)\tag{Lemma~\ref{lem:CSPGTopvsCSPG}}\\
    \Leftrightarrow{} & \prog P\text{ can derive $G$ on} \H'\\
    \Leftrightarrow{} & \prog P^{\top}\text{ can derive $G$ on} \H.\tag{$\ast$}
\end{align*}
Hence, $\prog P^\top$ recognizes $\cocsp(\G^\top)$.
\end{proof}

Note that $\prog P$ is (symmetric) linear if and only if $\prog P^{\top}$ and $\prog P^{\bot}$ are (symmetric) linear. Hence, as an immediate consequence of Lemma~\ref{lem:DatalogTopvsDatalog} we obtain the following corollary.
\begin{corollary}
Let $\G$ be a digraph. 
If $\cocsp(\G)$ is in ((symmetric) linear) Datalog then $\cocsp(\G^{\top})$ and $\cocsp(\G^{\bot})$ are in ((symmetric) linear) Datalog. 
\end{corollary}

It is easy to verify that $\C_i\models\NU 3$ and $\C_i\models\HM 1$. Hence, by Theorems~\ref{thm:NUimpliesInLinDL} and~\ref{thm:symLinDLiffHM}, we have that $\cocsp(\C_i)$ is in symmetric linear Datalog. By applying the previous corollary many times we obtain the following result. 

\begin{corollary}
Let $\G$ be a semicomplete graph. Then one of the following is true
\begin{enumerate}
    \item $\G\ppeq\T_n$ for some $n\geq1$ and $\cocsp(\G)$ is in symmetric linear Datalog; in this case $\csp(\G)$ is in $\Lclass$,
    \item $\G\ppeq\nCk nki$ for some $n,k\geq0, i\in\{2,3\}$ and $\cocsp(\G)$ is in symmetric linear Datalog; in this case $\csp(\G)$ is in $\Lclass$, or
    \item $\G\ppeq\K_3$; in this case $\csp(\G)$ is $\NP$-complete.
\end{enumerate}
\end{corollary}

\begin{figure}
    \centering
    \begin{tikzpicture}
    
    \node[var-b] (0) at (0:1) {};
    \node[var-b] (1) at (360/7:1) {};
    \node[var-b] (1a) at (360/7:{1.5}) {};
    \node[var-b] (1b) at (360/7:{2}) {};
    \node[var-b] (2) at (2*360/7:1) {};
    \node[var-b] (3) at (3*360/7:1) {};
    \node[var-b] (3a) at (3*360/7:{1.5}) {};
    \node[var-b] (3b) at (3*360/7:{2}) {};
    \node[var-b] (4) at (4*360/7:1) {};
    \node[var-b] (5) at (5*360/7:1) {};
    \node[var-b] (6) at (6*360/7:1) {};
    \path[->,>=stealth']
    (0) edge (1)
    (2) edge (1)
    (1) edge (1a)
    (1a) edge (1b)
    (2) edge (3)
    (4) edge (3)
    (5) edge (4)
    (6) edge (5)
    (3) edge (3a)
    (3a) edge (3b)
    (6) edge (0)
    ;
    \end{tikzpicture}
    \caption{Example of an obstruction for $\csp(\nCk 022)$.}
    \label{fig:ObstructionFor0C22}
\end{figure}
\begin{example}\label{exa:0C22Obstructions}
It is well known that a digraph is in $\csp(\C_2)$ if and only if it contains no orientation of a cycle of odd length.
A digraph $\H$ is not in $\csp(\nCk 022)$ if and only if there is a digraph $\F$ that is an orientation of a cycle of odd length where each vertex $u$ with in-degree two has an extra outgoing path $u\toEdge u_1\toEdge u_2$ such that $\F$ has a homomorphism into $\H$. 
See for example Figure~\ref{fig:ObstructionFor0C22}.
        Note that every vertex on the cycle of such a digraph $\F$ has an outgoing path of length at least two.  
\end{example}


\section{Transitive Tournaments and Cycles}
We understand $\T_1\ppstrictlyGreater\T_2\ppstrictlyGreater\cdots\ppstrictlyGreater\Ord$ and all disjoint unions of cycles. In this section look at products of these two types of digraphs and classify the pp-constructability order on these products. As one result we obtain lower covers of the submaximal directed cycles of $\DGPoset$, i.e., the directed cycles $\C_p$ where $p$ is prime. 

\begin{lemma}\label{lem:meetTnCisTnxC}
Let $n\geq 3$ and $\C$ be a disjoint union of cycles. 
Then the meet of $\T_n$ and $\C$ exists and can be represented by $\T_n\times\C$.
\end{lemma}
\begin{proof}
Assume without loss of generality that $\C$ is a core.
If $\A\ppleq\T_n$ and $\A\leq\C$, then $\A\ppleq\T_n\times\C$ for any finite structure $\A$.
Let $\H_0$, $\H_1$, and $\H_2$ be the first pp-power of $\T_n\times\C\ppleq\C$ respectively given by
\begin{center}
    \begin{tikzpicture}[scale=0.5]
    \node[var-f,label=left:$x$] (x) at (0,-1) {};
    \node[var-b] (a) at (1,-1) {};
    \node[var-b] (b) at (0,0) {};
    \node[var-f,label=right:$y$] (y) at (1,0) {};
    \path[->,>=stealth']
        (x) edge (b)
        (a) edge (b)
        (a) edge (y)
        ;
        \node at (0.5,-2) {$\Phi_0$};
    \end{tikzpicture}
    \hspace{2cm}
    \begin{tikzpicture}[scale=0.5]
    \node[var-f,label=left:$x$] (x) at (0,-1) {};
    \node[var-b] (a) at (0,1) {};
    \node[var-b] (b) at (0,0) {};
    \node[var-f,label=right:$y$] (y) at (1,0) {};
    \path[->,>=stealth']
        (x) edge (b)
        (b) edge (a)
        (y) edge (a)
        ;
        \node at (0.5,-2) {$\Phi_1$};
    \end{tikzpicture}
    \hspace{2cm}
    \begin{tikzpicture}[scale=0.5]
    \node[var-f,label=left:$x$] (x) at (0,-1) {};
    \node[var-b] (a) at (1,-1) {};
    \node[var-b] (b) at (1,-2) {};
    \node[var-f,label=right:$y$] (y) at (1,0) {};
    \path[->,>=stealth']
        (b) edge (x)
        (b) edge (a)
        (a) edge (y)
        ;
        \node at (0.5,-3) {$\Phi_2$};
    \end{tikzpicture}
\end{center}
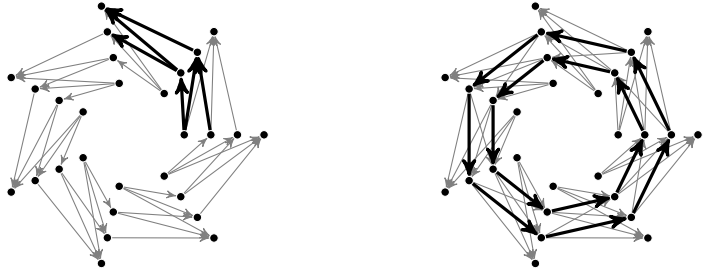
\begin{figure}
    \centering
    \begin{tikzpicture}[scale=0.7]
        \def \n {7}
        \def \radius {1}
        
        \foreach \s in {1,...,\n}
        {
        \node[var-b] (0\s) at ({360/\n * (\s - 1)}:\radius) {};
        \node[var-b] (1\s) at ({360/\n * (\s - 1)}:{\radius+0.5}) {};
        \node[var-b] (2\s) at ({360/\n * (\s - 1)}:{\radius+1}) {};
        \node[var-b] (3\s) at ({360/\n * (\s - 1)}:{\radius+1.5}) {};
        }
        \pgfmathtruncatemacro\nn{\n-1}
        \foreach \s in {1,...,\nn}
        {
         \pgfmathtruncatemacro\ss{\s+1}
        \path[gray,->,>=stealth']
            (0\s) edge (1\ss)
            (0\s) edge (2\ss)
            (0\s) edge (3\ss)
            (1\s) edge (2\ss)
            (1\s) edge (3\ss)
            (2\s) edge (3\ss)
            ;
        }
        
\path[->,>=stealth',gray]
            (0\n) edge (11)
            (0\n) edge (21)
            (0\n) edge (31)
            (1\n) edge (21)
            (1\n) edge (31)
            (2\n) edge (31)
            ;
            
\path[->,>=stealth',very thick]
            (01) edge (12)
            (01) edge (22)
            (11) edge (22)
            (12) edge (23)
            (12) edge (33)
            (22) edge (33)
            ;
    \end{tikzpicture}
    \hspace{2cm}
    \begin{tikzpicture}[scale=0.7]
        \def \n {7}
        \def \radius {1}
        
        \foreach \s in {1,...,\n}
        {
        \node[var-b] (0\s) at ({360/\n * (\s - 1)}:\radius) {};
        \node[var-b] (1\s) at ({360/\n * (\s - 1)}:{\radius+0.5}) {};
        \node[var-b] (2\s) at ({360/\n * (\s - 1)}:{\radius+1}) {};
        \node[var-b] (3\s) at ({360/\n * (\s - 1)}:{\radius+1.5}) {};
        }
        \pgfmathtruncatemacro\nn{\n-1}
        \foreach \s in {1,...,\nn}
        {
         \pgfmathtruncatemacro\ss{\s+1}
        \path[gray,->,>=stealth']
            (0\s) edge (1\ss)
            (0\s) edge (2\ss)
            (0\s) edge (3\ss)
            (1\s) edge (1\ss)
            (1\s) edge (2\ss)
            (1\s) edge (3\ss)
            (2\s) edge (1\ss)
            (2\s) edge (2\ss)
            (2\s) edge (3\ss)
            ;
        }
        
\path[gray,->,>=stealth']
            (0\n) edge (11)
            (0\n) edge (21)
            (0\n) edge (31)
            (1\n) edge (11)
            (1\n) edge (21)
            (1\n) edge (31)
            (2\n) edge (11)
            (2\n) edge (21)
            (2\n) edge (31)
            ;
            \path[->,>=stealth',very thick]
            (11) edge (12)
            (12) edge (13)
            (13) edge (14)
            (14) edge (15)
            (15) edge (16)
            (16) edge (17)
            (17) edge (11)
            
            (21) edge (22)
            (22) edge (23)
            (23) edge (24)
            (24) edge (25)
            (25) edge (26)
            (26) edge (27)
            (27) edge (21)
            ;

    \end{tikzpicture}
    \caption{$\T_4\times\C_7$, which is not 2-braided (left) the graph $\H_0$ pp-definable from $\T_4\times\C_7$ (right).}
    \label{fig:TnxC}
\end{figure}
For an example of $\H_0$ see Figure \ref{fig:TnxC} (right). The map $h_t\colon\C\to\H_0, c\mapsto (t,c)$ is an embedding for any $1\leq t\leq n-2$. The map $\pi_2\colon\H_0\to\C,(t,c)\mapsto c$ is a homomorphism. Hence, $\H_0\homeq\C$ and $\T_n\times\C\ppleq\C$.

Now let us show $\T_n\times\C\ppleq\T_n$. Observe that $\T_n\times\C$ is not $(n-2)$-braided, see Figure \ref{fig:TnxC} (left). To apply Theorem~\ref{thm:TnIsHMblocker} we still need to show that $\T_n\times\C$ is a core. 
Note that $h_0$ is an embedding of $\C$ into $\H_1$ and $h_{n-1}$ is an embedding of $\C$ into $\H_2$.
Let $f$ be an endomorphism of $\T_n\times\C\ppleq\T_n$.
Then $f$ must also be an endomorphism of $\H_0$, $\H_1$, and $\H_2$. Hence, $\pi_2\circ f\circ h_t\in\End(\C)=\Emb(\C)$ and $\pi_2\circ f\circ h_t$ is injective. Therefore $f$ restricted to the image of $h_t$ is injective for all $0\leq t< n$. 

Let $f(t_1,c_1)=f(t_2,c_2)$. 
Note that every $(t,c)\in \T_n\times\C$ is contained in exactly one path of length $n-1$. Therefore $\pi_1(t,c)\mapsto t$ is the only homomorphism from $\T_n\mathbin{\times}\C$ to $\T_n$ and $\pi_1\circ f=\pi_1$. Hence
\[t_1=\pi_1(t_1,c_1)=\pi_1(f(t_1,c_1))=\pi_1(f(t_2,c_2))=\pi_1(t_2,c_2)=t_2.\]
Therefore $(t_1,c_1)$ and $(t_2,c_2)$ are in the image of $h_{t_1}$ and $(t_1,c_1)=(t_2,c_2)$ since $f$ restricted to the image of $h_{t_1}$ is injective.
Therefore, $\T_n\times\C_p$ is a core. Since it is not $(n-2)$-braided we can apply Theorem~\ref{thm:TnIsHMblocker} to conclude that $\T_n\times\C_p\ppleq\T_n$.
\end{proof}

Combining the previous lemma with Theorem~\ref{thm:TnIsHMblocker} and Corollary~\ref{cor:PPvsLoopCaBlockerForSa} we obtain the following corollary.
\begin{corollary}
Let $n\geq 3$ and $p$ be a prime. Then $\T_n\times\C_p$ is a blocker for $\{\HM n,\Sigma_p\}$, i.e., any finite structure $\A$ such that $\A\not\models\HM n$ and $\A\not\models\Sigma_p$ can pp-construct $\T_n\times\C_p$.
\end{corollary}

Now we can classify the lower covers of directed cycles of prime lengths.
\begin{theorem}\label{thm:lowerCoversOfSubmaximalCycles}
Let $p$ be a prime. Then the distinct lower covers of $\C_p$ are  
$\T_3\mathbin{\times}\C_p$, $\C_{p\cdot q_1}$, $\C_{p\cdot q_2}$, $\dots$, where $q_1,q_2,\dots$ are the primes not equal to $p$.
\end{theorem}
\begin{proof}
Let $\A\ppstrictlyLess\C_p$. If $\A\not\models\Maltsev$, then $\A\ppleq\T_3$. Hence, by Lemma~\ref{lem:meetTnCisTnxC} $\A\ppleq\T_3\mathbin{\times}\C_p$ as desired. If $\A\models\Maltsev$, then we know from the classification of digraphs with a Maltsev polymorphism in Corollary~\ref{cor:classificationMaltsevPoset} that $\A$ must be below  $\C_{p\cdot q}$ for some prime $q\neq p$.

Again from Corollary~\ref{cor:classificationMaltsevPoset} we obtain that all the $\C_{p\cdot q}$ have different pp-constructability types. Lastly $\T_3\mathbin{\times}\C_p$ is not in the same pp-constructability class as  $\C_{p\cdot q}$ for any prime $q\neq p$.
\end{proof}

We can use our results to give a partial answer to Problem~\ref{pro:lowerCoversOfTthree}.
\begin{corollary}
Let $p$ be a prime. Then $\T_3\mathbin{\times}\C_p$ is a lower cover of $\T_3$.
\end{corollary}
Note that any other lower cover of $\T_3$ must satisfy $\Sigma_p$ for all $p$.
A natural candidate is $\T_4$. Unfortunately it is not a lower cover of $\T_3$ as shown in the following example.
\begin{example}\label{exa:T4isNotACoverOfT3}



Define the following two digraphs.
\begin{center}
\begin{tikzpicture}[scale=0.5]
\node[bullet] (0) at (0,0) {};
\node[bullet,label=right:$a$] (1) at (-1,1) {};
\node[bullet] (2) at (-1,2) {};
\node[bullet] (3) at (0,3) {};
\node[bullet] (4) at (1,1.5) {};
\node[bullet] (5) at (-2,1.5) {};
\node[bullet] (6) at (-2,2.5) {};
\path[->,>=stealth']
(0) edge (1)
(0) edge (4)
(0) edge (3)
(1) edge (2)
(1) edge (5)
(1) edge (6)
(2) edge (3)
(4) edge (3)
(5) edge (6)
;
\node at (-0.5,-1) {$\N_{123}\T_3$};
\end{tikzpicture}\hspace{20mm}
\begin{tikzpicture}[scale=0.5]
\node[bullet,label=left:0] (0) at (0,0) {};
\node[bullet,label=left:3] (1) at (0,3) {};
\node[bullet,label=left:2] (2) at (0,2) {};
\node[bullet,label=right:$4'$] (3) at (1,4) {};
\node[bullet,label=right:$3'$] (4) at (1,3) {};
\node[bullet,label=right:$2'$] (5) at (1,2) {};
\node[bullet,label=left:1] (6) at (-1,1) {};
\path[->,>=stealth']
(0) edge (2)
(0) edge (6)
(2) edge (1)
(4) edge (3)
(5) edge (1)
(5) edge (4)
(6) edge (2)
;
\node at (0.5,-1) {$\T_3^{+112}$};
\end{tikzpicture}
\end{center}



Using the program from Chapter~\ref{cha:4elements} I verified
\begin{itemize}
    \item $\T_4\not\models\HM2$ and $\N_{123}\T_3,\T_3^{+112},\T_3\models\HM2$ and
    \item $\T_4,\N_{123}\T_3,\T_3^{+112}\not\models\GFS2$ and $\T_3\models\GFS2$. Observe that the vertices $1111,1112,1123',1134'$ induce a path of length three in the indicator graph of $\T_3^{+112}$ for $\GFS2$ witnessing that $\T_3^{+112}\not\models\GFS2$.
\end{itemize}
Also using the program I found the three formulas
\begin{center}
\begin{tikzpicture}[scale=0.5]
    \node[var-f,label=below:$x_1$] (x1) at (0,0) {};
    \node[var-f,label=below:$x_2$] (x2) at (0+1,0) {};
    \node[var-f,label=below:$x_3$] (x3) at (0+2,0) {};
    \node[var-f,label=above:$y_1$] (y1) at (0,1) {};
    \node[var-f,label=above:$y_2$] (y2) at (0+1,1) {};
    \node[var-f,label=above:$y_3$] (y3) at (0+2,1) {};
    \path[>=stealth']
        (x1) edge[->] (y1)
        (x2) edge[dashed] (x1)
        (y2) edge[dashed] (x3)
        ;
    \node at (1,-2) {$\Phi_{1}$};
\end{tikzpicture}\hspace{20mm}
\begin{tikzpicture}[scale=0.5]
    \node[var-f,label=below:$a$] (x1) at (0,0) {};
    \node[var-f,label=below:$x_2$] (x2) at (0+1,0) {};
    \node[var-f,label=below:$x_3$] (x3) at (0+2,0) {};
    \node[var-f,label=above:$y_1$] (y1) at (0,1) {};
    \node[var-f,label=above:$y_2$] (y2) at (0+1,1) {};
    \node[var-f,label=above:$y_3$] (y3) at (0+2,1) {};
    \path[>=stealth']
        (x2) edge[dashed] (y1)
        (y2) edge[->] (x3)
        (y3) edge[->] (x3)
        ;
    \node at (1,-2) {$\Phi_{2}$};
\end{tikzpicture}\hspace{20mm}
\begin{tikzpicture}[scale=0.5]
    \node[var-f,label=left:$x$] (x1) at (0,0) {};
    \node[var-f,label=left:$y$] (y1) at (0,1) {};
    \node[var-b] (z) at (0,2) {};
    \path[>=stealth']
        (x1) edge[->] (y1)
        (y1) edge[->] (z)
        ;
    \node at (0,-2) {$\Phi_{3}$};
\end{tikzpicture}
\end{center}
that show 
\[\T_4\ppstrictlyLess \N_{123}\T_3\ppleq \T_3^{+112}\ppstrictlyLess\T_3.\]
Hence, $\T_4$ is not a lower cover of $\T_3$.
\end{example}

The example gives us a candidate for a lower cover of $\T_3$.
\begin{question}\label{que:candidateLowerCoverOfT3}
Is $\T_3^{+112}$ a lower cover of $\T_3$?
\end{question}

We now know that $\T_4$ is not a lower cover of $\T_3$, but it could still be that the following two questions have a positive answer.
\begin{question}
Are there finitely many classes of digraphs above $\T_4$?
\end{question}

\begin{question}
Is $\T_4, \C_2,\C_3,\C_5,\dots$ a maximal antichain?
\end{question}

After this excursion to the fourth topmost layer of $\DGPoset$ let us come back to $\OrdDG$. We have seen that transitive tournaments and disjoint unions of cycles play well together. We will now extend these results to $\OrdDG$ and disjoint unions of cycles.
\begin{lemma}\label{lem:meetOrdCisOrdxC}
Let $\C$ be a disjoint union of cycles. 
Then the meet of $\OrdDG$ and $\C$ exists and can be represented by $\OrdDG\mathbin{\times}\C$.
\end{lemma}
\begin{proof}
Assume without loss of generality that $\C$ is a core.
If $\A\ppleq\OrdDG$ and $\A\leq\C$, then $\A\ppleq\OrdDG\mathbin{\times}\C$ for any finite structure $\A$.
Analogously to the proof of Lemma~\ref{lem:meetTnCisTnxC} we can show that $\OrdDG\mathbin{\times}\C\ppleq\C$ and that $\OrdDG\mathbin{\times}\C$.
Next we show that $\OrdDG\mathbin{\times}\C$ is not $n$-braided for any $n$. Let $c\in\C$ and denote the  successor of a node $d\in\C$ by $d+1$. Consider the following $3n$-braid
\begin{center}
    \begin{tikzpicture}[scale = 1.2]
    \node[var-b,label=above:{\tiny$(1,c)$}] (10) at (0,1) {};
    \node[var-b,label=below:{\tiny$(0,c)$}] (00) at (0,0) {};
    \node[var-b,label=above:{\tiny$(2,c+1)$}] (11) at (1,1) {};
    \node[var-b,label=below:{\tiny$(1,c+1)$}] (01) at (1,0) {};
    \node[var-b,label=above:{\tiny$(3,c+2)$}] (12) at (2,1) {};
    \node[var-b,label=below:{\tiny$(2,c+2)$}] (02) at (2,0) {};
    \node[var-b,label=above:{\tiny$(1,c+3)$}] (13) at (3,1) {};
    \node[var-b,label=below:{\tiny$(0,c+3)$}] (03) at (3,0) {};
    \node[] (14) at (4,1) {$\dots$};
    \node[] (04) at (4,0) {$\dots$};
    \path[->,>=stealth']
        (00) edge (01)
        (00) edge (11)
        (10) edge (11)
        (01) edge (02)
        (01) edge (12)
        (11) edge (12)
        (02) edge[<-] (03)
        (02) edge[<-] (13)
        (12) edge[<-] (13)
        ;
    \node[var-b,label=above:{\tiny$(1,c+3n-4)$}] (10) at (5,1) {};
    \node[var-b,label=below:{\tiny$(0,c+3n-4)$}] (00) at (5,0) {};
    \node[var-b,label={[label distance=4mm]above:\tiny$(2,c+3n-3)$}] (11) at (6,1) {};
    \node[var-b,label={[label distance=4mm]below:\tiny$(1,c+3n-3)$}] (01) at (6,0) {};
    \node[var-b,label=above:{\tiny$(3,c+3n-2)$}] (12) at (7,1) {};
    \node[var-b,label=below:{\tiny$(2,c+3n-2)$}] (02) at (7,0) {};
    \node[var-b,label={[label distance=4mm]above:\tiny$(1,c+3n)$}] (13) at (8,1) {};
    \node[var-b,label={[label distance=4mm]below:\tiny$(0,c+3n)$}] (03) at (8,0) {};
    \path[->,>=stealth']
        (00) edge (01)
        (00) edge (11)
        (10) edge (11)
        (01) edge (02)
        (01) edge (12)
        (11) edge (12)
        (02) edge[<-] (03)
        (02) edge[<-] (13)
        (12) edge[<-] (13)
        ;
\end{tikzpicture}
\end{center}
and note that it witnesses that $\OrdDG$ is not $3n$-braided. Hence, by Theorem~\ref{thm:TnIsHMblocker}, Corollary \ref{cor:OrdisLeastLowerBoundForTn}, and Lemma~\ref{lem:OrdDGppeqOrd} we have that $\OrdDG\mathbin{\times}\C\ppleq\OrdDG$.
\end{proof}

The last goal of this section is to describe the smallest subposet of $\DGPoset$ that is closed under meets and joins and contains all disjoint unions of cycles, all $\T_n$'s, and $\OrdDG$. For this we still need the following well-known lemma.

\begin{lemma}\label{lem:AModelsSigmaBModelsSigmaImpliesAxBModlesSigma}
Let $\A,\B$ be finite structures and $\Sigma$ a minor condition. Then $\A\models\Sigma$ and $\B\models\Sigma$ implies $\A\times\B\models\Sigma$.
\end{lemma}
\begin{proof}
Let $h\in\Pol(\A)$ and $f\in\Pol(\B)$ of arity $n$. Define
\begin{align*}
    f\mathbin{\times}g\colon (\A\mathbin{\times}\B)^n&\to\A\mathbin{\times}\B\\ ((a_1,b_1),\dots,(a_n,b_n))&\mapsto(f(a_1,\dots,a_n),h(b_1,\dots,b_n))
\end{align*}
Clearly, $f\mathbin{\times}g\in\Pol(\A\times\B)$. Let $f_1,\dots,f_k\in\Pol(\A)$  and $g_1,\dots,g_k\in\Pol(\B)$ be witnesses for  $\A\models\Sigma$ and $\B\models\Sigma$, respectively. Then $f_1\mathbin{\times}g_1,\dots,f_k\mathbin{\times}g_k$ witness that $\A\times\B\models\Sigma$.
\end{proof}

Combining the results from this section with the results from Chapter~\ref{cha:cycles} we obtain the following corollary.
\begin{corollary}
Let $\T\in\{\T_3,\T_4,\T_5,\dots,\T_4^{\setminus(03)}\}$ and $\C$ be a disjoint union of cycles. Then
\begin{enumerate}
    \item $\T\times\C\models\Sigma_D$ if and only if $\C\models \Sigma_D$ for all finite $D\subset\N^+$ and
    \item $\T\times\C\models\HM k$ if and only if $\T\models\HM k$. 
\end{enumerate}
Let $\S,\T\in\{\T_3,\dots,\T_4^{\setminus(03)}\}$ and $\C$ and $\D$ be disjoint unions of cycles. Then
\begin{enumerate}\setcounter{enumi}{2}
    \item $\S\mathbin{\times}\C\ppleq\T\mathbin{\times}\D$ if and only if $\S\ppleq\T$ and $\C\ppleq\D$,
    \item $(\S\mathbin{\times}\C)\wedge(\T\mathbin{\times}\D)$ can be represented by $(\S\wedge\T)\mathbin{\times}(\C\wedge\D)$, and
    \item $(\S\mathbin{\times}\C)\vee(\T\mathbin{\times}\D)$ can be represented by $(\S\vee\T)\mathbin{\times}(\C\vee\D)$.
\end{enumerate}
\end{corollary}
In particular, the smallest subposet of $\DGPoset$ that is closed under meets and joins and contains all disjoint unions of cycles, all $\T_n$'s, and $\Ord$ looks like  the element $\P_1=\T_2$ together with infinitely many copies of $\SDPoset$ positioned in an infinite descending chain with a lower bound:
\[\P_1\cup\SDPoset\cup\T_3\mathbin{\times}\SDPoset\cup\dots\cup\OrdDG\mathbin{\times}\SDPoset.\]
Observe that we have seen in Example~\ref{exa:T4isNotACoverOfT3} that in $\DGPoset$ there are elements in between the elements of this subposet.
Interestingly this gives us infinitely many new greatest lower bounds of infinite descending chains.
\begin{corollary}
Let $\C$ be a disjoint union of cycles. Then the infinite descending chain $\T_3\mathbin{\times}\C\ppstrictlyGreater\T_4\mathbin{\times}\C\ppstrictlyGreater\T_5\mathbin{\times}\C\ppstrictlyGreater\cdots$  has the greatest lower bound $\OrdDG\mathbin{\times}\C$. 
\end{corollary}
\chapter{Orientations of Paths and Cycles}\label{cha:paths}
In Chapter~\ref{cha:cycles} we classified in particular  directed paths (in Lemma~\ref{lem:IdempConstructPaths}) and directed cycles (in Corollary~\ref{cor:PPvsLoopCaBlockerForSa}). In this chapter we try to go the next step and study orientations of paths and orientations of cycles. We will focus on orientations of paths but some of our results also extend to orientations of cycles. One motivation for this chapter was to find the lower covers of $\T_3$. We found a new candidate for a lower cover but could not proof that it is indeed a cover, see Problem~\ref{pro:lowerCoversOfTthreePathCandidate}.

\section{Paths of Height at most Three}
A digraph $G=(V,E)$ is \emph{an orientation of a path} if the digraph \[(V,E\cupdot\{(v,u)\mid (u,v)\in E\})\] is well defined and is an undirected path.
We have already seen $\P_n$ as the simplest examples of orientations of paths.
First we observe that orientations of paths are relatively high up in $\DGPoset$, i.e., they satisfy strong minor conditions.
\begin{lemma}
Let $\P$ be an orientation of a path. Then $\P\models\TS n$ for all $n\geq1$ and $\P\models\Majority$. In particular, $\P\models\NU n$ for all $n\geq3$.
\end{lemma}
\begin{proof}
Let $\P=(\{1,\dots,n\},E)$ such that for every $(u,v)\in E$ we have that $u+1=v$ or $u=v+1$. 
For every $n\in\N^+$ the map \[(x_1,\dots,x_n)\mapsto\min(x_1,\dots,x_n)\] is a totally symmetric operation on $\P$. We want to show that it is also a polymorphism. Let $(u_1,v_1),\dots,(u_n,v_n)\in E$,  $\min(u_1,\dots,u_n)=u_i$, and $\min(v_1,\dots,v_n)=v_j$. Note that since $|u_i-v_i|=1$ for all $i$ we have that $v_j\in\{u_i-1,u_i,u_i+1\}$.
If $v_i=u_i+1$, then $v_j\not=u_i$, a contradiction. If $v_i=u_i-1$, then $v_j=v_i$. Hence, $(u_i,v_j)=(u_i,v_i)\in E$, as desired. 
 Now we have two more cases to consider:
\begin{itemize}
    \item If $(u_i-1)\toEdge u_i \toEdge v_i$, then $v_j\not=u_i-1$ and $v_j=u_i+1=v_i$. Hence, $(u_i,v_i)\in E$.
    \item If $(u_i-1)\fromEdge u_i \toEdge v_i$, then  $(u_i,v_j)\in E$ since $v_j\in\{u_i-1,v_i\}$.
\end{itemize}
Therefore $\min$ is a polymorphism and $\P\models\TS n$ for all $n\geq1$.

The map $\median\colon\P^3\to\P$ is a majority operation on $\P$. Next we verify that it is also a polymorphism. 
Let $(u_1,v_1),(u_2,v_2),(u_3,v_3)\in E$, let $\median(u_1,u_2,u_3)=u_i$, and  let $\median(v_1,v_2,v_3)=v_j$. Since $\median$ is fully symmetric we can assume without loss of generality that $u_1\leq u_2\leq u_3$. Hence, $u_i=u_2$. There are four cases to consider:
\begin{itemize}
    \item If $u_1=u_2=u_3$, then $(u_2,v_j)=(u_j,v_j)\in E$.
    \item If $u_1=u_2<u_3$, then $v_1,v_2\leq v_3$. Hence, $v_j\in\{v_1,v_2\}$ and $(u_2,v_j)\in E$. 
    \item If $u_1<u_2=u_3$, then $v_1\leq v_2, v_3$. Hence, $v_j\in\{v_2,v_3\}$ and $(u_2,v_j)\in E$. 
    \item If $u_1<u_2<u_3$, then $v_1\leq v_2\leq v_3$. Hence, $v_j=v_2$ and $(u_2,v_j)\in E$. 
\end{itemize}
Therefore $\median$ is a polymorphism of $\P$ and $\P\models\Majority$.
Observe that from a majority operation we can construct near-unanimity operations of any arity by composition.
\end{proof}
\begin{figure}
    \centering
    \begin{tikzpicture}[scale=0.5]
    \node[var-b] (0) at (0,0) {};
    \node[var-b] (1) at (0,1) {};
    \node        (2) at (0,2) {$\vdots$};
    \node[var-b] (3) at (0,3) {};
    \node[var-b,label=left:$s$] (4) at (0,4) {};

    \node[var-b,label=right:$t$] (0') at (3,1) {};
    \node[var-b] (1') at (3,2) {};
    \node        (2') at (3,3) {$\vdots$};
    \node[var-b] (3') at (3,4) {};
    \node[var-b] (4') at (3,5) {};

    \draw plot [smooth] coordinates {(4) (0.8,2) (2.1,3)  (0')};
    \node at (1.6,2) {$\P$};
    
    \path[->,>=stealth']
        (0) edge (1)
        (3) edge (4)
        (0') edge (1')
        (3') edge (4')
        ;
    
    \end{tikzpicture}
    \caption{The core path $\P'$ obtained from $\P$.}
    \label{fig:corePathsEmbedAllPaths}
\end{figure}
Note that every orientation of a path that cannot be pp-constructed from $\P_1$ does not satisfy $\HM1$ and is therefore, by Theorem~\ref{thm:TnIsHMblocker}, below $\T_3$. 
A digraph $\G$ is \emph{balanced} if there is an $n\in\N$ such that $\G$ has a homomorphism to $\P_n$. 
If $\G$ is a balanced digraph and $n\in\N$ is minimal such that $\G$ has some homomorphism to $\P_n$. Then this homomorphism is unique and we denote it by $\height_{\G}$, or by $\height$ if $\G$ is clear from the context, and call $n$ the \emph{height} of $\G$. 
Since $\height_{\G}$ is unique we can make the following observation.
\begin{observation}
Let $\P$ be an orientation of a path such that the start and endpoint of $\P$ are the only nodes on their respective levels. Then $\P$ is a core.
\end{observation}
Let $\P$ be an orientation of a path of height $n$ with start vertex $s$ and end vertex $t$. Let $\P'$ be obtained from $\P$ by appending a directed incoming path of length $n+1$ to $s$ and a directed outgoing path of length $n+1$ to $t$, see Figure~\ref{fig:corePathsEmbedAllPaths}. 
Then, by the above observation, $\P'$ is a core that embeds $\P$.
Hence, cores of orientations of paths can embed any orientation of a path.
To simplify things let us for now restrict to orientations of paths with height at most three. 
Define the graphs presented in Figure~\ref{fig:ZigAndCycZig}
\begin{align*}
    \Zig_n&\coloneqq (\{0,\dots,2n+3\},E_n\cup\{(0,1),(1,2),(2n+2,2n+3)\})
    \text{, where}\\
        E_n&\coloneqq \{(2i+1,2i),(2i+1,2i+2)\mid i\in\{1,\dots,n\}\}
\text{ and}\\
    \ZigCyc_n&\coloneqq(\{0,\dots,2n\},E'_n\cup\{(2n,0)\})\text{, where}\\
        E'_n&\coloneqq \{(2i,2i+1),(2i+2,2i+1)\mid i\in\{0,\dots,n-1\}\}.
\end{align*}
\begin{figure}
    \centering
    \begin{tikzpicture}[scale=0.5]
    \node[var-b,label=left:0] (0) at (0,0) {};
    \node[var-b,label=left:1] (1) at (0,1) {};
    \node[var-b,label=left:2] (2) at (0,2) {};
    \node[var-b] (11) at (1,1) {};
    \node[var-b] (21) at (1,2) {};
    \node (12) at (2,1) {$\dots$};
    \node (22) at (2,2) {$\dots$};
    
    \node[var-b] (13) at (3,1) {};
    \node[var-b] (23) at (3,2) {};
    \node[var-b,label=right:\small $2n+1$] (14) at (4,1) {};
    \node[var-b,label=right:\small $2n+2$] (24) at (4,2) {};
    \node[var-b,label=right:\small $2n+3$] (34) at (4,3) {};
    \path[->,>=stealth']
        (0) edge (1)
        (1) edge (2)
        (0) edge (1)
        (11) edge (2)
        (11) edge (21)
        (13) edge (23)
        (14) edge (23)
        (14) edge (24)
        (24) edge (34)
        ;
    \end{tikzpicture}
    \hspace{2cm}
    \begin{tikzpicture}[scale=0.5]
    \node[var-b,label=left:0] (1) at (0,1) {};
    \node[var-b,label=left:1] (2) at (0,2) {};
    \node[var-b] (11) at (1,1) {};
    \node[var-b] (21) at (1,2) {};
    \node (12) at (2,1) {$\dots$};
    \node (22) at (2,2) {$\dots$};
    
    \node[var-b] (13) at (3,1) {};
    \node[var-b] (23) at (3,2) {};
    \node[var-b,label=right:\small $2n$] (14) at (4,1) {};
    \path[->,>=stealth']
        (1) edge (2)
        (11) edge (2)
        (11) edge (21)
        (13) edge (23)
        (14) edge (23)
        (14) edge[bend left] (1)
        ;
    \end{tikzpicture}
    \caption{The graphs $\Zig_n$ (left) and $\ZigCyc_n$ (right).}
    \label{fig:ZigAndCycZig}
\end{figure}
Note that $\Zig_0=\P_3$ and $\ZigCyc_1\ppeq\T_3$.

\begin{observation}
Orientations of paths of height at most three are homomorphically equivalent to $\P_0$, $\P_1$, $\P_2$, or $\Zig_n$ for some $n\in\N$. 
\end{observation}

\begin{lemma}
Let $n\in\N^+$\!. Then $\Zig_n\ppeq\ZigCyc_n\ppeq \T_{2n+1}$.
\end{lemma}
\begin{proof}
We show $\T_{2n+1}\ppleq \ZigCyc_{n}\ppleq \Zig_n\ppleq\T_{2n+1}$.
Consider the following formulas
\begin{center}
\begin{tikzpicture}[scale=0.5]
\node[var-f,label=below:{$x_1$}] (x1) at (0,0) {};
\node[var-f,label=below:{$x_2$}] (x2) at (1,0) {};
\node[var-f,label=below:{$x_3$}] (x3) at (2,0) {};
\node (x4) at (3,0) {$\dots$};
\node[var-f,label=below:{$x_n$}] (x5) at (4,0) {};

\node[var-f,label=above:{$y_1$}] (y1) at (0,1) {};
\node[var-f,label=above:{$y_2$}] (y2) at (1,1) {};
\node[var-f,label=above:{$y_3$}] (y3) at (2,1) {};
\node (y4) at (3,1) {$\dots$};
\node[var-f,label=above:{$y_n$}] (y5) at (4,1) {};
\path[>=stealth'] 
        (x5) edge[->] (y5)
        (x3) edge[->] (y3)
        (y2) edge[->] (x3)
        (x2) edge[->] (y2)
        (y1) edge[->] (x2)
        (x1) edge[->] (y1)
        ;
        \node at (2,-2) {$\Phi_1$};
\end{tikzpicture}
\hspace{1cm}
\begin{tikzpicture}[scale=0.5]
\node[var-f,label=below:{$x_1$}] (x1) at (0,0) {};
\node[var-f,label=below:{$x_2$}] (x2) at (1,0) {};
\node[var-f,label=below:{$x_3$}] (x3) at (2,0) {};
\node[var-f,label=below:{$0$}] (x4) at (3,0) {};

\node[var-f,label=above:{$y_1$}] (y1) at (0,1) {};
\node[var-f,label=above:{$1$}] (y2) at (1,1) {};
\node[var-f,label=above:{$y_3$}] (y3) at (2,1) {};
\node[var-f,label=above:{$y_4$}] (y4) at (3,1) {};
\path[>=stealth'] 
        (x1) edge[->] (y1)
        (x2) edge[dashed] (y3)
        (x3) edge[dashed] (y4)
        ;
        \node at (1.5,-2) {$\Phi_2$};
\end{tikzpicture}



\end{center}
and let $\ZigCyc_n'$ be the $n$-th pp-power of $\T_{2n+1}$ given by $\Phi_1$. 
Note that if $(u,v)$ is an edge in $\ZigCyc_n'$, then $u_0<v_0<u_1<\dots<u_n<v_n$ and there is a unique element $i$ of $\T_{2n+1}$ that does neither occur in $u$ nor in $v$. Denote the edge $(u,v)$ by $e_i$, $u$ by $s_i$, and $v$ by $t_i$. Clearly, there also exists an edge $e_i$ for every $i\in\T_{2n+1}$, see Figure~\ref{fig:CycZigPPConstruction} for an example. Hence, $\ZigCyc_n'$ has exactly $2n+1$ edges. Observe that
\begin{itemize}
    \item $s_{2i}=(0,2,\dots,2i-2,2i+1,2i+3,\dots,2n-1)=s_{2i-1}$ for $1\leq i\leq n$,
    \item $t_{2i}=(1,3,\dots,2i-1,2i+2,2i+4,\dots,2n)=t_{2i+1}$  for $0\leq i\leq n-1$, and
    \item $t_{2n}=(1,3,\dots,2n-1)=s_0$.
\end{itemize}
The map $h\colon \ZigCyc_n\to\ZigCyc_n'$,
\begin{align*}
    i\mapsto
    \begin{cases}
    s_i&\text{if $i$ is even}\\
    t_i&\text{if $i$ is odd}
    \end{cases}
\end{align*}
is clearly injective. Furthermore,
\begin{itemize}
    \item $(h(2i),h(2i+1))=(s_{2i},t_{2i+1})=(s_{2i},t_{2i})=e_{2i}$, 
    \item $(h(2i+2),h(2i+1))=(s_{2i+2},t_{2i+1})=(s_{2i+1},t_{2i+1})=e_{2i+1}$, and
    \item $(h(2n),h(0))=(s_{2n},s_{0})=(s_{2n},t_{2n})=e_{2n}$.
\end{itemize}
Hence, $h$ is an embedding. Since $\ZigCyc_n'$ and $\ZigCyc_n$ also have the same finite number of edges they are homomorphically equivalent. Hence $\T_{2n+1}\ppleq\ZigCyc_n$.

\begin{figure}
    \centering
    \begin{tikzpicture}[scale=1.0]
    \node[var-b,label=below:{\tiny 135}] (1) at (0,1) {};
    \node[var-b,label=above:{\tiny 246}] (2) at (0,2) {};
    \node[var-b,label=below:{\tiny 035}] (11) at (1,1) {};
    \node[var-b,label=above:{\tiny 146}] (21) at (1,2) {};
    \node[var-b,label=below:{\tiny 025}] (13) at (2,1) {};
    \node[var-b,label=above:{\tiny 136}] (23) at (2,2) {};
    \node[var-b,label=below:{\tiny 024}] (14) at (3,1) {};
    \path[->,>=stealth']
        (1)  edge node[left] {\small $e_0$} (2)
        (11) edge (2)
        (11) edge (21)
        (13) edge (21)
        (13) edge (23)
        (14) edge node[right] {\small $e_5$} (23)
        (14) edge[bend left=35] (1)
        ;
    \end{tikzpicture}
    
    \caption{Nonisolated vertices of $\ZigCyc_3'$.}
    \label{fig:CycZigPPConstruction}
\end{figure}

Let $\Zig_n'$ be the 4-th pp-power of $\ZigCyc_n$ given by $\Phi_2$. 
It is straightforward to verify that the maps
\begin{align*}
    h\colon \Zig_n&\to\Zig_n'\\
    0&\mapsto (2n,0,0,0)\\
    i&\mapsto (i-1,1,0,0)\text{ if $i$ is odd, $1\leq i\leq 2n+1$}\\
    i&\mapsto (i-1,1,1,0)\text{ if $i$ is even, $2\leq i\leq 2n+1$}\\
    2n+2&\mapsto (0,1,1,0)\\
    2n+3&\mapsto (1,1,1,1)
\intertext{and}
    h'\colon \Zig_n'&\to\Zig_n\\
    (2n,1,1,1)&\mapsto 0\\
    (0,1,1,0)&\mapsto 2n+2\\
    (1,1,1,1)&\mapsto 2n+3\\
    (i,j,k,\ell)&\mapsto i-1
\end{align*}
are homomorphisms. Hence, $\ZigCyc_n\ppleq\Zig_n$.
Observe that the following $(2n-1)$-braid
\begin{center}
    \begin{tikzpicture}[scale = 1.2]
    \node[var-b,label=above:{$0$}] (10) at (0,1) {};
    \node[var-b,label=below:{$2$}] (00) at (0,0) {};
    \node[var-b,label=above:{$1$}] (11) at (1,1) {};
    \node[var-b,label=below:{$3$}] (01) at (1,0) {};
    \node[var-b,label=above:{$2$}] (12) at (2,1) {};
    \node[var-b,label=below:{$4$}] (02) at (2,0) {};
    \node[var-b,label=above:{$3$}] (13) at (3,1) {};
    \node[var-b,label=below:{$5$}] (03) at (3,0) {};
    \node[] (14) at (4,1) {$\dots$};
    \node[] (04) at (4,0) {$\dots$};
    \path[->,>=stealth']
        (00) edge (01)
        (00) edge (11)
        (10) edge (11)
        (01) edge[<-] (02)
        (01) edge[<-] (12)
        (11) edge[<-] (12)
        (02) edge (03)
        (02) edge (13)
        (12) edge (13)
        ;
    \node[var-b,label=above:{\small$2n-3$}] (10) at (5,1) {};
    \node[var-b,label=below:{\small$2n-1$}] (00) at (5,0) {};
    \node[var-b,label={above:\small$2n-2$}] (11) at (6,1) {};
    \node[var-b,label={below:$2n$}] (01) at (6,0) {};
    \node[var-b,label=above:{\small$2n-1$}] (12) at (7,1) {};
    \node[var-b,label=below:{$0$}] (02) at (7,0) {};
    \path[->,>=stealth']
        (00) edge[<-] (01)
        (00) edge[<-] (11)
        (10) edge[<-] (11)
        (01) edge (02)
        (01) edge (12)
        (11) edge (12)
        ;
\end{tikzpicture}
\end{center}
witnesses that $\ZigCyc_n$ is not $(2n-1)$-braided. Hence, by Theorem~\ref{thm:TnIsHMblocker}, $\ZigCyc_n\ppleq\T_{2n+1}$ as desired.
\end{proof}

\begin{observation}
Let $\H$ be a digraph and $n,k\geq 1$. Note that 
\begin{align*}
        \H&\to\T_n&&\text{if and only if}&\P_n&\not\to\H.
\intertext{Guzm\'an-Pro and C\'esar Hern\'andez-Cruz showed in~\cite{SantiagoDualityPairs} in Theorem~12  that
}
\H&\to\ZigCyc_{k+1}&&\text{if and only if}&\Zig_k&\not\to\H.
\end{align*}
Hence, $\ZigCyc_{k+1}$ and $\T_n$ have a single \emph{dual} digraph. The dual graphs $\Zig_k$ and $\P_n$ are for no $n$ and $k$ homomorphically equivalent. However, $\ZigCyc_n\ppeq\T_{2n+1}$. 
To our knowledge this are the only known example of  digraphs with not homomorphically equivalent duals that have the same pp-constructability type.

\end{observation}

\begin{corollary}
The orientations of paths of height at most three form the following descending chain.
\[\P_0\ppstrictlyGreater\P_1\ppeq\P_2\ppeq\P_3\ppstrictlyGreater\Zig_1\ppstrictlyGreater\Zig_2\ppstrictlyGreater\cdots\]
\end{corollary}

\begin{figure}
    \centering
    \begin{tikzpicture}[scale=0.3]
    \node[var-b] (0) at (0,0) {};
    \node[var-b] (1) at (1,1) {};
    \node[var-b] (2) at (2,2) {};
    \node[var-b] (3) at (3,3) {};
    \node[var-b] (4) at (4,4) {};
    \node[var-b] (5) at (5,3) {};
    \node[var-b] (6) at (6,2) {};
    \node[var-b] (7) at (7,3) {};
    \node[var-b] (8) at (8,2) {};
    \node[var-b] (9) at (9,1) {};
    \node[var-b] (10) at (10,2) {};
    \node[var-b] (11) at (11,3) {};
    \node[var-b] (12) at (12,4) {};
    \node[var-b] (13) at (13,5) {};
    \path[->,>=stealth']
    (0) edge (1)
    (1) edge (2)
    (2) edge (3)
    (3) edge (4)
    (5) edge (4)
    (6) edge (5)
    (6) edge (7)
    (8) edge (7)
    (9) edge (8)
    (9) edge (10)
    (10) edge (11)
    (11) edge (12)
    (12) edge (13)
    ;
    \end{tikzpicture}\hspace{10mm}
    \begin{tikzpicture}[scale=0.3]
    \node[var-f,label=below:$x$] (1) at (1,1) {};
    \node[var-b] (2) at (2,2) {};
    \node[var-b] (3) at (3,3) {};
    \node[var-b] (4) at (4,4) {};
    \node[var-b] (5) at (5,3) {};
    \node[var-b] (6) at (6,2) {};
    \node[var-b] (7) at (7,3) {};
    \node[var-b] (8) at (8,2) {};
    \node[var-f,label=below:$y$] (9) at (9,1) {};
    \path[->,>=stealth']
    (1) edge (2)
    (2) edge (3)
    (3) edge (4)
    (5) edge (4)
    (6) edge (5)
    (6) edge (7)
    (8) edge (7)
    (9) edge (8)
    ;
    \end{tikzpicture}
    \caption{The smallest orientation of a path that can pp-construct $\Ord$ (left) and the pp-formula used in the construction (right).}
    \label{fig:smallestNLhardPath}
\end{figure}
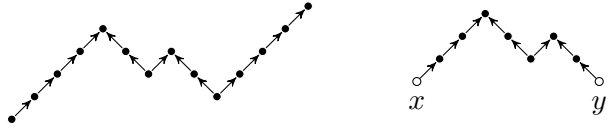

\section{NL-Hard Paths}
We have seen that for any $n$ there are orientations of paths that can pp-construct $\T_n$. Are there also orientations of paths that can pp-construct $\Ord$? This question was already answered by Egri in \cite{EgriNLhardPaths} in Theorem 28. In particular, he found the orientation of a path presented in Figure~\ref{fig:smallestNLhardPath} (left) that can pp-construct $\Ord$ using the formula on the right of Figure~\ref{fig:smallestNLhardPath}. 
Let $(n_1,\dots,n_k)$ be a tuple of natural numbers and $n=n_1+\dots+n_k$. Then $\OriP{n_1,\dots,n_k}$ is the orientation of $\P_n$ with the edges
\[0\toEdge\dots\toEdge n_1\fromEdge\dots\fromEdge (n_1+n_2)\toEdge\dots\toEdge (n_1+n_2+n_3)\fromEdge\dots\ n\]
and $\OriC{n_1,\dots,n_k}$ is the orientation of $\C_n$ with the edges 
\[0\toEdge\dots\toEdge n_1\fromEdge\dots\fromEdge (n_1+n_2)\toEdge\dots\toEdge (n_1+n_2+n_3)\fromEdge\dots\ 0.\]
For example 
\begin{align*}
    \C_n&=\OriC{n},\! &\Zig_1&=\OriP{2,1,2}, &\ZigCyc_2&=\OriC{1,1,1,1,1}=\OriC{2,1,1,1},\\
    \P_n&=\OriP{n}, &\Zig_2&=\OriP{2,1,1,1,2},\! &\ZigCyc_1&=\OriC{1,1,1}=\OriC{2,1},
\end{align*} the path from Figure~\ref{fig:smallestNLhardPath} is isomorphic to  $\OriP{4,2,1,2,4}$.
A  path $\OriP{n_1,\dots,n_k}$ is called \emph{symmetric} if $\OriP{n_1,\dots,n_k}$ is isomorphic to $\OriP{n_k,\dots,n_1}$. 

\begin{sidewaystable}
    \centering
    \begin{adjustbox}{width=\columnwidth,center}
    \begin{tabular}{l|rrrrrrrrrrrrrrrrrrrrrrrrrrr}
$n$& 1& 2& 3& 4& 5& 6& 7& 8& 9& 10& 11& 12& 13& 14& 15& 16& 17& 18& 19& 20& 21& 22& 23& 24& 25& 26& 27\\\hline
cores& 1& 1& 1& 1& 1& 2& 3& 5& 9& 17& 30& 61& 107& 226& 402& 852& 1544& 3256& 6001& 12556& 23512& 48771& 92547& 190538& 365382& 747654& 1445748\\\hline
$\operatorname{pathsHM}(n,1)$& 1& 1& 1& 1& 1& 1& 1& 1& 1& 1& 1& 1& 1& 1& 1& 1& 1& 1& 1& 1& 1& 1& 1& 1& 1& 1& 1\\
$\operatorname{pathsHM}(n,2)$& 0& 0& 0& 0& 0& 1& 2& 3& 6& 12& 19& 37& 59& 112& 187& 349& 588& 1087& 1855& 3400& 5873& 10685& 18605& 33684& 58977& 106373& 186802\\
$\operatorname{pathsHM}(n,3)$& 0& 0& 0& 0& 0& 0& 0& 0& 0& 0& 2& 6& 13& 34& 59& 145& 268& 613& 1147& 2486& 4680& 9773& 18497& 37625& 71485& 142744& 271629\\
$\operatorname{pathsHM}(n,4)$& 0& 0& 0& 0& 0& 0& 0& 1& 2& 3& 6& 13& 24& 55& 103& 226& 423& 901& 1707& 3568& 6825& 14060& 27069& 55081& 106487& 214600& 415960\\
$\operatorname{pathsHM}(n,5)$& 0& 0& 0& 0& 0& 0& 0& 0& 0& 0& 0& 0& 2& 6& 16& 40& 84& 206& 404& 952& 1846& 4210& 8237& 18178& 35875& 77148& 153017\\
$\operatorname{pathsHM}(n,6)$& 0& 0& 0& 0& 0& 0& 0& 0& 0& 1& 2& 3& 6& 13& 24& 56& 108& 244& 479& 1047& 2054& 4411& 8651& 18340& 36074& 75601& 149361\\
$\operatorname{pathsHM}(n,7)$& 0& 0& 0& 0& 0& 0& 0& 0& 0& 0& 0& 0& 0& 0& 2& 6& 14& 36& 78& 192& 407& 948& 1929& 4376& 8756& 19430& 38726\\
$\operatorname{pathsHM}(n,8)$& 0& 0& 0& 0& 0& 0& 0& 0& 0& 0& 0& 1& 2& 3& 6& 13& 24& 56& 104& 237& 456& 1008& 1987& 4285& 8512& 18065& 35933\\
$\operatorname{pathsHM}(n,9)$& 0& 0& 0& 0& 0& 0& 0& 0& 0& 0& 0& 0& 0& 0& 0& 0& 2& 6& 14& 36& 74& 184& 380& 890& 1819& 4126& 8381\\
$\operatorname{pathsHM}(n,10)$& 0& 0& 0& 0& 0& 0& 0& 0& 0& 0& 0& 0& 0& 1& 2& 3& 6& 13& 24& 56& 104& 237& 450& 997& 1941& 4200& 8292\\
$\operatorname{pathsHM}(n,11)$& 0& 0& 0& 0& 0& 0& 0& 0& 0& 0& 0& 0& 0& 0& 0& 0& 0& 0& 2& 6& 14& 36& 74& 184& 374& 878& 1772\\
$\operatorname{pathsHM}(n,12)$& 0& 0& 0& 0& 0& 0& 0& 0& 0& 0& 0& 0& 0& 0& 0& 1& 2& 3& 6& 13& 24& 56& 104& 237& 450& 997& 1930\\
$\operatorname{pathsHM}(n,13)$& 0& 0& 0& 0& 0& 0& 0& 0& 0& 0& 0& 0& 0& 0& 0& 0& 0& 0& 0& 0& 2& 6& 14& 36& 74& 184& 374\\
$\operatorname{pathsHM}(n,14)$& 0& 0& 0& 0& 0& 0& 0& 0& 0& 0& 0& 0& 0& 0& 0& 0& 0& 1& 2& 3& 6& 13& 24& 56& 104& 237& 450\\
$\operatorname{pathsHM}(n,15)$& 0& 0& 0& 0& 0& 0& 0& 0& 0& 0& 0& 0& 0& 0& 0& 0& 0& 0& 0& 0& 0& 0& 2& 6& 14& 36& 74\\
$\operatorname{pathsHM}(n,16)$& 0& 0& 0& 0& 0& 0& 0& 0& 0& 0& 0& 0& 0& 0& 0& 0& 0& 0& 0& 1& 2& 3& 6& 13& 24& 56& 104\\
$\operatorname{pathsHM}(n,17)$& 0& 0& 0& 0& 0& 0& 0& 0& 0& 0& 0& 0& 0& 0& 0& 0& 0& 0& 0& 0& 0& 0& 0& 0& 2& 6& 14\\
$\operatorname{pathsHM}(n,18)$& 0& 0& 0& 0& 0& 0& 0& 0& 0& 0& 0& 0& 0& 0& 0& 0& 0& 0& 0& 0& 0& 1& 2& 3& 6& 13& 24\\
$\operatorname{pathsHM}(n,19)$& 0& 0& 0& 0& 0& 0& 0& 0& 0& 0& 0& 0& 0& 0& 0& 0& 0& 0& 0& 0& 0& 0& 0& 0& 0& 0& 2\\
$\operatorname{pathsHM}(n,20)$& 0& 0& 0& 0& 0& 0& 0& 0& 0& 0& 0& 0& 0& 0& 0& 0& 0& 0& 0& 0& 0& 0& 0& 1& 2& 3& 6\\
$\operatorname{pathsHM}(n,21)$& 0& 0& 0& 0& 0& 0& 0& 0& 0& 0& 0& 0& 0& 0& 0& 0& 0& 0& 0& 0& 0& 0& 0& 0& 0& 0& 0\\
$\operatorname{pathsHM}(n,22)$& 0& 0& 0& 0& 0& 0& 0& 0& 0& 0& 0& 0& 0& 0& 0& 0& 0& 0& 0& 0& 0& 0& 0& 0& 0& 1& 2\\
$\operatorname{pathsHM}(n,23)$& 0& 0& 0& 0& 0& 0& 0& 0& 0& 0& 0& 0& 0& 0& 0& 0& 0& 0& 0& 0& 0& 0& 0& 0& 0& 0& 0\\
$\operatorname{pathsHM}(n,24)$& 0& 0& 0& 0& 0& 0& 0& 0& 0& 0& 0& 0& 0& 0& 0& 0& 0& 0& 0& 0& 0& 0& 0& 0& 0& 0& 0\\
$\operatorname{pathsHM}(n,\infty)$& 0& 0& 0& 0& 0& 0& 0& 0& 0& 0& 0& 0& 0& 1& 2& 12& 24& 89& 178& 558& $\leq1144$& $\leq3139$& $\leq6515$& $\leq16545$& $\leq34405$& $\leq82955$& $\leq172894$
    \end{tabular}
    \end{adjustbox}
    \caption{The number of core orientations of paths with $n$ vertices partitioned by their satisfaction of $\HM k$.}
    \label{tab:HMpathsNumbers}
\end{sidewaystable}
Let $\operatorname{pathsHM}(n,k)$ denote the number of core orientations of paths with $n$ vertices that satisfy $\HM k$ but not $\HM {k-1}$ and let $\operatorname{pathsHM}(n,\infty)$ be the number of core orientations of paths with $n$ vertices that do not satisfy $\HM k$ for any $k$.
Using the program we also used in Chapter~\ref{cha:hardTrees} I computed the numbers $\operatorname{pathsHM}(n,k)$ for $n\leq 21$ and $k\in\N\cup\{\infty\}$ and for $n\leq 27$ and $k\leq 30$.
These numbers are presented in Table~\ref{tab:HMpathsNumbers} and there is a visual representation in Table~\ref{tab:HMpathsTable}. 
One challenge in computing these numbers was to verify that an orientation of a path with  $n$ vertices does not satisfy $\HM k$ for any $k$.
Recall from Lemma~\ref{lem:noHMimpliesACanppDefineOrder} that a structure $\A$ with $n$ vertices satisfies: 
\begin{align*}
    &\A\ppleq\Ord &&\text{iff}  &&\A\not\models\HM {n\cdot(n-1)^2} &&\text{iff} &&\A\not\models\HM k\text{ for any }k.
\end{align*} 
The first condition was easier to verify than the second one. Still the program could only compute $\operatorname{pathsHM}(n,\infty)$ for $n\leq 21$. 
The core orientations of paths with at most 16 vertices that can pp-construct $\Ord$ are presented in Table~\ref{tab:NLHardPaths} together with the formula used for the pp-construction. As in Figure~\ref{fig:smallestNLhardPath}, the free variables of the formula are the endpoints of the path. The core orientation of a path with the pp-formula with the most variables I encountered when computing $\operatorname{pathsHM}(n,\infty)$ for $n\leq 21$ is the path 
\begin{align*}
&\OriP{2,1,3,3,3,1,1,2,1,3}\text{ with the pp-formula}\\
&\OriP{2,1,2,1,2,1,1,2,1,2,2,1,2,2,1,2,1,2}.
\end{align*}

\begin{table}
    \centering
    \pgfplotstableread{ 
Label 1 2 3 4 5 6 7 8 9 10 11 12 13 14 15 16 17 18 19 20 21 22 23 24 25 26 27 28 29 30 
1 1.00000 0.00000 0.00000 0.00000 0.00000 0.00000 0.00000 0.00000 0.00000 0.00000 0.00000 0.00000 0.00000 0.00000 0.00000 0.00000 0.00000 0.00000 0.00000 0.00000 0.00000 0.00000 0.00000 0.00000 0.00000 0.00000 0.00000 0.00000 0.00000 0.00000 
2 1.00000 0.00000 0.00000 0.00000 0.00000 0.00000 0.00000 0.00000 0.00000 0.00000 0.00000 0.00000 0.00000 0.00000 0.00000 0.00000 0.00000 0.00000 0.00000 0.00000 0.00000 0.00000 0.00000 0.00000 0.00000 0.00000 0.00000 0.00000 0.00000 0.00000 
3 1.00000 0.00000 0.00000 0.00000 0.00000 0.00000 0.00000 0.00000 0.00000 0.00000 0.00000 0.00000 0.00000 0.00000 0.00000 0.00000 0.00000 0.00000 0.00000 0.00000 0.00000 0.00000 0.00000 0.00000 0.00000 0.00000 0.00000 0.00000 0.00000 0.00000 
4 1.00000 0.00000 0.00000 0.00000 0.00000 0.00000 0.00000 0.00000 0.00000 0.00000 0.00000 0.00000 0.00000 0.00000 0.00000 0.00000 0.00000 0.00000 0.00000 0.00000 0.00000 0.00000 0.00000 0.00000 0.00000 0.00000 0.00000 0.00000 0.00000 0.00000 
5 1.00000 0.00000 0.00000 0.00000 0.00000 0.00000 0.00000 0.00000 0.00000 0.00000 0.00000 0.00000 0.00000 0.00000 0.00000 0.00000 0.00000 0.00000 0.00000 0.00000 0.00000 0.00000 0.00000 0.00000 0.00000 0.00000 0.00000 0.00000 0.00000 0.00000 
6 0.50000 0.50000 0.00000 0.00000 0.00000 0.00000 0.00000 0.00000 0.00000 0.00000 0.00000 0.00000 0.00000 0.00000 0.00000 0.00000 0.00000 0.00000 0.00000 0.00000 0.00000 0.00000 0.00000 0.00000 0.00000 0.00000 0.00000 0.00000 0.00000 0.00000 
7 0.33333 0.66666 0.00000 0.00000 0.00000 0.00000 0.00000 0.00000 0.00000 0.00000 0.00000 0.00000 0.00000 0.00000 0.00000 0.00000 0.00000 0.00000 0.00000 0.00000 0.00000 0.00000 0.00000 0.00000 0.00000 0.00000 0.00000 0.00000 0.00000 0.00000 
8 0.20000 0.60000 0.00000 0.20000 0.00000 0.00000 0.00000 0.00000 0.00000 0.00000 0.00000 0.00000 0.00000 0.00000 0.00000 0.00000 0.00000 0.00000 0.00000 0.00000 0.00000 0.00000 0.00000 0.00000 0.00000 0.00000 0.00000 0.00000 0.00000 0.00000 
9 0.11111 0.66666 0.00000 0.22222 0.00000 0.00000 0.00000 0.00000 0.00000 0.00000 0.00000 0.00000 0.00000 0.00000 0.00000 0.00000 0.00000 0.00000 0.00000 0.00000 0.00000 0.00000 0.00000 0.00000 0.00000 0.00000 0.00000 0.00000 0.00000 0.00000 
10 0.05882 0.70588 0.00000 0.17647 0.00000 0.05882 0.00000 0.00000 0.00000 0.00000 0.00000 0.00000 0.00000 0.00000 0.00000 0.00000 0.00000 0.00000 0.00000 0.00000 0.00000 0.00000 0.00000 0.00000 0.00000 0.00000 0.00000 0.00000 0.00000 0.00000 
11 0.03333 0.63333 0.06666 0.20000 0.00000 0.06666 0.00000 0.00000 0.00000 0.00000 0.00000 0.00000 0.00000 0.00000 0.00000 0.00000 0.00000 0.00000 0.00000 0.00000 0.00000 0.00000 0.00000 0.00000 0.00000 0.00000 0.00000 0.00000 0.00000 0.00000 
12 0.01639 0.60655 0.09836 0.21311 0.00000 0.04918 0.00000 0.01639 0.00000 0.00000 0.00000 0.00000 0.00000 0.00000 0.00000 0.00000 0.00000 0.00000 0.00000 0.00000 0.00000 0.00000 0.00000 0.00000 0.00000 0.00000 0.00000 0.00000 0.00000 0.00000 
13 0.00934 0.55140 0.12149 0.22429 0.01869 0.05607 0.00000 0.01869 0.00000 0.00000 0.00000 0.00000 0.00000 0.00000 0.00000 0.00000 0.00000 0.00000 0.00000 0.00000 0.00000 0.00000 0.00000 0.00000 0.00000 0.00000 0.00000 0.00000 0.00000 0.00000 
14 0.00442 0.49557 0.15044 0.24336 0.02654 0.05752 0.00000 0.01327 0.00000 0.00442 0.00000 0.00000 0.00000 0.00000 0.00000 0.00000 0.00000 0.00000 0.00000 0.00000 0.00000 0.00000 0.00000 0.00000 0.00000 0.00000 0.00000 0.00000 0.00000 0.00442 
15 0.00248 0.46517 0.14676 0.25621 0.03980 0.05970 0.00497 0.01492 0.00000 0.00497 0.00000 0.00000 0.00000 0.00000 0.00000 0.00000 0.00000 0.00000 0.00000 0.00000 0.00000 0.00000 0.00000 0.00000 0.00000 0.00000 0.00000 0.00000 0.00000 0.00497 
16 0.00117 0.40962 0.17018 0.26525 0.04694 0.06572 0.00704 0.01525 0.00000 0.00352 0.00000 0.00117 0.00000 0.00000 0.00000 0.00000 0.00000 0.00000 0.00000 0.00000 0.00000 0.00000 0.00000 0.00000 0.00000 0.00000 0.00000 0.00000 0.00000 0.01408 
17 0.00064 0.38082 0.17357 0.27396 0.05440 0.06994 0.00906 0.01554 0.00129 0.00388 0.00000 0.00129 0.00000 0.00000 0.00000 0.00000 0.00000 0.00000 0.00000 0.00000 0.00000 0.00000 0.00000 0.00000 0.00000 0.00000 0.00000 0.00000 0.00000 0.01554 
18 0.00030 0.33384 0.18826 0.27671 0.06326 0.07493 0.01105 0.01719 0.00184 0.00399 0.00000 0.00092 0.00000 0.00030 0.00000 0.00000 0.00000 0.00000 0.00000 0.00000 0.00000 0.00000 0.00000 0.00000 0.00000 0.00000 0.00000 0.00000 0.00000 0.02733 
19 0.00016 0.30911 0.19113 0.28445 0.06732 0.07982 0.01299 0.01733 0.00233 0.00399 0.00033 0.00099 0.00000 0.00033 0.00000 0.00000 0.00000 0.00000 0.00000 0.00000 0.00000 0.00000 0.00000 0.00000 0.00000 0.00000 0.00000 0.00000 0.00000 0.02966 
20 0.00007 0.27078 0.19799 0.28416 0.07582 0.08338 0.01529 0.01887 0.00286 0.00446 0.00047 0.00103 0.00000 0.00023 0.00000 0.00007 0.00000 0.00000 0.00000 0.00000 0.00000 0.00000 0.00000 0.00000 0.00000 0.00000 0.00000 0.00000 0.00000 0.04444 
21 0.00004 0.24978 0.19904 0.29027 0.07851 0.08735 0.01731 0.01939 0.00314 0.00442 0.00059 0.00102 0.00008 0.00025 0.00000 0.00008 0.00000 0.00000 0.00000 0.00000 0.00000 0.00000 0.00000 0.00000 0.00000 0.00000 0.00000 0.00000 0.00000 0.04865 
22 0.00002 0.21908 0.20038 0.28828 0.08632 0.09044 0.01943 0.02066 0.00377 0.00485 0.00073 0.00114 0.00012 0.00026 0.00000 0.00006 0.00000 0.00002 0.00000 0.00000 0.00000 0.00000 0.00000 0.00000 0.00000 0.00000 0.00000 0.00000 0.00000 0.06436 
23 0.00001 0.20103 0.19986 0.29248 0.08900 0.09347 0.02084 0.02147 0.00410 0.00486 0.00079 0.00112 0.00015 0.00025 0.00002 0.00006 0.00000 0.00002 0.00000 0.00000 0.00000 0.00000 0.00000 0.00000 0.00000 0.00000 0.00000 0.00000 0.00000 0.07039 
24 0.00000 0.17678 0.19746 0.28908 0.09540 0.09625 0.02296 0.02248 0.00467 0.00523 0.00096 0.00124 0.00018 0.00029 0.00003 0.00006 0.00000 0.00001 0.00000 0.00000 0.00000 0.00000 0.00000 0.00000 0.00000 0.00000 0.00000 0.00000 0.00000 0.08683 
25 0.00000 0.16141 0.19564 0.29144 0.09818 0.09872 0.02396 0.02329 0.00497 0.00531 0.00102 0.00123 0.00020 0.00028 0.00003 0.00006 0.00000 0.00001 0.00000 0.00000 0.00000 0.00000 0.00000 0.00000 0.00000 0.00000 0.00000 0.00000 0.00000 0.09416 
26 0.00000 0.14227 0.19092 0.28703 0.10318 0.10111 0.02598 0.02416 0.00551 0.00561 0.00117 0.00133 0.00024 0.00031 0.00004 0.00007 0.00000 0.00001 0.00000 0.00000 0.00000 0.00000 0.00000 0.00000 0.00000 0.00000 0.00000 0.00000 0.00000 0.11095 
27 0.00000 0.12924 0.18793 0.28779 0.10587 0.10334 0.02679 0.02486 0.00579 0.00573 0.00122 0.00133 0.00025 0.00031 0.00005 0.00007 0.00000 0.00001 0.00000 0.00000 0.00000 0.00000 0.00000 0.00000 0.00000 0.00000 0.00000 0.00000 0.00000 0.11933 
}\testdata
\begin{tikzpicture}
\begin{axis}[
            ybar stacked,   
            ymin=0,         
            ymax = 1,
            xtick=data,     
            enlarge x limits={0.05},
            xticklabels from table={\testdata}{Label},
            tick label style={font=\tiny},
            width=11.6cm, height = 7cm,line width=0.7pt,bar width=2.5mm
]
\def\plotcommand#1{
    \addplot [fill=black!#1!blue!90!white] table [y=\s, meta=Label,x expr=\coordindex] {\testdata};
}
\foreach \s in {1,...,14}
{
\pgfmathparse{ln(\s)*100/3}
\expandafter\plotcommand\expandafter{\pgfmathresult}
}
\addplot [fill=orange!80!white!80!red] table [y=30, meta=Label,x expr=\coordindex] {\testdata};
\end{axis}
\node[align=left] (2) at (11.3,0.3) {\small have $\HM2$\\\small but no $\HM1$};
\node[align=left] (4) at (11.3,2.4) {\small have $\HM4$\\\small but no $\HM3$};
\node[align=left] (6) at (11.3,3.5) {\small have $\HM5$\\\small but no $\HM4$};
\node[align=left] (NL) at (11,5.1) {likely \\ NL-hard};
\draw (9.55,0.3) -- (2);
\draw (9.55,2.4) -- (4);
\draw (9.55,3.5) -- (6);
\draw (9.55,5.1) -- (NL);
\end{tikzpicture}
    \caption{Distribution of core orientations of paths in L.}
    \label{tab:HMpathsTable}
\end{table}

There are  a lot of interesting observations about Table~\ref{tab:HMpathsNumbers}. 
Firstly, note that the path in Figure~\ref{fig:smallestNLhardPath} has 14 vertices. As seen in the table all orientations of paths with at most 13 vertices do have a Hagemann Mitschke chain. Hence, the path from Figure~\ref{fig:smallestNLhardPath} is the smallest orientation of a path that can pp-construct $\Ord$ 
Note that $\Zig_n$ has $2n+4$ vertices, $\Zig_n\models\HM{2n}$, and $\Zig_n\not\models\HM{2n-1}$.
Observe from Table~\ref{tab:HMpathsNumbers} that for even $n=2k\geq 4$ there is exactly one core orientation of a path with $n$ vertices that satisfies $\HM{n}$ but not $\HM{n-1}$. This path is $\Zig_k$. The program could verify that the following is true for orientations of paths with at most 21 vertices.

\begin{conjecture}\label{con:pathsWithnoHMnAreNLHard}
Let $\P(n_1,\dots,n_k)$ be an orientation of a path with $k\geq2$. If $\P\not\models \HM {k-1}$ then $\P\ppleq\Ord$. In particular: let $\P$ be an orientation of a path with $n$ vertices. If $\P\not\models\HM{n-1}$ then $\P\ppleq\Ord$.
\end{conjecture}
All orientations of paths with at most 21 vertices where this bound is tight are of the form $\OriP{n,k,\dots,k,n'}$ for some $n,n'>k>0$. 
Note that if the conjecture is true then for any $\P$, an orientation of a path with $n$ vertices, it would hold that  
\begin{align*}
    &\P\ppleq\Ord &&\text{iff}  &&\P\not\models\HM {n-1} &&\text{iff} &&\P\not\models\HM k\text{ for any }k.
\end{align*}
In particular, all $\leq$ in Table~\ref{tab:HMpathsNumbers} could be removed. 
The distribution  of core orientations of paths satisfying $\HM n$ for $n\leq 30$ for core orientations of paths with up to 27 vertices is presented in Table~\ref{tab:HMpathsTable}. Observe that the fraction of core orientations of paths that can pp-construct $\Ord$, displayed in orange, is steadily increasing (assuming Conjecture~\ref{con:pathsWithnoHMnAreNLHard}). 

\begin{conjecture}
As $n$ tends to infinity the fraction of core orientations of paths with $n$ vertices that can pp-construct $\Ord$ tends to one.
\end{conjecture}

\begin{table}
    \centering
    \begin{tabular}{l|l|l}
    path & pp-formula & symmetric \\\hline
    $\OriP{4, 2, 1, 2, 4}$& $\OriP{3, 2, 1, 2}$ &yes\\\hline 
    $\OriP{5, 2, 1, 2, 4}$& $\OriP{3, 2, 1, 2}$ &no\\\hline 
    $\OriP{5, 2, 1, 2, 5}$& $\OriP{3, 2, 1, 2}$ &yes\\ 
    $\OriP{6, 2, 1, 2, 5}$& $\OriP{3, 2, 1, 2}$ &no\\
    $\OriP{4, 2, 1, 2, 3, 1, 2}$& $\OriP{3, 2, 1, 2}$ &no\\
    $\OriP{3, 1, 2, 1, 1, 2, 2, 3}$& $\OriP{2, 1, 2, 1, 1, 2, 2, 2, 1, 2}$ &no\\
    $\OriP{3, 1, 2, 2, 2, 1, 1, 3}$& $\OriP{2, 1, 2, 2, 2, 1, 1, 2, 1, 2}$ &no\\
    $\OriP{4, 2, 1, 1, 1, 2, 4}$& $\OriP{3, 2, 1, 1, 1, 2}$ &yes\\
    $\OriP{4, 1, 1, 2, 1, 2, 4}$& $\OriP{3, 1, 1, 2, 1, 2}$ &no
    \end{tabular}
    \caption{Smallest paths that can pp-construct $\Ord$ and a pp-formula for the construction.}
    \label{tab:NLHardPaths}
\end{table}

When looking at every other column of Table~\ref{tab:HMpathsNumbers} we can observe that the following surprising statement seems to be true.
\begin{conjecture}
For every $K\geq 0$ there exists $N\geq0$ such that for all $n\geq N$, $k\geq n-K$ we have $\operatorname{pathsHM}(n,k)=\operatorname{pathsHM}(n-2,k-2)$.
\end{conjecture}

I did not investigate these conjectures thoroughly, hence this could be an easy entrance point for future research. 
One observation I used to compute Table~\ref{tab:HMpathsNumbers} is that many orientations of paths are in the same pp-constructability class as shorter orientations of paths. To obtain a formal result we first need some more definitions. 
Let $\G=(V,E)$ be a digraph. Then the \emph{line graph} of $\G$, denoted $L(\G)$, is the digraph 
\[(E,\{((u,v),(v,w))\mid (u,v),(v,w)\in E\}).\] 
There are some initial results about line graphs of digraphs in~\cite{Aigner1967}, for example that for digraphs the map $\G\mapsto L(\G)$ is not injective.
We define $\G^\ast$ to be the digraph obtained from $\G$ by adding two vertices $\bot$ and $\top$ and adding an edge from $\bot$ to each source in $\G$ and an edge from each sink in $\G$ to $\top$. 
Note that the map $\G\mapsto L(\G^\ast)$ is injective. Observe that $\G^\ast$ cannot always be pp-constructed from $\G$. For example, $\T_3^\ast$ is isomorphic to $\T_4$, which cannot pp-construct $\T_3$. However, we can show that $\G$ and the line graph of $\G^\ast$ have the same pp-constructability type.

\begin{lemma}\label{lem:betterLineGraphOfGIsPPeqToG}
Let $\G$ be a graph with at least one edge. Then $\G\ppeq L(\G^\ast)$.
\end{lemma}
\begin{proof}
First we show $\G\ppleq L(\G^\ast)$.
Let $\H$ be the third pp-power of $\G$ given by
\begin{center}
\begin{tikzpicture}[scale=0.5]
\node[var-f,label=below:{$x_1$}] (x1) at (0,0) {};
\node[var-f,label=below:{$x_3$}] (x2) at (2,0) {};
\node[var-f,label=below:{$x_2$}] (x3) at (1,0) {};

\node[var-f,label=above:{$y_1$}] (y1) at (0,1) {};
\node[var-f,label=above:{$y_3$}] (y2) at (2,1) {};
\node[var-f,label=above:{$y_2$}] (y3) at (1,1) {};
\path[>=stealth'] 
        (x1) edge[->] (x3)
        (x2) edge[dashed] (y1)
        (y2) edge[dashed] (y3)
        ;
\end{tikzpicture}
\end{center}
and let 
$(u',v')$ be an edge in $\G$. It is easy to verify that the map 
\begin{align*}
    L(\G^\ast)&\to\H\\
    (u,v)&\mapsto
    \begin{cases}
    (u,v,v)&\text{if $u,v\in V(\G)$}\\
    (u,u,u)&\text{if $v=\top$}\\
    (u',v',v)&\text{if $u=\bot$}
    \end{cases}
\intertext{is an embedding and that the map}
    \H&\to L(\G^\ast)\\
    (u,v,w)&\mapsto 
    \begin{cases}
    (u,v)&\text{if } (u,v)\in E\text{ and } v=w\\
    (u,\top)&\text{if }(u,v)\notin E\text{ and } v=w\\
    (\bot,v)&\text{if }(u,v)\in E\text{ and } v\not=w\\
    (u',v')&\text{if }(u,v)\notin E\text{ and } v\not=w
    \end{cases}
\end{align*}
is a homomorphism. Hence, $\G\ppleq L(\G^\ast)$. To show $\G\ppgeq L(\G^\ast)$ let $\H'$ be the first pp-power of $L(\G^\ast)$ given by the formula
\begin{center}
\begin{tikzpicture}[scale=0.5]
\node[var-f,label=right:{$x$}] (x1) at (0,0) {};
\node[var-b,] (a) at (-1,1) {};
\node[var-b,] (b) at (-1,2) {};

\node[var-f,label=right:{$y$}] (y1) at (0,1) {};
\path[>=stealth',->] 
        (x1) edge (a)
        (a) edge (b)
        (y1) edge (b)
        ;
\end{tikzpicture}
\end{center}
and consider the map $\H'\to\G, (u,v)\mapsto v$. To show that this is a homomorphism let $((u,v),(u',v')$ be an edge in $\H'$. Then there are $(u_1,v_1),(u_2,v_2)$ in $\H'$ such that 
\[(u,v)\toEdge (u_1,v_1)\toEdge(u_2,v_2)\fromEdge(u',v')\]
is a path in $L(\G^\ast)$. Hence, $v=u_1$ and $v_1=u_2=v'$. Therefore, $(v,v')$ is equal to $(u_1,v_1)$, which is an edge in $\G$ as desired.
Now consider a map $\G\to\H'$ that sends each vertex $v$ to some edge $(u,v)$ in $\G^\ast$. To verify that such a map is a homomorphism let $(v,v')$ be an edge in $\G$ and let $(u,v)$ and $(u',v')$ be the images of $v$ and $v'$, respectively. Let $w$ be a vertex of $\G^\ast$ such that $(v',w)$ is an edge in $\G^\ast$. Note that if $v'$ has no successors in $\G$, then $w=\top$. Then
\[(u,v)\toEdge (v,v')\toEdge(v',w)\fromEdge(u',v')\]
is a path in $L(\G^\ast)$. Therefore, $((u,v),(u',v'))$ is an edge in $\H'$ as desired. Hence, both maps are homomorphisms and $\G\ppleq L(\G^\ast)$.
\end{proof}

We can apply this lemma to orientations of paths and orientations of cycles. 

\begin{theorem}\label{thm:pathsAdd1Sub1}
Let $(n_1,\dots,n_k)$ be a tuple of positive natural numbers. Then
\begin{enumerate}
    \item $\OriP{n_1,\dots,n_k}\ppeq\OriP{n_1+1,\dots,n_k+1}$ and
    \item if $k$ is even, then $\OriC{n_1,\dots,n_k}\ppeq\OriC{n_1+1,\dots,n_k+1}$. 
\end{enumerate}
\end{theorem}
\begin{proof}
Observe that
\begin{align*}
    L(\OriP{n_1,\dots,n_k}^\ast)&\text{ is hom. eq. to }\OriP{n_1+1,\dots,n_k+1}\text{ and that}\\
L(\OriC{n_1,\dots,n_k}^\ast)&\text{ is hom. eq. to } \OriC{n_1+1,\dots,n_k+1}\text{ if $k$ is even}.
\end{align*} 
Hence, both statements follow from Lemma~\ref{lem:betterLineGraphOfGIsPPeqToG}.
\end{proof}
In particular, we have shown that $\P_n=\OriP{n}\ppeq\OriP{n+1}=\P_{n+1}$ and that $\T_3=\OriC{2,1}\ppeq\OriC{3,2}\ppeq\OriC{k+1,k}$ for all $n,k\geq1$.
Somewhat surprisingly it seems to have been unnoticed that we can combine the following theorems to show that the L-NL dichotomy conjecture is true for orientations of paths.
A theorem of Kazda.
\begin{theorem}[Theorem 2.7 in \cite{Kazda-n-permute}]
Let $\A$ be a finite relational structure such that $\cocsp(\A)$ is in linear Datalog  and $\A\models\HM n$ for some $n$. Then $\cocsp(\A)$ is in symmetric linear Datalog, in particular $\csp(\A)$ is in $\Lclass$.
\end{theorem}

A theorem of Barto, Kozik, and Willard.
\begin{theorem}[Theorem 7 in \cite{BartoKozikWillard}]
Let $\A$ be a finite relational structure such that $\A\models\NU n$ for some $n$. Then $\csp(\A)$ has bounded pathwidth duality and $\csp(\A)$ is in $\NL$. 
\end{theorem}

A theorem of Dalmau.
\begin{theorem}[Theorem 5 in \cite{Dalmau_2005LinearDatalog}]\label{thm:DalmauBoundedPathwidthDualityIffLinearDL}
Let $\A$ be a relational structure. Then $\csp(\A)$ has bounded pathwidth duality if and only if $\cocsp(\A)$ is in linear Datalog.
\end{theorem}

Combining these theorems with Lemma~\ref{lem:noHMimpliesACanppDefineOrder} we obtain that the NL-dichotomy conjecture is true for orientations of paths.
\begin{corollary}\label{cor:NLDichotomyHoldForPaths}
Let $\A$ be a finite relational structure such that $\cocsp(\A)$ is in linear Datalog. This is in particular true if $\A$ is an orientation of a path. Then exactly one of the following is true
\begin{enumerate}
    \item $\A\models\HM n$ for some $n\geq1$ and $\cocsp(\A)$ is in symmetric linear Datalog; in this case $\csp(\A)$ is in $\Lclass$ or
    \item $\A\not\models\HM n$ for any $n\geq1$ and $\A\ppleq\Ord$; in this case $\csp(\A)$ is $\NL$-complete. 
\end{enumerate}
\end{corollary}




\section{Paths with Two Changes in Direction}
Let us now consider orientations of paths $\P=\OriP{n_1,\dots,n_k}$, where $k$ is small. If $k\leq3$ then, by Theorem~\ref{thm:pathsAdd1Sub1}, $\P$ is in the same pp-constructability class as $\P_0$, $\P_1$ or $\OriP{n,1,k}$ for some $n\geq k\geq2$. 


\begin{lemma}
Let $n\geq k\geq 2$. Then 
 $\OriP{n',1,k'}\ppleq\OriP{n,1,k}$ for all $n'\geq n$ $k'\geq k$.
\end{lemma}
\begin{proof}
Consider the following two formulas 
\begin{center}
\begin{tikzpicture}[scale=0.5]
\node[var-f,label=left:{$x$}] (x1) at (0,0) {};
\node[var-b] (a) at (1,0) {};
\node[var-b] (b) at (1,1) {};
\node[var-b] (c) at (1,2) {};

\node[var-f,label=left:{$y$}] (y1) at (0,1) {};
\path[>=stealth',->] 
        (x1) edge (y1)
        (a) edge (y1)
        (a) edge (b)
        (b) edge (c)
        ;
\end{tikzpicture}\hspace{20mm}
\begin{tikzpicture}[scale=0.5]
\node[var-f,label=left:{$x$}] (x1) at (0,0) {};
\node[var-b] (a) at (1,-1) {};
\node[var-b] (b) at (1,0) {};
\node[var-b] (c) at (1,1) {};

\node[var-f,label=left:{$y$}] (y1) at (0,1) {};
\path[>=stealth',->] 
        (x1) edge (y1)
        (x1) edge (c)
        (a) edge (b)
        (b) edge (c)
        ;
\end{tikzpicture}
\end{center}
The first one shows $\OriP{n,1,k+1}\ppleq\OriP{n,1,k}$ and the second one shows $\OriP{n+1,1,k}\ppleq\OriP{n,1,k}$.
\end{proof}

To separate orientations of paths of the form $\OriP{n,1,k}$ we introduce a new minor condition.

\begin{definition}\label{def:elevatorGFS}
Let $n\geq 1$. The \emph{$n$-elevator condition}, denoted $\GFS n$, consists of the following identities
\begin{align*}
    f(x,x,y_1,y_2,\dots,y_n)&\approx f(z,y_1,z,y_2,\dots,y_n)\\
    &\approx f(z,y_2,y_1,z,\dots,y_n)\\
    &\phantom{{}\approx{}}\vdots\\
    &\approx f(z,y_n,y_1,y_2,\dots,z).
\end{align*}
\end{definition}
Note that $\GFS 1$ is equal to $\HM 1$.
Similar to $n$-braids we can define the following.
Let $\A$ be a structure. An \emph{$n$-staircase of $\A$} consists of pp-formulas $\phi_0(x_0,y_0),\dots,\phi_{n+1}(x_{n+1},y_{n+1})$ together with $(n+3)$-tuples $t_1,\dots,t_n$ such that 
\begin{align*}
    \A\models{} &\phi_{0}(t_{k,0},t_{k,1})\AND\dots\AND \phi_{k-1}(t_{k,k-1},t_{k,k})\AND\\
    &\exists x,y\colon \phi_k(t_{k,k},x)\AND\phi_{k}(y,x)\AND\phi_k(y,t_{k,k+1})\AND\\
    &\phi_{k+1}(t_{k,k+1},t_{k,k+2})\AND\dots\AND \phi_{n+1}(t_{k,n+1},t_{k,n+2})
\end{align*}
for all $1\leq k\leq n$.
This $n$-staircase can be seen as the structure $\B$ that has a homomorphism into $\A$ expanded by constants and $\phi_0^{\A},\dots,\phi_{n+1}^{\A}$, where $\B$ is
\begin{center}
    \begin{tikzpicture}[scale = 1.2]
    \node[var-b,label=above:$t_{1,0}$] (1-1) at (-1,1) {};
    \node[var-b,label=above:$t_{1,1}$] (10) at (0,1) {};
    \node[var-b] (00) at (0,0) {};
    \node[var-b] (11) at (1,1) {};
    \node[var-b,label=below:$t_{1,2}$] (01) at (1,0) {};
    \node[var-b,label=below:$t_{1,3}$] (02) at (2,0) {};
    \node at (3,0) {$\dots$};
    \node[var-b,label=below:$t_{1,n}$] (03) at (4,0) {};
    \node[var-b,label=below:$t_{1,n+1}$] (04) at (5,0) {};
    \node[var-b,label=below:$t_{1,n+2}$] (05) at (6,0) {};
\tikzset{decoration={snake,amplitude=.25mm,segment length=1.3mm,post length=0.7mm,pre length=0.7mm}}
    \draw[decorate,->,>=stealth'] (1-1) -- node[above] {\small$0$} (10);
    \draw[decorate,->,>=stealth'] (00) -- node[above] {\small$1$} (01);
    \draw[decorate,->,>=stealth'] (00) -- node[above] {\small$1$} (11);
    \draw[decorate,->,>=stealth'] (10) -- node[above] {\small$1$} (11);
    \draw[decorate,->,>=stealth'] (01) -- node[above] {\small$2$} (02);
    
    \draw[decorate,->,>=stealth'] (03) -- node[above] {\small$n$} (04);
    \draw[decorate,->,>=stealth'] (04) -- node[above] {\small$n+1$} (05);

\node[rotate=90] at (-0.5,-0.9) {$\dots$};
\node[rotate=90] at (0.5,-0.9) {$\dots$};
\node[rotate=90] at (5.5,-0.9) {$\dots$};

    \node[var-b,label=above:$t_{k,0}$] (1-1) at (-1,-2) {};
    \node[var-b,label=above:$t_{k,1}$] (10) at (0,-2) {};
    \node at (1,-2) {$\dots$};
    \node[var-b,label=above:$t_{k,k}$] (01) at (2,-2) {};
    \node[var-b] (00) at (2,-3) {};
    \node[var-b] (11) at (3,-2) {};
    \node[var-b,label=below:$t_{k,k+1}$] (02) at (3,-3) {};
    \node at (4,-3) {$\dots$};
    \node[var-b,label=below:$t_{k,n+1}$] (04) at (5,-3) {};
    \node[var-b,label=below:$t_{k,n+2}$] (05) at (6,-3) {};
\tikzset{decoration={snake,amplitude=.25mm,segment length=1.3mm,post length=0.7mm,pre length=0.7mm}}
    \draw[decorate,->,>=stealth'] (1-1) -- node[above] {\small$0$} (10);
    \draw[decorate,->,>=stealth'] (00) -- node[above] {\small$k$} (02);
    \draw[decorate,->,>=stealth'] (00) -- node[above] {\small$k$} (11);
    \draw[decorate,->,>=stealth'] (01) -- node[above] {\small$k$} (11);
    
    \draw[decorate,->,>=stealth'] (04) -- node[above] {\small$n+1$} (05);

\node[rotate=90] at (-0.5,-3.9) {$\dots$};
\node[rotate=90] at (0.5,-3.9) {$\dots$};
\node[rotate=90] at (5.5,-3.9) {$\dots$};
    \node[var-b,label=above:$t_{n,0}$] (1-1) at (-1,-5.0) {};
    \node[var-b,label=above:$t_{n,1}$] (10) at (0,-5.0) {};
    \node[var-b,label=above:$t_{n,2}$] (01) at (1,-5.0) {};
    \node[var-b,label=above:$t_{n,3}$] (02) at (2,-5.0) {};
    \node at (3,-5) {$\dots$};
    \node[var-b,label=above:$t_{n,n}$] (03) at (4,-5) {};
    \node[var-b] (00) at (4,-6) {};
    \node[var-b] (11) at (5,-5) {};
    \node[var-b,label=below:$t_{n,n+1}$] (04) at (5,-6) {};
    \node[var-b,label=below:$t_{n,n+2}$] (05) at (6,-6) {};
\tikzset{decoration={snake,amplitude=.25mm,segment length=1.3mm,post length=0.7mm,pre length=0.7mm}}
    \draw[decorate,->,>=stealth'] (1-1) -- node[above] {\small$0$} (10);
    \draw[decorate,->,>=stealth'] (10) -- node[above] {\small$1$} (01);
    \draw[decorate,->,>=stealth'] (01) -- node[above] {\small$2$} (02);
    
    \draw[decorate,->,>=stealth'] (00) -- node[above] {\small$n$} (04);
    \draw[decorate,->,>=stealth'] (00) -- node[above] {\small$n$} (11);
    \draw[decorate,->,>=stealth'] (03) -- node[above] {\small$n$} (11);
    
    \draw[decorate,->,>=stealth'] (04) -- node[above] {\small$n+1$} (05);
    
\end{tikzpicture}
\end{center}
and \begin{tikzpicture}[scale=0.7]
    \node[var-f,label=above:$x$] (x1) at (0,0) {};
    \node[var-f,label=above:$y$] (x2) at (1,0) {};
\tikzset{decoration={snake,amplitude=.25mm,segment length=1.3mm,post length=0.7mm,pre length=0.7mm}}
    \draw[decorate,->,>=stealth'] (x1) -- node[above] {\small$i$} (x2);
\end{tikzpicture} 
abbreviates 
the relation $\{(a,b)\mid \A\models\phi_i(a,b)\}$.
A structure $\A$ is called \emph{$n$-elevated} if for any $n$-staircase $\phi_1(x_1,y_1),\dots,\phi_n(x_n,y_n)$ with $t_1,\dots,t_n$ there are $a_0,\dots,a_{n+2}$ such that
\begin{enumerate}
    \item $\A\models\phi_0(a_0,a_1)\AND\dots\AND \phi_{n+1}(a_{n+1},a_{n+2})$,
    \item $t_{1,i}=\dots=t_{n,i}$ implies $a_i=t_{1,i}$ for all $0\leq i\leq n+2$, and 
    \item $t_{1,i}=t_{1,k},\dots,t_{n,i}=t_{n,k}$ implies $a_i=a_k$ for all $0\leq i<k\leq n+2$.
\end{enumerate}
Observe that a structure is $1$-elevated if and only if it is $1$-braided. 

\begin{lemma}\label{lem:GFSBraidedness}
Let $\A$ be a finite relational core structure. If $\A\models\GFS n$, then $\A$ is $n$-elevated.
\end{lemma}
\begin{proof}
Let $a\in A$ and $f\colon A^{n+2}\to A$ be a polymorphism of $\A$ satisfying $\GFS n$. Since $\A$ is a finite core we can assume without loss of generality that $f$ is idempotent. Consider an $n$-staircase $\phi_0(x_0,y_0),\dots,\phi_{n+1}(x_{n+1},y_{n+1})$ together with $(n+3)$-tuples $t_1,\dots,t_n$ of $\A$. For $0\leq i\leq n+2$ define 
\[a_i\coloneqq f(a,a,t_{1,i},\dots,t_{n,i}).\] 
If $t_{1,i}=\dots=t_{n,i}$ then 
\[f(a,a,t_{1,i},\dots,t_{n,i})=f(t_{1,i},t_{1,i},t_{1,i},\dots,t_{1,i})=t_{1,i}\]
for any $0\leq i \leq n+2$. Similarly, if
$t_{1,i}=t_{1,k},\dots,t_{n,i}=t_{n,k}$, then $a_i=a_k$ for all $0\leq i<k\leq n+2$.
Since $f$ is a polymorphism we have  $\A\models\Phi_0(a_0,a_1)$ and $\A\models\Phi_{n+1}(a_{n+1},a_{n+2})$. 
Let $1\leq i\leq n+1$. We show $\A\models\Phi_i(a_i,a_{i+1})$. By assumption there must be $b,c\in A$ such that 
\[\A\models\phi_i(t_{i,i},c)\AND\phi_{i}(b,c)\AND\phi_k(b,t_{i,i+1}).\]
Observe that
\begin{align*}
    a_i&=f(a,a,t_{1,i},\dots,t_{n,i})& a_{i+1}&=f(a,a,t_{1,i+1},\dots,t_{n,i+1})\\
    &=f(b,b,t_{1,i},\dots,t_{i,i},\dots,t_{n,i})& &=f(c,c,t_{1,i+1},\dots,t_{i,i+1},\dots,t_{n,i+1})\\
&&&=f(c,t_{i,i+1},t_{1,i+1},\dots,c,\dots,t_{n,i+1})\\
\end{align*}
and that $\A$ satisfies the formulas $\Phi_i(b,c),\Phi_i(b,t_{i,i+1})$, $\Phi_i(t_{1,i},t_{1,i+1})$, $\dots$, $\Phi_i(t_{i,i},c)$, $\dots$, and $\Phi_i(t_{n,i},t_{n,i+1})$. 
Since $f$ is a polymorphism $\A$ must also satisfy $\Phi_i(a_i,a_{i+1})$. Hence, $\A$ is $n$-elevated.
\end{proof}
Whether the other implication holds is still open.

\begin{question}
Let $\A$ be a finite relational core structure. If $\A$ is $n$-elevated, then $\A\models\GFS n$.
\end{question}

As $\GFS n$ and $\HM n$ are both generalizations of $\Maltsev$ it is interesting to see how the $\HM n$ blocker structures  with respect to $\GFS n$.
\begin{lemma}
Let $n\geq 3$. Then $\T_n\not\models\GFS {n-2}$.
\end{lemma}
\begin{proof}
For $1\leq i\leq n-2$ define $t_i\coloneqq(0,\dots,i,1,\dots,n-i)$ and observe that
\begin{center}
    \begin{tikzpicture}[scale = 1.0]
    \node[var-b,label=above:$0$] (1-2) at (-3,1) {};
    \node[var-b,label=above:$1$] (1-1) at (-2,1) {};
    \node at (-1,1) {$\dots$};
    \node[var-b,label=above:$i$] (10) at (0,1) {};
    \node[var-b] (00) at (0,0) {};
    \node[var-b] (11) at (1,1) {};
    \node[var-b,label=below:$1$] (01) at (1,0) {};
    \node at (2,0) {$\dots$};
    \node[var-b,label=above:$n-i-1$] (04) at (3,0) {};
    \node[var-b,label=below:$n-i$] (05) at (4,0) {};

    \draw[->,>=stealth'] (1-2) -- (1-1);
    \draw[->,>=stealth'] (00) --  (01);
    \draw[->,>=stealth'] (00) --  (11);
    \draw[->,>=stealth'] (10) --  (11);
    
    \draw[->,>=stealth'] (04) --  (05);
    
\end{tikzpicture}
\end{center}
is a path in $\T_n$. Note that there is no directed path of length $n$ in $\T_n$. Hence, $(t_1,\dots,t_{n-2})$ shows that $\T_n$ is not $(n-2)$-elevated. By Lemma~\ref{lem:GFSBraidedness} we have that $\T_n\not\models\GFS{n-2}$.
\end{proof}

Combining the previous lemma with Theorem~\ref{thm:TnIsHMblocker} we obtain the following corollary.
\begin{corollary}\label{cor:GFSimpliesHM}
Let $n\geq1$ and $\A$ be a finite structure. Then 
\[\A\not\models\HM{n}\Leftrightarrow \A\ppleq\T_{n+2}\Rightarrow \A\not\models\GFS n,\]
in particular $\A\models\GFS n$ implies $\A\models \HM n$.
\end{corollary}
Observe that we could have also shown this corollary by proving that $n$-elevatedness implies $n$-braidedness. 
Next we study the behaviour of $\OriP{n,1,k}$ with respect to $\GFS m$.

\begin{lemma}\label{lem:Pn1kNotModelsGFSkM1}
Let $n\geq k\geq2$. Then $\OriP{n,1,k}\not\models\GFS {k-1}$.
\end{lemma}
\begin{proof}
Since $\OriP{n,1,k}\ppleq\OriP{k,1,k}$ we can assume without loss of generality that $n=k$. Then the vertices of $\OriP{k,1,k}$ are $0,\dots,2k+1$. 
For $1\leq i\leq k-1$ define $t_i\coloneqq(i-1,\dots,k-1,k+2,\dots,k+i+2)$ and observe that
\begin{center}
    \begin{tikzpicture}[scale = 1.0]
    \node[var-b,label=above:$i-1$] (1-2) at (-3,1) {};
    \node[var-b,label=above:$i$] (1-1) at (-2,1) {};
    \node at (-1,1) {$\dots$};
    \node[var-b,label=above:$k-1$] (10) at (0,1) {};
    \node[var-b] (00) at (0,0) {};
    \node[var-b] (11) at (1,1) {};
    \node[var-b,label=below:$k+2$] (01) at (1,0) {};
    \node at (2,0) {$\dots$};
    \node[var-b,label=above:$k+i+1$] (04) at (3,0) {};
    \node[var-b,label=below:$k+i+2$] (05) at (4,0) {};

    \draw[->,>=stealth'] (1-2) -- (1-1);
    \draw[->,>=stealth'] (00) --  (01);
    \draw[->,>=stealth'] (00) --  (11);
    \draw[->,>=stealth'] (10) --  (11);
    
    \draw[->,>=stealth'] (04) --  (05);
    
\end{tikzpicture}
\end{center}
is a path in $\OriP{k,1,k}$. Note that there is no directed path of length $k+1$ in $\OriP{k,1,k}$. Hence, $(t_1,\dots,t_{k-1})$ shows that $\OriP{k,1,k}$ is not $(n-1)$-elevated. By Lemma~\ref{lem:GFSBraidedness} we have that $\OriP{k,1,k}\not\models\GFS{k-1}$.
\end{proof}
For $n,k\leq 4$ my program verified the following conjecture.
\begin{conjecture}
Let $n\geq k\geq2$. Then $\OriP{n,1,k}\models\GFS k$.
\end{conjecture}

If this conjecture is true then it would separate the paths $\OriP{2,1,2}$, $\OriP{3,1,3}$, $\OriP{4,1,4}$, $\dots$ from one another. What happens with $\OriP{n,1,k}$ for $n> k$? It could have the same pp-constructability type as $\OriP{k,1,k}$. I failed to prove this but I discovered the following. 
\begin{lemma}
We have that $\T_3\ppeq\OriP{2,1,2}\ppeq\OriP{3,1,2}\ppeq\OriP{n,1,2,n-1}$ for all $n\geq 3$.
\end{lemma}
\begin{proof}
We show 
\[\T_3\ppeq\T_3^{-+n}\ppleq \OriP{0,n,1,2,n-1}\ppeq \OriP{n,1,2,n-1}\ppleq\T_3.\] 
Clearly, $\OriP{0,n,1,2,n-1}\ppeq \OriP{n,1,2,n-1}\ppleq\T_3$. The equality $\T_3\ppeq\T_3^{-+n}$ was shown in Lemma~\ref{lem:4ElementsT3leqT3n}.  
Let $\H$ be the $(n+1)$-st pp-power of $\T_3^{-+n}$ given by 
\begin{center}
    \begin{tikzpicture}[scale=0.7]
    \node[var-f,label=above:$3$] (x1) at (0,0) {};
    \node[var-f,label=above:$y_n$] (x2) at (-1,0) {};
    \node[var-f,label=above:$y_{n-1}$] (x3) at (-2,0) {};
    \node at (-3.5,0) {$\dots$};
    \node[var-f,label=above:$y_2$] (x5) at (-5,0) {};
    \node[var-f,label=above:$y_1$] (x6) at (-6,0) {};
    
    \node[var-f,label=below:$x_{n+1}$] (y1) at (0,-1) {};
    \node[var-f,label=below:$x_n$] (y2) at (-1,-1) {};
    \node at (-2.5,-1) {$\dots$};
    \node[var-f,label=below:$x_3$] (y4) at (-4,-1) {};
    \node[var-f,label=below:$x_2$] (y5) at (-5,-1) {};
    \node[var-f,label=below:$x_{1}$] (y6) at (-6,-1) {};
    \path[>=stealth']
        (y1) edge[dashed] (x2)
        (x3) edge[<-] (y2)
        (x5) edge[<-] (y4)
        (x6) edge[<-] (y5)
        (y5) edge[dashed] (y6)
        ;
    \end{tikzpicture}
\end{center}
The core of $\H$ for $n=4$ is presented in Figure~\ref{fig:CoreOfHIsP4123}.
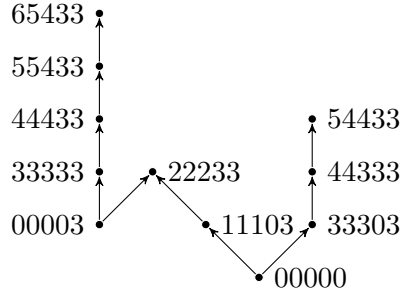
\begin{figure}
    \centering
\begin{tikzpicture}[scale=0.7]
\node[var-b, label={right:$00000$}] (00000) at (1,0) {};
\node[var-b, label={right:$11103$}] (11103) at (0,1) {};
\node[var-b, label={left:$00003$}] (00003) at (-2,1) {};
\node[var-b, label={right:$33303$}] (33303) at (2,1) {};
\node[var-b, label={right:$22233$}] (22233) at (-1,2) {};
\node[var-b, label={left:$33333$}] (33333) at (-2,2) {};
\node[var-b, label={right:$44333$}] (44133) at (2,2) {};
\node[var-b, label={left:$44433$}] (44433) at (-2,3) {};
\node[var-b, label={right:$54433$}] (52433) at (2,3) {};
\node[var-b, label={left:$55433$}] (55433) at (-2,4) {};
\node[var-b, label={left:$65433$}] (65433) at (-2,5) {};
\path[->,>=stealth']
(00000) edge (11103)
(00000) edge (33303)
(11103) edge (22233)
(00003) edge (22233)
(00003) edge (33333)
(33303) edge (44133)
(33333) edge (44433)
(44133) edge (52433)
(44433) edge (55433)
(55433) edge (65433)
;
\end{tikzpicture}
    \caption{The core of $\H$ for $n=4$.}
    \label{fig:CoreOfHIsP4123}
\end{figure}
Observe the following
\begin{itemize}
    \item $(00\dots003)$ has an outgoing path of length $n$,
    \item $(00\dots000)$ has an outgoing path of length $n-1$,
    \item there is an edge from $(00\dots003)$ to $(22\dots233)$, and
    \item there is a path of length two from $(00\dots000)$ over $(11\dots103)$ to $(22\dots233)$. 
    \end{itemize}
Hence, $\OriP{0,n,1,2,n-1}$ has an injective homomorphism into $\H$.
Let $v$ be a vertex in $\H$.
\begin{itemize}
    \item every directed path in $\H$ has length at most $n$,
    \item if $v$ has in-degree at least two then $v$ is of the form $(22\dots2a3)$ and $v$ has no outgoing edges,
    \item if $v$ has in-degree at least two and an incoming path of length at least two, then $v=(22\dots233)$, and
    \item if $v$ has out-degree at least two then $v$ is of the form $(00\dots00a)$ and $v$ has no incoming edges.
\end{itemize}
Combining these observation we conclude that the core of $\H$ is isomorphic to $\OriP{0,n,1,2,n-1}$. 
\end{proof}

By replacing the constant 3 in the pp-formula of the previous proof by 0 we can show that $\T_3\ppeq \T_3^{-+n}\cupdot \P_{n+1}$.
Recall from Problem~\ref{pro:lowerCoversOfTthree} that we are still looking for lover covers of $\T_3$. The digraph $\OriP{3,1,3}$ lies strictly below $\T_3$ and still satisfies lots of strong identities. Hence, it is a good candidate.
\begin{question}\label{pro:lowerCoversOfTthreePathCandidate}
Is $\OriP{3,1,3}$ a lower cover of $\T_3$?
\end{question}






\section{Level-wise Hagemann-Mitschke Chains}\label{sec:levelWiseSatisfyability}
The next goal is to show that orientations of paths of the form $\OriP{n,1,k}$ satisfy $\HM2$. In particular this shows that we cannot use $\HM n$'s to separate paths of the form $\OriP{n,1,k}$.
We first have to introduce a new concept for balanced digraphs.
If $\G$ is a balanced digraph (in particular, if $\G$ is an orientation of a path), the test whether $\G\models\Sigma$ can sometimes be significantly simplified. 
Define $\G^k_{\height}$ as the subgraph of $\G^k$ induced by the set 
\[\{(v_1,\dots,v_k)\in V(\G^k)\mid \height(v_1)=\dots=\height(v_K)\}.\] 
A \emph{level-wise polymorphism} is a homomorphism from $\G^k_{\height}$ to $\G$.
We say that a minor condition $\Sigma$ is \emph{level-wise satisfied} by $\G$, denoted $\G\models\Sigma$ level-wise, if $\G$ has level-wise polymorphisms that satisfy $\Sigma$.



While we do not have a general construction, for many minor conditions we can show that if a linear condition is level-wise satisfied by $\G$, then it is satisfied by $\G$. The idea is to start with polymorphisms satisfying the identities level-wise, and then extend those polymorphisms for tuples of vertices that are not all on the same level, in such a way as to satisfy the identities. Similar constructions have appeared in \cite{SpecialTriads,BartoB13,Bulin18,BulinDelicJacksonNiven}. 

\begin{lemma}\label{lem:HMLevelwise}
Let $\G$ be a balanced digraph and $n\geq1$. Then $\G\models\HM n$ level-wise if and only if $\G\models\HM n$.
\end{lemma}
\begin{proof}
If $\G\models\HM n$ then $\G\models\HM n$ level-wise. To show the other direction let $p'_1,\dots,p_n'$ be level-wise polymorphisms of $\G$ that satisfy $\HM n$. Note that the map $h\colon \G\to\G, v\mapsto p'_i(v,v,v)$ does not depend on the choice of $i$. Observe that $h$ is an endomorphism of $\G$.
Define maps $p_1,\dots,p_n$ from $\G^3$ to $\G$ by the following table: 
\begin{center}
\begin{adjustbox}{width=\columnwidth,center}
\begin{tabular}{l|llcl}
 &  $p_1(v_1,v_2,v_3)$ & $p_2(v_1,v_2,v_3)$ & $\dots$ & $p_n(v_1,v_2,v_3)$\\\hline
$\height(v_1)=\height(v_2)=\height(v_3)$ & $p'_1(v_1,v_2,v_3)$ & $p'_2(v_1,v_2,v_3)$ & $\dots$ & $p'_n(v_1,v_2,v_3)$\\
$\height(v_1)=\height(v_2)\neq\height(v_3)$ & $h(v_2)$ & $h(v_2)$& & $h(v_3)$\\
$\height(v_1)\neq\height(v_2)=\height(v_3)$ & $h(v_1)$ & $h(v_2)$& & $h(v_2)$\\
$\height(v_1)=\height(v_3)\neq\height(v_2)$ & $h(v_2)$ & $h(v_2)$& & $h(v_2)$\\
all levels different & $h(v)$ & $h(v)$& & $h(v)$ 
\end{tabular}
\end{adjustbox}
\end{center}
where $v$ is the vertex in $\{v_1,v_2,v_3\}$ with the smallest level.
Note that $p_1,\dots,p_n$ are polymorphisms.
Let $u,v$ be vertices of $\G$. If $\height(u)=\height(v)$, then
\begin{align*}
p_1(u,u,u) &=p_1'(u,u,u) =p_1'(u,v,v)= p_1(u,v,v) \\
p_{i}(u,u,v) & = p'_{i}(u,u,v) = p'_{i+1}(u,v,v)= p_{i+1}(u,v,v) && \text{ for all } i \in [n-1]  \\
p_n(u,u,v) & =p'_n(u,u,v)= p'_n(v,v,v)= p_n(v,v,v).
\end{align*}
If $\height(u)\neq\height(v)$, then
\begin{align*}
p_1(u,u,u) =p_1'(u,u,u) &=h(u)= p_1(u,v,v) \\
p_{i}(u,u,v) & = h(u)= p_{i+1}(u,v,v) && \!\text{for all } i \in [n-1] \\
p_n(u,u,v) & =h(v)= p'_n(v,v,v) = p_n(v,v,v).
\end{align*}
Hence, $p_1,\dots,p_n$ are polymorphisms of $\G$ that satisfy $\HM n$.
\end{proof}

\begin{theorem}
Let $\P$ be an orientation of a path with at most two vertices on each level, i.e., $\P$ is homomorphically equivalent to a path of the form $\OriP{n_1,1,n_2,1,\dots n_{k-1},1,n_k}$ or $\OriP{0,n_1,1,n_2,1,\dots n_{k-1},1,n_k}$ for some $n_1,\dots,n_k\geq 2$. Then $\P\models\HM2$.
\end{theorem}
\begin{proof}
By Lemma~\ref{lem:HMLevelwise} it suffices to show that $\P\models\HM2$ level-wise.
Observe that if two vertices $u,v$ in $\P$ are on the same level, then exactly one of these vertices must have in- or out-degree two.
Define maps $p_1$ and $p_2$ on $u,v,w$ with $\height(v_1)=\height(v_2)=\height(v_3)$. Note that $|\{v_1,v_2.v_3\}|\leq 2$. 
\begin{align*}
    p_1(v_1,v_2,v_3)&\coloneqq
    \begin{cases}
    v_1&\text{if }v_2=v_3\\
    v&\text{if }v_2\neq v_3\text{ and $v\in\{v_2,v_3\}$ has in- or out-degree two}
    \end{cases}\\
    p_1(v_1,v_2,v_3)&\coloneqq
    \begin{cases}
    v_3&\text{if }v_1=v_2\\
    v&\text{if }v_1\neq v_2\text{ and $v\in\{v_1,v_2\}$ has in- or out-degree two}
    \end{cases}
\end{align*}

First we verify that $p_1,p_2$ satisfy $\HM2$ level-wise. Let $u,v$ be vertices in $\P$ on the same level. If $u=v$, then all identities hold, since $p_1$ and $p_2$ are idempotent. 
If $u\neq v$, then let  $w\in\{u,v\}$ be the vertex with in- or out-degree two: 
\begin{align*}
p_1(u,u,u) &=u= p_1(u,v,v) \\
p_{1}(u,u,v) & = w= p_{2}(u,v,v) && \text{ for all } i \in \{1,\dots,n-1\} \\
p_2(u,u,v) & =v= p_2(v,v,v)\\
\end{align*}

To show that $p_1$ is a level-wise polymorphism let $(u_1,v_1), (u_2,v_2), (u_3,v_3)$ be edges in $\P$ with $\height(u_1)=\height(u_2)=\height(u_3)$. Then $v_1,v_2,v_3$ are also on the same level. Let $u=p_1(u_1,u_2,u_3)$ and $v=p_1(v_1,v_2,v_3)$.
If $u_1,u_2,u_3$ have out-degree one, then $u_1=u_2=u_3=u$ and $v_1=v_2=v_3=v$.  Hence, $(u,v)$ is an edge.
Similarly, if $v_1,v_2,v_3$ have in-degree one, then $(u,v)$ is an edge. 
Assume that there are $u'$ in $\{u_1,u_2,u_3\}$ and $v'$ in $\{v_1,v_2,v_3\}$  with in- or out-degree two. Note that $(u_i,v')$ and $(u',v_i)$ are edges in $\P$ for all $i\in\{1,2,3\}$. Hence, we only need to consider the case that $u\neq u'$ and $v\neq v'$. Then $u_2=u_3=u'$ and $u=u_1$. Analogously $v_2=v_3=v'$ and $v=v_1$. Hence, $(u,v)=(u_1,v_1)$ is an edge in $\P$ and $p_1$ is a level-wise polymorphism.
Analogously we can show that $p_2$ is a level-wise polymorphism.
Hence, $p_1,p_2$ are level-wise polymorphisms of $\P$ that satisfy $\HM2$.
\end{proof}
We will make use of level-wise polymorphisms again in Chapter~\ref{cha:hardTrees}.
Note that this theorem together with Lemma~\ref{lem:Pn1kNotModelsGFSkM1} shows in particular that the backwards implication of Corollary~\ref{cor:GFSimpliesHM} is false since $\OriP{n,1,n}\models\HM 2$ and $\OriP{n,1,n}\not\models\GFS {n-1}$ for all $n\geq2$. Hence, $\GFS n$ is stronger then $\HM n$ for all $n\geq2$.

\section{Cycles with Two Changes in Direction}\label{sec:CyclesWithTwoChangesInDirection}
Let us make some initial observations about orientations of cycles $\OriC{n,k}$ with $1\leq n<k$. By Theorem~\ref{thm:pathsAdd1Sub1} we have that $\OriC{n,k}\ppeq\OriC{1,k+1-n}$. 

\begin{lemma}
Let $k>1$. Then $\OriC{1,k}\ppleq \T_3$ and $\OriC{1,k}\ppleq\C_{k-1}$.
\end{lemma}
\begin{proof}
Clearly, $\OriC{1,k}$ is not rectangular. Hence, it is not 1-braided. Since $\OriC{1,k}$ is a core we can apply Theorem~\ref{thm:TnIsHMblocker} to obtain $\OriC{1,k}\ppleq\T_3$.
To see why $\OriC{1,k}\ppleq\C_{k-1}$ observe that the digraph pp-defined from $\OriC{1,k}$ by the formula
\begin{center}
\begin{tikzpicture}[scale=0.5]
\node[var-f,label=left:{$x$}] (x1) at (0,0) {};
\node[var-b] (a) at (1,0) {};
\node[var-f,label=right:{$y$}] (b) at (1,1) {};

\node[var-b] (y1) at (0,1) {};
\path[>=stealth',->] 
        (x1) edge (y1)
        (a) edge (y1)
        (a) edge (b)
        ;
\end{tikzpicture}
\end{center}
is homomorphically equivalent to $\C_{k-1}$.  
\end{proof}

With the help of my program I verified  $\Pol(\T_3\mathbin{\times}\C_2)\models\Sigma(\OriC{1,3},\P_4)$ and $\Pol(\OriC{1,3})\not\models\Sigma(\OriC{1,3},\P_4)$ (see Definition~\ref{def:freeCondition}). Therefore $\OriC{1,3}$ lies strictly below $\T_3\mathbin{\times}\C_2$. 
\begin{question}
Some open questions about these orientations of cycles:
\begin{enumerate}
    \item Do $\OriC{1,3},\OriC{1,4},\OriC{1,5},\dots$ form an anti-chain?
    \item Does $\Pol(\OriC{1,k})\models\Sigma_p$ hold if and only if $p$ divides $k-1?$
    \item Does $\Pol(\OriC{1,k})\models\HM 2$ hold?
\end{enumerate}\end{question}

\chapter{Digraphs of Order Four}\label{cha:4elements}

By now it is probably clear to you, dear reader, that describing all of $\DGPoset$ is beyond the capabilities of the author of this thesis. Hence, we look again at a natural simplification of the problem. This time we describe the subposet of $\DGPoset$ consisting of digraphs with at most four vertices. 
Vucaj and Bodirsky described the subposet of $\PPPoset$ consisting of all structures with a two-element domain \cite{PPPoset}. There is also ongoing work on the subposet consisting of all structures with a three-element domain \cite{VucajZhukBodirsky21ThreeElements,vucajzhuk2023submaximal}. 

Most of the results from this chapter have been found with the help of a computer program that can 
\begin{itemize}
    \item generate random minor conditions,
    \item test whether a digraph satisfies a given minor condition,
    \item generate random pp-formulas, and
    \item test whether a digraph can pp-construct another digraph with a given pp-formula.
\end{itemize}
The program was also used in Chapter~\ref{cha:hardTrees}. 

There are 3161 digraphs with at most four vertices up to isomorphism (compare to OEIS A000595). Restricting to cores decreases the number of graphs we have to consider significantly. 
There are 100 core digraphs with at most four vertices up to isomorphism. Only 28 of those have a Siggers polymorphism. They are displayed in Figure~\ref{fig:4ElementPGraphs}. The other 72 digraphs are displayed in Figure~\ref{fig:4ElementNPGraphs}. If $\G$ is a digraph, then the \emph{reverse} of $\G$ is the digraph $\G^R = (V;E^R)$ where $E^R = \{(y,x) \mid (x,y) \in E\}$. 
The operation that obtains $\G^R$ from $\G$ is called \emph{edge reversal}. Note that $\G$ and $\G^R$ have the same pp-constructability type. The digraphs in Figure~\ref{fig:4ElementNPGraphs} that are isomorphic to their reverse are marked with an $\ast$.

\tikzstyle{bullet}=[circle,fill,inner sep=0pt,minimum size=3pt]
\usetikzlibrary{arrows}
\def\scale{0.48}
\def\hdist{0mm}
\def\vdist{3mm}
\def\GAA{
\begin{tikzpicture}[scale=\scale]
\clip(-0.8,-1.5) rectangle (1.8,1.2);
\node[bullet] (0) at (1,0){};
\node[bullet] (1) at (0,0){};
\node[bullet] (3) at (0,1){};
\path[->,>=stealth']
(0) edge (1)
(0) edge (3)
(1) edge (3)
;
\node at (0.5,-1) {\small $\T_3$};
\end{tikzpicture}
}
\def\GAB{
\begin{tikzpicture}[scale=\scale]
\clip(-0.8,-1.5) rectangle (1.8,1.2);
\node[bullet] (0) at (1,0){};
\node[bullet] (1) at (0,0){};
\node[bullet] (2) at (0,1){};
\node[bullet] (3) at (1,1){};
\path[->,>=stealth']
(0) edge (1)
(0) edge (2)
(0) edge (3)
(1) edge (2)
(1) edge (3)
(2) edge (3)
;
\node at (0.5,-1) {\small $\T_4$};
\end{tikzpicture}
}
\def\GAC{
\begin{tikzpicture}[scale=\scale]
\clip(-0.8,-1.5) rectangle (1.8,1.2);
\node[bullet] (0) at (0,0){};
\path[->,>=stealth']
(0) edge[loop] (0)
;
\node at (0.5,-1) {\small $\C_1$};
\end{tikzpicture}
}
\def\GAD{
\begin{tikzpicture}[scale=\scale]
\clip(-0.8,-1.5) rectangle (1.8,1.2);
\node[bullet] (0) at (0,0){};
\node[bullet] (2) at (1,1){};
\path[->,>=stealth']
(0) edge (2)
(2) edge (0)
;
\node at (0.5,-1) {\small $\C_2$};
\end{tikzpicture}
}
\def\GAE{
\begin{tikzpicture}[scale=\scale]
\clip(-0.8,-1.5) rectangle (1.8,1.2);
\node[bullet] (0) at (0,0){};
\node[bullet] (2) at (1,1){};
\node[bullet] (3) at (1,0){};
\path[->,>=stealth']
(0) edge (2)
(2) edge (3)
(3) edge (0)
;
\node at (0.5,-1) {\small $\C_3$};
\end{tikzpicture}
}
\def\GAF{
\begin{tikzpicture}[scale=\scale]
\clip(-0.8,-1.5) rectangle (1.8,1.2);
\node[bullet] (0) at (0,0){};
\node[bullet] (1) at (0,1){};
\node[bullet] (2) at (1,1){};
\node[bullet] (3) at (1,0){};
\path[->,>=stealth']
(0) edge (1)
(1) edge (2)
(2) edge (3)
(3) edge (0)
;
\node at (0.5,-1) {\small $\C_4$};
\end{tikzpicture}
}
\def\GAG{
\begin{tikzpicture}[scale=\scale]
\clip(-0.8,-1.5) rectangle (1.8,1.2);
\node[bullet] (0) at (0,0){};
\node[bullet] (2) at (1,1){};
\path[->,>=stealth']
(0) edge (2)
;
\node at (0.5,-1) {\small $\P_1$};
\end{tikzpicture}
}
\def\GAH{
\begin{tikzpicture}[scale=\scale]
\clip(-0.8,-1.5) rectangle (1.8,1.2);
\node[bullet] (0) at (0,0){};
\node[bullet] (2) at (0,1){};
\node[bullet] (3) at (1,1){};
\path[->,>=stealth']
(0) edge (2)
(2) edge (3)
;
\node at (0.5,-1) {\small $\P_2$};
\end{tikzpicture}
}
\def\GAI{
\begin{tikzpicture}[scale=\scale]
\clip(-0.8,-1.5) rectangle (1.8,1.2);
\node[bullet] (0) at (0,0){};
\node[bullet] (1) at (0,1){};
\node[bullet] (2) at (1,1){};
\node[bullet] (3) at (1,0){};
\path[->,>=stealth']
(0) edge (1)
(1) edge (2)
(2) edge (3)
;
\node at (0.5,-1) {\small $\P_3$};
\end{tikzpicture}
}
\def\GAJ{
\begin{tikzpicture}[scale=\scale]
\clip(-0.8,-1.5) rectangle (1.8,1.2);
\node[bullet] (0) at (1,0){};
\node[bullet] (1) at (0,0){};
\node[bullet] (2) at (0,1){};
\node[bullet] (3) at (1,1){};
\path[->,>=stealth']
(0) edge (1)
(0) edge (2)
(1) edge (2)
(2) edge (3)
;
\node at (0.5,-1) {\small $\T_3^+$};
\end{tikzpicture}
}
\def\GBA{
\begin{tikzpicture}[scale=\scale]
\clip(-0.8,-1.5) rectangle (1.8,1.2);
\node[bullet] (0) at (0,1){};
\node[bullet] (1) at (0,0){};
\node[bullet] (2) at (1,0){};
\node[bullet] (3) at (1,1){};
\path[->,>=stealth']
(0) edge (1)
(0) edge (2)
(1) edge (2)
(3) edge (0)
;
\node at (0.5,-1) {\small $\T_3^-$};
\end{tikzpicture}
}
\def\GBB{
\begin{tikzpicture}[scale=\scale]
\clip(-0.8,-1.5) rectangle (1.8,1.2);
\node[bullet] (0) at (1,0){};
\node[bullet] (1) at (0,0){};
\node[bullet] (2) at (0,1){};
\node[bullet] (3) at (1,1){};
\path[->,>=stealth']
(0) edge (1)
(0) edge (2)
(0) edge (3)
(1) edge (2)
(2) edge (3)
;
\node at (0.5,-1) {\small $\T_4^{\setminus(13)}$};
\end{tikzpicture}
}
\def\GBC{
\begin{tikzpicture}[scale=\scale]
\clip(-0.8,-1.5) rectangle (1.8,1.2);
\node[bullet] (0) at (1,1){};
\node[bullet] (1) at (0,1){};
\node[bullet] (2) at (1,0){};
\node[bullet] (3) at (0,0){};
\path[->,>=stealth']
(0) edge (1)
(0) edge (2)
(1) edge (2)
(1) edge (3)
(3) edge (2)
;
\node at (0.5,-1) {\small $\T_4^{\setminus(02)}$};
\end{tikzpicture}
}
\def\GBD{
\begin{tikzpicture}[scale=\scale]
\clip(-0.8,-1.5) rectangle (1.8,1.2);
\node[bullet] (0) at (1,0){};
\node[bullet] (1) at (0,0){};
\node[bullet] (2) at (0,1){};
\node[bullet] (3) at (1,1){};
\path[->,>=stealth']
(0) edge (1)
(0) edge (2)
(1) edge (2)
(2) edge (3)
(3) edge (2)
;
\node at (0.5,-1) {\small $\T_3\hspace{-0.2mm}\cup_2\hspace{-0.2mm}\C_2$};
\end{tikzpicture}
}
\def\GBE{
\begin{tikzpicture}[scale=\scale]
\clip(-0.8,-1.5) rectangle (1.8,1.2);
\node[bullet] (0) at (0,1){};
\node[bullet] (1) at (1,1){};
\node[bullet] (2) at (0,0){};
\node[bullet] (3) at (1,0){};
\path[->,>=stealth']
(0) edge (1)
(0) edge (2)
(0) edge (3)
(1) edge (0)
(2) edge (3)
;
\node at (0.5,-1) {\small $\T_3\hspace{-0.2mm}\cup_0\hspace{-0.2mm}\C_2$};
\end{tikzpicture}
}
\def\GBF{
\begin{tikzpicture}[scale=\scale]
\clip(-0.8,-1.5) rectangle (1.8,1.2);
\node[bullet] (0) at (1,0){};
\node[bullet] (1) at (0,0){};
\node[bullet] (2) at (0,1){};
\node[bullet] (3) at (1,1){};
\path[->,>=stealth']
(0) edge (1)
(0) edge (2)
(1) edge (2)
(1) edge (3)
(2) edge (3)
;
\node at (0.5,-1) {\small $\T_4^{\setminus(03)}$};
\end{tikzpicture}
}
\def\GBG{
\begin{tikzpicture}[scale=\scale]
\clip(-0.8,-1.5) rectangle (1.8,1.2);
\node[bullet] (0) at (1,1){};
\node[bullet] (1) at (0,0){};
\node[bullet] (2) at (0,1){};
\node[bullet] (3) at (1,0){};
\path[<-,>=stealth']
(0) edge (1)
(0) edge (2)
(0) edge (3)
(1) edge (2)
(2) edge (3)
(3) edge (1)
;
\node at (0.5,-1) {\small $\C_3^+$};
\end{tikzpicture}
}
\def\GBH{
\begin{tikzpicture}[scale=\scale]
\clip(-0.8,-1.5) rectangle (1.8,1.2);
\node[bullet] (0) at (0,0){};
\node[bullet] (1) at (1,0){};
\node[bullet] (2) at (1,1){};
\node[bullet] (3) at (0,1){};
\path[<-,>=stealth']
(0) edge (1)
(0) edge (2)
(1) edge (2)
(1) edge (3)
(3) edge (0)
(3) edge (2)
;
\node at (0.5,-1) {\small $\C_3^-$};
\end{tikzpicture}
}
\def\GBI{
\begin{tikzpicture}[scale=\scale]
\clip(-0.8,-1.5) rectangle (1.8,1.2);
\node[bullet] (0) at (0,0){};
\node[bullet] (1) at (1,1){};
\node[bullet] (2) at (0,1){};
\node[bullet] (3) at (1,0){};
\path[->,>=stealth']
(0) edge (1)
(0) edge (2)
(1) edge (0)
(3) edge (0)
(3) edge (2)
;
\node at (0.5,-1) {\small $\T_3\hspace{-0.2mm}\cup_1\hspace{-0.2mm}\C_2$};
\end{tikzpicture}
}
\def\GBJ{
\begin{tikzpicture}[scale=\scale]
\clip(-0.8,-1.5) rectangle (1.8,1.2);
\node[bullet] (0) at (0,0){};
\node[bullet] (1) at (0,1){};
\node[bullet] (3) at (1,1){};
\path[->,>=stealth']
(0) edge (1)
(0) edge (3)
(1) edge (0)
(1) edge (3)
;
\node at (0.5,-1) {\small $\C_2^+$};
\end{tikzpicture}
}
\def\GCA{
\begin{tikzpicture}[scale=\scale]
\clip(-0.8,-1.5) rectangle (1.8,1.2);
\node[bullet] (0) at (1,1){};
\node[bullet] (1) at (0,1){};
\node[bullet] (3) at (0,0){};
\path[->,>=stealth']
(0) edge (1)
(0) edge (3)
(1) edge (3)
(3) edge (1)
;
\node at (0.5,-1) {\small $\C_2^-$};
\end{tikzpicture}
}
\def\GCB{
\begin{tikzpicture}[scale=\scale]
\clip(-0.8,-1.5) rectangle (1.8,1.2);
\node[bullet] (0) at (0,0){};
\node[bullet] (1) at (0,1){};
\node[bullet] (2) at (1,1){};
\node[bullet] (3) at (1,0){};
\path[->,>=stealth']
(0) edge (1)
(0) edge (2)
(1) edge (0)
(1) edge (2)
(2) edge (3)
;
\node at (0.5,-1) {\small $\C_2^{++}$};
\end{tikzpicture}
}
\def\GCC{
\begin{tikzpicture}[scale=\scale]
\clip(-0.8,-1.5) rectangle (1.8,1.2);
\node[bullet] (0) at (1,1){};
\node[bullet] (1) at (0,1){};
\node[bullet] (2) at (0,0){};
\node[bullet] (3) at (1,0){};
\path[->,>=stealth']
(0) edge (1)
(0) edge (2)
(1) edge (2)
(2) edge (1)
(3) edge (0)
;
\node at (0.5,-1) {\small $\C_2^{--}$};
\end{tikzpicture}
}
\def\GCD{
\begin{tikzpicture}[scale=\scale]
\clip(-0.8,-1.5) rectangle (1.8,1.2);
\node[bullet] (0) at (1,0){};
\node[bullet] (1) at (0,0){};
\node[bullet] (2) at (0,1){};
\node[bullet] (3) at (1,1){};
\path[->,>=stealth']
(0) edge (1)
(0) edge (2)
(0) edge (3)
(1) edge (2)
(1) edge (3)
(2) edge (1)
(2) edge (3)
;
\node at (0.5,-1) {\small $\nCk112$};
\end{tikzpicture}
}
\def\GCE{
\begin{tikzpicture}[scale=\scale]
\clip(-0.8,-1.5) rectangle (1.8,1.2);
\node[bullet] (0) at (1,0){};
\node[bullet] (1) at (0,0){};
\node[bullet] (2) at (0,1){};
\node[bullet] (3) at (1,1){};
\path[->,>=stealth']
(0) edge (1)
(0) edge (2)
(0) edge (3)
(1) edge (2)
(1) edge (3)
(2) edge (3)
(3) edge (2)
;
\node at (0.5,-1) {\small $\nCk202$};
\end{tikzpicture}
}
\def\GCF{
\begin{tikzpicture}[scale=\scale]
\clip(-0.8,-1.5) rectangle (1.8,1.2);
\node[bullet] (0) at (1,1){};
\node[bullet] (1) at (0,1){};
\node[bullet] (2) at (0,0){};
\node[bullet] (3) at (1,0){};
\path[->,>=stealth']
(0) edge (1)
(0) edge (2)
(0) edge (3)
(1) edge (0)
(1) edge (2)
(1) edge (3)
(2) edge (3)
;
\node at (0.5,-1) {\small $\nCk022$};
\end{tikzpicture}
}
\def\GCG{
\begin{tikzpicture}[scale=\scale]
\clip(-0.8,-1.5) rectangle (1.8,1.2);
\node[bullet] (0) at (0,0){};
\node[bullet] (1) at (0,1){};
\node[bullet] (2) at (1,0){};
\node[bullet] (3) at (1,1){};
\path[->,>=stealth']
(0) edge (1)
(0) edge (2)
(1) edge (3)
(3) edge (2)
;
\node at (0.5,-1) {\small $\C(1,3)$};
\end{tikzpicture}
}
\def\GCH{
\begin{tikzpicture}[scale=\scale]
\clip(-0.8,-1.5) rectangle (1.8,1.2);
\node[bullet] (0) at (0,0){};
\node[bullet] (1) at (0,1){};
\node[bullet] (2) at (1,0){};
\node[bullet,gray] (3) at (1,1){};
\path[->,>=stealth']
(0) edge (1)
(0) edge (2)
(0) edge[gray] (3)
(1) edge (0)
(1) edge (2)
(1) edge[gray] (3)
(2) edge (0)
(2) edge (1)
(2) edge[gray] (3)
;
\node at (0.5,-1) {\small G71};
\end{tikzpicture}
}
\def\GCI{
\begin{tikzpicture}[scale=\scale]
\clip(-0.8,-1.5) rectangle (1.8,1.2);
\node[bullet,gray] (0) at (1,1){};
\node[bullet] (1) at (0,1){};
\node[bullet] (2) at (0,0){};
\node[bullet] (3) at (1,0){};
\path[->,>=stealth']
(0) edge[gray] (1)
(0) edge[gray] (2)
(0) edge[gray] (3)
(1) edge (2)
(1) edge (3)
(2) edge (1)
(2) edge (3)
(3) edge (1)
(3) edge (2)
;
\node at (0.5,-1) {\small G72};
\end{tikzpicture}
}
\def\GCJ{
\begin{tikzpicture}[scale=\scale]
\clip(-0.8,-1.5) rectangle (1.8,1.2);
\node[bullet] (0) at (0,0){};
\node[bullet] (1) at (1,0){};
\node[bullet] (2) at (0,1){};
\node[bullet,gray] (3) at (1,1){};
\path[->,>=stealth']
(0) edge (1)
(0) edge (2)
(0) edge[gray] (3)
(1) edge (0)
(1) edge (2)
(1) edge[gray] (3)
(2) edge (0)
(2) edge[gray] (3)
;
\node at (0.5,-1) {\small G67};
\end{tikzpicture}
}
\def\GDA{
\begin{tikzpicture}[scale=\scale]
\clip(-0.8,-1.5) rectangle (1.8,1.2);
\node[bullet,gray] (0) at (1,1){};
\node[bullet] (1) at (0,0){};
\node[bullet] (2) at (0,1){};
\node[bullet] (3) at (1,0){};
\path[->,>=stealth']
(0) edge[gray] (1)
(0) edge[gray] (2)
(0) edge[gray] (3)
(1) edge (2)
(1) edge (3)
(2) edge (1)
(2) edge (3)
(3) edge (1)
;
\node at (0.5,-1) {\small G68};
\end{tikzpicture}
}
\def\GDB{
\begin{tikzpicture}[scale=\scale]
\clip(-0.8,-1.5) rectangle (1.8,1.2);
\node[bullet] (0) at (0,0){};
\node[bullet] (1) at (1,0){};
\node[bullet] (2) at (0,1){};
\node[bullet,gray] (3) at (1,1){};
\path[->,>=stealth']
(0) edge (1)
(0) edge (2)
(0) edge[gray] (3)
(1) edge (0)
(1) edge (2)
(2) edge (0)
(2) edge[gray] (3)
;
\node at (0.5,-1) {\small G69};
\end{tikzpicture}
}
\def\GDC{
\begin{tikzpicture}[scale=\scale]
\clip(-0.8,-1.5) rectangle (1.8,1.2);
\node[bullet] (0) at (0,0){};
\node[bullet] (1) at (0,1){};
\node[bullet] (2) at (1,0){};
\node[bullet,gray] (3) at (1,1){};
\path[->,>=stealth']
(0) edge (1)
(0) edge (2)
(1) edge (0)
(1) edge (2)
(2) edge (0)
(3) edge[gray] (0)
(3) edge[gray] (1)
;
\node at (0.5,-1) {\small G70};
\end{tikzpicture}
}
\def\GDD{
\begin{tikzpicture}[scale=\scale]
\clip(-0.8,-1.5) rectangle (1.8,1.2);
\node[bullet] (0) at (0,0){};
\node[bullet] (1) at (0,1){};
\node[bullet] (2) at (1,1){};
\node[bullet] (3) at (1,0){};
\path[->,>=stealth']
(0) edge (1)
(0) edge (2)
(0) edge (3)
(1) edge (0)
(1) edge (2)
(2) edge (3)
;
\node at (0.5,-1) {\small G58};
\end{tikzpicture}
}
\def\GDE{
\begin{tikzpicture}[scale=\scale]
\clip(-0.8,-1.5) rectangle (1.8,1.2);
\node[bullet] (0) at (1,0){};
\node[bullet] (1) at (1,1){};
\node[bullet] (2) at (0,0){};
\node[bullet] (3) at (0,1){};
\path[->,>=stealth']
(0) edge (1)
(0) edge (2)
(1) edge (2)
(1) edge (3)
(2) edge (3)
(3) edge (2)
;
\node at (0.5,-1) {\small G59};
\end{tikzpicture}
}
\def\GDF{
\begin{tikzpicture}[scale=\scale]
\clip(-0.8,-1.5) rectangle (1.8,1.2);
\node[bullet] (0) at (0,0){};
\node[bullet] (1) at (0,1){};
\node[bullet,gray] (2) at (1,1){};
\node[bullet] (3) at (1,0){};
\path[->,>=stealth']
(0) edge (1)
(0) edge[gray] (2)
(0) edge (3)
(1) edge (0)
(1) edge[gray] (2)
(3) edge (1)
(3) edge[gray] (2)
;
\node at (0.5,-1) {\small G63};
\end{tikzpicture}
}
\def\GDG{
\begin{tikzpicture}[scale=\scale]
\clip(-0.8,-1.5) rectangle (1.8,1.2);
\node[bullet,gray] (0) at (1,1){};
\node[bullet] (1) at (0,1){};
\node[bullet] (2) at (0,0){};
\node[bullet] (3) at (1,0){};
\path[->,>=stealth']
(0) edge[gray] (1)
(0) edge[gray] (2)
(0) edge[gray] (3)
(1) edge (2)
(1) edge (3)
(2) edge (1)
(3) edge (2)
;
\node at (0.5,-1) {\small G64};
\end{tikzpicture}
}
\def\GDH{
\begin{tikzpicture}[scale=\scale]
\clip(-0.8,-1.5) rectangle (1.8,1.2);
\node[bullet] (0) at (0,1){};
\node[bullet] (1) at (0,0){};
\node[bullet,gray] (2) at (1,1){};
\node[bullet] (3) at (1,0){};
\path[->,>=stealth']
(0) edge (1)
(0) edge[gray] (2)
(0) edge (3)
(1) edge (0)
(1) edge[gray] (2)
(3) edge (1)
;
\node at (0.5,-1) {\small G66};
\end{tikzpicture}
}
\def\GDI{
\begin{tikzpicture}[scale=\scale]
\clip(-0.8,-1.5) rectangle (1.8,1.2);
\node[bullet] (0) at (0,0){};
\node[bullet] (1) at (0,1){};
\node[bullet] (2) at (1,0){};
\node[bullet,gray] (3) at (1,1){};
\path[->,>=stealth']
(0) edge (1)
(0) edge (2)
(1) edge (0)
(2) edge (1)
(3) edge[gray] (0)
(3) edge[gray] (1)
;
\node at (0.5,-1) {\small G65};
\end{tikzpicture}
}
\def\GDJ{
\begin{tikzpicture}[scale=\scale]
\clip(-0.8,-1.5) rectangle (1.8,1.2);
\node[bullet] (0) at (0,0){};
\node[bullet] (1) at (0,1){};
\node[bullet] (2) at (1,1){};
\node[bullet] (3) at (1,0){};
\path[->,>=stealth']
(0) edge (1)
(0) edge (2)
(1) edge (0)
(1) edge (2)
(3) edge (0)
(3) edge (1)
;
\node at (0.5,-1) {\small G60$\ast$};
\end{tikzpicture}
}
\def\GEA{
\begin{tikzpicture}[scale=\scale]
\clip(-0.8,-1.5) rectangle (1.8,1.2);
\node[bullet] (0) at (0,0){};
\node[bullet] (1) at (0,1){};
\node[bullet] (2) at (1,0){};
\node[bullet,gray] (3) at (1,1){};
\path[->,>=stealth']
(0) edge (1)
(0) edge (2)
(1) edge (0)
(1) edge[gray] (3)
(2) edge (1)
(2) edge[gray] (3)
;
\node at (0.5,-1) {\small G61};
\end{tikzpicture}
}
\def\GEB{
\begin{tikzpicture}[scale=\scale]
\clip(-0.8,-1.5) rectangle (1.8,1.2);
\node[bullet] (0) at (0,1){};
\node[bullet] (1) at (0,0){};
\node[bullet] (2) at (1,0){};
\node[bullet,gray] (3) at (1,1){};
\path[->,>=stealth']
(0) edge (1)
(0) edge (2)
(1) edge (0)
(2) edge (1)
(3) edge[gray] (0)
(3) edge[gray] (2)
;
\node at (0.5,-1) {\small G62};
\end{tikzpicture}
}
\def\GEC{
\begin{tikzpicture}[scale=\scale]
\clip(-0.8,-1.5) rectangle (1.8,1.2);
\node[bullet] (0) at (0,0){};
\node[bullet] (1) at (1,1){};
\node[bullet] (2) at (0,1){};
\node[bullet] (3) at (1,0){};
\path[->,>=stealth']
(0) edge (1)
(0) edge (2)
(1) edge (2)
(1) edge (3)
(3) edge (0)
;
\node at (0.5,-1) {\small G56};
\end{tikzpicture}
}
\def\GED{
\begin{tikzpicture}[scale=\scale]
\clip(-0.8,-1.5) rectangle (1.8,1.2);
\node[bullet] (0) at (0,1){};
\node[bullet] (1) at (1,1){};
\node[bullet] (2) at (0,0){};
\node[bullet] (3) at (1,0){};
\path[->,>=stealth']
(0) edge (1)
(0) edge (2)
(1) edge (2)
(2) edge (3)
(3) edge (1)
;
\node at (0.5,-1) {\small G57};
\end{tikzpicture}
}
\def\GEE{
\begin{tikzpicture}[scale=\scale]
\clip(-0.8,-1.5) rectangle (1.8,1.2);
\node[bullet] (0) at (0,0){};
\node[bullet] (1) at (0,1){};
\node[bullet] (2) at (1,1){};
\node[bullet] (3) at (1,0){};
\path[->,>=stealth']
(0) edge (1)
(0) edge (2)
(0) edge (3)
(1) edge (0)
(1) edge (2)
(1) edge (3)
(2) edge (0)
(2) edge (1)
(2) edge (3)
(3) edge (0)
(3) edge (1)
(3) edge (2)
;
\node at (0.5,-1) {\small G55$\ast$};
\end{tikzpicture}
}
\def\GEF{
\begin{tikzpicture}[scale=\scale]
\clip(-0.8,-1.5) rectangle (1.8,1.2);
\node[bullet] (0) at (0,0){};
\node[bullet] (1) at (0,1){};
\node[bullet] (2) at (1,1){};
\node[bullet] (3) at (1,0){};
\path[->,>=stealth']
(0) edge (1)
(0) edge (2)
(0) edge (3)
(1) edge (0)
(1) edge (2)
(1) edge (3)
(2) edge (0)
(2) edge (1)
(2) edge (3)
(3) edge (0)
(3) edge (1)
;
\node at (0.5,-1) {\small G54$\ast$};
\end{tikzpicture}
}
\def\GEG{
\begin{tikzpicture}[scale=\scale]
\clip(-0.8,-1.5) rectangle (1.8,1.2);
\node[bullet] (0) at (0,0){};
\node[bullet] (1) at (0,1){};
\node[bullet] (2) at (1,0){};
\node[bullet] (3) at (1,1){};
\path[->,>=stealth']
(0) edge (1)
(0) edge (2)
(0) edge (3)
(1) edge (0)
(1) edge (2)
(1) edge (3)
(2) edge (0)
(2) edge (1)
(2) edge (3)
(3) edge (0)
;
\node at (0.5,-1) {\small G52};
\end{tikzpicture}
}
\def\GEH{
\begin{tikzpicture}[scale=\scale]
\clip(-0.8,-1.5) rectangle (1.8,1.2);
\node[bullet] (0) at (0,0){};
\node[bullet] (1) at (1,1){};
\node[bullet] (2) at (1,0){};
\node[bullet] (3) at (0,1){};
\path[->,>=stealth']
(0) edge (1)
(0) edge (2)
(0) edge (3)
(1) edge (0)
(1) edge (2)
(1) edge (3)
(2) edge (0)
(2) edge (3)
(3) edge (0)
(3) edge (2)
;
\node at (0.5,-1) {\small G53};
\end{tikzpicture}
}
\def\GEI{
\begin{tikzpicture}[scale=\scale]
\clip(-0.8,-1.5) rectangle (1.8,1.2);
\node[bullet] (0) at (0,0){};
\node[bullet] (1) at (0,1){};
\node[bullet] (3) at (1,0){};
\path[->,>=stealth']
(0) edge (1)
(0) edge (3)
(1) edge (0)
(1) edge (3)
(3) edge (0)
(3) edge (1)
;
\node at (0.5,-1) {\small G37$\ast$};
\end{tikzpicture}
}
\def\GEJ{
\begin{tikzpicture}[scale=\scale]
\clip(-0.8,-1.5) rectangle (1.8,1.2);
\node[bullet] (0) at (0,0){};
\node[bullet] (1) at (0,1){};
\node[bullet] (2) at (1,0){};
\node[bullet] (3) at (1,1){};
\path[->,>=stealth']
(0) edge (1)
(0) edge (2)
(0) edge (3)
(1) edge (0)
(1) edge (2)
(1) edge (3)
(2) edge (0)
(2) edge (1)
(3) edge (0)
(3) edge (2)
;
\node at (0.5,-1) {\small G51$\ast$};
\end{tikzpicture}
}
\def\GFA{
\begin{tikzpicture}[scale=\scale]
\clip(-0.8,-1.5) rectangle (1.8,1.2);
\node[bullet] (0) at (0,1){};
\node[bullet] (1) at (0,0){};
\node[bullet] (2) at (1,0){};
\node[bullet] (3) at (1,1){};
\path[->,>=stealth']
(0) edge (1)
(0) edge (2)
(0) edge (3)
(1) edge (0)
(1) edge (2)
(1) edge (3)
(2) edge (0)
(2) edge (1)
(3) edge (2)
;
\node at (0.5,-1) {\small G38};
\end{tikzpicture}
}
\def\GFB{
\begin{tikzpicture}[scale=\scale]
\clip(-0.8,-1.5) rectangle (1.8,1.2);
\node[bullet] (0) at (1,0){};
\node[bullet] (1) at (0,0){};
\node[bullet] (2) at (0,1){};
\node[bullet] (3) at (1,1){};
\path[->,>=stealth']
(0) edge (1)
(0) edge (2)
(0) edge (3)
(1) edge (0)
(1) edge (2)
(2) edge (0)
(2) edge (1)
(3) edge (1)
(3) edge (2)
;
\node at (0.5,-1) {\small G39};
\end{tikzpicture}
}
\def\GFC{
\begin{tikzpicture}[scale=\scale]
\clip(-0.8,-1.5) rectangle (1.8,1.2);
\node[bullet] (0) at (0,0){};
\node[bullet] (1) at (0,1){};
\node[bullet] (2) at (1,1){};
\node[bullet] (3) at (1,0){};
\path[->,>=stealth']
(0) edge (1)
(0) edge (2)
(0) edge (3)
(1) edge (0)
(1) edge (2)
(1) edge (3)
(2) edge (0)
(2) edge (3)
(3) edge (0)
;
\node at (0.5,-1) {\small G49$\ast$};
\end{tikzpicture}
}
\def\GFD{
\begin{tikzpicture}[scale=\scale]
\clip(-0.8,-1.5) rectangle (1.8,1.2);
\node[bullet] (0) at (0,0){};
\node[bullet] (1) at (0,1){};
\node[bullet] (2) at (1,0){};
\node[bullet] (3) at (1,1){};
\path[->,>=stealth']
(0) edge (1)
(0) edge (2)
(0) edge (3)
(1) edge (0)
(1) edge (2)
(1) edge (3)
(2) edge (0)
(2) edge (3)
(3) edge (1)
(3) edge (2)
;
\node at (0.5,-1) {\small G50$\ast$};
\end{tikzpicture}
}
\def\GFE{
\begin{tikzpicture}[scale=\scale]
\clip(-0.8,-1.5) rectangle (1.8,1.2);
\node[bullet] (0) at (0,0){};
\node[bullet] (1) at (0,1){};
\node[bullet] (2) at (1,0){};
\node[bullet] (3) at (1,1){};
\path[->,>=stealth']
(0) edge (1)
(0) edge (2)
(0) edge (3)
(1) edge (0)
(1) edge (2)
(1) edge (3)
(2) edge (0)
(2) edge (3)
(3) edge (1)
;
\node at (0.5,-1) {\small G43};
\end{tikzpicture}
}
\def\GFF{
\begin{tikzpicture}[scale=\scale]
\clip(-0.8,-1.5) rectangle (1.8,1.2);
\node[bullet] (0) at (1,1){};
\node[bullet] (1) at (0,1){};
\node[bullet] (2) at (0,0){};
\node[bullet] (3) at (1,0){};
\path[->,>=stealth']
(0) edge (1)
(0) edge (2)
(0) edge (3)
(1) edge (0)
(1) edge (2)
(2) edge (1)
(2) edge (3)
(3) edge (1)
(3) edge (2)
;
\node at (0.5,-1) {\small G44};
\end{tikzpicture}
}
\def\GFG{
\begin{tikzpicture}[scale=\scale]
\clip(-0.8,-1.5) rectangle (1.8,1.2);
\node[bullet] (0) at (0,0){};
\node[bullet] (1) at (1,0){};
\node[bullet] (2) at (0,1){};
\node[bullet] (3) at (1,1){};
\path[->,>=stealth']
(0) edge (1)
(0) edge (2)
(0) edge (3)
(1) edge (0)
(1) edge (2)
(1) edge (3)
(2) edge (0)
(2) edge (3)
(3) edge (2)
;
\node at (0.5,-1) {\small G45$\ast$};
\end{tikzpicture}
}
\def\GFH{
\begin{tikzpicture}[scale=\scale]
\clip(-0.8,-1.5) rectangle (1.8,1.2);
\node[bullet] (0) at (0,0){};
\node[bullet] (1) at (1,0){};
\node[bullet] (3) at (0,1){};
\path[->,>=stealth']
(0) edge (1)
(0) edge (3)
(1) edge (0)
(1) edge (3)
(3) edge (0)
;
\node at (0.5,-1) {\small G15$\ast$};
\end{tikzpicture}
}
\def\GFI{
\begin{tikzpicture}[scale=\scale]
\clip(-0.8,-1.5) rectangle (1.8,1.2);
\node[bullet] (0) at (0,1){};
\node[bullet] (1) at (0,0){};
\node[bullet] (2) at (1,1){};
\node[bullet] (3) at (1,0){};
\path[->,>=stealth']
(0) edge (1)
(0) edge (2)
(0) edge (3)
(1) edge (0)
(1) edge (2)
(1) edge (3)
(2) edge (0)
(3) edge (1)
;
\node at (0.5,-1) {\small G41};
\end{tikzpicture}
}
\def\GFJ{
\begin{tikzpicture}[scale=\scale]
\clip(-0.8,-1.5) rectangle (1.8,1.2);
\node[bullet] (0) at (0,1){};
\node[bullet] (1) at (1,1){};
\node[bullet] (2) at (0,0){};
\node[bullet] (3) at (1,0){};
\path[->,>=stealth']
(0) edge (1)
(0) edge (2)
(1) edge (0)
(1) edge (2)
(2) edge (0)
(2) edge (3)
(3) edge (0)
(3) edge (2)
;
\node at (0.5,-1) {\small G42};
\end{tikzpicture}
}
\def\GGA{
\begin{tikzpicture}[scale=\scale]
\clip(-0.8,-1.5) rectangle (1.8,1.2);
\node[bullet] (0) at (0,0){};
\node[bullet] (1) at (0,1){};
\node[bullet] (2) at (1,0){};
\node[bullet] (3) at (1,1){};
\path[<-,>=stealth']
(0) edge (1)
(0) edge (2)
(0) edge (3)
(1) edge (0)
(1) edge (2)
(1) edge (3)
(2) edge (0)
(3) edge (2)
;
\node at (0.5,-1) {\small G20};
\end{tikzpicture}
}
\def\GGB{
\begin{tikzpicture}[scale=\scale]
\clip(-0.8,-1.5) rectangle (1.8,1.2);
\node[bullet] (0) at (1,0){};
\node[bullet] (1) at (0,0){};
\node[bullet] (2) at (0,1){};
\node[bullet] (3) at (1,1){};
\path[<-,>=stealth']
(0) edge (1)
(0) edge (2)
(0) edge (3)
(1) edge (0)
(1) edge (2)
(2) edge (1)
(3) edge (1)
(3) edge (2)
;
\node at (0.5,-1) {\small G21};
\end{tikzpicture}
}
\def\GGC{
\begin{tikzpicture}[scale=\scale]
\clip(-0.8,-1.5) rectangle (1.8,1.2);
\node[bullet] (0) at (0,0){};
\node[bullet] (1) at (0,1){};
\node[bullet] (2) at (1,1){};
\node[bullet] (3) at (1,0){};
\path[->,>=stealth']
(0) edge (1)
(0) edge (2)
(0) edge (3)
(1) edge (0)
(1) edge (2)
(1) edge (3)
(2) edge (3)
(3) edge (2)
;
\node at (0.5,-1) {\small G33$\ast$};
\end{tikzpicture}
}
\def\GGD{
\begin{tikzpicture}[scale=\scale]
\clip(-0.8,-1.5) rectangle (1.8,1.2);
\node[bullet] (0) at (0,0){};
\node[bullet] (1) at (0,1){};
\node[bullet] (2) at (1,1){};
\node[bullet] (3) at (1,0){};
\path[->,>=stealth']
(0) edge (1)
(0) edge (2)
(0) edge (3)
(1) edge (0)
(1) edge (2)
(2) edge (0)
(2) edge (3)
(3) edge (0)
(3) edge (1)
;
\node at (0.5,-1) {\small G48$\ast$};
\end{tikzpicture}
}
\def\GGE{
\begin{tikzpicture}[scale=\scale]
\clip(-0.8,-1.5) rectangle (1.8,1.2);
\node[bullet] (0) at (0,0){};
\node[bullet] (1) at (0,1){};
\node[bullet] (2) at (1,1){};
\node[bullet] (3) at (1,0){};
\path[->,>=stealth']
(0) edge (1)
(0) edge (2)
(0) edge (3)
(1) edge (0)
(1) edge (2)
(2) edge (0)
(2) edge (3)
(3) edge (0)
;
\node at (0.5,-1) {\small G47$\ast$};
\end{tikzpicture}
}
\def\GGF{
\begin{tikzpicture}[scale=\scale]
\clip(-0.8,-1.5) rectangle (1.8,1.2);
\node[bullet] (0) at (0,0){};
\node[bullet] (1) at (1,0){};
\node[bullet] (2) at (0,1){};
\node[bullet] (3) at (1,1){};
\path[->,>=stealth']
(0) edge (1)
(0) edge (2)
(0) edge (3)
(1) edge (0)
(1) edge (2)
(2) edge (0)
(2) edge (3)
(3) edge (1)
(3) edge (2)
;
\node at (0.5,-1) {\small G46$\ast$};
\end{tikzpicture}
}
\def\GGG{
\begin{tikzpicture}[scale=\scale]
\clip(-0.8,-1.5) rectangle (1.8,1.2);
\node[bullet] (0) at (0,0){};
\node[bullet] (1) at (1,0){};
\node[bullet] (2) at (0,1){};
\node[bullet] (3) at (1,1){};
\path[->,>=stealth']
(0) edge (1)
(0) edge (2)
(0) edge (3)
(1) edge (0)
(1) edge (2)
(2) edge (0)
(2) edge (3)
(3) edge (1)
;
\node at (0.5,-1) {\small G22};
\end{tikzpicture}
}
\def\GGH{
\begin{tikzpicture}[scale=\scale]
\clip(-0.8,-1.5) rectangle (1.8,1.2);
\node[bullet] (0) at (0,0){};
\node[bullet] (1) at (0,1){};
\node[bullet] (2) at (1,0){};
\node[bullet] (3) at (1,1){};
\path[->,>=stealth']
(0) edge (1)
(0) edge (2)
(1) edge (0)
(1) edge (2)
(2) edge (0)
(2) edge (3)
(3) edge (0)
(3) edge (1)
;
\node at (0.5,-1) {\small G23};
\end{tikzpicture}
}
\def\GGI{
\begin{tikzpicture}[scale=\scale]
\clip(-0.8,-1.5) rectangle (1.8,1.2);
\node[bullet] (0) at (0,1){};
\node[bullet] (1) at (1,1){};
\node[bullet] (2) at (0,0){};
\node[bullet] (3) at (1,0){};
\path[->,>=stealth']
(0) edge (1)
(0) edge (2)
(0) edge (3)
(1) edge (0)
(1) edge (2)
(2) edge (0)
(2) edge (3)
(3) edge (2)
;
\node at (0.5,-1) {\small G40$\ast$};
\end{tikzpicture}
}
\def\GGJ{
\begin{tikzpicture}[scale=\scale]
\clip(-0.8,-1.5) rectangle (1.8,1.2);
\node[bullet] (0) at (0,0){};
\node[bullet] (1) at (1,0){};
\node[bullet] (2) at (0,1){};
\node[bullet] (3) at (1,1){};
\path[->,>=stealth']
(0) edge (1)
(0) edge (2)
(0) edge (3)
(1) edge (0)
(1) edge (2)
(2) edge (0)
(3) edge (1)
(3) edge (2)
;
\node at (0.5,-1) {\small G24$\ast$};
\end{tikzpicture}
}
\def\GHA{
\begin{tikzpicture}[scale=\scale]
\clip(-0.8,-1.5) rectangle (1.8,1.2);
\node[bullet] (0) at (0,1){};
\node[bullet] (1) at (0,0){};
\node[bullet] (2) at (1,0){};
\node[bullet] (3) at (1,1){};
\path[->,>=stealth']
(0) edge (1)
(0) edge (2)
(0) edge (3)
(1) edge (0)
(1) edge (2)
(2) edge (1)
(2) edge (3)
(3) edge (1)
;
\node at (0.5,-1) {\small G25$\ast$};
\end{tikzpicture}
}
\def\GHB{
\begin{tikzpicture}[scale=\scale]
\clip(-0.8,-1.5) rectangle (1.8,1.2);
\node[bullet] (0) at (0,0){};
\node[bullet] (1) at (1,0){};
\node[bullet] (2) at (0,1){};
\node[bullet] (3) at (1,1){};
\path[->,>=stealth']
(0) edge (1)
(0) edge (2)
(0) edge (3)
(1) edge (0)
(1) edge (2)
(2) edge (0)
(3) edge (1)
;
\node at (0.5,-1) {\small G16};
\end{tikzpicture}
}
\def\GHC{
\begin{tikzpicture}[scale=\scale]
\clip(-0.8,-1.5) rectangle (1.8,1.2);
\node[bullet] (0) at (0,0){};
\node[bullet] (1) at (0,1){};
\node[bullet] (2) at (1,0){};
\node[bullet] (3) at (1,1){};
\path[->,>=stealth']
(0) edge (1)
(0) edge (2)
(1) edge (0)
(1) edge (2)
(2) edge (0)
(2) edge (3)
(3) edge (0)
;
\node at (0.5,-1) {\small G17};
\end{tikzpicture}
}
\def\GHD{
\begin{tikzpicture}[scale=\scale]
\clip(-0.8,-1.5) rectangle (1.8,1.2);
\node[bullet] (0) at (0,0){};
\node[bullet] (1) at (0,1){};
\node[bullet] (2) at (1,1){};
\node[bullet] (3) at (1,0){};
\path[->,>=stealth']
(0) edge (1)
(0) edge (2)
(0) edge (3)
(1) edge (0)
(1) edge (2)
(2) edge (3)
(3) edge (0)
;
\node at (0.5,-1) {\small G18};
\end{tikzpicture}
}
\def\GHE{
\begin{tikzpicture}[scale=\scale]
\clip(-0.8,-1.5) rectangle (1.8,1.2);
\node[bullet] (0) at (0,0){};
\node[bullet] (1) at (1,0){};
\node[bullet] (2) at (0,1){};
\node[bullet] (3) at (1,1){};
\path[->,>=stealth']
(0) edge (1)
(0) edge (2)
(1) edge (0)
(1) edge (3)
(2) edge (0)
(3) edge (0)
(3) edge (2)
;
\node at (0.5,-1) {\small G19};
\end{tikzpicture}
}
\def\GHF{
\begin{tikzpicture}[scale=\scale]
\clip(-0.8,-1.5) rectangle (1.8,1.2);
\node[bullet] (0) at (0,0){};
\node[bullet] (1) at (0,1){};
\node[bullet] (2) at (1,1){};
\node[bullet] (3) at (1,0){};
\path[->,>=stealth']
(0) edge (1)
(0) edge (2)
(0) edge (3)
(1) edge (0)
(1) edge (2)
(2) edge (3)
(3) edge (1)
(3) edge (2)
;
\node at (0.5,-1) {\small G34$\ast$};
\end{tikzpicture}
}
\def\GHG{
\begin{tikzpicture}[scale=\scale]
\clip(-0.8,-1.5) rectangle (1.8,1.2);
\node[bullet] (0) at (0,0){};
\node[bullet] (1) at (0,1){};
\node[bullet] (2) at (1,1){};
\node[bullet] (3) at (1,0){};
\path[->,>=stealth']
(0) edge (1)
(0) edge (2)
(0) edge (3)
(1) edge (0)
(1) edge (2)
(2) edge (3)
(3) edge (1)
;
\node at (0.5,-1) {\small G10};
\end{tikzpicture}
}
\def\GHH{
\begin{tikzpicture}[scale=\scale]
\clip(-0.8,-1.5) rectangle (1.8,1.2);
\node[bullet] (0) at (0,1){};
\node[bullet] (1) at (0,0){};
\node[bullet] (2) at (1,0){};
\node[bullet] (3) at (1,1){};
\path[->,>=stealth']
(0) edge (1)
(0) edge (2)
(1) edge (0)
(2) edge (1)
(2) edge (3)
(3) edge (0)
(3) edge (1)
;
\node at (0.5,-1) {\small G11};
\end{tikzpicture}
}
\def\GHI{
\begin{tikzpicture}[scale=\scale]
\clip(-0.8,-1.5) rectangle (1.8,1.2);
\node[bullet] (0) at (0,0){};
\node[bullet] (1) at (0,1){};
\node[bullet] (2) at (1,1){};
\node[bullet] (3) at (1,0){};
\path[->,>=stealth']
(0) edge (1)
(0) edge (2)
(0) edge (3)
(1) edge (0)
(1) edge (2)
(2) edge (3)
(3) edge (2)
;
\node at (0.5,-1) {\small G29$\ast$};
\end{tikzpicture}
}
\def\GHJ{
\begin{tikzpicture}[scale=\scale]
\clip(-0.8,-1.5) rectangle (1.8,1.2);
\node[bullet] (0) at (0,0){};
\node[bullet] (2) at (0,1){};
\node[bullet] (3) at (1,0){};
\path[->,>=stealth']
(0) edge (2)
(0) edge (3)
(2) edge (0)
(3) edge (2)
;
\node at (0.5,-1) {\small G3$\ast$};
\end{tikzpicture}
}
\def\GIA{
\begin{tikzpicture}[scale=\scale]
\clip(-0.8,-1.5) rectangle (1.8,1.2);
\node[bullet] (0) at (0,0){};
\node[bullet] (1) at (0,1){};
\node[bullet] (2) at (1,1){};
\node[bullet] (3) at (1,0){};
\path[->,>=stealth']
(0) edge (1)
(0) edge (2)
(0) edge (3)
(1) edge (0)
(2) edge (1)
(2) edge (3)
(3) edge (1)
(3) edge (2)
;
\node at (0.5,-1) {\small G35};
\end{tikzpicture}
}
\def\GIB{
\begin{tikzpicture}[scale=\scale]
\clip(-0.8,-1.5) rectangle (1.8,1.2);
\node[bullet] (0) at (0,0){};
\node[bullet] (1) at (0,1){};
\node[bullet] (2) at (1,1){};
\node[bullet] (3) at (1,0){};
\path[->,>=stealth']
(0) edge (1)
(0) edge (2)
(0) edge (3)
(1) edge (0)
(2) edge (1)
(2) edge (3)
(3) edge (1)
;
\node at (0.5,-1) {\small G12$\ast$};
\end{tikzpicture}
}
\def\GIC{
\begin{tikzpicture}[scale=\scale]
\clip(-0.8,-1.5) rectangle (1.8,1.2);
\node[bullet] (0) at (0,0){};
\node[bullet] (1) at (0,1){};
\node[bullet] (2) at (1,1){};
\node[bullet] (3) at (1,0){};
\path[->,>=stealth']
(0) edge (1)
(0) edge (2)
(0) edge (3)
(1) edge (0)
(2) edge (1)
(2) edge (3)
(3) edge (2)
;
\node at (0.5,-1) {\small G30};
\end{tikzpicture}
}
\def\GID{
\begin{tikzpicture}[scale=\scale]
\clip(-0.8,-1.5) rectangle (1.8,1.2);
\node[bullet] (0) at (1,1){};
\node[bullet] (1) at (1,0){};
\node[bullet] (2) at (0,0){};
\node[bullet] (3) at (0,1){};
\path[->,>=stealth']
(0) edge (1)
(0) edge (2)
(1) edge (0)
(1) edge (2)
(2) edge (3)
(3) edge (0)
(3) edge (2)
;
\node at (0.5,-1) {\small G31};
\end{tikzpicture}
}
\def\GIE{
\begin{tikzpicture}[scale=\scale]
\clip(-0.8,-1.5) rectangle (1.8,1.2);
\node[bullet] (0) at (0,0){};
\node[bullet] (1) at (0,1){};
\node[bullet] (2) at (1,1){};
\node[bullet] (3) at (1,0){};
\path[->,>=stealth']
(0) edge (1)
(0) edge (2)
(0) edge (3)
(1) edge (0)
(2) edge (1)
(3) edge (2)
;
\node at (0.5,-1) {\small G5};
\end{tikzpicture}
}
\def\GIF{
\begin{tikzpicture}[scale=\scale]
\clip(-0.8,-1.5) rectangle (1.8,1.2);
\node[bullet] (0) at (0,1){};
\node[bullet] (1) at (0,0){};
\node[bullet] (2) at (1,1){};
\node[bullet] (3) at (1,0){};
\path[->,>=stealth']
(0) edge (1)
(0) edge (2)
(1) edge (0)
(2) edge (1)
(2) edge (3)
(3) edge (1)
;
\node at (0.5,-1) {\small G6};
\end{tikzpicture}
}
\def\GIG{
\begin{tikzpicture}[scale=\scale]
\clip(-0.8,-1.5) rectangle (1.8,1.2);
\node[bullet] (0) at (0,0){};
\node[bullet] (1) at (0,1){};
\node[bullet] (2) at (1,1){};
\node[bullet] (3) at (1,0){};
\path[->,>=stealth']
(0) edge (1)
(0) edge (2)
(0) edge (3)
(1) edge (0)
(2) edge (3)
(3) edge (2)
;
\node at (0.5,-1) {\small G26};
\end{tikzpicture}
}
\def\GIH{
\begin{tikzpicture}[scale=\scale]
\clip(-0.8,-1.5) rectangle (1.8,1.2);
\node[bullet] (0) at (1,1){};
\node[bullet] (1) at (1,0){};
\node[bullet] (2) at (0,0){};
\node[bullet] (3) at (0,1){};
\path[->,>=stealth']
(0) edge (1)
(0) edge (2)
(1) edge (0)
(1) edge (2)
(2) edge (3)
(3) edge (2)
;
\node at (0.5,-1) {\small G27};
\end{tikzpicture}
}
\def\GII{
\begin{tikzpicture}[scale=\scale]
\clip(-0.8,-1.5) rectangle (1.8,1.2);
\node[bullet] (0) at (0,0){};
\node[bullet] (1) at (0,1){};
\node[bullet] (2) at (1,1){};
\node[bullet] (3) at (1,0){};
\path[->,>=stealth']
(0) edge (1)
(0) edge (2)
(1) edge (0)
(1) edge (2)
(2) edge (3)
(3) edge (0)
(3) edge (1)
;
\node at (0.5,-1) {\small G13$\ast$};
\end{tikzpicture}
}
\def\GIJ{
\begin{tikzpicture}[scale=\scale]
\clip(-0.8,-1.5) rectangle (1.8,1.2);
\node[bullet] (0) at (0,0){};
\node[bullet] (1) at (0,1){};
\node[bullet] (2) at (1,1){};
\node[bullet] (3) at (1,0){};
\path[->,>=stealth']
(0) edge (1)
(0) edge (2)
(1) edge (0)
(1) edge (2)
(2) edge (3)
(3) edge (0)
;
\node at (0.5,-1) {\small G7};
\end{tikzpicture}
}
\def\GJA{
\begin{tikzpicture}[scale=\scale]
\clip(-0.8,-1.5) rectangle (1.8,1.2);
\node[bullet] (0) at (0,0){};
\node[bullet] (1) at (0,1){};
\node[bullet] (2) at (1,0){};
\node[bullet] (3) at (1,1){};
\path[->,>=stealth']
(0) edge (1)
(0) edge (2)
(1) edge (0)
(2) edge (3)
(3) edge (0)
(3) edge (1)
;
\node at (0.5,-1) {\small G8};
\end{tikzpicture}
}
\def\GJB{
\begin{tikzpicture}[scale=\scale]
\clip(-0.8,-1.5) rectangle (1.8,1.2);
\node[bullet] (0) at (0,0){};
\node[bullet] (1) at (0,1){};
\node[bullet] (2) at (1,1){};
\node[bullet] (3) at (1,0){};
\path[->,>=stealth']
(0) edge (1)
(0) edge (2)
(1) edge (0)
(1) edge (3)
(2) edge (1)
(2) edge (3)
(3) edge (0)
(3) edge (2)
;
\node at (0.5,-1) {\small G36$\ast$};
\end{tikzpicture}
}
\def\GJC{
\begin{tikzpicture}[scale=\scale]
\clip(-0.8,-1.5) rectangle (1.8,1.2);
\node[bullet] (0) at (0,0){};
\node[bullet] (1) at (0,1){};
\node[bullet] (2) at (1,1){};
\node[bullet] (3) at (1,0){};
\path[->,>=stealth']
(0) edge (1)
(0) edge (2)
(1) edge (0)
(1) edge (3)
(2) edge (1)
(2) edge (3)
(3) edge (0)
;
\node at (0.5,-1) {\small G14$\ast$};
\end{tikzpicture}
}
\def\GJD{
\begin{tikzpicture}[scale=\scale]
\clip(-0.8,-1.5) rectangle (1.8,1.2);
\node[bullet] (0) at (1,0){};
\node[bullet] (1) at (1,1){};
\node[bullet] (2) at (0,0){};
\node[bullet] (3) at (0,1){};
\path[->,>=stealth']
(0) edge (1)
(0) edge (2)
(1) edge (0)
(1) edge (3)
(2) edge (1)
(2) edge (3)
(3) edge (2)
;
\node at (0.5,-1) {\small G32$\ast$};
\end{tikzpicture}
}
\def\GJE{
\begin{tikzpicture}[scale=\scale]
\clip(-0.8,-1.5) rectangle (1.8,1.2);
\node[bullet] (0) at (0,0){};
\node[bullet] (1) at (0,1){};
\node[bullet] (2) at (1,1){};
\node[bullet] (3) at (1,0){};
\path[->,>=stealth']
(0) edge (1)
(0) edge (2)
(1) edge (0)
(1) edge (3)
(2) edge (1)
(3) edge (0)
;
\node at (0.5,-1) {\small G9$\ast$};
\end{tikzpicture}
}
\def\GJF{
\begin{tikzpicture}[scale=\scale]
\clip(-0.8,-1.5) rectangle (1.8,1.2);
\node[bullet] (0) at (1,1){};
\node[bullet] (1) at (1,0){};
\node[bullet] (2) at (0,0){};
\node[bullet] (3) at (0,1){};
\path[->,>=stealth']
(0) edge (1)
(0) edge (2)
(1) edge (0)
(2) edge (1)
(2) edge (3)
(3) edge (2)
;
\node at (0.5,-1) {\small G28$\ast$};
\end{tikzpicture}
}
\def\GJG{
\begin{tikzpicture}[scale=\scale]
\clip(-0.8,-1.5) rectangle (1.8,1.2);
\node[bullet] (0) at (0,0){};
\node[bullet] (1) at (0,1){};
\node[bullet] (2) at (1,1){};
\node[bullet] (3) at (1,0){};
\path[->,>=stealth']
(0) edge (1)
(0) edge (2)
(1) edge (0)
(2) edge (3)
(3) edge (0)
;
\node at (0.5,-1) {\small G4$\ast$};
\end{tikzpicture}
}
\def\GJH{
\begin{tikzpicture}[scale=\scale]
\clip(-0.8,-1.5) rectangle (1.8,1.2);
\node[bullet] (0) at (0,0){};
\node[bullet] (1) at (0,1){};
\node[bullet] (2) at (1,1){};
\node[bullet] (3) at (1,0){};
\path[->,>=stealth']
(0) edge (1)
(0) edge (2)
(1) edge (2)
(1) edge (3)
(2) edge (3)
(3) edge (0)
;
\node at (0.5,-1) {\small G2$\ast$};
\end{tikzpicture}
}
\def\GJI{
\begin{tikzpicture}[scale=\scale]
\clip(-0.8,-1.5) rectangle (1.8,1.2);
\node[bullet] (0) at (0,0){};
\node[bullet] (1) at (0,1){};
\node[bullet] (2) at (1,1){};
\node[bullet] (3) at (1,0){};
\path[->,>=stealth']
(0) edge (1)
(0) edge (2)
(1) edge (2)
(2) edge (3)
(3) edge (0)
;
\node at (0.5,-1) {\small G1$\ast$};
\end{tikzpicture}
}
\begin{figure}
\centering
\GAC
\hspace{\hdist}
\begin{tikzpicture}[scale=\scale]
\clip(-0.8,-1.5) rectangle (1.8,1.2);
\node[bullet] (0) at (0,0){};
\node at (0.5,-1) {\small $\P_0$};
\end{tikzpicture}
\hspace{\hdist}
\GAG
\hspace{\hdist}
\GAH
\hspace{\hdist}
\GAI
\hspace{\hdist}
\GAD
\hspace{\hdist}
\GAE
\hspace{\hdist}
\GAF
 
\vspace{\vdist}
\GAA
\hspace{\hdist}
\GAJ
\hspace{\hdist}
\GBA
\hspace{\hdist}
\GBJ
\hspace{\hdist}
\GCA
\hspace{\hdist}
\GCB
\hspace{\hdist}
\GCC
\hspace{\hdist}
\GBG

\vspace{\vdist}
\GBH
\hspace{\hdist}
\GBI
\hspace{\hdist}
\GBD
\hspace{\hdist}
\GBE
\hspace{\hdist}
\GCG
\hspace{\hdist}
\GAB
\hspace{\hdist}
\GBB
\hspace{\hdist}
\GBC
 
\vspace{\vdist}
\GBF
\hspace{\hdist}
\GCE
\hspace{\hdist}
\GCF
\hspace{\hdist}
\GCD

\caption{All core digraphs with at most 4 vertices with a Siggers polymorphism.}
\label{fig:4ElementPGraphs}
\end{figure}

\begin{figure}
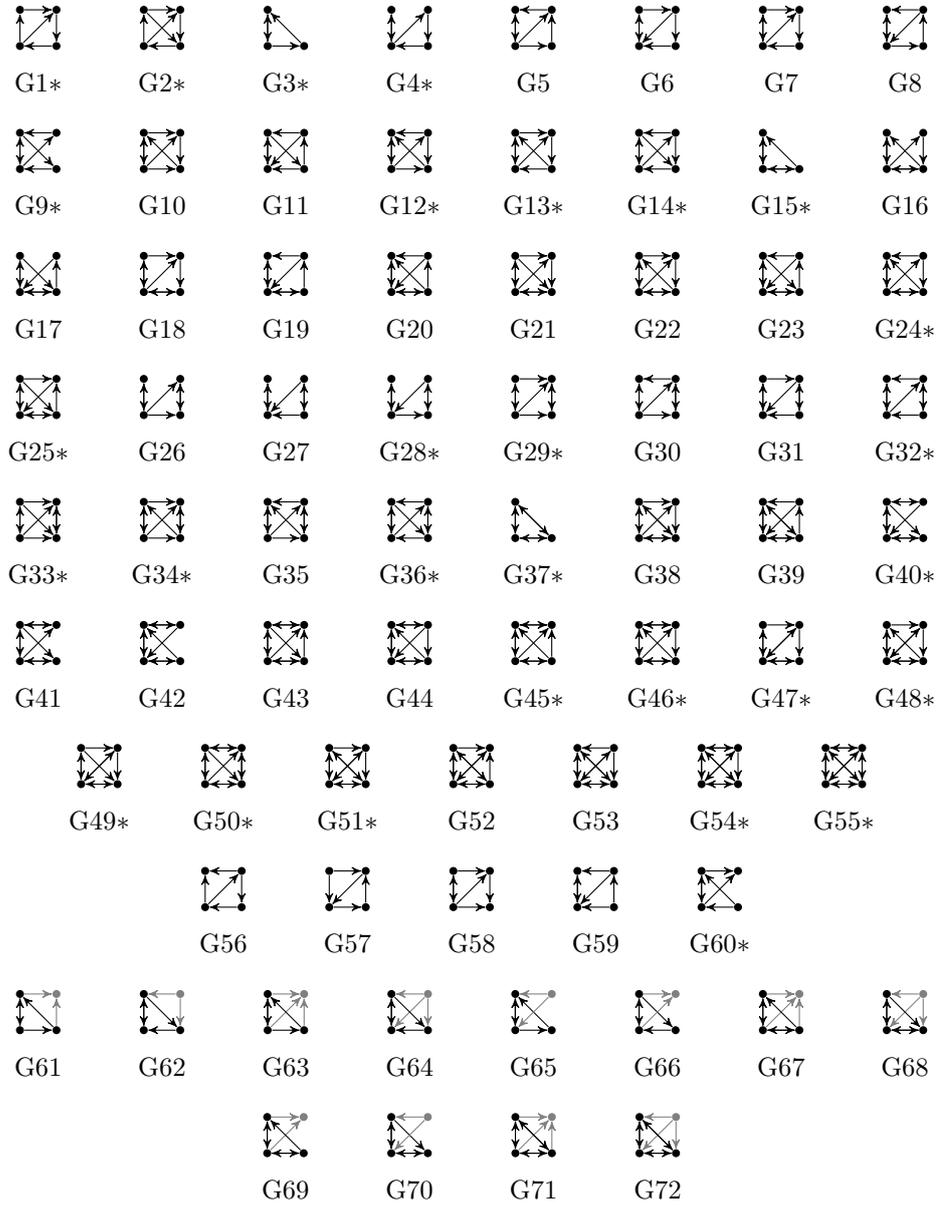

\centering
 
\vspace{\vdist}
\GJI
\hspace{\hdist}
\GJH
\hspace{\hdist}
\GHJ
\hspace{\hdist}
\GJG
\hspace{\hdist}
\GIE
\hspace{\hdist}
\GIF
\hspace{\hdist}
\GIJ
\hspace{\hdist}
\GJA

\vspace{\vdist}
\GJE
\hspace{\hdist}
\GHG
\hspace{\hdist}
\GHH
\hspace{\hdist}
\GIB
\hspace{\hdist}
\GII
\hspace{\hdist}
\GJC
\hspace{\hdist}
\GFH
\hspace{\hdist}
\GHB

\vspace{\vdist}
\GHC
\hspace{\hdist}
\GHD
\hspace{\hdist}
\GHE
\hspace{\hdist}
\GGA
\hspace{\hdist}
\GGB
\hspace{\hdist}
\GGG
\hspace{\hdist}
\GGH
\hspace{\hdist}
\GGJ

\vspace{\vdist}
\GHA
\hspace{\hdist}
\GIG
\hspace{\hdist}
\GIH
\hspace{\hdist}
\GJF
\hspace{\hdist}
\GHI
\hspace{\hdist}
\GIC
\hspace{\hdist}
\GID
\hspace{\hdist}
\GJD

\vspace{\vdist}
\GGC
\hspace{\hdist}
\GHF
\hspace{\hdist}
\GIA
\hspace{\hdist}
\GJB
\hspace{\hdist}
\GEI
\hspace{\hdist}
\GFA
\hspace{\hdist}
\GFB
\hspace{\hdist}
\GGI

\vspace{\vdist}
\GFI
\hspace{\hdist}
\GFJ
\hspace{\hdist}
\GFE
\hspace{\hdist}
\GFF
\hspace{\hdist}
\GFG
\hspace{\hdist}
\GGF
\hspace{\hdist}
\GGE
\hspace{\hdist}
\GGD
 
\vspace{\vdist}
\GFC
\hspace{\hdist}
\GFD
\hspace{\hdist}
\GEJ
\hspace{\hdist}
\GEG
\hspace{\hdist}
\GEH
\hspace{\hdist}
\GEF
\hspace{\hdist}
\GEE

\vspace{\vdist}
\GEC
\hspace{\hdist}
\GED
\hspace{\hdist}
\GDD
\hspace{\hdist}
\GDE
\hspace{\hdist}
\GDJ

\vspace{\vdist}
\GEA
\hspace{\hdist}
\GEB
\hspace{\hdist}
\GDF
\hspace{\hdist}
\GDG
\hspace{\hdist}
\GDI
\hspace{\hdist}
\GDH
\hspace{\hdist}
\GCJ
\hspace{\hdist}
\GDA

\vspace{\vdist}
\GDB
\hspace{\hdist}
\GDC
\hspace{\hdist}
\GCH
\hspace{\hdist}
\GCI

\caption{All core digraphs with at most 4 vertices without a Siggers polymorphism.}
\label{fig:4ElementNPGraphs}
\end{figure}

\begin{figure}
    \centering
    \begin{tikzpicture}
    \def\CTwoCenter{4}
    \def\CThreeCenter{-4}
    
    
    \node (C1) at (0,0) {$\C_1,\P_0$};
    \node (P1) at (0,-1) {$\P_1,\P_2,\P_3$};
    \node (T3) at (0,-2) {$\T_3,\T_3^+,\T_3^-$};
    \node (C2) at (\CTwoCenter,-2) {$\C_2,\C_4$};
    \node (C3) at (\CThreeCenter,-2) {$\C_3$};
    
    \node (C13) at (\CTwoCenter,-3) {$\C(1,3)$};
    \node (T4) at (0,-3) {$\T_4$};
    
    \node (T4M02) at (0,-4) {$\T_4^{\setminus(02)},\T_4^{\setminus(13)}$};
    \node (C2P1) at (\CTwoCenter,-5) {$\C_2^+,\C_2^-,\nCk112$};
    \node (C3P1) at (\CThreeCenter,-5) {$\C_3^+,\C_3^-$};
    
    \node (C2P2) at (\CTwoCenter-2,-7) {$\C_2^{++},\C_2^{--}$};

    \node (T4P32) at (\CTwoCenter,-8) {$\nCk022,\nCk202$};
    
    \node (T4M03) at (0,-5) {$\T_4^{\setminus(03)}\ppeq\Ord$};
    \node (T3PC2at1) at (\CTwoCenter-4,-6) {$\T_3\cup_0\C_2, \T_3\cup_1\C_2, \T_3\cup_2\C_2$};
    
    \node (K3) at (\CThreeCenter,-9) {$\K_3,\K_4,\dots$};
    
    \path
        (P1) edge (C1)
        (T3) edge (P1)
        (C2) edge (P1)
        (C3) edge (P1)
        (T4) edge (T3)
        (C13) edge (T3)
        (C13) edge (C2)
        (T4M02) edge (T4)
        (C2P1) edge (C13)
        (C2P1) edge (T4M02)
        (C3P1) edge (C3)
        (C3P1) edge (T4M02)
        
        (T3PC2at1) edge (C2P1)
        (C2P2) edge (C2P1)
        (T4P32) edge (C2P1)
        
        (T3PC2at1) edge (T4M03)
        (T4M03) edge (T4M02)
        
        (K3) edge (T4P32)
        (K3) edge (T3PC2at1)
        (K3) edge (C2P2)
        (K3) edge (C3P1)
        ;
        \path[->,dashed]
            (T4P32) edge (C2P2)
            ;
    \end{tikzpicture}
    \caption{Poset of digraphs with at most 4 vertices.}
    \label{fig:4ElementPoset}
\end{figure}
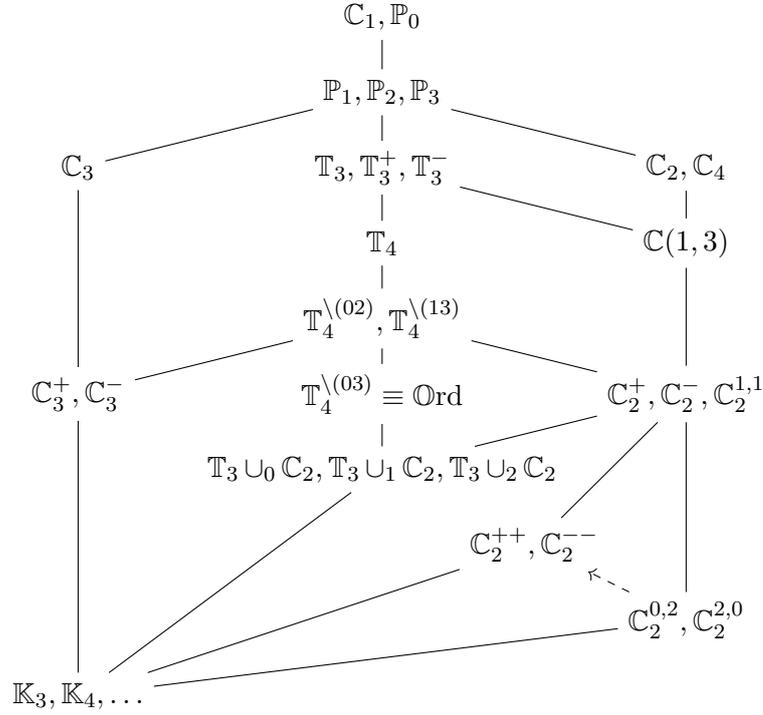

The poset of all digraphs with at most four vertices ordered by pp-constructions is presented in Figure~\ref{fig:4ElementPoset}. The dashed edge indicates an open problem. This problem is possibly the simplest one presented in this thesis. Note that, by Theorem~\ref{thm:freestructure}, it suffices to check whether
the free structure of $\Pol(\nCk022)$ generated by $\C_2^{++}$ 
has a homomorphism into $\C_2^{++}$. 
However, the vertices of this free structure are the homomorphisms from  $(\C_2^{++})^4$ to $\C_2^{++}$.
The number of these homomorphisms is bounded by the number of 4-ary operations on a four element set: \[4^{4^4}=2^{512}\approx 10^{154}\] Most of these operations are not homomorphisms, and there are faster ways to determine the homomorphisms than to iterate over all 4-ary operations. Nevertheless, the sheer size of this upper bound gives an idea of why the problem is difficult to solve with a computer.

\begin{question}\label{pro:4Elements}
    Can the digraph $\nCk022$  pp-construct $\C_2^{++}$?
\end{question}
Note that, by Lemma~\ref{lem:nCkNoNUs}, the digraph $\nCk022$ has no quasi near unanimity polymorphism of any arity. Hence, for no $n$ can $\NU n$ be used to disprove $\nCk022\ppleq\C_2^{++}$.
Next, we will look in detail at all pp-constructions and minor identities used to obtain this poset. 

\section{Collapses}

First we will look at the pp-constructions. Most of the pp-formulas we will need are presented in Figure~\ref{fig:4ElementGadgets}. 
Many of them were found by my computer program and were then simplified by hand. Most of them can be checked easily and it is left to the reader to do so.

\def\scale{0.5}

\def\GA{    
\begin{tikzpicture}[scale=\scale]
\node[var-f,label=left:1] at (0,0) {};
    \node[var-f,label=left:0] at (0,-1) {};
    \node at (0,-2) {$\Phi_{1}$};
\end{tikzpicture}
}
\def\GB{ 
\begin{tikzpicture}[scale=\scale]
    \node[var-f,label=left:$y$] (y) at (0,0) {};
    \node[var-f,label=left:$x$] (x) at (0,-1) {};
    \path[->,>=stealth'] (y) edge (x);
    \node at (0,-2) {$\Phi_{2}$};
\end{tikzpicture}
}
\def\GC{ 
\begin{tikzpicture}[scale=\scale]
    \node[var-f,label=left:$x$] (x) at (0,-1) {};
    \node[var-b] (a) at (1,-1) {};
    \node[var-b] (b) at (0,0) {};
    \node[var-f,label=right:$y$] (y) at (1,0) {};
    \path[->,>=stealth']
        (x) edge (b)
        (a) edge (b)
        (a) edge (y)
        ;
    \node at (0.5,-2) {$\Phi_{3}$};
\end{tikzpicture}
}

\def\GD{ 
\begin{tikzpicture}[scale=\scale]
    \node[var-f,label=left:$x$] (x) at (0,-1) {};
    \node[var-b] (a) at (0,1) {};
    \node[var-f,label=left:$y$] (y) at (0,0) {};
    \path[->,>=stealth']
        (x) edge (y)
        (y) edge (a)
        ;
    \node at (0,-2) {$\Phi_{4}$};
    
\end{tikzpicture}
}

\def\GE{ 
\begin{tikzpicture}[scale=\scale]
\node[var-f,label=left:$x$] (x) at (0,-1) {};
    \node[var-b] (a) at (0,0) {};
    \node[var-b] (b) at (0,1) {};
    \node[var-f,label=left:$y$] (y) at (0,2) {};
    \path[->,>=stealth']
        (x) edge (a)
        (a) edge (b)
        (b) edge (y)
        ;
    \node at (0,-2) {$\Phi_{5}$};

\end{tikzpicture}
}

\def\GF{ 
\begin{tikzpicture}[scale=\scale]

    \node[var-f,label=left:$x$] (x) at (0,0) {};
    \node[var-b] (a) at (0,1) {};
    \node[var-b] (b) at (0,2) {};
    \node[var-f,label=right:$y$] (y) at (1,1) {};
    \path[->,>=stealth']
        (x) edge (a)
        (a) edge (b)
        (y) edge (b)
        ;
    \node at (0.5,-1) {$\Phi_{6}$};

\end{tikzpicture}
}

\def\GI{ 
\begin{tikzpicture}[scale=\scale]

    \node[var-f,label=left:$x$] (x) at (0,0) {};
    \node[var-b] (a) at (0,1) {};
    \node[var-b] (b) at (1,0) {};
    \node[var-b] (c) at (1,1) {};
    \node[var-f,label=right:$y$] (y) at (1,2) {};
    
    \path[->,>=stealth']
        (x) edge (a)
        (b) edge (a)
        (b) edge (c)
        (c) edge (y)
        ;
    \node at (0.5,-1) {$\Phi_{7}$};

\end{tikzpicture}
}

\def\GG{ 
\begin{tikzpicture}[scale=\scale]
\node[var-f,label=below:$x_1$] (x1) at (0,0) {};
    \node[var-f,label=below:$x_2$] (x2) at (0+1,0) {};
    \node[var-b] (a) at (0+0.5,1) {};
    \node[var-f,label=above:$y_1$] (y1) at (0,2) {};
    \node[var-f,label=above:$y_2$] (y2) at (0+1,2) {};
    \path[<-,>=stealth']
        (x1) edge (x2)
        (x1) edge (a)
        (y2) edge (y1)
        (y2) edge (a)
        ;
    \node at (0.5,-2) {$\Phi_{8}$};
    
\end{tikzpicture}
}

\def\GH{ 
\begin{tikzpicture}[scale=\scale]

    \node[var-f,label=below:$x_1$] (x1) at (0,0) {};
    \node[var-f,label=below:$x_2$] (x2) at (0+1,0) {};
    \node[var-f,label=below:$x_3$] (x3) at (0+2,0) {};
    \node[var-f,label=below:$0$] (x4) at (0+3,0) {};
    \node[var-f,label=above:$1$] (y1) at (0,1) {};
    \node[var-f,label=above:$y_2$] (y2) at (0+1,1) {};
    \node[var-f,label=above:$y_3$] (y3) at (0+2,1) {};
    \node[var-f,label=above:$y_4$] (y4) at (0+3,1) {};
    \path[>=stealth']
        (x1) edge[->] (y2)
        (x2) edge[->] (y2)
        (x2) edge[->] (y4)
        (y4) edge[->] (x3)
        (y2) edge[dashed] (y3)
        ;
    \node at (1.5,-2) {$\Phi_{16}$};
\end{tikzpicture}
}

\def\GK{ 
\begin{tikzpicture}[scale=\scale]

    \node[var-f,label=below:$x_1$] (x1) at (0,0) {};
    \node[var-f,label=below:$x_2$] (x2) at (0+1,0) {};
    \node[var-f,label=below:$x_3$] (x3) at (0+2,0) {};
    \node[var-f,label=above:$0$] (y1) at (0,1) {};
    \node[var-f,label=above:$y_2$] (y2) at (0+1,1) {};
    \node[var-f,label=above:$y_3$] (y3) at (0+2,1) {};
    \path[>=stealth']
        (y2) edge[->] (x1)
        (x2) edge[dashed] (y3)
        (x3) edge[<->] (x2)
        ;
    \node at (1,-2) {$\Phi_{13}$};
\end{tikzpicture}
}

\def\GL{ 
\begin{tikzpicture}[scale=\scale]

    \node[var-f,label=below:$x_1$] (x1) at (0,0) {};
    \node[var-f,label=below:$2$] (x2) at (0+1,0) {};
    \node[var-f,label=above:$1$] (y1) at (0,1) {};
    \node[var-f,label=above:$y_2$] (y2) at (0+1,1) {};
    \path[>=stealth']
        (x1) edge[->] (y2)
        ;
    \node at (0.5,-2) {$\Phi_{9}$};
\end{tikzpicture}
}

\def\GQ{ 
\begin{tikzpicture}[scale=\scale]

    \node[var-f,label=below:$x_1$] (x1) at (0,0) {};
    \node[var-f,label=below:$x_2$] (x2) at (1,0) {};
    \node[var-f,label=above:$y_1$] (y1) at (0,1) {};
    \node[var-f,label=above:$y_2$] (y2) at (1,1) {};
    \path[>=stealth']
        (x1) edge[->] (y2)
        (x2) edge[dashed] (y1)
        ;
    \node at (0.5,-2) {$\Phi_{11}$};
\end{tikzpicture}
}

\def\GR{ 
\begin{tikzpicture}[scale=\scale]

    \node[var-f,label=below:$x_1$] (x1) at (0,0) {};
    \node[var-f,label=below:$x_2$] (x2) at (1,0) {};
    \node[var-f,label=below:$x_3$] (x3) at (2,0) {};
    \node[var-f,label=above:$y_1$] (y1) at (0,1) {};
    \node[var-f,label=above:$y_2$] (y2) at (1,1) {};
    \node[var-f,label=above:$y_3$] (y3) at (2,1) {};
    \path[>=stealth']
        (x2) edge[dashed] (y1)
        (x3) edge[dashed] (y2)
        (x1) edge[->] (y3)
        ;
    \node at (1,-2) {$\Phi_{12}$};
\end{tikzpicture}
}

\def\GM{ 
\begin{tikzpicture}[scale=\scale]

    \node[var-f,label=below:$x_1$] (x1) at (0,0) {};
    \node[var-f,label=below:$x_2$] (x2) at (1,0) {};
    \node[var-f,label=above:$1$] (y1) at (0,1) {};
    \node[var-f,label=above:$y_2$] (y2) at (1,1) {};
    \path[>=stealth']
        (x1) edge[->] (y2)
        (x2) edge[->] (y2)
        ;
    \node at (0.5,-2) {$\Phi_{10}$};
\end{tikzpicture}
}

\def\GO{ 
\begin{tikzpicture}[scale=\scale]

    \node[var-f,label=below:$x_1$] (x1) at (0,0) {};
    \node[var-f,label=below:$x_2$] (x2) at (1,0) {};
    \node[var-f,label=below:$1$] (x3) at (2,0) {};
    \node[var-f,label=above:$3$] (y1) at (0,1) {};
    \node[var-f,label=above:$y_2$] (y2) at (1,1) {};
    \node[var-f,label=above:$y_3$] (y3) at (2,1) {};
    \path[>=stealth']
        (y2) edge[->] (x1)
        (x2) edge[->] (y3)
        ;
    \node at (1,-2) {$\Phi_{14}$};
\end{tikzpicture}
}


\def\GP{ 
\begin{tikzpicture}[scale=\scale]

    \node[var-f,label=below:$1$] (x1) at (0,0) {};
    \node[var-f,label=below:$x_2$] (x2) at (0+1,0) {};
    \node[var-f,label=below:$x_3$] (x3) at (0+2,0) {};
    \node[var-f,label=above:$y_1$] (y1) at (0,1) {};
    \node[var-f,label=above:$y_2$] (y2) at (0+1,1) {};
    \node[var-f,label=above:$0$] (y3) at (0+2,1) {};
    \path[>=stealth']
        (y2) edge[->] (x3)
        (x3) edge[->] (y1)
        (y1) edge[->] (x2)
        ;
    \node at (1,-2) {$\Phi_{15}$};
\end{tikzpicture}
}

\begin{figure}
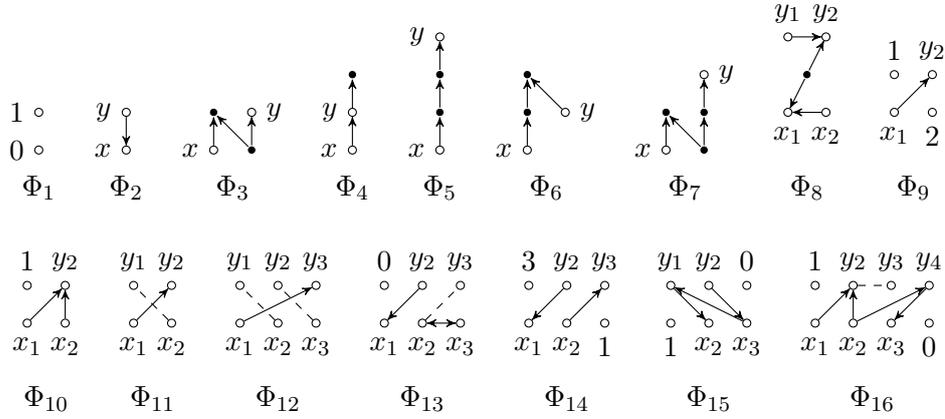

\centering

\GA
\GB
\GC
\GD
\GE
\GF
\GI
\GG
\GL

\vspace{4mm}
\GM
\GQ
\GR
\GK
\GO
\GP
\GH

    \caption{List of pp-formulas used in  pp-constructions in this chapter.}
    \label{fig:4ElementGadgets}
\end{figure}

\begin{itemize}
    \item $\P_0\ppleq \C_1$ using $x=y$.
    \item $\C_1\ppleq \P_0$ using $\bot$.
    \item $\C_2\ppleq\P_1$, $\C_3\ppleq\P_1$, and $\T_3\ppleq\P_1$ using $\Phi_1$.
    \item $\C_2^+\ppeq\C_2^-$, $\C_3^+\ppeq\C_3^-$, $\T_3^+\ppeq\T_3^-$,  $\T_3\cup_0\C_2\ppeq \T_3\cup_2\C_2$, $\T_4^{\setminus(01)}\ppeq\T_4^{\setminus(13)}$,  $\nCk022\ppeq\nCk202$,   and $\C_2^{++}\ppeq\C_2^{--}$ using $\Phi_2$.
    
    \item $\C(1,3)\ppleq\C_2$ and $\T_4^{\setminus(03)}\ppleq\Ord$ using $\Phi_3$.
    \item $\P_3\ppleq\P_2\ppleq\P_1\ppleq\P_0$,\ \   $\C_3^+\ppleq\C_3$,\ \  $\T_4\ppleq\T_3\ppleq\P_1$,\ \ $\T_3^+\ppleq\T_3$,\ \ $\nCk112\ppleq\C_2^-$, $\nCk022\ppleq\C_2^+$, and $\C_2^{++}\ppleq\C_2^+\ppleq\C_2$ using $\Phi_4$.
    \item $\T_3\ppleq\T_3^+$ by Lemma~\ref{lem:4ElementsT3leqT3n}.
    \item $\C_4\ppleq\C_2$ using $\Phi_5$($=x\stackrel{2}{\to}y$ see Chapter~\ref{cha:cycles}).
    \item $\C_2\ppleq\C_4$ and $\P_1\ppleq\P_2$ using $\Phi_{11}$ (compare to Example~\ref{ex.jakubconstruction}).
    \item $\P_1\ppleq\P_3$ using $\Phi_{12}$ (see Lemma~\ref{lem:IdempConstructPaths}).
    \item $\T_3\cup_0\C_2\ppleq\C_2^+$ and $\T_3\cup_1\C_2\ppleq\C_2^+$ using $\Phi_6$.
    \item $\T_3\cup_1\C_2\ppleq\Ord$ and $\T_3\cup_2\C_2\ppleq\Ord$ using $\Phi_7$. 
    
    \item $\T_3\cup_1\C_2\ppeq \T_3\cup_2\C_2$ by Lemma~\ref{lem:4ElementsT3u1C2EqualsT3u2C2}.
    
    \item $\C_2^+\ppleq\nCk112$ using $\Phi_8$ (see Lemma~\ref{lem:0C1iPPeq1C1i}).
    
    \item $\C(1,3)\ppleq\T_3$ using $\Phi_9$, where $(1,2)$ is the edge in between the two vertices in  $\C(1,3)$ that have in- and out-degree one (compare to proof of Theorem \ref{thm:TnIsHMblocker}).
    \item $\C_2^+\ppleq\C(1,3)$ using $\Phi_{13}$.
    \item $\C_2^+\ppleq \T_4^{\setminus(02)}$ and $\C_3^+\ppleq \T_4^{\setminus(02)}$ using $\Phi_{10}$ (see Lemma~\ref{lem:4ElementsCnPleqT4M02}).
    \item $\T_4^{\setminus(02)}\ppleq\T_4$ using $\Phi_{14}$ (compare to proof of Theorem \ref{thm:TnIsHMblocker}).
    \item $\Ord\ppleq \T_4^{\setminus(13)}$ using $\Phi_{16}$  (hint $0010\toEdge1000\toEdge1110\toEdge1111$).
    
    \item $\Ord\ppleq\T_4^{\setminus(03)}$ using $\Phi_{15}$ (hint $111\toEdge110\toEdge100\toEdge000$). 

    
    \item $\K_3$ can pp-construct all other finite digraphs, by Theorem~\ref{cor:K3IsBottom}.
    
    \item G1-G55 can pp-construct $\K_3$, by Theorem~\ref{thm:BartoKozikNiven} since they are all smooth cores that are not disjoint unions of cycles.
    
    \item G61-G72 can pp-construct $\K_3$ since they they have G3, G15, or G37 as a pp-definable subgraph (indicated by the gray vertices).
    
    \item G56-G60 can pp-construct $\K_3$ since they do not have a Siggers polymorphism (verified by computer).
    

\end{itemize}


Let $\C_n^+\coloneqq(\{a,0,\dots,n-1\},\{(v,v+1)\mid v\in \Z_n\}\cup\{v,a\mid v \in \Z_n\})$.

\begin{lemma}\label{lem:4ElementsCnPleqT4M02}
Let $n\geq2$. Then $\C_n^+\ppleq \T_4^{\setminus(02)}$ using $\Phi_{10}$.
\end{lemma}
\begin{proof}
Let $\H$ be the second pp-power of $\C_n^+$ defined by  $\Phi_{10}$. Note that $(1,a)$ is the only vertex in $\H$ which has in-degree greater one and out-degree zero. Furthermore, all edges that do not end in $(1,a)$ are of the form $((i,i),(1,i+1))$. The digraph $\H$ can be drawn like this: 
\[\begin{tikzpicture}
\node[] (3) at (0,3) {$(1,a)$};
\node[] (2) at (0,2) {$(1,2)$};
\node[] (1) at (0,1) {$(1,1)$};
\node[] (0) at (0,0) {$(0,0)$};

\node[] (2b) at (2,2) {$(2,3)$};
\node[] (2c) at (3,2) {$\dots$};

\node[] (1b) at (-2,1) {$(1,3)$};
\node[] (0b) at (-2,0) {$(2,2)$};
\node[] (1c) at (-3,1) {$\dots$};
\node[] (0c) at (-3,0) {$\dots$};

\path[>=stealth',->]
    (0) edge (1)
    (1) edge (2)
    (2) edge (3)
    (0) edge[bend right=40] (3)
    (1) edge[bend right=40] (3)
    
    (2b) edge (3)
    (1b) edge (3)
    (0b) edge (1b)
    (0b) edge (3)
    ;
\end{tikzpicture}\]
Hence, $\H$ is homomorphically equivalent to $\T_4^{\setminus(02)}$ and therefore $\C_n^+\ppleq \T_4^{\setminus(02)}$ as desired.
\end{proof}
For $n\in\N$ define the digraphs
\begin{align*}
    \T_3^{+n}&\coloneqq (\Z_{n+3},\{(v,v+1)\mid v\in \Z_{n+3}\}\cup\{(0,2)\})\text{ and}\\
    \T_3^{-+n}&\coloneqq (\Z_{n+3},\{(0,2),(0,3)\}\cup\{(v,v+1)\mid v\in \{0,1,3,4\dots,n,n+1\}\}).
\end{align*} 
The last pp-construction missing from the list is presented in the following lemma.
\begin{lemma}\label{lem:4ElementsT3leqT3n}
Let $n\in\N$. Then $\T_3\ppeq\T_3^{+n}\ppeq\T_3^{-+(n+2)}$.
\end{lemma}
\begin{proof}
Define the pp-formulas
\begin{center}
\begin{tikzpicture}[scale=0.6]
\node[var-b] (x1) at (7,0) {};
\node[var-f,label=below:{\small $x_1$}] (x2) at (1,0) {};
\node[var-f,label=below:{\small $x_2$}] (x3) at (2,0) {};
\node (x5) at (3.5,0) {$\dots$};
\node[var-f,label=below:{\small $x_{n-1}$}] (x6) at (5,0) {};
\node[var-f,label=below:{\small $x_n$}] (x7) at (6,0) {};

\node[var-f,label=above:{\small $y_1$}] (y2) at (1,1) {};
\node[var-f,label=above:{\small $y_2$}] (y3) at (2,1) {};
\node[var-f,label=above:{\small $y_3$}] (y4) at (3,1) {};
\node (y5) at (4.5,1) {$\dots$};
\node[var-f,label=above:{\small $y_n$}] (y7) at (6,1) {};
\path[>=stealth'] 
        (x1) edge[->] (y7)
        (x2) edge[dashed] (y3)
        (x3) edge[dashed] (y4)
        (x6) edge[dashed] (y7)
        (x7) edge[->] (y2)
        ;
    \node at (3.5,-2) {$\Phi_{1}$};
\end{tikzpicture}
\hspace{20mm}
\begin{tikzpicture}[scale=0.6]

    \node[var-f,label=below:$x_1$] (x1) at (0,0) {};
    \node[var-f,label=below:$x_2$] (x2) at (0+1,0) {};
    \node[var-f,label=below:$x_3$] (x3) at (0+2,0) {};
    \node[var-f,label=above:$y_1$] (y1) at (0,1) {};
    \node[var-f,label=above:$y_2$] (y2) at (0+1,1) {};
    \node[var-f,label=above:$y_3$] (y3) at (0+2,1) {};
    \path[>=stealth']
        (x1) edge[->] (y1)
        (x2) edge[dashed] (x1)
        (x3) edge[dashed] (y2)
        ;
    \node at (1,-2) {$\Phi_{2}$};
\end{tikzpicture}

\end{center}
and let $\H$ be the $n$-th pp-power of $\T_3$ given by $\Phi_1$ and $\H'$ be the third pp-power $\T_3^{+n}$ given by $\Phi_2$.
Note that there are the following induced substructures in $\H$ and $\H'$ respectively:
\begin{center}
\begin{tikzpicture}[scale=0.7]
\node[var-b,label=right:$111\dots10$] (0) at (0,0) {};
\node[var-b,label=right:$111\dots11$] (1) at (0,1) {};
\node[var-b,label=right:$211\dots11$] (2) at (0,2) {};
\node[rotate=90] at (0,3) {$\dots$};
\node[var-b,label=right:$222\dots21$] (4) at (0,4) {};
\node[var-b,label=right:$222\dots22$] (5) at (0,5) {};
\path[->,>=stealth']
(0) edge (1)
(1) edge (2)
(0) edge[bend left] (2)
(4) edge (5)
;
\end{tikzpicture}
\hspace{20mm}
\begin{tikzpicture}[scale=0.7]
\node[var-b,label=left:$001$] (0) at (0,0) {};
\node[var-b,label=left:$111$] (1) at (0,1) {};
\node[var-b,label=left:$211$] (2) at (0,2) {};

\node[var-b,label=right:$112$] (3) at (1,1) {};
\node[var-b,label=right:$223$] (4) at (1,2) {};
\node[rotate=90] at (1,3) {$\dots$};
\node[var-b,label={right:$(n+1)(n+1)(n+2)$}] (5) at (1,4) {};
\node[var-b,label={right:$(n+2)(n+2)(n+2)$}] (6) at (1,5) {};
\path[->,>=stealth']
(0) edge (1)
(1) edge (2)
(0) edge[bend right] (2)
(0) edge (3)
(3) edge (4)
(5) edge (6)
;
\end{tikzpicture}
\end{center}
Hence, $\T_3^{+n}$ embeds into $\H$ and $\T_3^{-+(n+2)}$ embeds into $\H'$.
The map 
\begin{align*}
    \H&\to\T_3^{+n}\\
    (v_1,\dots,v_n) &\mapsto
    \begin{cases}
    i+1&\text{if $v_i=2$ and $v_{i+1},\dots,v_n\in\{0,1\}$}\\
    v_n&\text{if $v_{1},\dots,v_n\in\{0,1\}$}
    \end{cases}
\end{align*}
is a homomorphism. Hence, $\T_3\ppleq\T_3^{+n}$.
The map 
\begin{align*}
\H'&\to\T_3^{-+(n+2)}\\
(u,v,w)&\mapsto     
\begin{cases}
u&\text{if }uvw\in\{001,111\}\text{ or }uv=21\\
u+2&\text{otherwise}
\end{cases}
\end{align*}
is a homomorphism. Hence, $\T_3^{+n}\ppleq\T_3^{-+(n+2)}$. Clearly $\T_3^{-+(n+2)}\ppleq\T_3$.
\end{proof}

\begin{lemma}\label{lem:4ElementsT3u1C2EqualsT3u2C2}
We have $\T_3\cup_1\C_2\ppeq\T_3\cup_2\C_2$.
\end{lemma}
\begin{proof}
Let 
\begin{align*}
    \T_3\cup_1\C_2&=(V,E_1)=(\{0,1,2,a\},\{(0,1),(0,2),(1,2),(1,a),(a,1)\}) \\ 
    \text{and }\T_3\cup_2\C_2&=(V,E_2)=(\{0,1,2,a\},\{(0,1),(0,2),(1,2),(2,a),(a,2)\}).
\end{align*}
First we show $\T_3\cup_1\C_2\ppleq\T_3\cup_2\C_2$.
Let $\H$ be the 2-nd pp-power of $\T_3\cup_1\C_2$ given by the formula
\begin{center}
\begin{tikzpicture}[scale=0.7]

    \node[var-f,label=below:$x_1$] (x1) at (0,0) {};
    \node[var-f,label=below:$x_2$] (x2) at (1,0) {};
    \node[var-f,label=above:$y_1$] (y1) at (0,1) {};
    \node[var-f,label=above:$y_2$] (y2) at (1,1) {};
    \path[>=stealth']
        (y1) edge[->] (x2)
        ;
    \tikzset{decoration={snake,amplitude=.25mm,segment length=1.3mm,post length=0.6mm,pre length=0.6mm}}
    \draw[decorate] (y1) -- (x1);
    \draw[decorate] (y2) -- (x1);
    \draw[decorate,->,>=stealth'] (y2) -- (x2);
\end{tikzpicture}
\end{center}
where \begin{tikzpicture}[scale=0.7]
    \node[var-f,label=above:$x$] (x1) at (0,0) {};
    \node[var-f,label=above:$y$] (x2) at (1,0) {};
\tikzset{decoration={snake,amplitude=.25mm,segment length=1.3mm,post length=0.6mm,pre length=0.6mm}}
    \draw[decorate,->,>=stealth'] (x1) -- (x2);
\end{tikzpicture} 
is an abbreviation for 
the pp-definable relation 
\[\{(x,y)\mid \T_3\cup_1\C_2\models
\begin{tikzpicture}[scale=0.6]
    \node[var-f,label=above:$x$] (0) at (0,0) {};
    \node[var-b] (1) at (1,0) {};
    \node[var-b] (2) at (2,0) {};
    \node[var-f,label=above:$y$] (3) at (3,0) {};
    \path[>=stealth']
        (0) edge[->] (1)
        (1) edge[<-] (2)
        (2) edge[->] (3)
        ;
\end{tikzpicture}
\}=E_1\cup\{(1,1),(a,2)\}.\]
Note that 
\begin{tikzpicture}[scale=0.7]
    \node[var-f,label=above:$x$] (x1) at (0,0) {};
    \node[var-f,label=above:$y$] (x2) at (1,0) {};
\tikzset{decoration={snake,amplitude=.25mm,segment length=1.3mm,post length=0.6mm,pre length=0.6mm}}
    \draw[decorate] (x1) -- (x2);
\end{tikzpicture}
represents the relation $\{1,a\}^2\setminus\{(a,a)\}$.
The subgraph of $\H$ consisting of the non-isolated vertices is the following:
\begin{center}
    \begin{tikzpicture}
     \node (0) at (0,0) {$(1,2)$};
     \node (1) at (0,1) {$(1,a)$};
     \node (2) at (0,2) {$(1,1)$};
     \node (a) at (1.5,2) {$(a,a)$};
     \node (a1) at (0,3) {$(a,1)$};
     \node (a2) at (-1.5,2) {$(a,2)$};
     
     \path[>=stealth']
        (0) edge[->] (1)
        (0) edge[->,bend left=40] (2)
        (1) edge[->] (2)
        (2) edge (a)
        (2) edge[->] (a1)
        (a2) edge[->] (2)
        ;
        
    \end{tikzpicture}
\end{center}
This graph is clearly homomorphically equivalent to $\T_3\cup_2\C_2$.


Now we show $\T_3\cup_2\C_2\ppleq\T_3\cup_1\C_2$.
Let $\H'$ be the 2-nd pp-power of $\T_3\cup_2\C_2$ given by the formula
\begin{center}
\begin{tikzpicture}[scale=0.7]
    \node[var-f,label=below:$x_1$] (x1) at (0,0) {};
    \node[var-f,label=below:$x_2$] (x2) at (1,0) {};
    \node[var-f,label=above:$y_1$] (y1) at (0,1) {};
    \node[var-f,label=above:$y_2$] (y2) at (1,1) {};
    \node[var-b,label=above:$0$] (0b) at (2,0.5) {};
    \node[var-b,label=above:$0$] (0a) at (-1,0.5) {};
    \tikzset{decoration={snake,amplitude=.25mm,segment length=1.3mm,post length=0.6mm,pre length=0.6mm}}
    \draw[decorate] (y1) -- (x1);
    \draw[decorate] (y2) -- (x2);
    \draw[decorate] (x2) -- (y1);
    \draw[decorate] (x2) -- (0b);
    \draw[decorate] (y2) -- (0b);
    \draw[decorate] (x1) -- (0a);
    \draw[decorate] (y1) -- (0a);
    \tikzset{decoration={snake,amplitude=.25mm,segment length=0.7mm,post length=0.6mm,pre length=0.6mm}}
    \draw[decorate,very thick] (x1) -- (y2);
\end{tikzpicture}
\end{center}
where
\begin{align*}
&\begin{tikzpicture}[scale=0.7]
\node[var-f,label=above:$x$] (0) at (0,0) {};
\node[var-f,label=above:$y$] (1) at (1,0) {};
\tikzset{decoration={snake,amplitude=.25mm,segment length=1.3mm,post length=0.6mm,pre length=0.6mm}}
    \draw[decorate,->,>=stealth'] (0) -- (1);
\end{tikzpicture}
\text{ abbreviates }
\begin{tikzpicture}[scale=0.6]
    \node[var-f,label=above:$x$] (0) at (0,0) {};
    \node[var-b] (1) at (1,0) {};
    \node[var-b] (2) at (2,0) {};
    \node[var-f,label=above:$y$] (3) at (3,0) {};
    \path[>=stealth']
        (0) edge[->] (1)
        (1) edge[<-] (2)
        (2) edge[<-] (3)
        ;
\end{tikzpicture}\text{, and}\\
&\begin{tikzpicture}[scale=0.7]
\node[var-f,label=above:$x$] (0) at (0,0) {};
\node[var-f,label=above:$y$] (1) at (1,0) {};
\tikzset{decoration={snake,amplitude=.25mm,segment length=0.7mm,post length=0.6mm,pre length=0.6mm}}
\draw[decorate,very thick,->,>=stealth'] (0) -- (1);
\end{tikzpicture}
\text{ abbreviates }
\begin{tikzpicture}[scale=0.6]
    \node[var-f,label=above:$x$] (0) at (0,0) {};
    \node[var-b] (1) at (1,0) {};
    \node[var-b] (2) at (2,0) {};
    \node[var-b] (3) at (3,0) {};
    \node[var-b] (4) at (4,0) {};
    \node[var-b] (5) at (5,0) {};
    \node[var-b] (6) at (6,0) {};
    \node[var-b] (7) at (7,0) {};
    \node[var-b] (8) at (8,0) {};
    \node[var-f,label=above:$y$] (9) at (9,0) {};
    \node[var-b] (a) at (7,1) {};
    \path[>=stealth']
        (0) edge[->] (1)
        (1) edge[<-] (2)
        (2) edge[<-] (3)
        (3) edge[<-] (4)
        (4) edge[->] (5)
        (5) edge[<-] (6)
        (6) edge[<-] (7)
        (7) edge[->] (8)
        (8) edge[<-] (9)
        (a) edge[->] (7)
        ;
\end{tikzpicture}
.
\end{align*}
Note that the formula for $\H'$ uses undirected edges everywhere.
Observe that on $\{0,2\}$ the thin edge defines the relation $\{0,2\}^2\setminus\{(2,2)\}$ and the thick edge
 defines the relation $\{0,2\}^2\setminus\{(0,0)\}$. Furthermore neither $(0,1)$ nor $(0,a)$ are in the relation defined by the thin edge. Hence, \begin{tikzpicture}[scale=0.7]
\node[var-b,label=above:$0$] (0) at (0,0) {};
\node[var-f,label=above:$x$] (1) at (1,0) {};
\tikzset{decoration={snake,amplitude=.25mm,segment length=1.3mm,post length=0.6mm,pre length=0.6mm}}
    \draw[decorate] (0) -- (1);
\end{tikzpicture} is only satisfied by 0 and 2. Therefore all edges in $\H'$ are between the four vertices $00,02,20$, and $22$. 

The map $h$ that sends $0,1,2,a$ to $20,00,02,22$ respectively is an embedding from $\T_3\cup_1\C_2$ to $\H'$. Since all vertices that are not in the image of $h$ are isolated $\H'$ is homomorphically equivalent to $\T_3\cup_1\C_2$.
\end{proof}

Note that we proved in particular that $\T_3\cup_2\C_2$ is pp-equivalent to the two element structure $(\{0,1\},\{0,1\}^2\setminus\{(0,0)\},\{0,1\}^2\setminus\{(1,1)\})$, whose polymorphism-clone is generated by majority (its the clone $[d_3]$ in \cite{PPPoset}).



\section{Separations}
To prove the separations we use the following minor conditions:
$\Const$, $\HM1$, $\HM2$, $\HM4$, $\HM5$, $\Sigma_2$, $\Sigma_3$ $\Majority$, 
and $\GFS3$. 
We will discuss each of these conditions in detail.
Before we do so, recall that 
\begin{itemize}
    \item if $\A\ppleq\B$ and $\B\models\Sigma$ then $\A\models\Sigma$ and
    \item if $\A\ppleq\B$ and $\A\not\models\Sigma$ then $\B\not\models\Sigma$.
\end{itemize}
This greatly reduces the number of statements to verify. For each condition $\Sigma$ it suffices to determine some digraphs $\G_1,\dots,\G_n$ that satisfy $\Sigma$ and some digraphs $\G'_1,\dots,\G'_m$ that do not satisfy $\Sigma$, such that for each digraph $\H$ we have that either $\H$ can be pp-constructed by some $\G_i$ or $\H$ can pp-construct some $\G'_i$. In the first case $\H\models\Sigma$ and in the latter case $\H\not\models\Sigma$. 

Though most of the following results have been found by my program, whenever a digraph $\H$ does not satisfy a condition $\Sigma$ my program could compute a fairly small subgraph of the indicator structure $\Ind(\H,\Sigma)$ preventing a homomorphism to $\H$. These subgraphs are presented here and it is left to the reader to verify that they are indeed subgraphs of the corresponding indicator structure and that they do not have a homomorphism to corresponding digraph. Recall that, since we can restrict to idempotent polymorphisms, we have that vertices of the indicator structure (introduced in Definition~\ref{def:indicatorStructure}) of the form $\{f(u,u,\dots,u),\dots\}$ must be mapped to $u$. To ease readability, if the condition involves only one function symbol $f$, then we replace $f(u_1,u_2,\dots,u_n)$ by $u_1u_2\dots u_n$.

\begin{enumerate}
    \item[$\Const$] 
    \begin{itemize}
        \item is satisfied by $\C_1$, witnessed by $0\mapsto0$.
        \item is not satisfied by $\P_1$, since $\{0,1\}$ is a loop in $\Ind(\P_2,\Const)$.
    \end{itemize}
    \item [$\HM1$]
    \begin{itemize}
        \item is satisfied by $\C_2$, witnessed by $(x,y,z)\mapsto x-y+z \pmod 2$.
        \item is satisfied by $\C_3$, witnessed by $(x,y,z)\mapsto x-y+z \pmod 3$.
        \item is not satisfied by $\T_3$, since $\{001,122\}$ is a loop in $\Ind(\T_3,\HM1)$. 
    \end{itemize}
    \item [$\HM2$]
    \begin{itemize}
        \item is satisfied by $\C(1,3)$, verified by computer.
        \item is satisfied by $\C_3$, since $\C_3\models\HM1$.
        \item is not satisfied by $\T_4$, 
        since it is not 2-braided (see Theorem~\ref{thm:TnIsHMblocker}), witnessed by the 2-braid
        \[\begin{tikzpicture}[scale=0.5]
        \node[var-b,label=below:0] (00) at (0,0) {};
        \node[var-b,label=above:1] (10) at (0,1) {};
        \node[var-b,label=below:1] (01) at (1,0) {};
        \node[var-b,label=above:2] (11) at (1,1) {};
        \node[var-b,label=below:2] (02) at (2,0) {};
        \node[var-b,label=above:3] (12) at (2,1) {};
        \path[>=stealth',->]
        (00) edge (01)
        (00) edge (11)
        (10) edge (11)
        (01) edge (02)
        (01) edge (12)
        (11) edge (12)
        ;
        \end{tikzpicture}\]
    \end{itemize}
    \item [$\HM4$]
    \begin{itemize}
        \item is satisfied by $\C_2^{++}$, verified by computer.
        \item is satisfied by $\C_3^+$, verified by computer.
        \item is not satisfied by $\nCk022$, since it is not 4-braided, witnessed by 
        \[\begin{tikzpicture}[scale=0.5]
        \node[var-b,label=below:2] (00) at (0,0) {};
        \node[var-b,label=above:$a$] (10) at (0,1) {};
        \node[var-b,label=below:1] (01) at (1,0) {};
        \node[var-b,label=above:$b$] (11) at (1,1) {};
        \node[var-b,label=below:$b$] (02) at (2,0) {};
        \node[var-b,label=above:$a$] (12) at (2,1) {};
        \node[var-b,label=below:$a$] (03) at (3,0) {};
        \node[var-b,label=above:1] (13) at (3,1) {};
        \node[var-b,label=below:$b$] (04) at (4,0) {};
        \node[var-b,label=above:2] (14) at (4,1) {};
        \path[>=stealth']
        (00) edge[<-] (01)
        (00) edge[<-] (11)
        (10) edge[<-] (11)
        (01) edge[<-] (02)
        (01) edge[<-] (12)
        (11) edge[<-] (12)
        (02) edge[->] (03)
        (02) edge[->] (13)
        (12) edge[->] (13)
        (03) edge[->] (04)
        (03) edge[->] (14)
        (13) edge[->] (14)
        ;
        \end{tikzpicture}\]

        \item is not satisfied by $\Ord$, since $\Ord\not\models\HM5$
    \end{itemize}
    \item [$\HM5$]
    \begin{itemize}
        \item is satisfied by $\C_2^{++}$, since $\C_2^{++}\models\HM4$.
        \item is satisfied by $\C_3^+$, since $\C_3^+\models\HM4$.
        \item is satisfied by $\nCk022$, verified by computer.
        \item is not satisfied by $\Ord$, since it is not 5-braided, witnessed by
        \[\begin{tikzpicture}[scale=0.5]
        \node[var-b,label=below:0] (00) at (0,0) {};
        \node[var-b,label=above:1] (10) at (0,1) {};
        \node[var-b,label=below:0] (01) at (1,0) {};
        \node[var-b,label=above:1] (11) at (1,1) {};
        \node[var-b,label=below:0] (02) at (2,0) {};
        \node[var-b,label=above:1] (12) at (2,1) {};
        \node[var-b,label=below:0] (03) at (3,0) {};
        \node[var-b,label=above:1] (13) at (3,1) {};
        \node[var-b,label=below:0] (04) at (4,0) {};
        \node[var-b,label=above:1] (14) at (4,1) {};
        \node[var-b,label=below:0] (05) at (5,0) {};
        \node[var-b,label=above:1] (15) at (5,1) {};
        \path[>=stealth',->]
        (00) edge (01)
        (00) edge (11)
        (10) edge (11)
        (01) edge (02)
        (01) edge (12)
        (11) edge (12)
        (02) edge (03)
        (02) edge (13)
        (12) edge (13)
        (03) edge (04)
        (03) edge (14)
        (13) edge (14)
        (04) edge (05)
        (04) edge (15)
        (14) edge (15)
        ;
        \end{tikzpicture}\]
    \end{itemize}
    \item[$\Sigma_2$] 
    \begin{itemize}
        \item is satisfied by $\C_3^+$, verified by computer.
        \item is satisfied by $\Ord$, witnessed by $\min^{(2)}$.
        \item is not satisfied by $\C_2$, $\{01,10\}$ is a loop in $\Ind(\C_2,\Sigma_2)$.
    \end{itemize}
    \item[$\Sigma_3$] 
    \begin{itemize}
        \item is satisfied by $\T_3\cup_1\C_2$, verified by computer.
        \item is satisfied by $\T_3\cup_0\C_2$, verified by computer.
        \item is satisfied by $\nCk022$, verified by computer.
        \item is satisfied by $\C_2^{++}$, verified by computer.
        \item is not satisfied by $\C_3$, $\{012,120\}$ is a loop in $\Ind(\C_3,\Sigma_3)$.
    \end{itemize}
    \item [$\Majority$]
\begin{itemize}

    \item satisfied by $\C_3^+$, witnessed by the polymorphism $f$ defined by: $f$ is a majority operation, $f(x,y,z)=x$ if $\{x,y,z\}=\{0,1,2\}$, and $f(x,y,z)=a$ if $|\{x,y,z\}|=3$ and $a \in \{x,y,z\}$.
    
    \item satisfied by $\T_3\cup_1\C_2$, verified by computer.
    \item satisfied by $\T_3\cup_0\C_2$, verified by computer.
    \item not satisfied by $\C_2^{++}$, let
    \begin{align*}
    &\C_2^{++}=
    \begin{tikzpicture}[scale=0.5,baseline=10]
    \node[var-b,label=below:$a$] (0) at (0,0) {};
    \node[var-b,label=below:$b$] (1) at (1,0) {};
    \node[var-b,label=right:$1$] (a) at (0.5,1) {};
    \node[var-b,label=right:$2$] (b) at (0.5,2) {};
    \path[>=stealth',->]
        (0) edge[<->] (1)
        (0) edge (a)
        (1) edge (a)
        (a) edge (b)
        ;
    \end{tikzpicture}
    &&\text{then}
    &\begin{tikzpicture}[scale=0.5, baseline=10]
    \node[var-b,label=left:{$\{aaa,aba,\dots\}$}] (0) at (0,0) {};
    \node[var-b,label=right:{$\{bbb,bba,1bb\dots\}$}] (1) at (1,0) {};
    \node[var-b,label=left:{$\{1ab\}$}] (a) at (0.5,1) {};
    \node[var-b,label=left:{$\{21a\}$}] (b) at (0.5,2) {};
    \path[>=stealth',->]
        (0) edge[<->] (1)
        (0) edge (a)
        (1) edge (a)
        (a) edge (b)
        (1) edge[bend right] (b)
        ;
    \end{tikzpicture}
    \end{align*}
    is a subgraph of $\Ind(\C_2^{++},\Majority)$.

    \item not satisfied by $\nCk022$, let
    \begin{align*}
    &\nCk022=
    \begin{tikzpicture}[scale=0.5,baseline=10]
    \node[var-b,label=below:$a$] (0) at (0,0) {};
    \node[var-b,label=below:$b$] (1) at (1,0) {};
    \node[var-b,label=right:$1$] (a) at (0.5,1) {};
    \node[var-b,label=right:$2$] (b) at (0.5,2) {};
    \path[>=stealth',->]
        (0) edge[<->] (1)
        (0) edge (a)
        (1) edge (a)
        (0) edge[bend left=5] (b)
        (1) edge[bend right=5] (b)
        (a) edge (b)
        ;
    \end{tikzpicture}
    &&\text{then}
    &\begin{tikzpicture}[scale=0.5,baseline=10]
    \node[var-b,label=left:{$\{aaa,aba\dots\}$}] (0) at (0,0) {};
    \node[var-b,label=right:{$\{bbb,bba\dots\}$}] (1) at (1,0) {};
    \node[var-b,label=right:{$\{1ab\}$}] (a) at (0.5,1) {};
    \node[var-b,label=right:{$\{111,211,\dots\}$}] (b) at (0.5,2) {};
    \path[>=stealth',->]
        (0) edge[<->] (1)
        (0) edge (a)
        (1) edge (a)
        (a) edge (b)
        (0) edge[bend left=5] (b)
        (1) edge[bend right=5] (b)
        ;
    \end{tikzpicture}
    \end{align*}
    is a subgraph of $\Ind(\nCk022,\Majority)$.
\end{itemize}
\item [$\GFS3$]
\begin{itemize}
    \item satisfied by $\T_4$, verified by computer.
    \item satisfied by $\C(1,3)$, verified by computer.
    \item not satisfied by $\C_3$, since $\{00012,00120,00201,\dots\}$  is a loop in $\Ind(\C_3,\GFS3)$.
    \item not satisfied by $\T_4^{\setminus(02)}$, since 
    \[
    \begin{tikzpicture}[scale=0.6]
    \node[var-b,label=right:{$\{10010,01111,\dots\}$}] (0) at (0,0) {};
    \node[var-b,label=right:{$\{33121,21111,\dots\}$}] (1) at (0,1) {};
    \node[var-b,label=right:{$\{22222,33222,\dots\}$}] (2) at (0,2) {};
    \path[>=stealth',->]
        (0) edge (1)
        (1) edge (2)
        (0) edge[bend left] (2)
        ;
    \end{tikzpicture}
    \]is a subgraph of $\Ind(\T_4^{\setminus(02)},\GFS3)$.
    Though the vertices of this subgraph look incomprehensible, we can read of the following $3$-staircase witnessing that $\T_4^{\setminus(02)}\not\models\GFS3$.
    \begin{center}
    \begin{tikzpicture}[scale = 0.7]
    
    \node[var-b,label=above:$1$] (1-2) at (-3,3.5) {};
    \node[var-b,label=above:$1$] (1-1) at (-2,3.5) {};
    \node[var-b] (00) at (-2,2.5) {};
    \node[var-b] (11) at (-1,3.5) {};
    \node[var-b,label=below:$1$] (10) at (-1,2.5) {};
    \node[var-b,label=below:$2$] (01) at (0,2.5) {};
    \node[var-b,label=below:$1$] (04) at (1,2.5) {};
    \node[var-b,label=below:$1$] (05) at (2,2.5) {};

    \draw[dashed] (1-2) -- (1-1);
    \draw[->,>=stealth'] (00) --  (10);
    \draw[->,>=stealth'] (00) --  (11);
    \draw[->,>=stealth'] (10) --  (01);
    \draw[->,>=stealth'] (1-1) --  (11);
    
    \draw[<-,>=stealth'] (01) --  (04);
    \draw[dashed] (04) --  (05);
    \node[var-b,label=above:$1$] (1-2) at (-3,1) {};
    \node[var-b,label=above:$1$] (1-1) at (-2,1) {};
    \node[var-b,label=above:$2$] (10) at (-1,1) {};
    \node[var-b] (00) at (-1,0) {};
    \node[var-b] (11) at (0,1) {};
    \node[var-b,label=below:$2$] (01) at (0,0) {};
    \node[var-b,label=below:$1$] (04) at (1,0) {};
    \node[var-b,label=below:$1$] (05) at (2,0) {};

    \draw[dashed] (1-2) -- (1-1);
    \draw[->,>=stealth'] (00) --  (01);
    \draw[->,>=stealth'] (00) --  (11);
    \draw[->,>=stealth'] (10) --  (11);
    \draw[->,>=stealth'] (1-1) --  (10);
    
    \draw[<-,>=stealth'] (01) --  (04);
    \draw[dashed] (04) --  (05);
    
    \node[var-b,label=above:$0$] (1-2) at (-3,-1.5) {};
    \node[var-b,label=above:$0$] (1-1) at (-2,-1.5) {};
    \node[var-b,label=above:$1$] (10) at (-1,-1.5) {};
    \node[var-b,label=above:$2$] (01) at (0,-1.5) {};
    \node[var-b] (00) at (0,-2.5) {};
    \node[var-b] (11) at (1,-1.5) {};
    \node[var-b,label=below:$0$] (04) at (1,-2.5) {};
    \node[var-b,label=below:$0$] (05) at (2,-2.5) {};

    \draw[dashed] (1-2) -- (1-1);
    \draw[<-,>=stealth'] (00) --  (04);
    \draw[<-,>=stealth'] (00) --  (11);
    \draw[->,>=stealth'] (10) --  (01);
    \draw[->,>=stealth'] (1-1) --  (10);
    
    \draw[<-,>=stealth'] (01) --  (11);
    \draw[dashed] (04) --  (05);
\end{tikzpicture}
\end{center}

 
\end{itemize}
\end{enumerate}

The reader may use Table~\ref{tab:4ElementTable} to verify that the only comparison, by $\ppleq$, of digraphs with at most four vertices not covered by our analysis is the one in Problem~\ref{pro:4Elements}.

\begin{sidewaystable}

    \centering
    \resizebox{\columnwidth}{!}{
    \begin{tabular}{l|c|c|c|c|c|c|c|c|c|c|c|c|c|c|c|c|c}
    \tiny
        $\ppleq$ & $\C_1$ & $\P_1$ & $\C_2$ & $\C_3$ & $\T_3$ & $\C(1,3)$ & $\T_4$ & $\T_4^{\setminus(02)}$ & $\C_3^+$ & $\C_2^+$ & $\C_2^{++}$ & $\nCk022$ & $\Ord$ & $\T_3\cup_1\C_2$ & $\T_3\cup_0\C_2$ & $\K_4$ \\\hline
        $\C_1$   &  R      & $\Const$ & $\Const$ & $\Const$ & $\Const$ & $\Const$ & $\Const$ & $\Const$ & $\Const$ & $\Const$ & $\Const$ & $\Const$ & $\Const$ & $\Const$ & $\Const$ & $\Const$  \\\hline
        $\P_1$   & $x=y$  & R & $\Sigma_2$ & $\Sigma_3$&$\HM1$ &$\Sigma_2$&$\HM2$&$\HM2$     &$\Sigma_3$&$\Sigma_2$&$\Sigma_2$&$\Sigma_2$&$\HM5$  &$\Sigma_2$&$\Sigma_2$&$\Sigma_2$\\\hline
        $\C_2$ & T & $\Phi_1$ & R & $\Sigma_3$ & $\HM1$       & $\HM1$ & $\HM2$ &$\HM2$      &$\Sigma_3$&$\HM1$    &$\HM1$    &$\HM1$    &$\HM1$&$\HM1$&$\HM1$  &$\HM1$  \\\hline
        $\C_3$ & T & $\Phi_1$ & $\Sigma_2$ & R &$\HM1$    & $\Sigma_2$ & $\HM2$ &$\HM2$      &$\HM1$    &$\Sigma_2$&$\Sigma_2$&$\Sigma_2$&$\HM5$  &$\Sigma_2$&$\Sigma_2$&$\Sigma_2$\\\hline
        $\T_3$ & T & $\Phi_1$ & $\Sigma_2$ & $\Sigma_3$ & R&$\Sigma_2$ & $\HM2$ &$\HM2$      &$\Sigma_3$&$\Sigma_2$&$\Sigma_2$&$\Sigma_2$&$\HM5$  &$\Sigma_2$&$\Sigma_2$&$\Sigma_2$\\\hline
        $\C(1,3)$ & T & T & $\Phi_3$ & $\Sigma_3$& $\Phi_9$&R          & $\HM2$ &$\HM2$      &$\Sigma_3$&$\HM2$    &$\HM2$    &$\HM2$    &$\HM5$  &$\HM2$  &$\HM2$  &$\HM2$  \\\hline
        $\T_4$ & T &T &$\Sigma_2$ & $\Sigma_3$&$\Phi_4$    &$\Sigma_2$ & R      &$\GFS3$     &$\Sigma_3$&$\Sigma_2$&$\Sigma_2$&$\Sigma_2$&$\HM5$  &$\Sigma_2$&$\Sigma_2$&$\Sigma_2$\\\hline
        $\T_4^{\setminus(02)}$ & T &T &$\Sigma_2$ & $\Sigma_3$&T&$\Sigma_2$ & $\Phi_{14}$ & R&$\Sigma_3$&$\Sigma_2$&$\Sigma_2$&$\Sigma_2$&$\HM5$  &$\Sigma_2$&$\Sigma_2$&$\Sigma_2$\\\hline
        $\C_3^+$ & T &T & $\Sigma_2$ & $\Phi_4$ &T         & $\Sigma_2$ & T     & $\Phi_{10}$& R        &$\Sigma_2$&$\Sigma_2$&$\Sigma_2$&$\HM5$  &$\Sigma_2$&$\Sigma_2$&$\Sigma_2$\\\hline
        $\C_2^+$ & T &T & $\Phi_4$ & $\Sigma_3$ & T        & $\Phi_{13}$& T     & $\Phi_{10}$&$\Sigma_3$&R        &$\Majority$&$\HM4$    &$\HM5$  &$\HM4$  &$\HM4$  &$\HM4$  \\\hline
        $\C_2^{++}$ & T &T & T & $\Sigma_3$ & T            &T           & T     & T          &$\Sigma_3$&$\Phi_4$  &R         &$\HM4$    &$\HM5$  &$\HM5$  &$\HM5$  &$\HM5$  \\\hline
        $\nCk022$ & T &T &T & $\Sigma_3$&T  &T           & T     & T          &$\Sigma_3$&$\Phi_4$  &$\boldsymbol ?$         & R        &$\HM5$  &$\HM5$  &$\HM5$  &$\HM5$  \\\hline
        $\Ord$ & T &T &$\Sigma_2$ & $\Sigma_3$&T           & $\Sigma_2$ & T     & $\Phi_{16}$&$\Sigma_3$&$\Sigma_2$&$\Sigma_2$&$\Sigma_2$&R       &$\Sigma_2$&$\Sigma_2$&$\Sigma_2$\\\hline
        $\T_3\cup_1\C_2$ & T &T &T&$\Sigma_3$&T&            T           & T     & T          &$\Sigma_3$&$\Phi_6$&$\Majority$  &$\Majority$&$\Phi_7$&R& L.~\ref{lem:4ElementsT3u1C2EqualsT3u2C2} &$\Sigma_3$\\\hline
        $\T_3\cup_0\C_2$ & T &T &T&$\Sigma_3$&T&            T           & T     & T          &$\Sigma_3$&$\Phi_6$&$\Majority$  &$\Majority$&$\Phi_7$& L.~\ref{lem:4ElementsT3u1C2EqualsT3u2C2} &R&$\Sigma_3$\\\hline
        $\K_4$ & T& T& T& T& T& T& T& T& L.~\ref{lem:KnppdefinesAllGraphs}& T&  L.~\ref{lem:KnppdefinesAllGraphs} &L.~\ref{lem:KnppdefinesAllGraphs} & T & L.~\ref{lem:KnppdefinesAllGraphs} & L.~\ref{lem:KnppdefinesAllGraphs} & R\\\hline
    \end{tabular}
    }
    \caption{
    The order $\ppleq$ on digraphs with at most four vertices.  
    The entry in the row of $\H$ and the column of $\G$ is either a minor condition $\Sigma$, a pp-formula $\Phi$, a lemma, an R, or a T. In the first case $\G\not\ppleq \H$, $\G\models\Sigma$ and $\H\not\models\Sigma$. In the remaining four cases $\G\ppleq\H$ holds either by using $\Phi$, as proved in the lemma, by reflexivity, or by transitivity, respectively.  
    }
    \label{tab:4ElementTable}

\end{sidewaystable}

At the end of this chapter I would like to mention some more observations about digraphs of order four.

\begin{itemize}
    \item Whenever one of the minor conditions I tested for this chapter was satisfied by some digraph $\G$ with at most four vertices it was was already satisfied by conservative polymorphisms of $\G$.
    \item For the CSPs of some of the digraphs in this chapter we the obstruction set (see Example~\ref{exa:0C22Obstructions}) has a simple description. 
    \begin{itemize}
        \item A digraph is not in $\csp(\T_n)$  if and only if $\P_n$ has a homomorphism into it.
        \item A digraph is not in $\csp(\C_2)$ if and only if there is an orientation of a cycle of odd length that has a homomorphism into it.
        \item We determined that a digraph is not in $\csp(\C_2^+)$ if and only if there is no digraph $\F$ that has a homomorphism into it, where $\F$ is some orientation of a cycle of odd length such that each vertex with in-degree two has an extra out edge.
        \[\begin{tikzpicture}
        
        \node[var-b,scale=0.5]  at (-30:1) {};
        \node[var-b,scale=0.5]  at (-20:1) {};
        \node[var-b,scale=0.5]  at (-10:1) {};
        \node[var-b,scale=0.5]  at (50:1) {};
        \node[var-b,scale=0.5]  at (60:1) {};
        \node[var-b,scale=0.5]  at (70:1) {};
        \node[var-b] (0) at (0:1) {};
        \node[var-b] (1) at (20:1) {};
        \node[var-b] (2) at (40:1) {};
        \node[var-b] (3) at (20:{1.3}) {};
        \path[->,>=stealth']
        (0) edge (1)
        (2) edge (1)
        (1) edge (3)
        ;
        \end{tikzpicture}\]
        \item We determined that a digraph is not in $\csp(\nCk022)$ if and only if there is no digraph $\F$ that has a homomorphism into it, where $\F$ is some orientation of a cycle of odd length such that each vertex $u$ with in-degree two has an extra outgoing path $u\to u_1\to u_2$.
        \[\begin{tikzpicture}
        
        \node[var-b,scale=0.5]  at (-30:1) {};
        \node[var-b,scale=0.5]  at (-20:1) {};
        \node[var-b,scale=0.5]  at (-10:1) {};
        \node[var-b,scale=0.5]  at (50:1) {};
        \node[var-b,scale=0.5]  at (60:1) {};
        \node[var-b,scale=0.5]  at (70:1) {};
        \node[var-b] (0) at (0:1) {};
        \node[var-b] (1) at (20:1) {};
        \node[var-b] (2) at (40:1) {};
        \node[var-b] (3) at (20:{1.3}) {};
        \node[var-b] (4) at (20:{1.6}) {};
        \path[->,>=stealth']
        (0) edge (1)
        (2) edge (1)
        (1) edge (3)
        (3) edge (4)
        ;
        \end{tikzpicture}\]
        \item Sebastian Meyer determined (personal communication) that a digraph is not in $\csp(\C_2^{++})$ if and only if there is no digraph $\F$ that has a homomorphism into it, where $\F$ is some orientation of a cycle of odd length such that each vertex $u$ with in-degree two has an extra outgoing path of net-length two.
        
        \[\begin{tikzpicture}
        
        \node[var-b,scale=0.5]  at (-30:1) {};
        \node[var-b,scale=0.5]  at (-20:1) {};
        \node[var-b,scale=0.5]  at (-10:1) {};
        \node[var-b,scale=0.5]  at (50:1) {};
        \node[var-b,scale=0.5]  at (60:1) {};
        \node[var-b,scale=0.5]  at (70:1) {};
        \node[var-b] (0) at (0:1) {};
        \node[var-b] (1) at (20:1) {};
        \node[var-b] (2) at (40:1) {};
        \node[var-b] (2b1) at (20:{1.3}) {};
        \node[var-b] (2a1) at ($(20:{1})+(-0.09,0.4)$) {};
        \node[var-b] (2b2) at ($(20:{1.3})+(-0.09,0.4)$) {};
        \node[var-b] (2a2) at ($(20:{1})+2*(-0.09,0.4)$) {};
        \node[var-b] (2b3) at ($(20:{1.3})+2*(-0.09,0.4)$) {};
        \node[var-b] (2c3) at ($(20:{1.6})+2*(-0.09,0.4)$) {};
        \path[->,>=stealth']
        (0) edge (1)
        (2) edge (1)
        (1) edge (2b1)
        (2a1) edge (2b1)
        (2a1) edge (2b2)
        (2a2) edge (2b2)
        (2a2) edge (2b3)
        (2b3) edge (2c3)

        ;
        \end{tikzpicture}\]
    \end{itemize}
\end{itemize}
\chapter{The Smallest Hard Trees}\label{cha:hardTrees}

In the previous chapter we have seen that there are orientations of paths whoose CSPs is NL-hard. What if we allow trees instead of paths? Gutjahr, Welzl, and Woeginger showed that there are orientations of trees whoose CSP is NP-hard~\cite{GutjahrWW92}.
In the final chapter of this thesis we go on a hunt for the smallest orientation of a tree whose CSP is NP-hard. This is joint work with Manuel Bodirsky, Jakub Bul\'in, and Michael Wernthaler and has been published in \cite{BodirskyBulinStarkeWernthalerSmallestHardTrees}. The Rust code has been improved since the publication. Hence, the claims about the implementation in this chapter differ slightly from the ones made in the original publication. 


Unfortunately, there is no graph theoretic characterization of which orientations of finite trees 
have NP-hard CSPs (assuming P${}\neq{}$NP). The first NP-hard tree  was found by Gutjahr, Welzl, and Woeginger and had 287 vertices~\cite{GutjahrWW92}. This was later improved by Gutjahr to a smaller NP-hard tree with 81 vertices~\cite{Gutjahr}, and then to an NP-hard tree with just 45 vertices by Hell, Ne\v{s}et\v{r}il, and Zhu~\cite{HellNZ96}. The tree $\T$ constructed there is even a \emph{triad}, i.e., a tree with exactly one vertex of degree three and all other vertices of degree one or two.
An NP-hard triad with 39 vertices was found by Barto, Kozik, Mar\'oti, and Niven~\cite{SpecialTriads,SpecialTriadsErratum} using an in-depth analysis of the polymorphisms of triads; they conjectured that their triad is the smallest NP-hard tree (assuming $\operatorname{P}\neq\operatorname{NP}$).
This approach lead to a study of certain classes of trees~\cite{BartoB13,Bulin18}.
Fischer~\cite{Jana} used a computer search and found an NP-hard tree with just 30 vertices (refuting the conjecture of Barto et al.\ mentioned above). Later, independently, Tatarko constructed a 26-vertex NP-hard triad, by manual analysis of polymorphisms. 
In 2022 Bodirsky, Bul\'in, Wernthaler, and I used a computer program to search all orientations of trees with up to 20 vertices and found the smallest orientation of a tree with NP-hard CSP.
See Table~\ref{tab:history}. One of the goals of this chapter is to explain how we computed this tree. 
The basic idea of our approach is very simple: enumerate all trees and check for the existence of a Siggers polymorphism using the indicator graph (see Definition~\ref{def:indicatorStructure}). Hence, we were also able to study other properties of small trees related to minor conditions: We compute the smallest tree that is NL-hard (assuming $\operatorname{L}\neq\operatorname{NL}$), the smallest tree that cannot be solved by arc consistency, and the smallest tree that cannot be solved by Datalog. 
Our experimental results also support a conjecture of Bul\'in concerning a question of Hell, Ne\v{s}et\v{r}il and Zhu, namely that `easy trees lack the ability to count'.

\renewcommand{\arraystretch}{1.1}
\begin{table}
\begin{center}
\begin{tabular}{llll}
author & year & size & comment \\
\hline
Gutjahr, Welzl, and Woeginger~\cite{GutjahrWW92} & 1992 & 287 & First published\\
Gutjahr~\cite{Gutjahr} & 1991 & 81 & PhD thesis \\
Hell, Ne\v{s}et\v{r}il, and Zhu~\cite{HellNZ96} & 1996 & 45 & Triad \\
Barto, Kozik, Mar\'oti, and Niven~\cite{SpecialTriads} & 2009 & 39 & Triad \\
Fischer~\cite{Jana} & 2015 & 30 & Master thesis \\
Tatarko~\cite{Tatarko} & 2019 & 26 & Triad\\
Bodirsky, Bul\'in, Starke, Wernthaler~\cite{BodirskyBulinStarkeWernthalerSmallestHardTrees} & 2022 & 22 & Smallest triad \\
Bodirsky, Bul\'in, Starke, Wernthaler~\cite{BodirskyBulinStarkeWernthalerSmallestHardTrees} & 2022 & 20 & Smallest tree
\end{tabular}
\end{center}
\caption{A time-line of the history of the smallest known NP-hard orientation of a tree.} 
\label{tab:history}
\end{table}






Several important conjectures about the computational complexity of CSPs   remain open: most notably the question for which finite structures $\A$ the problem $\csp(\A)$ is in NL (non-deterministic log-space), and for which finite structures $\A$ it is in L (deterministic log-space). 
As in the case of P versus NP-hard, it appears that these questions are closely linked to central dividing lines in universal algebra, as illustrated by the following conjectures.

\begin{conjecture}[Larose and Tesson~\cite{LaroseTesson}]\label{conj:NL} If the polymorphisms of a finite structure $\Hb$ with finite relational signature contain a Kearnes-Kiss chain (defined in Definition~\ref{def:kk}), then $\csp(\Hb)$ is in $\NL$. \end{conjecture}

It is known that if $\Hb$ does not satisfy the condition from Conjecture~\ref{conj:NL}, then $\Hb$ is hard for complexity classes that are not believed to be in NL (more details can be found in Theorem~\ref{thm:KKiffNoHorn3SATorNo3Linp}).
Conjecture~\ref{conj:NL}
is wide open and we believe it to be one of the most difficult research problems in the theory of finite-domain constraint satisfaction that remains open.

\begin{conjecture}
[Egri, Larose, and Tesson~\cite{EgriLaroseTessonLogspace}]\label{conj:L}
If the polymorphisms of a finite structure $\Hb$ with finite relational signature contain a Noname chain  
 (defined in  Definition~\ref{def:nn})
 then $\csp(\Hb)$ is in $\Lclass$.
 \end{conjecture} 
 
Also here it is known that if $\Hb$ does not satisfy the condition from Conjecture~\ref{conj:L}, then $\Hb$ is hard for complexity classes that are not believed to be in L. 
 
We mention that both conjectures can equivalently be phrased by the inability to  pp-construct certain finite structures that are known to be L-hard, Mod$_p$L-hard, or NL-hard (see Section~\ref{sect:specific-polymorphisms}).
Kazda proved a conditional result that states that resolving the first conjecture would also provide a solution to the second~\cite{Kazda-n-permute} (see Theorem~\ref{thm:symLinDLiffHM}). 

Again, these conjectures are already open if $\Hb$ is a finite digraph, or even if $\Hb$ is an orientation of a finite tree. It is also known that answering the question of containment in NL for finite digraphs would also answer the question for general finite structures~\cite{BulinDelicJacksonNiven}. For orientations of finite trees, however, the question might be easier to resolve.
For brevity, an orientation of a finite tree is simply called a \emph{tree} in this chapter and we adopt the following terminology: a digraph $\Hb$ is \emph{NP-hard} if $\csp(\Hb)$ is NP-hard, and \emph{tractable} if $\csp(\Hb)$ is in P. Similarly, we say that $\Hb$ is
\emph{P-hard}, \emph{NL-hard}, \emph{in NL}, \emph{in L}, \emph{NP-complete},  \emph{NL-complete}, etc.\ if $\csp(\Hb)$ has that property.

\section{Descriptive Complexity} 
Besides the computational complexity of CSPs, the \emph{descriptive complexity} of CSPs has been studied intensively, and leads to a fruitful interplay of finite model theory, graph theory, and universal algebra. Since the results obtained in this context are highly relevant for the open conjectures mentioned above, 
we provide a brief introduction to the  most prominent concepts. A digraph $\Hb$ has \emph{tree duality} if for all finite digraphs $\G$, if whenever all trees that map homomorphically to $\G$ also map to $\Hb$, then $\G$ maps to $\Hb$. It is well known that a finite digraph $\Hb$ has tree duality if and only if the so-called \emph{arc-consistency procedure} solves $\csp(\Hb)$~\cite{FederVardi}. This procedure is of central importance to our work, for many independent reasons that we mention later, and will be introduced in detail in Section~\ref{sect:AC}.

For every finite digraph $\Hb$, the arc-consistency procedure for $\csp(\Hb)$ can be formulated as a Datalog program~\cite{FederVardi}. 
Feder and Vardi proved that $\cocsp(\Hb)$ is in Datalog if and only if $\Hb$ has so-called \emph{bounded treewidth duality}; the definition of this concept is similar to the concept of tree duality but we omit it since it is not needed in this chapter. Bounded treewidth duality can be strengthened to \emph{bounded pathwidth duality}, which corresponds precisely to recognizability by a linear Datalog program~\cite{Dalmau}. 

A structure $\A$ has the \emph{ability to count} \cite{AbilityToCount} if $\A$ can encode solving systems of linear equations over $\mathbb Z_p$ (for some prime $p$), i.e., $\A$ can pp-construct $\TLinP$; thus making $\csp(\A)$ Mod$_p$L-hard. The ability to count adds substantial complexity to the CSP. Structures that \emph{cannot count} (i.e., that lack the ability to count) are all tractable and even in Datalog (see Theorem~\ref{thm:DLiff34WNUiff3Linp})
and this result, known as the \emph{bounded width theorem}, was an important intermediate step towards the resolution of the Feder-Vardi CSP dichotomy conjecture. Based on this theorem, the lack of the ability to count has a number of equivalent characterizations: \emph{bounded width}, bounded treewidth duality,  definability in Datalog, solvability by \emph{Singleton Arc Consistency}~\cite{Kozik-SLAC}. 

Several important classes of structures exhibit a dichotomy between NP-hardness and the lack of the ability to count (assuming $\operatorname{P}\neq\operatorname{NP}$), which we will refer to as ``easy structures cannot count''. 
Examples include undirected graphs \cite{HellNesetril}, \emph{smooth digraphs} (digraphs without sources and sinks)~\cite{BartoKozikNiven}, \emph{conservative digraphs} (digraphs expanded with all subsets of vertices as unary relations) \cite{HellRafiey-list-homomorphism-digraphs}, \emph{binary conservative structures} (even 3-conservative) \cite{Kazda-binary-conservative}.
We note that this phenomenon also occurs for many large classes of infinite structures: for example for all first-order expansions of the basic relations of RCC5~\cite{BodirskyBodorUIP};
see~\cite{Qualitative-Survey} for a survey on the question of which infinite-domain CSPs are in Datalog. In this chapter, however, we only consider finite structures. Additionally, for classes of finite structures, 
if the easy structures in the class cannot count, then the algebraic tractability condition 
for that class can be tested in polynomial time~\cite{MetaChenLarose}.

In \cite{Bulin18} Bul\'in conjectured that `easy trees cannot count' (i.e., only NP-hard trees possess the ability to count) and showed that this conjecture held for a large yet structurally limited subclass of trees.
\begin{conjecture}\label{conj:bulin}
Let $\T$ be a tree. 
If $\T\not\models\WNU{3,4}$, then $\T\ppleq\K_3$. 
\end{conjecture}

This conjecture (which is a rephrasing of Conjecture 2 in~\cite{Bulin18}) would provide a negative answer to an open question posed by Hell, Ne\v{s}et\v{r}il, and Zhu~\cite{HNZ} (Open Problem~1 at the end of the article): they asked whether there exists a tractable tree which does not have bounded treewidth duality.

\section{Contributions}
\label{sect:results}

In this chapter, we continue the research program of systematic, computer-based investigation and classification of the complexity of CSPs by encoding conditions ensuring tractability as constraint problems via the indicator structure (see Definition~\ref{def:indicatorStructure}) \cite{JeavonsCohenGyssensATest}. This program was initiated in \cite{GaultJeavonsImplementing}; our approach is more directly influenced by \cite{MetaChenLarose} but builds upon recent developments of the algebraic theory. We obtain the following results.

\begin{enumerate}
    \item We find 36 NP-hard trees with 20 vertices; moreover, we prove that all smaller trees and all other trees with 20 vertices are tractable. 
        \item We find four NP-hard triads with 22 vertices, and prove that all smaller triads and all other triads with 22 vertices are tractable.
\item We show that all the trees with at most 20 vertices that are not NP-hard can be solved by Datalog, confirming Conjecture~\ref{conj:bulin} for trees with at most 20 vertices. 
    \item We find a tree with 19 vertices that cannot be solved by arc consistency, and prove that all smaller trees and all other trees with 19 vertices can be solved by arc consistency.
    \item We find 8 NL-hard trees with 12 vertices; moreover, we prove that all smaller trees and all other trees with 12 vertices are in L. 
\end{enumerate}

Even though we draw from the results of the universal-algebraic approach to the CSP which led to the theorems of Bulatov and of Zhuk, and use state-of-the-art computers for our computations, these tasks remain challenging due to the huge number of trees: for example, even considered up to isomorphism, there are
139,354,922,608 trees with 20 vertices (see Table~\ref{table:otrees}), which is prohibitive even if we could test the algebraic tractability condition, i.e., the existence of a Siggers polymorphism, within milliseconds. 
Several further contributions of this article are related to the way we managed to overcome these difficulties.

A well-known key simplification is to only consider trees that are cores. Note that, if $\T$ is a tree, then the core $\T'$ of $\T$ is a tree as well (and its size is at most the size of $\T$). Hence, it suffices to work with trees that are cores.
Hell and Ne\v{s}et\v{r}il proved that deciding whether a given digraph is a core is coNP-complete~\cite{cores}. However,
it follows from a result of Chen and Mengel~\cite{ChenMengel} (Lemma 25) that for trees this problem can be decided in polynomial time. Our next result provides a more efficient algorithm for the same task and remarkably seems to be unnoticed in the literature.

\begin{enumerate}[resume]
    \item A single execution of the arc consistency procedure can be used to decide whether a given tree is a core.
\end{enumerate}

There are far too many trees with at most 20 vertices to run the core test on each of them. Our next contribution is a method to generate the core trees more directly, rather than generating all trees and then discarding the non-cores (the details can be found in Section~\ref{sect:gen}). 
Applying our method we were able to construct all trees that are cores up to size 20, and all triads that are cores up to size 22, which was essential to achieve the results 1.-5.\ above.

\begin{enumerate}[resume]
    \item 
    We computed the number of trees that are cores up to isomorphism for sizes up to 20 (see Table~\ref{table:otrees}). In particular, there are 779268 core trees of size 20.
\end{enumerate}
These are still too many to be tested for the algebraic tractability condition 
if this is implemented naively. We therefore use results from the universal algebraic approach to first run more efficient tests for certain sufficient conditions,
such as the existence of a binary symmetric polymorphism; 
this will be explained in Section~\ref{sect:specific-polymorphisms}.

Finally, we identify trees that are important `test-cases' for the open problems that have been mentioned earlier. 
\begin{enumerate}[resume]
    \item We computed the two smallest trees that are not known to be in  NL; they have 16 vertices, and they are the smallest trees that do not have a majority polymorphism. 
    \item We computed 28 smallest trees that are candidates for failing the condition in Conjecture~\ref{conj:NL} and hence 
    might be P-hard (and that are thus candidates for not being in NL, unless $\operatorname{NL}=\operatorname{P}$); they have 18 vertices. 
\end{enumerate}


\section{Notation for Trees}
\label{sect:prelims}

An \emph{undirected tree} is a connected  undirected graph without cycles. If $u, v \in T$ and $\T$ is an undirected tree, then  there exists a unique path $\mathbb P$ from $u$ to $v$ in $\T$;
the number of edges of $\mathbb P$ is denoted by $\dist(u,v)$. 
A vertex $v\in T$ is called a \emph{center} of $\T$ if $v$ lies in the middle of a longest path in $\T$.
An edge $e\in \edges{\T}$ is called a \emph{bicenter} of $\T$ if $e$ is the middle edge of a longest path in $\T$. 
We will use the following classical result.

\begin{theorem}[Jordan (1869)]\label{thm:jordan} 
An undirected tree $\T$ has exactly one center or one bicenter.
\end{theorem}

If $\mathbb O$ is an orientation of a tree and $u,v \in O$,  then $\dist(u,v)$ (center of $\mathbb O$,  bicenter $\mathbb O$)
is meant with respect to the underlying undirected tree. 
As mentioned in the introduction: an orientation of a finite tree will simply be called a \emph{tree}. 
A \emph{rooted tree} is a tuple $(\T,r)$, where  $\T$ is a tree and $r\in T$; $r$ is then called the \emph{root} of $\T$. 
A rooted tree $(\T,r)$ is called a \emph{rooted core} if every endomorphism of $\T$ that fixes $r$ is injective. 
The \emph{depth} of a rooted tree $(\T,r)$ is $\max\{\operatorname{dist}(r,v)\mid v \in T\}$.

\section{The Arc-consistency Procedure}
\label{sect:AC}
One of the most efficient algorithms employed by constraint solvers to reduce the search space 
is the \emph{arc-consistency} procedure.
In the graph homomorphism literature, the algorithm is sometimes
called the \emph{consistency check algorithm}.
The arc-consistency procedure is important for us for several reasons:
\begin{itemize}
    \item It plays a crucial role for efficiently deciding whether a given tree is a core (Section~\ref{sect:cores}). 
    \item It is well suited for combination with exhaustive search to prune the search space, and this will be relevant in Section~\ref{subsect:testing-linear}.
    \item It is an important fragment of Datalog of independent interest from the point of view of the CSP theory (see Section~\ref{sect:AC-Pol}), and we will later perform experiments to compute the smallest tree that cannot be solved by arc consistency (Section~\ref{sect:NAC}). 
\end{itemize}
We will give a short description of the procedure: 
Let $\G$ and $\Hb$ be finite digraphs. We would like to determine whether there exists a homomorphism from $\G$ to $\Hb$. 
The idea of the arc-consistency procedure is to maintain for each $x \in G$ a set $L(x) \subseteq H$.
Informally, each element of $L(x)$ represents a candidate for an image of $x$ under a homomorphism from $\G$ to $\Hb$.
The algorithm initializes each list $L(x)$ with $H$ and successively removes vertices from these lists; it only removes 
a vertex $u \in H$ from $L(x)$ if there is no homomorphism from $\G$ to $\Hb$ that maps $x$ to $u$. 
To detect vertices $x,u$ such that $u$ can be removed from $L(x)$, the algorithm uses two rules (in fact, one rule and
a symmetric version of the same rule): if $(x,y) \in \edges{\G}$, then
\begin{align*}
\text{remove $u$ from $L(x)$ if there is no $v \in L(y)$
with $(u,v) \in \edges{\Hb}$;}\\
\text{remove $v$ from $L(y)$ if there is no $u \in L(x)$
with $(u,v) \in \edges{\Hb}$.}
\end{align*}
If eventually we cannot remove any vertex from any list with these rules any more, the digraph $\G$ together with 
the lists for each vertex is called \emph{arc-consistent}.
Note that formally we may view $L$ as a function $L \colon G \to 2^H$  from the vertices of $\G$ to sets of vertices of $\Hb$.

Note that we may run the algorithm also on digraphs $\G$ where for some $x \in G$ the list $L(x)$ is already set to some subset of $H$. 
In this setting, the input consists of $\G$ and the given lists, and we are looking for a homomorphism $h$ from $\G$ to $\Hb$ such that 
for every $x \in G$ we have $h(x) \in L(x)$. 
The pseudocode of the entire arc-consistency procedure is displayed in Algorithm~\ref{alg:ac}. The standard 
arc-consistency procedure $\AC_{\Hb}(\G)$ is then obtained by calling $\AC_{\Hb}(\G,L)$ with $L(x) \coloneqq H$ for all $x \in G$.

\RestyleAlgo{ruled}
\SetAlgoVlined{}
\begin{algorithm}
\DontPrintSemicolon{}
\SetKwInOut{Input}{input}\SetKwInOut{Data}{data}
\SetKw{Remove}{remove}
\Input{a finite digraph $\G$}
\Data{a list $L(x) \subseteq H$ for each vertex $x \in G$}
\caption{$\AC_{\Hb}(\G,L)$ }\label{alg:ac}
\Repeat{no list changes}{
  \ForEach{$(x,y) \in \edges{\G}$}{
    \If{there is no $v \in L(y)$ with $(u,v) \in \edges{\Hb}$}{\Remove{$u$ from $L(x)$}}
    \If{there is no $u \in L(x)$ with $(u,v) \in \edges{\Hb}$}{\Remove{$v$ from $L(y)$}}
  }
  \If{$L(x)$ is empty for some vertex $x \in G$}{\bf reject}
}
\end{algorithm}

Clearly, if the algorithm removes all vertices from one of the lists, then there is no homomorphism from $\G$ to $\Hb$. 
It follows that if $\AC_{\Hb}$ rejects $\G$, then there is no homomorphism from $\G$ to $\Hb$.  
The converse implication does not hold in general.
For instance, let $\Hb$ be the loopless digraph with two vertices and two edges, denoted $\mathbb K_2$, and let $\G$ 
be $\mathbb K_3 = (\{0,1,2\};\neq)$. In this case, $\AC_{\Hb}$ does not remove any vertex from any list, but obviously 
there is no homomorphism from $\mathbb K_3$ to $\mathbb K_2$.

The arc-consistency procedure can be implemented so that it runs in $O(|\edges{\G}| \cdot |H|^3)$, e.g.~by Mackworth's AC-3 algorithm~\cite{Mackworth}.

\section{Cores of Trees}
\label{sect:cores}
Recall that the problem of deciding whether a given digraph is a core is coNP-complete~\cite{cores}. 
The following theorem implies that 
whether a given  finite tree is a core can be tested in polynomial time.

\begin{theorem}\label{thm:ac} 
Let $\T$ be a finite tree. Then the following are equivalent.
\begin{enumerate}
    \item $\T$ is a core;
\item $\End(\T) = \{\id_{T}\}$;
\item $\AC_{\T}(\T)$
terminates such that the list for each vertex contains a single element.
\end{enumerate}
\end{theorem}

\begin{corollary}\label{cor:treecoreP}
There is a polynomial-time algorithm to decide whether a given finite tree is a core. 
\end{corollary}

The equivalence of 1 and 2 in Theorem~\ref{thm:ac} was already shown in greater generality in Lemma~4.1 in~\cite{LLT}. 
To prove Theorem~\ref{thm:ac} we first show the following two useful lemmata.

\begin{lemma}\label{lem:ac-new}
Let $\T$ be a finite tree and let $\Hb$ be a finite digraph such that $\AC_{\Hb}(\T)$ does not reject. Let $t\in T$, and 
let $a\in H$ be such that $a\in L(t)$ after running $\AC_{\Hb}(\T)$. Then there is a homomorphism $h\colon \T\to \Hb$ such that $h(t)=a$. 
\end{lemma}

\begin{proof}
Let $S$ be a maximal subtree of $\T$ such that $t\in S$ and there exists a partial homomorphism $h\colon S\to\Hb$ with $h(t)=a$. 
If $S\neq T$, then there exists $x\in S$ and $y\in T\setminus S$ such that either $(x,y)$ or $(y,x)$ is an edge in $\T$; without 
loss of generality assume that $(x,y)\in\edges{\T}$. Because the value $u=h(x)$ was not removed from $L(x)$ when running $\AC_{\Hb}(\T)$, 
it follows that there exists $v\in L(y)$ such that $(u,v)\in\edges{\mathbb H}$. But then setting $h(y)=v$ extends $h$ to a 
partial homomorphism from $S\cup\{y\}$ to $\Hb$ contradicting maximality of the subtree $S$.
\end{proof}

\begin{lemma}\label{lem:rootedTreeAutoImpliesEndo}
Let $(\mathbb T,r)$ be a rooted tree  with an automorphism that is not the identity. Then $(\T,r)$ has a non-injective endomorphism.
\end{lemma}
\begin{proof}
Let $h$ be an automorphism of $(\mathbb T,r)$ that is not the identity. We prove the statement by induction on the number of vertices of $\mathbb T$.
Consider the components of the graph obtained from $\T$ by deleting $r$. If there is a component $C$ such that $h$ does not map $C$ into itself, then the mapping which agrees with $h$ on $C$ and which fixes all other vertices of $\T$ is a non-injective endomorphism of $(\T,r)$.

If each component $C$ is mapped by $h$ into itself, then each $h|_C$ is an automorphism of $(C,r_C)$, where $r_C$ is the unique neighbor of $r$ that lies in $C$. Since $h$ is not $\operatorname{id}_T$ there must be some $C$ such that $h|_C$ is not $\operatorname{id}_C$ and by the induction hypothesis there exists a non-injective endomorphism $h'$ of $(C,r_C)$. Since $h'(r_C)=r_C$ the mapping which extends $h'$ to $T$ by fixing all other vertices of $\T$ is a non-injective endomorphism of $(\T,r)$.
\end{proof}

\begin{proof}[Proof of Theorem~\ref{thm:ac}]
We prove the equivalence of 1.\ and 2., and then the equivalence of 2.\ and 3.

Clearly, 2.\ implies 1. Conversely, suppose that $\T$ has an endomorphism $h$ which is not the identity map. If $h$ is not injective, then $\mathbb T$ is not a core and we are done. 
Hence, suppose that $h$ is an automorphism. Note that by Theorem~\ref{thm:jordan}, if $\T$ has a center, then $h(c)=c$ and if $\T$ has a bicenter $(x,y) \in E(\T)$, then $h(\{x,y\})=\{x,y\}$. In the latter case, since $(y,x)\notin E(\T)$, we must have $h(x)=x$ and $h(y)=y$.
In both cases $h$ has a fixed point $r$ and $h$ is an automorphism of $(\T,r)$. By Lemma~\ref{lem:rootedTreeAutoImpliesEndo}, $\T$ has a non-injective endomorphism and is therefore not a core.

To prove that 2.\ implies 3., we prove the contrapositive. 
Suppose that $\AC_{\T}(\T)$ terminates with $|L(x)|>1$ for some $x \in T$. Then there exists $y\in L(x)$ with $y\neq x$. Lemma~\ref{lem:ac-new} implies that there is an endomorphism $h$ of $\T$ such that $h(x) = y$ and thus $h\neq\id_{T}$.

To see that 3.\ implies 2., note that $L(x)=\{x\}$ since because of the identity endomorphism $\id_{T}$, $x$ cannot be removed from $L(x)$. Therefore, for any endomorphism $h\colon\T\to\T$ and any $x\in T$ we must have $h(x)\in L(x)$, and so $h(x)=x$ and $h=\id_{T}$.
\end{proof}

\section{Generating all Core Trees}
\label{sect:gen} 
In this section we present an algorithm to generate all core trees with $n$ vertices up to isomorphism. 
To this end, we first present a known algorithm that generates all trees with $n$ vertices up to isomorphism \cite{103287}.
Later we explain how to modify this algorithm to directly generate core trees.

We also refer to the isomorphism classes of trees as \emph{unlabeled trees}, as opposed to \emph{labeled trees}, which are trees with vertex set $\{1,\dots,n\}$ for some $n \in {\mathbb N}$. 
The difference between the enumeration of labelled and  unlabeled trees is significant:
while the number of labelled trees is Sloane's integer sequence A097629, given by $2 (2n)^{n-2}$, 
the number of unlabeled trees  is Sloane's integer sequence A000238, which grows
asymptotically as $c d^n / n^{5/2}$
where $d \approx 5.6465$ and $c \approx 0.2257$ are constants;
the initial terms are shown in Table~\ref{table:otrees}. 
However, these numbers are still too large to apply the core test to all the unlabeled trees separately. 
The number of unlabeled trees that are cores is again much smaller. 
We therefore present a modification of the generation algorithms that allows us to generate  unlabeled trees that are cores directly without enumerating all unlabeled trees.

Let $\geq$ be some total order on all rooted trees that linearly extends the order by depth. 
The idea of the algorithm is to generate all unlabeled rooted  trees with at most $n-1$ vertices and then use Theorem~\ref{thm:jordan}.

\begin{algorithm}[ht]
\SetKwInOut{Input}{input}\SetKwInOut{Output}{output}
\caption{$\operatorname{GenerateTrees}$}\label{alg:genTrees}
\DontPrintSemicolon{}
\Input{a positive integer $n$}
\Output{a list of trees with $n$ vertices}
\Begin{
$\mathrm{Trees} \leftarrow \emptyset$\;
\tcp{bicenter}
\ForEach{$(\T_1,r_1),(\T_2,r_2)$ rooted trees where
  $|T_1| + |T_2| = n$ and $\operatorname{depth}(\T_1,r_1)=\operatorname{depth}(\T_2,r_2)$
}{
  $T \coloneqq T_1\uplus T_2$\\
  $E \coloneqq \{(r_1,r_2)\}\uplus \edges{\T_1}\uplus \edges{\T_2}$\\
  $\T \coloneqq (T;E)$\\ 
  $\mathrm{Trees} \leftarrow \mathrm{Trees} \cup \{\T\}$
}

\hspace{0.5cm}

\tcp{center}
\ForEach{ $(\T_1,r_1)\geq(\T_2,r_2)\geq\dots\geq(\T_m,r_m)$ rooted trees where
$|T_1|+\dots+|T_m|=n-1$ and $\operatorname{depth}(\T_1,r_1)=\operatorname{depth}(\T_2,r_2)$
}{
  \ForEach{ $s\in\{0,1\}^m$ where $i<j$ and $(\T_i,r_i)=(\T_j,r_j)$ imply $s_i\leq s_j$}{
  $T \coloneqq \{r\}\uplus T_1\uplus\dots\uplus T_m$\\
  $E \coloneqq \{(r_i,r)\mid r_i=1\}\uplus\{(r,r_i)\mid r_i=0\}\uplus E(\T_1)\uplus\dots\uplus E(\T_m)$\\
  $\T \coloneqq (T;E)$\\
  $\mathrm{Trees} \leftarrow \mathrm{Trees} \cup \{\T\}$
 }
}
\Return{$\mathrm{Trees}$}
}
\end{algorithm}

\begin{algorithm}
\caption{$\operatorname{GenerateRootedTrees}$}\label{alg:genRootedTrees}
\SetKwInOut{Input}{input}\SetKwInOut{Output}{output}
\Input{two positive integers $n$, $d$}
\Output{a list of rooted trees with $n$ vertices and depth
$d$}
\Begin{
\If{$n=0$}{
\Return{$\emptyset$}}
\If{$n=1$}{
\Return{$\{((\{r\};\emptyset),r)\}$}}
$\mathrm{RootedTrees} \leftarrow \emptyset$\\
\ForEach{$(\T_1,r_1)\geq\dots\geq(\T_m,r_m)$ rooted trees where
 $|T_1|+\dots+|T_m|=n-1$  and $\operatorname{depth}(\T_1,r_1)=d-1$}{
    \ForEach{$s\in\{0,1\}^m$ where $i<j$ and $(\T_i,r_i)=(\T_j,r_j)$ imply $s_i\leq s_j$}{
  $T \coloneqq \{r\}\uplus T_1\uplus\dots\uplus T_m$\\
  $E \coloneqq \{(r_i,r)\mid r_i=1\}\uplus\{(r,r_i)\mid r_i=0\}\uplus \edges{\T_1}\uplus\dots\uplus \edges{\T_m}$\\
  $\T \coloneqq (T;E)$\\
  $\mathrm{RootedTrees} \leftarrow \mathrm{RootedTrees} \cup \{(\T,r)\}$
  }
    }
\Return{$\mathrm{RootedTrees}$}
}
\end{algorithm}

It is easy to verify that Algorithm~\ref{alg:genRootedTrees} produces all unlabeled rooted trees with $n$ vertices and depth~$d$.
Analogously,  Algorithm~\ref{alg:genTrees} generates all unlabeled trees with $n$ vertices.
Remarkably, there are no isomorphism checks necessary and Algorithm~\ref{alg:genTrees} runs in linear time in the number of unlabeled trees with $n$ vertices plus the number of unlabeled rooted trees with at most $n-1$ vertices. 

Let us make some observations.
Let $(\T_1,r_1)\geq\dots\geq(\T_m,r_m)$ be rooted trees, $s\in\{0,1\}^m$, $T \coloneqq \{r\}\uplus T_1\uplus\dots\uplus T_m$, and 
\[E \coloneqq \{(r_i,r)\mid s_i=1\}\uplus\{(r,r_i)\mid s_i=0\}\uplus E(\T_1)\uplus\dots\uplus E(\T_m).\]
\begin{itemize}
    \item A rooted tree $(\T,r)$ is a rooted core if and only if $\AC_{\T}(\T,L)$, where $L(r)=\{r\}$ and $L(x)=T$ for $x\in T\setminus\{r\}$, terminates such that the list for each vertex contains a single element.
    \item By Corollary \ref{cor:treecoreP}, testing whether a (rooted) tree is a (rooted) core can be checked in polynomial time using the arc-consistency procedure. 
    \item If $(\T,r)$ is a rooted core, then  $(\T_i,r_i)$ is a rooted core for every $i$.
    \item If $\T$ is a core and $r$ is its center, then  $(\T_i,r_i)$ is a rooted core for every $i$.
    \item If $(T_1\uplus T_2;\{(r_1,r_2)\}\uplus \edges{\T_1}\uplus \edges{\T_2})$ is a core and $(r_1,r_2)$ is its bicenter, then $(\T_1,r_1)$ and $(\T_2,r_2)$ are rooted cores.
    \item If two trees that are cores are homomorphically equivalent, then they are isomorphic. 
\end{itemize}

To generate all oriented trees that are cores we slightly modify both algorithms. In both functions we only add trees to the output if they are cores or rooted cores, respectively.
By the above observations, these modified algorithms generate each tree with $n$ vertices that is a core exactly once. We do not know whether our algorithm is a polynomial-delay generation procedure for unlabeled core trees. In practice, it is fast enough to generate all core trees with at most 20 vertices within  reasonable time (see Section~\ref{sect:exp}).

\section{Testing Minor Conditions}\label{subsect:testing-linear}


To test for the existence of polymorphisms satisfying a given minor condition $\Sigma$,
we run the arc-consistency procedure for $\Hb$ on the indicator digraph $\Ind(\Sigma,\Hb)$ and then perform an exhaustive search. While this procedure is not (provably) in P, it is very efficient in practice. Observe that the size of the indicator graph grows exponentially in the arity of the function symbols and only linearly in the number of function symbols. Hence, we will always try to use conditions where all function symbols have a small arity.

We initialize the lists for vertices of $\Ind(\Sigma,\Hb)$ with preset values dictated by non-height-one identities, as explained above. Additionally, for every $u \in H$, we initialize the list for every vertex of $\Ind(\Sigma,\Hb)$ of the form $f(u,\dots,u)$ with $\{u\}$ (since it suffices to look for idempotent polymorphisms and this reduces the search space). For the remaining vertices of $\Ind(\Sigma,\Hb)$, the lists are initialized to $H$.

If $\AC_{\Hb}$ detects an inconsistency, we can be sure that no polymorphisms of $\H$ satisfy $\Sigma$. Otherwise, we select some vertex $x \in \Ind(\Sigma,\Hb)$, and set $L(x)$ to $\{u\}$ for some $u \in L(x)$. 
Then we proceed recursively with the resulting lists. If $\AC_{\Hb}$ now detects an empty list, we backtrack, but remove $u$ from $L(x)$. Finally, if the algorithm does not detect an empty list at the first level of the recursion, we end up with singleton lists for each vertex $x \in \Ind(\Sigma,\Hb)$, which defines a homomorphism from $\Ind(\Sigma,\Hb)$ to $\Hb$. The restriction of this homomorphism to the vertices of $\Ind(\Sigma,\Hb)$ for a specific function symbol can then be interpreted as (the function table of) a polymorphism of $\Hb$, and these polymorphisms satisfy $\Sigma$.

There are numerous heuristics that often help to speed up this backtracking procedure. One of the best known is called 
Maintaining Arc Consistency (MAC) \cite{Sabin1994ContradictingCW}. This family of algorithms has the arc-consistency procedure
at its core and takes advantage of the incremental design of the backtracking procedure by maintaining data structures which help to 
reduce the number of consistency checks. Another common way to speed up the search procedure is to choose the vertex $x \in \Ind(\Sigma,\Hb)$ 
that has a list of smallest size.

\subsection{Level-wise Satisfiability}\label{sect:level-wise}

Let $\G$ be a balanced digraph. Recall from Section~\ref{sec:levelWiseSatisfyability} the definition of level-wise satisfiability. 
When testing whether a minor condition is level-wise satisfied, we do not need to construct the full indicator digraph. Instead, for every function symbol (say of arity $k$) we construct only the subgraph $\G^k_{\height}$ of $\G^k$ consisting of all same-level $k$-tuples (i.e., tuples in which all vertices are from the same level). 
This optimization is particularly useful when testing for polymorphisms with large arities such as the condition $\TS n$ for all $n$; see Section~\ref{sect:NAC}. Even if we have not proved that satisfiability of $\Sigma$ is equivalent  level-wise satisfiability of $\Sigma$, we can still run the the level-wise test first and if it fails we can be sure that there are no polymorphisms satisfying $\Sigma$, not even level-wise ones.

\section{Specific Polymorphism Conditions}\label{sect:specific-polymorphisms}
In this section we focus on certain concrete 
minor conditions that are relevant for studying the membership of CSPs in the most prominent complexity classes in the subsequent sections. 
An overview of the classes and the respective minor conditions 
is given in Figure~\ref{fig:landscape}. Solid arrows indicate implications, dotted arrows indicate conjectures. 
Figure~\ref{fig:tax} shows
the relationships between relevant minor conditions that are defined throughout the section. The left side shows the general case of finite-domain structures with possibly infinite signatures and the right side shows the case for trees assuming Conjecture \ref{conj:bulin} (and $\operatorname{P}\neq\operatorname{NP}$). The implications are either immediate or from the literature~\cite{HobbyMcKenzie} (Chapter 9), ~\cite{BartoKozikWillard,KearnesMarkovicMcKenzie,Barto-cd}.

\begin{figure}
    \centering
    \begin{tikzpicture}[scale=0.71]
        \node[align=left] at (0.3,0) {Computational\\ Complexity};
        \node[align=left] at (4.5,0) {Descriptive\\ Complexity};
        \node[align=left] at (9,0) {pp-Constructions};
        \node[align=left] at (13.5,0) {Minor\\ Conditions};
        
        \node (L1) at (0,-2) {$\in \operatorname{L}$};
        \node[align=left] (LSD1) at (4.5,-2) {Symmetric\\linear\\Datalog};
        \node (Leq1) at (9,-2) {\parbox{1.5cm}{\centering no $\Ord$ no $\TLinP$}};
        \node[align=left] (HM1) at (13.5,-2) {Noname};
        
        \node (L2) at (0,-5) {$\in \operatorname{NL}$};
        \node[align=left] (LSD2) at (4.5,-5) {Linear\\Datalog};
        \node (Leq2) at (9,-5) {\parbox{2.25cm}{\centering no Horn-3SAT no $\TLinP$}};
        \node[align=left] (HM2) at (13.5,-5) {Kearnes- \\Kiss};
        
        \node (L3) at (0,-8) {$\in \operatorname{P}$};
        \node (Leq3) at (9,-8) {\centering no $\mathbb K_3$};
        \node[align=left] (HM3) at (13.5,-8) {Kearnes-\\Markovi\'{c}-\\McKenzie}; 
        
        \draw (-1.7,-0.7) -- (16,-0.7);
        
        \path
        (HM1) edge[<->] node[above] {$\operatorname{HMcK}$} (Leq1)
        (LSD1) edge[->] node[below] {E} (Leq1)
        (LSD1) edge[->] node[above] {LT} (L1)
        
        (HM2) edge[<->] node[above] {$\operatorname{HMcK}$} (Leq2)
        (LSD2) edge[->] node[below] {AC} (Leq2)
        (LSD2) edge[->] node[above] {D} (L2)
        
        (HM3) edge[<->] node[above] {$\KMM$} (Leq3)
        
        (L1) edge[->] (L2)
        (L2) edge[->] (L3)
        (LSD1) edge[->] (LSD2)
        (Leq1) edge[->] (Leq2)
        (Leq2) edge[->] (Leq3)
        (HM1) edge[->] (HM2)
        (HM2) edge[->] (HM3)
        
        (L3) edge[->,bend right=10,dashed] node[below] {$\operatorname{P}\neq\operatorname{NP}$} (Leq3)        
        (Leq3) edge[->,bend right=10] node[above] {BZ} (L3)

        (L2) edge[->,bend right=18,dashed] node[below] {$\operatorname{NL}\neq\operatorname{P}$, $\operatorname{NL}\neq\operatorname{Mod}_p\text{L}$} (Leq2)
        (L1) edge[->,bend right=18,dashed] node[below] {$\operatorname{L}\neq\operatorname{NL}$, $\operatorname{L}\neq\operatorname{Mod}_p\text{L}$} (Leq1)
        ;
        \path[dashed]
        
        (Leq1) edge[->,bend right=15,in=190] node[above] {Conj.~\ref{conj:L}} (LSD1)
        (Leq2) edge[->,bend right=15] node[above] {Conj.~\ref{conj:NL}} (LSD2)
        
        ;
    \end{tikzpicture}
    \caption{An overview of (computational and descriptive) complexity classes that are relevant for finite-domain CSPs, of important pp-constructions, and of the respective polymorphism conditions.
    LT stands for Larose and Tesson~\cite{LaroseTesson}, E stands for Egri \cite{EgristConNotinSymLinDL}, $\operatorname{HMcK}$ stands for Hobby-McKenzie~\cite{HobbyMcKenzie} (Chapter 9), D stands for Dalmau~\cite{Dalmau_2005LinearDatalog}, 
    AC stands for Afrati and Cosmadakis~\cite{AfratiCosmadakis}, BZ stands for Bulatov~\cite{BulatovFVConjecture} and Zhuk~\cite{ZhukFVConjecture}, and KMM stands for Kearnes, Markovi\'c, and McKenzie~\cite{KearnesMarkovicMcKenzie}. 
    }
    \label{fig:landscape}
\end{figure}

\begin{figure}
    \centering
    \begin{tikzpicture}[align=center, node distance=2.8cm]
      \tikzset{class/.style={fill=black!10, minimum size=0.6cm, rounded corners=3mm}}
        \node (34WNU) at (0,0) {$\WNU{3,4}$};
        \node[below left of = 34WNU] (KMM) {$\KMM$};
        \node[above left of = 34WNU] (KK) {$(\exists k)\KK k$};
        \node[above left of = KMM] (HMcK) {$(\exists k)\HMcK k$};
        \node[above left of = KK] (NN) {$(\exists k)\NN k$};
        \node at ($(NN)!0.5!(HMcK)$) (HM) {$(\exists k)\HM k$};
        \node[class, above of = KK] (NU) {$(\exists k)\J k$ \\ $(\exists k)\NU k$};
        \node[above right of = NN] (Maj) {Majority};
        \node[above of = 34WNU] (TS) {$(\forall k)\mathrm{TS}(k)$};
        
        \path[->]
        (Maj) edge (NU)
        (NU) edge (KK)
        (NN) edge (KK)
        (KK) edge (34WNU)
        (NN) edge (HM)
        (34WNU) edge (KMM)
        (HMcK) edge (KMM)
        (HM) edge (HMcK)
        (KK) edge (HMcK)
        (TS) edge (34WNU)
        ;

        \node[class] (34WNU2) at (6.1,0) {$\KMM$ \\ $\WNU{3,4}$};
        \node[class, above left of = 34WNU2] (KK2) {$(\exists k)\HMcK k$ \\ $(\exists k)\KK k$};
        \node[class, above left of = KK2] (NN2) {$(\exists k)\HM k$ \\ $(\exists k)\NN k$};
        \node[ above right of = NN2] (Maj2) {Majority};
        \node[class, above of = KK2] (NU2) {$(\exists k)\J k$ \\ $(\exists k)\NU k$};
        \node[ above of = 34WNU2] (TS2) {$(\forall k)\mathrm{TS}(k)$};
        
        \path[->]
        (Maj2) edge (NU2)
        (NU2) edge (KK2)
        (NN2) edge (KK2)
        (KK2) edge (34WNU2)
        (TS2) edge (34WNU2)
        ;
    \end{tikzpicture}
    
    \caption{
    Taxonomy of polymorphism conditions of structures with a finite domain (and possibly infinite signature) that are relevant for the computational complexity of CSPs, ordered by strength; the arrows point from stronger conditions to weaker ones.
    The right picture shows the situation for trees assuming Conjecture~\ref{conj:bulin}.
    $\operatorname{NN}$ stands for \emph{Noname} (Definition \ref{def:nn}),
    $\operatorname{HM}$ stands for \emph{Hagemann-Mitschke} (Definition \ref{def:hm}), 
    $\operatorname{KK}$ stands for \emph{Kearnes-Kiss} (Definition \ref{def:kk}),
    $\operatorname{HMcK}$ stands for \emph{Hobby-McKenzie} (Definition \ref{def:hmck}),
    $\operatorname{TS}$ stands for \emph{totally symmetric} (Definition \ref{def:ts}),
    $\WNU{3,4}$ stands for \emph{3-4 weak near-unanimity pair} (Definition \ref{def:wnu34}),
    $\KMM$ stands for \emph{Kearnes-Markovi\'{c}-McKenzie} (Definition \ref{def:kmm}).
    }
    \label{fig:tax}
\end{figure}

\subsection{Containment in P}
\label{sect:Siggers}
The characterization of the algebraic condition for tractability which is the most suitable for testing with a computer 
consists of a pair of ternary operations~\cite{KearnesMarkovicMcKenzie}. 

\begin{definition}\label{def:kmm}
The \emph{Kearnes-Markovi\'{c}-McKenzie condition}, denoted by $\KMM$, consists of the identities 
\begin{align*}
p(x,y,y) &\approx q(y,x,x)\approx q(x,x,y)\\
p(x,y,x)&\approx q(x,y,x).
\end{align*}
\end{definition}

Using this characterization, the CSP dichotomy can be stated as follows.


\begin{theorem}[\cite{ZhukFVConjecture,BulatovFVConjecture,KearnesMarkovicMcKenzie}]
Let $\A$ be a finite relational structure. Then
\[\A\ppleq\K_3\text{ if and only if }\A\models\Siggers\text{ if and only if }\A\models\KMM.\]
\end{theorem}

This characterization is optimal in the following sense: every minor condition equivalent to $\KMM$ involves either an operation of arity at least 4 or at least two operations of arity 3 \cite{KearnesMarkovicMcKenzie}. This means it is usually the condition that yields the smallest indicator structure. 
However, there are several minor conditions that imply 
the existence of Kearnes-Markovi\'{c}-McKenzie polymorphisms
and that are easier to test. In particular, we use the following. 

\begin{definition}\label{def:wnu}
For $k\geq2$ the \emph{$k$-ary weak near-unanimity condition}, denoted $\WNU k$, consists of the identities
\[f(y,x,\dots,x) \approx f(x,y,x,\dots,x) \approx \cdots \approx f(x,\dots,x,y).\]
\end{definition}

Note that $\WNU2=\Sigma_2$. It is known that the existence of a $k$-ary weak near-unanimity polymorphism with $k\geq 2$ implies the existence of Kearnes-Markovi\'{c}-McKenzie polymorphisms
\cite{Siggers, KearnesMarkovicMcKenzie}, and that the existence of Kearnes-Markovi\'{c}-McKenzie polymorphisms 
implies the existence of a  $k$-ary weak near-unanimity polymorphism for some $k \geq 2$~\cite{wnuf}.
Hence, in particular, if a finite digraph $\Hb$ satisfies $\Sigma_2$ then $\csp(\Hb)$ can be solved in polynomial time. 

For both $\KMM$ and $\WNU k$, it is enough to test for level-wise satisfiability as discussed in Section \ref{sect:level-wise}. We prove the following more general claim.

\begin{lemma}\label{lemma:level-wise-two-variable}
Let $\Sigma$ be a minor condition in two variables such that both the variables appear on each side in every identity from $\Sigma$. Then a balanced digraph level-wise satisfies $\Sigma$ if and only if it satisfies $\Sigma$.
\end{lemma}
\begin{proof}
Fix some level-wise polymorphisms $f,\dots$ that satisfy $\Sigma$, and define polymorphisms $f',\dots$ in the following way (say $f$ is $k$-ary): for $x_1,\dots,x_k$ let $\ell=\min\{\level{x_i}\mid 1\leq i\leq k\}$ and let $j\in\{1,2,\dots,k\}$ be the smallest index such that $\level{x_j}=\ell$, define
\begin{align*}
f'(x_1,x_2,\dots,x_k)\coloneqq
\begin{cases}
f(x_1,x_2,\dots,x_k)&\text{if }\level{x_1}=\level{x_2}=\dots=\level{x_k}\\
x_j&\text{otherwise.}
\end{cases}
\end{align*}
To verify that the $f'$s are polymorphisms, note that if $(x_i,y_i)$ is an edge for $i\in\{1,2,\dots,k\}$, then $f'(x_1,x_2,\dots,x_k)$ and $f'(y_1,y_2,\dots,y_k)$ fall under the same case of the definition. 
If it is the second case, $x_j$ lies on the smallest level out of $\{\level{x_i}\mid 1\leq i\leq k\}$ if and only if $y_j$ lies on the smallest level out of $\{\level{y_i}\mid 1\leq i\leq k\}$. Hence, the selected coordinate $j$ is the same. 
 
Now let $x,y$ be the two variables appearing in $\Sigma$. To see that every identity from $\Sigma$ is satisfied, note that the only interesting case is when $\level{x}\neq\level{y}$, and $f'$ then chooses the variable on the lower level. 
\end{proof}

\subsection{Containment in Datalog}
Recall the following theorem from Chapter~\ref{cha:DatalogIntro}.
\begin{manualtheorem}{\ref{thm:DLiff34WNUiff3Linp}}[\cite{BoundedWidthJournal,Maltsev-Cond}, Theorem 47 in \cite{Pol}]
Let $\A$ be a finite relational structure over a finite signature. Then the following are equivalent. 
\begin{enumerate}
    \item $\cocsp(\A)$ is in Datalog, in particular $\csp(\A)$ is in P,
\item $\A\not\ppleq\TLinP$ for all primes $p$, i.e., $\A$ lacks the ability to count, and
\item $\A\models\WNU{3,4}$.
\end{enumerate}
\end{manualtheorem}

Note that by the above theorem, if Conjecture~\ref{conj:bulin} is true (and assuming P $\neq$ NP), then every tree with Kearnes-Markovi\'{c}-McKenzie polymorphisms has $\WNU{3,4}$ polymorphisms. In particular, it implies that every tree with $\KMM$ polymorphisms has a $\WNU 3$ polymorphism: a claim that is weaker but still open. Also note that by Lemma \ref{lemma:level-wise-two-variable}, it is enough to test for level-wise $\WNU{3,4}$.

\subsection{Containment in NL}
\label{sect:NLTriads}

Recall from Theorem~\ref{thm:NUimpliesInLinDL} that if $\A\models\NU n$ for some $n\geq3$, then $\cocsp(\A)$ is in linear Datalog and $\csp(\A)$ is in NL. However, testing for $\NU n$ becomes very expensive for large $n$. Fortunately for us there is an equivalent condition using only ternary function symbols.

\begin{definition}\label{def:j} 
For $n \geq 0$, the \emph{quasi J\'onsson condition of length $n$}, denoted $\J n$, 
consists of the identities
\begin{align*}
j_1(x,x,x) & \approx j_1(x,x,y) \\
j_{2i-1}(x,y,y) & \approx j_{2i}(x,y,y) && \text{ for all } i \in \{1,\dots,n\} \\
j_i(x,y,x) & \approx j_i(x,x,x) && \text{ for all } i \in \{1,\dots,2n+1\} \\
j_{2i}(x,x,y) & \approx j_{2i+1}(x,x,y) && \text{ for all } i  \in \{1,\dots,n\} \\
j_{2n+1}(x,y,y) & \approx j_{2n+1}(y,y,y).
\end{align*}
\end{definition}

Note that $\J0=\Majority=\NU3$ and that $\J n$ implies $\J{n+1}$ for every $n \geq 0$. 
Barto~\cite{Barto-cd} proved that if a finite structure with a finite relational signature has polymorphisms that form a J\'onsson chain, then it also has a near-unanimity polymorphism (albeit its arity in the proof is doubly exponential in the size of the domain). In the other direction, it is well known that $\NU n$ implies $\J{n-2}$, syntactically. Therefore, we do not test for near-unanimities of arities higher than 3; it is more efficient to test for a J{\'o}nsson chain.

\begin{theorem}
Let $\A$ be a finite relational structure. Then $\A\models\NU n$ for some $n\geq3$ if and only if $\A\models\J n$ for some $n\geq0$.
\end{theorem}





Note that the existence of polymorphisms of $\Hb$ that form a J\'onsson chain is only a sufficient condition for the containment of $\csp(\Hb)$ in NL. An incomparable sufficient condition for the containment of $\csp(\Hb)$ in NL was identified in~\cite{CarvalhoDalmauKrokhin}. The condition presented there also has a characterization via minor identities, but the arities of the operations are prohibitively large so that we did not implement this test for trees. A necessary condition is given by Corollary~\ref{cor:linDLimpliesNotHornSat}.


It is widely believed that NL is a proper subclass of P. 
Another complexity class which is believed to be a proper subclass of P is the class 
Mod$_p$L, for some prime $p$: this is defined to be the class of problems such that there exists a non-deterministic log-space machine $M$ such that an instance is in the class if and only if the number of accepting paths of $M$ on the instance is divisible by $p$; see~\cite{LaroseTesson}. It is well known that $\csp(\TLinP)$ is Mod$_p$L-complete (see the discussion in Section 1.3 of~\cite{LaroseTesson}). 
If NL would contain Mod$_p$L then this would be a considerable breakthrough in complexity theory.

\begin{definition}[from Theorem 9.11 in~\cite{HobbyMcKenzie}] \label{def:kk}
For $n\geq2$ the \emph{quasi Kearnes-Kiss condition of length
$n$}, denoted $\KK n$, consists of the identities 
\begin{align*}
d_1(x,y,y) & \approx d_1(x,x,x) \\
d_i(x,y,y) & \approx d_{i+1}(x,y,y) && \text{for even } i \in \{0,1,\dots,n-1\} 
\\ 
d_i(x,y,x) & \approx d_{i+1}(x,y,x) 
&& \text{for even } i \in \{0,1,\dots,n-1\} \\
d_i(x,x,y) & \approx d_{i+1}(x,x,y) && \text{for odd } i  \in \{1,\dots,n-1\} \\
d_{n-1}(x,y,y) & \approx d_{n-1}(y,y,y) && \text{if $n$ odd} \\
d_{n-1}(x,y,x) & \approx d_{n-1}(x,x,x) && \text{if $n$ odd} \\
d_{n-1}(x,x,y) & \approx d_{n-1}(y,y,y) && \text{if $n$ even.} 
\end{align*}
\end{definition}
Observe that our definition varies slightly from the one Hobby and McKenzie give in Theorem 9.11 in~\cite{HobbyMcKenzie}. We have replaced 
\begin{itemize}
    \item the identities $x\approx d_0(x,y,z)$ and $d_0(x,y,y)\approx d_1(x,y,y)$ by the identity 
$d_1(x,y,y)\approx d_1(x,x,x)$,
\item the identities  $d_n(x,y,z)\approx z$ and $d_{n-1}(x,x,y)\approx d_n(x,x,y)$ by the identity $d_{n-1}(x,x,y) \approx d_{n-1}(y,y,y)$ if $n$ is even, and 
\item the three identities 
\[\text{$d_n(x,y,z)\approx z$, $d_{n-1}(x,y,y)\approx d_n(x,y,y)$, $d_{n-1}(x,y,x)\approx d_n(x,y,x)$}\] by $d_{n-1}(x,y,y) \approx d_{n-1}(y,y,y)$ and $d_{n-1}(x,y,x) \approx d_{n-1}(x,x,x)$ if $n$ is odd. 
\end{itemize}
Clearly, both conditions are satisfied by the same finite core structures. Observe that our minor condition satisfies (ii) in Theorem~\ref{thm:StructurToDigraph}.
Note that $\KK n$ implies $\KK {n+1}$ for every $n \geq 2$.  
Also note that the existence of a J{\'o}nsson chain implies the existence of a Kearnes-Kiss chain \cite{HobbyMcKenzie}; namely $
\J n$ trivially implies $\KK{2n+4}$. 

\begin{theorem}[see~\cite{Pol} and~\cite{HobbyMcKenzie}]\label{thm:KKiffNoHorn3SATorNo3Linp}
A finite structure $\A$ does not satisfy $\KK n$ for any $n \geq 2$ if and only if $\A$ can pp-construct $\HornSAT$ or $\TLinP$ for some prime $p$. 
\end{theorem}

Since $\HornSAT$ is P-hard and $\TLinP$ is Mod$_p$L-hard we have in both cases that the structure $\A$ is not in NL (unless $\operatorname{NL}=\operatorname{P}$ or $\operatorname{NL}={}${Mod$_p$L}). 
Note that in both cases $\cocsp(\H)$ is not in linear Datalog.
If the conjecture that `easy trees cannot count' (Conjecture~\ref{conj:bulin}) is true, then, by Theorems~\ref{thm:KKiffNoHorn3SATorNo3Linp}, \ref{thm:DLiff34WNUiff3Linp}, and~\ref{thm:HMcKiffnoHorn3SAT}, the existence of a Kearnes-Kiss chain is equivalent to the existence of a Hobby-McKenzie chain for trees.

\subsection{Containment in L}
One of the strongest known sufficient conditions for containment in L is Theorem~\ref{thm:symLinDLiffHM} from Kazda, which states that if $\cocsp(\Hb)$ is in linear Datalog, then either 
$\Hb\models\HM n$ for some $n \geq 1$, 
and $\cocsp(\Hb)$ is in symmetric linear Datalog (and hence is in L) or $\Hb$ can pp-construct $\Ord$.
It is well known that $\Ord$ is NL-complete (see, e.g.,~\cite{Immerman}). Hence, in the latter case $\H$ would also be NL-complete. 
We now present a polymorphism condition that characterizes the finite structures that can pp-construct 
$\Ord$ \emph{or} $\TLinP$ for some prime $p$.

\begin{definition}[Theorem 9.15 in \cite{HobbyMcKenzie}]\label{def:nn} 
For $n \geq 0$, the \emph{quasi Noname condition of length $n$}, denoted $\NN n$, consists of the identities
\begin{align*}
f_0(x,y,y,z) & \approx f_0(x,x,x,x) \\
f_i(x,x,y,x) & \approx f_{i+1}(x,y,y,x) && \text{for all } i \in \{0,\dots,n-1\} \\
f_i(x,x,y,y) & \approx f_{i+1}(x,y,y,y)
&& \text{for all } i \in \{0,\dots,n-1\} \\
f_n(x,x,y,z) & \approx f_n(z,z,z,z). 
\end{align*}
\end{definition}
Note that $\NN n$ implies $\NN {n+1}$ for every $n \geq 0$.

\begin{theorem}[\cite{HobbyMcKenzie, Pol}]\label{thm:NNiffNoOrdOrNo3Linp}
A finite structure does not satisfy $\NN n$ for any $n \geq 1$ if and only if it can pp-construct the structure $\Ord$ or the structure $\TLinP$ for some prime $p$. 
\end{theorem}

It follows that if a finite digraph $\Hb$ does not satisfy $\NN n$ for any $n \geq 1$, then $\csp(\Hb)$ is NL-hard or Mod$_p$L-hard.
Hence, $\Hb$ is in this case not in L, unless $\operatorname{L}=\operatorname{NL}$ or $\operatorname{L}={}$ Mod$_p$L.
Note that Conjecture~\ref{conj:bulin} together with Theorems~\ref{thm:NNiffNoOrdOrNo3Linp}, \ref{thm:DLiff34WNUiff3Linp}, and~\ref{thm:symLinDLiffHM} implies that $\NN n$ for some $n$ and $\HM n$ for some $n$ are equivalent for trees.

\subsection{Solvability by Arc Consistency}
\label{sect:AC-Pol}
Solvability by Arc-Consistency
(and tree duality) can be characterized in terms of height one conditions as well. 

\begin{definition}\label{def:ts}
Let $n\geq2$. The \emph{$n$-ary totally symmetric condition}, denoted $\TS n$, consists of all identities  
\begin{align*}
s_n(x_1,\dots,x_n) &\approx s_n(y_1,\dots,y_n),
\intertext{where
$x_1,\dots,x_n$ and $y_1,\dots,y_n$ are variables with}  
\{x_1,\dots,x_n\} &= \{y_1,\dots,y_n\}.    
\end{align*}

\end{definition}

Note that $\TS 2=\Sigma_2=\WNU 2$.
The digraph $\Hb$ can be solved by arc consistency if and only if $\Hb\models\TS n$ for all $n\geq 2$~\cite{FederVardi,DalmauPearson}. Note that $\TS 4$ implies $\WNU{3,4}$. Also note that a finite digraph $\Hb$ satisfies $\TS n$ for all $n>2$ if and only if it satisfies $\TS{2\cdot |V(\Hb)|}$ (see the proof given in~\cite{DalmauPearson}). The arity $2\cdot |V(\Hb)|$ is still fairly large; therefore it is particularly useful that the level-wise test is sufficient.

\begin{lemma}[see proof of Lemma 4.1 in {\cite{SpecialTriads}}]\label{lemma:level-wise-ts}
For any balanced digraph $\Hb$ and $n>0$, $\Hb$ level-wise satisfies $\TS n$  if and only if $\Hb$ satisfies $\TS n$.
\end{lemma}
\begin{proof}
Let $s_n$ be a level-wise $n$-ary polymorphism of $\Hb$ that satisfies the condition $\TS n$. 
We can construct an $n$-ary totally symmetric polymorphism $s'_n$ of $\Hb$ by applying $s_n$ to the set of vertices on the smallest level. That is, for an input tuple $(x_1,\dots,x_n)$ let $\ell=\min\{\level{x_i}\mid 1\leq i\leq n\}$, $\{x_i\mid\level{x_i}=\ell\}=\{x_{i_1},\dots,x_{i_k}\}$, and set 
\[
s'_n(x_1,\dots,x_n)=s_n(x_{i_1},\dots,x_{i_k},\underbrace{x_{i_k},\dots,x_{i_k}}_{(n-k)\text{ times}}).
\]
Clearly, the definition of $s'_n(x_1,\dots,x_n)$ depends only on the set $\{x_1,\dots,x_n\}$. To see that $s'_n$ is a polymorphism, note that similarly as in Lemma \ref{lemma:level-wise-two-variable} if $(x_i,y_i)$ is an edge for $i\in\{1,2,\dots,n\}$, then $x_{i}$ lies on the smallest level among the levels of the arguments if and only if $y_{i}$ does. The rest follows from the fact that $s_n$ is totally symmetric on each level. 
\end{proof}

\section{Experimental Results}
\label{sect:exp}

\begin{table}
    \centering
    \begin{tabular}{l|cc}
         &  \thead{Implementation\\ in Rust} & \thead{Implementation\\ in Python} \\\hline
         generate trees with up to 14 vertices & x&x\\
         generate trees with up to 20 vertices & x&\\
         test for majority & x&x\\
         test for $\TS2$ & x&x\\
         test for $\TS n$ with $n>2$ & &x\\
         test for $\KMM$ & x&x \\
         test for other conditions & x & x
    \end{tabular}
    \caption{
    Comparing the usage of our Rust and Python implementations.
    }
    \label{tab:RustVsPython}
\end{table}

We implemented the AC-3 algorithm for establishing arc-consistency and used its adaptation, known as the MAC-3 algorithm, 
for maintaining arc-consistency during the backtracking procedure described in Section~\ref{subsect:testing-linear}. 
The lists and related operations were implemented by doubly-linked lists. The code is written in 
Rust 
and the experiments were run on a Intel(R) Xeon(R) CPU E5-2680 v3 (12 cores) @ 2.50GHz with Linux.
We also used another implementation written in  Python. 
All tests for chains of polymorphisms and for totally symmetric polymorphisms were performed using this implementation on a AMD Ryzen 5 4500U (with 8 cores) @ 2.38 GHz with Windows. Some results were verified with both implementations, see Table~\ref{tab:RustVsPython}.
An efficient implementation was essential to obtain our results.
\footnote{
Using the constraint modeling language MiniZinc~\cite{minizinc} and the solver Gecode~\cite{gecode} we are able to verify the polymorphism conditions for concrete trees. Also, together with the program Nauty~\cite{nauty} used to generate trees up to isomorphism we can verify the number of core trees for reasonably small sizes, but this approach was orders of magnitude slower than required.}
Both implementations and all trees presented in the figures in this section can be found at \codelink.
\begin{table}
\hfuzz=1000pt
\begin{center}
\begin{tabular}{rrrrrc S[table-format=3.5]}
\toprule
\thead{$n$} & \thead{trees} & \thead{cores} & \thead{rooted \\ cores} & \thead{AC \\ calls} & \hspace{-2mm}\thead{time per \\ AC call}\hspace{-2mm} & {\thead{total time}}\\
\toprule
1   & 1            & 1       & 1          & 1       &  10\,\si{\micro\second}      &        1\,\si{\milli\second}  \\
2   & 1            & 1       & 2          & 1       &  10\,\si{\micro\second}      &        1\,\si{\milli\second}  \\
3   & 3            & 1       & 3          & 3       &  11\,\si{\micro\second}      &        1\,\si{\milli\second}  \\
4   & 8            & 1       & 6          & 8       &  12\,\si{\micro\second}      &        1\,\si{\milli\second}  \\
5   & 27           & 1       & 11         & 19      &  13\,\si{\micro\second}      &        1\,\si{\milli\second}  \\
6   & 91           & 2       & 28         & 39      &  15\,\si{\micro\second}       &        1\,\si{\milli\second}  \\
7   & 350          & 3       & 63         & 94      &  24\,\si{\micro\second}       &        1\,\si{\milli\second}  \\
8   & 1376         & 7       & 170        & 198     &  25\,\si{\micro\second}       &        2\,\si{\milli\second}  \\
9   & 5743         & 15      & 439        & 439     &  30\,\si{\micro\second}       &       5\,\si{\milli\second}  \\
10  & 24635        & 36      & 1200       & 953     &  36\,\si{\micro\second}       &       13\,\si{\milli\second}  \\
11  & 108968       & 85      & 3307       & 2180    &  40\,\si{\micro\second}       &       33\,\si{\milli\second}  \\
12  & 492180       & 226     & 9380       & 5050    &  48\,\si{\micro\second}       &      84\,\si{\milli\second}  \\
13  & 2266502      & 578     & 26731      & 12218   &  55\,\si{\micro\second}       &      236\,\si{\milli\second}  \\
14  & 10598452     & 1569    & 77508      & 29785   &  66\,\si{\micro\second}       &        657\,\si{\milli\second}  \\
15  & 50235931     & 4243    & 226399     & 74902   &  71\,\si{\micro\second}      &    2.0\,\si{\second}  \\
16  & 240872654    & 11848   & 668228     & 190632  &  84\,\si{\micro\second}      &   5.7\,\si{\second}  \\
17  & 1166732814   & 33104   & 1984592    & 496373  &  98\,\si{\micro\second}      &   16.6\,\si{\second}  \\
18  & 5702001435   & 94221   & 5937276    & 1308847 &  121\,\si{\micro\second}      &    49.3\,\si{\second}  \\
19  & 28088787314  & 269455  & 17856807   & 3512229 &  129\,\si{\micro\second}      &    2.5\,\si{\minute}  \\
20  & 139354922608 & 779268  & 53996424   & 9538804 &  142\,\si{\micro\second}        &   7.4\,\si{\minute}
\end{tabular}
\end{center}
\caption{The number of unlabeled trees with $n$ vertices together with results of Algorithm~\ref{alg:genTrees}.}
\label{table:otrees}
\end{table}
Table~\ref{table:otrees} shows the number of unlabeled trees with $n$ vertices and the number of those that are cores.
The table suggests that the
fraction of trees that are cores quickly goes to 0. The next columns contain the number of unlabeled rooted cores with $n$ vertices, the number of AC calls, 
and the mean cpu time per AC call on a tree with $n$ vertices. The final column in the table shows the computation time needed to generate all the unlabeled 
core trees with $n$ vertices with Algorithm \ref{alg:genTrees}.

Now we present the results of testing the discussed minor conditions on these trees to classify them with respect to their computational complexity. In some cases we manage to compute all the minimal trees in the respective complexity class; the corresponding results are presented in Section~\ref{sect:sht}. In Section~\ref{sect:opentrees} we present trees whose precise complexity status is open. 
While the numbers of trees given in the text are up to isomorphism, the corresponding figures make a further restriction based on the following fact. 

\begin{remark}
\label{rem:reverse}
An operation is a polymorphism of $\Hb$ if and only if it is a polymorphism of $\Hb^R$.
\end{remark}

The remark justifies that 
our figures contain exactly one of the trees $\T$, $\T^R$.
It turns out that the trees in our figures that do not satisfy a certain minor condition have a unique minimal subtree which has no idempotent polymorphisms satisfying the respective condition.
The vertices and edges drawn in gray do not belong to this minimal subtree of $\T$.

\subsection{The Smallest Hard Trees}
\label{sect:sht}
In this section we present the smallest trees that are NP-hard and that are NL-hard, under standard assumptions from complexity theory. 
We also compute the smallest tree that cannot be solved by arc consistency, the smallest trees that cannot be solved by Datalog,
and the smallest trees that cannot be solved by  symmetric linear Datalog; these results hold without any assumptions from complexity theory.

\subsubsection{The Smallest NP-Hard Trees}
Our algorithm found that all trees with at most 19 vertices have 
Kearnes-Markovi\'{c}-McKenzie polymorphisms 
and hence are tractable. 
It also found that there exist exactly 36 trees with 20 vertices that have no Kearnes-Markovi\'{c}-McKenzie polymorphisms
and hence are NP-hard. 
For such an NP-hard tree $\T$ with 20 vertices it takes our algorithm about 0.07 seconds to construct the indicator digraph $\Ind(\KMM,\T)$ for the Kearnes-Markovi\'{c}-McKenzie polymorphisms and about 0.03 seconds to verify that this indicator digraph does not have a homomorphism to $\T$. When applying the level trick, it takes about 0.01 seconds to construct the indicator digraph and 0.03 seconds to verify that $\Ind(\KMM,\T)$ does not have a homomorphism to $\T$.

\input{trees}

The trees with 20 vertices that have no Kearnes-Markovi\'{c}-McKenzie polymorphisms 
are displayed in Figure~\ref{fig:hard-trees}.
Each of these trees has a unique smallest subtree without idempotent Kearnes-Markovi\'{c}-McKenzie polymorphisms. These subtrees are clearly not cores and have Kearnes-Markovi\'{c}-McKenzie polymorphisms. Note that these subtrees are the same for the trees A1--A8 and for A10--A18.
Since there are way too many trees up to isomorphism we were not able to answer the following open problem.
\begin{question}
Is it true that the smallest trees without idempotent $\KMM$ polymorphisms have 18 vertices and are displayed in Figure~\ref{fig:hard-trees}?
\end{question}

Moreover, there are 4 smallest triads with 22 vertices that have no Kearnes-Markovi\'{c}-McKenzie polymorphisms;
these are shown in Figure~\ref{fig:hardTriads}. 
All smaller triads have a binary symmetric polymorphism.

\def\scale{0.4}
\def\hdist{1cm}
\def\vdist{1cm}
\begin{figure}
\centering
\begin{tikzpicture}[scale=\scale]
\node[bullet,gray] (21) at (4,0) {};
\node[bullet] (4) at (0.5,1) {};
\node[bullet] (14) at (6,1) {};
\node[bullet] (20) at (4,1) {};
\node[bullet] (1) at (2,2) {};
\node[bullet] (5) at (0,2) {};
\node[bullet] (3) at (1,2) {};
\node[bullet] (9) at (5,2) {};
\node[bullet] (13) at (6,2) {};
\node[bullet] (19) at (4,2) {};
\node[bullet] (2) at (1,3) {};
\node[bullet] (6) at (0,3) {};
\node[bullet] (0) at (2,3) {};
\node[bullet] (10) at (5,3) {};
\node[bullet] (12) at (6,3) {};
\node[bullet] (16) at (3,3) {};
\node[bullet] (18) at (4,3) {};
\node[bullet] (7) at (0,4) {};
\node[bullet] (11) at (5.5,4) {};
\node[bullet] (15) at (2.5,4) {};
\node[bullet] (17) at (3.5,4) {};
\node[bullet,gray] (8) at (0,5) {};
\path[->,>=stealth']
(1) edge (2)
(1) edge (0)
(4) edge (5)
(4) edge (3)
(5) edge (6)
(6) edge (7)
(7) edge[gray] (8)
(0) edge (15)
(3) edge (2)
(9) edge (0)
(9) edge (10)
(10) edge (11)
(12) edge (11)
(13) edge (12)
(14) edge (13)
(16) edge (15)
(16) edge (17)
(18) edge (17)
(19) edge (18)
(20) edge (19)
(21) edge[gray] (20)
;
\node at (3,-1) {Triad 1};
\end{tikzpicture}
\hspace{\hdist}
\begin{tikzpicture}[scale=\scale]
\node[bullet,gray] (21) at (6,0) {};
\node[bullet] (4) at (0.5,1) {};
\node[bullet] (14) at (6,1) {};
\node[bullet] (20) at (4,1) {};
\node[bullet] (1) at (2,2) {};
\node[bullet] (5) at (0,2) {};
\node[bullet] (3) at (1,2) {};
\node[bullet] (9) at (5,2) {};
\node[bullet] (13) at (6,2) {};
\node[bullet] (19) at (4,2) {};
\node[bullet] (2) at (1,3) {};
\node[bullet] (6) at (0,3) {};
\node[bullet] (0) at (2,3) {};
\node[bullet] (10) at (5,3) {};
\node[bullet] (12) at (6,3) {};
\node[bullet] (16) at (3,3) {};
\node[bullet] (18) at (4,3) {};
\node[bullet] (7) at (0,4) {};
\node[bullet] (11) at (5.5,4) {};
\node[bullet] (15) at (2.5,4) {};
\node[bullet] (17) at (3.5,4) {};
\node[bullet,gray] (8) at (0,5) {};
\path[->,>=stealth']
(1) edge (2)
(1) edge (0)
(4) edge (5)
(4) edge (3)
(5) edge (6)
(6) edge (7)
(7) edge[gray] (8)
(0) edge (15)
(3) edge (2)
(9) edge (0)
(9) edge (10)
(10) edge (11)
(12) edge (11)
(13) edge (12)
(14) edge (13)
(16) edge (15)
(16) edge (17)
(18) edge (17)
(19) edge (18)
(20) edge (19)
(21) edge[gray] (14)
;
\node at (3,-1) {Triad 2};
\end{tikzpicture}
\caption{The smallest NP-hard triads (up to edge reversal and assuming P $\neq$ NP), they have 22 vertices.}
\label{fig:hardTriads}
\end{figure}
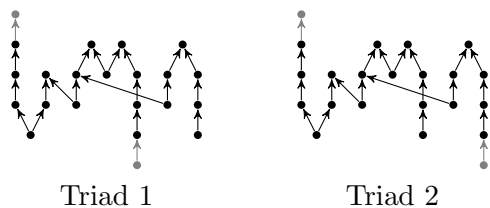

\subsubsection{The Smallest NL-hard Trees}
There are 8 trees with 12 vertices that are NL-hard. 
Two of them are isomorphic to their reverse, so we only 
 display 5 trees in Figure~\ref{fig:NLHardTrees}, called B1, B2, B3, B4, and B5. The proof that they are NL-hard can be found below.
All other trees with at most 12 vertices satisfy  $\HM 8$.

\def\scale{0.4}
\def\hdist{5mm}
\def\vdist{1cm}
\begin{figure}
\centering\begin{tikzpicture}[scale=\scale]
\node[bullet] (11) at (0,0) {};
\node[bullet] (6) at (2,0) {};
\node[bullet] (7) at (1,-1) {};
\node[bullet] (5) at (2,-1) {};
\node[bullet] (10) at (0,-1) {};
\node[bullet] (8) at (0,-2) {};
\node[bullet] (0) at (1,-2) {};
\node[bullet] (4) at (2,-2) {};
\node[bullet] (3) at (1,-3) {};
\node[bullet] (9) at (0,-3) {};
\node[bullet,gray] (1) at (2,-3) {};
\node[bullet,gray] (2) at (2,-4) {};
\path[<-,>=stealth']
(7) edge (8)
(7) edge (0)
(8) edge (9)
(0) edge (3)
(0) edge[gray] (1)
(5) edge (4)
(4) edge (3)
(1) edge[gray] (2)
(11) edge (10)
(10) edge (8)
(6) edge (5)
;
\node at (1,-5) {Tree B1};
\end{tikzpicture}
\hspace{\hdist}
\begin{tikzpicture}[scale=\scale]
\node[bullet,gray] (9) at (0,0) {};
\node[bullet] (0) at (1,1) {};
\node[bullet] (8) at (0,1) {};
\node[bullet] (7) at (0,2) {};
\node[bullet] (5) at (1,2) {};
\node[bullet] (1) at (2,2) {};
\node[bullet] (2) at (2,3) {};
\node[bullet] (10) at (0,3) {};
\node[bullet] (6) at (1,3) {};
\node[bullet] (3) at (2,4) {};
\node[bullet] (11) at (0,4) {};
\node[bullet,gray] (4) at (2,5) {};
\path[->,>=stealth']
(2) edge (3)
(3) edge[gray] (4)
(7) edge (10)
(7) edge (6)
(10) edge (11)
(0) edge (5)
(0) edge (1)
(5) edge (6)
(1) edge (2)
(9) edge[gray] (8)
(8) edge (7)
;
\node at (1,-1) {Tree B2};
\end{tikzpicture}
\hspace{\hdist}
\begin{tikzpicture}[scale=\scale]
\node[bullet,gray] (5) at (0,0) {};
\node[bullet] (4) at (0,-1) {};
\node[bullet] (11) at (2,-1) {};
\node[bullet] (0) at (1,-2) {};
\node[bullet] (3) at (0,-2) {};
\node[bullet] (10) at (2,-2) {};
\node[bullet] (6) at (1,-3) {};
\node[bullet] (9) at (2,-3) {};
\node[bullet] (1) at (0,-3) {};
\node[bullet] (7) at (1,-4) {};
\node[bullet] (2) at (0,-4) {};
\node[bullet,gray] (8) at (1,-5) {};
\path[<-,>=stealth']
(0) edge (6)
(0) edge (1)
(6) edge (7)
(4) edge (3)
(3) edge (1)
(7) edge[gray] (8)
(9) edge (7)
(1) edge (2)
(11) edge (10)
(10) edge (9)
(5) edge[gray] (4)
;
\node at (1,-6) {Tree B3};
\end{tikzpicture}
\hspace{\hdist}
\begin{tikzpicture}[scale=\scale]
\node[bullet,gray] (3) at (0,0) {};
\node[bullet] (2) at (0,1) {};
\node[bullet] (9) at (2,1) {};
\node[bullet] (6) at (1,2) {};
\node[bullet] (1) at (0,2) {};
\node[bullet] (8) at (2,2) {};
\node[bullet] (7) at (2,3) {};
\node[bullet] (4) at (0,3) {};
\node[bullet] (0) at (1,3) {};
\node[bullet] (10) at (2,4) {};
\node[bullet] (5) at (0,4) {};
\node[bullet,gray] (11) at (2,5) {};
\path[->,>=stealth']
(6) edge (7)
(6) edge (0)
(7) edge (10)
(1) edge (4)
(1) edge (0)
(4) edge (5)
(8) edge (7)
(2) edge (1)
(10) edge[gray] (11)
(9) edge (8)
(3) edge[gray] (2)
;
\node at (1,-1) {Tree B4};
\end{tikzpicture}
\hspace{\hdist}
\begin{tikzpicture}[scale=\scale]
\node[bullet,gray] (11) at (0,0) {};
\node[bullet] (10) at (0,-1) {};
\node[bullet] (2) at (2,-1) {};
\node[bullet] (1) at (2,-2) {};
\node[bullet] (9) at (0,-2) {};
\node[bullet] (6) at (1,-2) {};
\node[bullet] (3) at (2,-3) {};
\node[bullet] (0) at (1,-3) {};
\node[bullet] (7) at (0,-3) {};
\node[bullet] (8) at (0,-4) {};
\node[bullet] (4) at (2,-4) {};
\node[bullet,gray] (5) at (2,-5) {};
\path[<-,>=stealth']
(1) edge (3)
(1) edge (0)
(3) edge (4)
(7) edge (8)
(9) edge (7)
(11) edge[gray] (10)
(10) edge (9)
(4) edge[gray] (5)
(6) edge (0)
(6) edge (7)
(2) edge (1)
;
\node at (1,-6) {Tree B5};
\end{tikzpicture}
\caption{The smallest NL-hard trees (up to edge reversal and assuming L $\neq$ NL).}
\label{fig:NLHardTrees}
\end{figure}
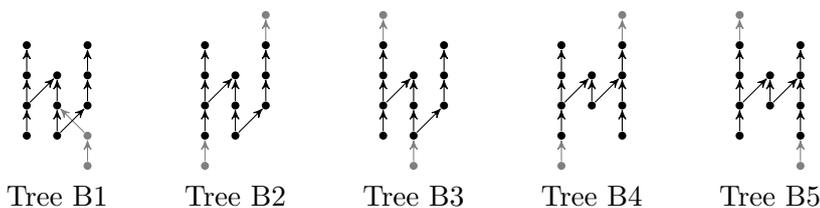

Since the trees B1-B5 have a majority polymorphism they are in NL. 
To prove that B1-B5 are NL-hard, we show that they can pp-construct $\Ord$. Hence, we need to construct the three relations \[\text{$\{0\}, \{1\},$ and $\{(0,0),(0,1),(1,1)\}$.}\] First 
note that in a core tree $\T$ any singleton set is pp-definable from $\T$, since  $\End(\T)=\{\operatorname{id}_{T}\}$ by Theorem~\ref{thm:ac}. The following two graphs represent two pp-formulas $\phi_1(x,y)$ and $\phi_2(x,y)$. The filled vertices stand for existentially quantified variables.

\begin{center}
\begin{tikzpicture}[scale=0.5]
    \node[var-f,draw,label=below:$x$] (0) at (0,0) {};
    \node[var-b] (1) at (1,1) {};
    \node[var-b] (2) at (1,0) {};
    \node[var-b] (4) at  (1,-1) {};
    \node[var-f,label=below:$y$] (6) at (2,0) {};
    \node[var-b] (5) at (0,2) {};
    \node[var-b] (9) at (0,1) {};
    \node[var-b] (7) at (2,2) {};
    \node[var-b] (8) at (2,1) {};
    
    \path[->,>=stealth']
        (0) edge (1)
        (2) edge (1)
        (4) edge (2)
        (4) edge (6)
        (9) edge (5)
        (0) edge (5)
        (8) edge (7)
        (6) edge (7)
        ;
\end{tikzpicture}
\hspace{3cm}
\begin{tikzpicture}[scale=0.5]
    \node[var-f,draw,label=below:$x$] (0) at (0,0) {};
    \node[var-b] (1) at (1,1) {};
    \node[var-b] (2) at (1,0) {};
    \node[var-b] (3) at (2,1) {};
    \node[var-b] (4) at  (2,2) {};
    \node[var-f,label=right:$y$] (6) at (2,0) {};
    \node[var-b] (7) at (2,-1) {};
    \path[->,>=stealth']
        (0) edge (1)
        (2) edge (1)
        (2) edge (3)
        (6) edge (3)
        (3) edge (4)
        (7) edge (6)
        ;
\end{tikzpicture}
\end{center}

The trees B1, B2, B3 can pp-define a structure that is homomorphically equivalent to $\Ord$ using  $\phi_1(x,y)$ for $\edges{\Ord}$. The trees B4 and B5 can do the same using  $\phi_2(x,y)$ for $\edges{\Ord}$. Since $\Ord$ is NL-hard, B1-B5 are NL-hard as well.

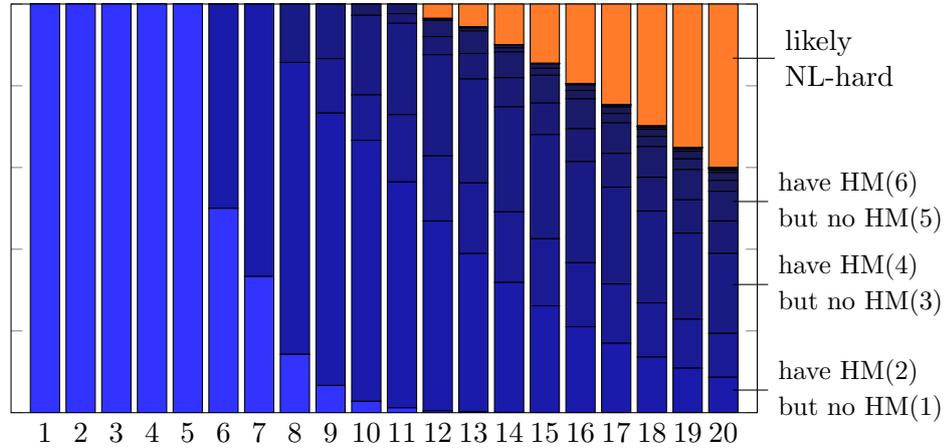
\begin{figure}
    \centering
\begin{tikzpicture}
\begin{axis}[
            ybar stacked,   
            ymin=0,         
            ymax = 1,
            xtick=data,     
            xticklabels from table={\testdata}{Label}, width=11.4cm, height = 7cm,bar width=0.385cm, enlarge x limits=0.05
]

\def\plotcommand#1{
    \addplot [fill=black!#1!blue!90!white] table [y=\s, meta=Label,x expr=\coordindex] {\testdata};
}

\addplot [fill=black!0!blue!80!white] table [y=1, meta=Label,x expr=\coordindex] {\testdata};

\foreach \s in {2,...,14}
{
\pgfmathparse{ln(\s+1)*100/3}

\expandafter\plotcommand\expandafter{\pgfmathresult}
}

\addplot [fill=orange!80!white!80!red] table [y=29, meta=Label,x expr=\coordindex] {\testdata};
\end{axis}
\node[align=left] (2) at (11.18,0.3) {\small have $\HM2$\\\small but no $\HM1$};
\node[align=left] (4) at (11.18,1.7) {\small have $\HM4$\\\small but no $\HM3$};
\node[align=left] (6) at (11.18,2.8) {\small have $\HM6$\\\small but no $\HM5$};
\node[align=left] (NL) at (10.9,4.7) {likely \\ NL-hard};
\draw (9.5,0.3) -- (2);
\draw (9.5,1.7) -- (4);
\draw (9.5,2.8) -- (6);
\draw (9.5,4.7) -- (NL);
\end{tikzpicture}
    \caption{Distribution of core trees in L.}
    \label{fig:HMtrees}
\end{figure}

So for 12 vertices, 8 out of 226 trees are NL-hard (assuming L $\neq$ NL). In Figure~\ref{fig:HMtrees} we present how this distribution in core trees changes with increasing number of vertices. Every tree with at most 20 vertices falls into one of two cases:
\begin{itemize}
    \item it satisfies $\HM{16}$ and has a majority polymorphism, hence it is in symmetric linear Datalog; in this case it is in L, or
    \item it has no $\HM{30}$.
\end{itemize}

We strongly suspect that in the latter case the trees have no $\HM n $ for any $n$, can pp-construct $\Ord$, and are therefore NL-hard.

\pgfplotstableread{ 
Label	1	2	3	4	5	6	7	8	9	10	11	12	13	14	15	16	17	18	19	20	21	22	23	24	25	26	27	28	29
1	1.00000	0.00000	0.00000	0.00000	0.00000	0.00000	0.00000	0.00000	0.00000	0.00000	0.00000	0.00000	0.00000	0.00000	0.00000	0.00000	0.00000	0.00000	0.00000	0.00000	0.00000	0.00000	0.00000	0.00000	0.00000	0.00000	0.00000	0.00000	0.00000
2	1.00000	0.00000	0.00000	0.00000	0.00000	0.00000	0.00000	0.00000	0.00000	0.00000	0.00000	0.00000	0.00000	0.00000	0.00000	0.00000	0.00000	0.00000	0.00000	0.00000	0.00000	0.00000	0.00000	0.00000	0.00000	0.00000	0.00000	0.00000	0.00000
3	1.00000	0.00000	0.00000	0.00000	0.00000	0.00000	0.00000	0.00000	0.00000	0.00000	0.00000	0.00000	0.00000	0.00000	0.00000	0.00000	0.00000	0.00000	0.00000	0.00000	0.00000	0.00000	0.00000	0.00000	0.00000	0.00000	0.00000	0.00000	0.00000
4	1.00000	0.00000	0.00000	0.00000	0.00000	0.00000	0.00000	0.00000	0.00000	0.00000	0.00000	0.00000	0.00000	0.00000	0.00000	0.00000	0.00000	0.00000	0.00000	0.00000	0.00000	0.00000	0.00000	0.00000	0.00000	0.00000	0.00000	0.00000	0.00000
5	1.00000	0.00000	0.00000	0.00000	0.00000	0.00000	0.00000	0.00000	0.00000	0.00000	0.00000	0.00000	0.00000	0.00000	0.00000	0.00000	0.00000	0.00000	0.00000	0.00000	0.00000	0.00000	0.00000	0.00000	0.00000	0.00000	0.00000	0.00000	0.00000
6	0.50000	0.50000	0.00000	0.00000	0.00000	0.00000	0.00000	0.00000	0.00000	0.00000	0.00000	0.00000	0.00000	0.00000	0.00000	0.00000	0.00000	0.00000	0.00000	0.00000	0.00000	0.00000	0.00000	0.00000	0.00000	0.00000	0.00000	0.00000	0.00000
7	0.33333	0.66667	0.00000	0.00000	0.00000	0.00000	0.00000	0.00000	0.00000	0.00000	0.00000	0.00000	0.00000	0.00000	0.00000	0.00000	0.00000	0.00000	0.00000	0.00000	0.00000	0.00000	0.00000	0.00000	0.00000	0.00000	0.00000	0.00000	0.00000
8	0.14286	0.71429	0.00000	0.14286	0.00000	0.00000	0.00000	0.00000	0.00000	0.00000	0.00000	0.00000	0.00000	0.00000	0.00000	0.00000	0.00000	0.00000	0.00000	0.00000	0.00000	0.00000	0.00000	0.00000	0.00000	0.00000	0.00000	0.00000	0.00000
9	0.06667	0.66667	0.13333	0.13333	0.00000	0.00000	0.00000	0.00000	0.00000	0.00000	0.00000	0.00000	0.00000	0.00000	0.00000	0.00000	0.00000	0.00000	0.00000	0.00000	0.00000	0.00000	0.00000	0.00000	0.00000	0.00000	0.00000	0.00000	0.00000
10	0.02778	0.63889	0.11111	0.19444	0.00000	0.02778	0.00000	0.00000	0.00000	0.00000	0.00000	0.00000	0.00000	0.00000	0.00000	0.00000	0.00000	0.00000	0.00000	0.00000	0.00000	0.00000	0.00000	0.00000	0.00000	0.00000	0.00000	0.00000	0.00000
11	0.01176	0.55294	0.16471	0.22353	0.02353	0.02353	0.00000	0.00000	0.00000	0.00000	0.00000	0.00000	0.00000	0.00000	0.00000	0.00000	0.00000	0.00000	0.00000	0.00000	0.00000	0.00000	0.00000	0.00000	0.00000	0.00000	0.00000	0.00000	0.00000
12	0.00442	0.46460	0.15929	0.24779	0.04425	0.03982	0.00000	0.00442	0.00000	0.00000	0.00000	0.00000	0.00000	0.00000	0.00000	0.00000	0.00000	0.00000	0.00000	0.00000	0.00000	0.00000	0.00000	0.00000	0.00000	0.00000	0.00000	0.00000	0.03540
13	0.00173	0.38754	0.17301	0.25433	0.06228	0.05536	0.00692	0.00346	0.00000	0.00000	0.00000	0.00000	0.00000	0.00000	0.00000	0.00000	0.00000	0.00000	0.00000	0.00000	0.00000	0.00000	0.00000	0.00000	0.00000	0.00000	0.00000	0.00000	0.05536
14	0.00064	0.31804	0.17272	0.25685	0.07138	0.06310	0.01020	0.00701	0.00000	0.00064	0.00000	0.00000	0.00000	0.00000	0.00000	0.00000	0.00000	0.00000	0.00000	0.00000	0.00000	0.00000	0.00000	0.00000	0.00000	0.00000	0.00000	0.00000	0.09943
15	0.00024	0.26137	0.16403	0.25477	0.07730	0.06811	0.01697	0.01061	0.00094	0.00047	0.00000	0.00000	0.00000	0.00000	0.00000	0.00000	0.00000	0.00000	0.00000	0.00000	0.00000	0.00000	0.00000	0.00000	0.00000	0.00000	0.00000	0.00000	0.14518
16	0.00008	0.21042	0.15648	0.24755	0.08018	0.07335	0.02043	0.01367	0.00219	0.00110	0.00000	0.00008	0.00000	0.00000	0.00000	0.00000	0.00000	0.00000	0.00000	0.00000	0.00000	0.00000	0.00000	0.00000	0.00000	0.00000	0.00000	0.00000	0.19446
17	0.00003	0.16992	0.14515	0.23683	0.08265	0.07461	0.02308	0.01577	0.00369	0.00187	0.00018	0.00006	0.00000	0.00000	0.00000	0.00000	0.00000	0.00000	0.00000	0.00000	0.00000	0.00000	0.00000	0.00000	0.00000	0.00000	0.00000	0.00000	0.24616
18	0.00001	0.13609	0.13274	0.22465	0.08262	0.07493	0.02479	0.01767	0.00467	0.00261	0.00036	0.00016	0.00000	0.00001	0.00000	0.00000	0.00000	0.00000	0.00000	0.00000	0.00000	0.00000	0.00000	0.00000	0.00000	0.00000	0.00000	0.00000	0.29867
19	0.00000	0.10889	0.11986	0.21063	0.08154	0.07410	0.02593	0.01865	0.00557	0.00327	0.00068	0.00030	0.00002	0.00001	0.00000	0.00000	0.00000	0.00000	0.00000	0.00000	0.00000	0.00000	0.00000	0.00000	0.00000	0.00000	0.00000	0.00000	0.35054
20	0.00000	0.08666	0.10727	0.19592	0.07930	0.07275	0.02641	0.01934	0.00618	0.00384	0.00092	0.00045	0.00006	0.00002	0.00000	0.00000	0.00000	0.00000	0.00000	0.00000	0.00000	0.00000	0.00000	0.00000	0.00000	0.00000	0.00000	0.00000	0.40089
}\testdata

\subsubsection{The Smallest Tree not Solved by Datalog} 
It turns out that every tree with at most 20 vertices which is not NP-hard can be solved by Datalog, thus confirming Conjecture~\ref{conj:bulin}. 
In fact, up to 20 vertices all trees that have 
Kearnes-Markovi\'{c}-McKenzie polymorphisms
either have a majority polymorphism or totally symmetric polymorphisms of all arities. The picture is however more complex for larger trees: there exists a tree which can be solved by Datalog but does not have a near-unanimity polymorphism (of any arity) and does not have totally symmetric polymorphisms of all arities, see \cite[Proposition 5.5]{BartoB13} for an example and \cite{Bulin18} for its solvability in Datalog.

\subsubsection{The Smallest Tree not Solved by Arc Consistency} 
\label{sect:NAC}
\begin{figure}
\centering\begin{tikzpicture}[scale=0.4]
\node[var-b] (4) at (0,0) {};
\node[var-b] (18) at (1,0) {};
\node[var-b] (3) at (0,1) {};
\node[var-b] (5) at (4,1) {};
\node[var-b] (11) at (3,1) {};
\node[var-b] (14) at (2,1) {};
\node[var-b] (17) at (1,1) {};
\node[var-b] (2) at (0,2) {};
\node[var-b] (15) at (1,2) {};
\node[var-b] (0) at (2,2) {};
\node[var-b] (6) at (4,2) {};
\node[var-b] (10) at (3,2) {};
\node[var-b] (1) at (0,3) {};
\node[var-b] (7) at (4,3) {};
\node[var-b] (9) at (2,3) {};
\node[var-b] (12) at (3,3) {};
\node[var-b] (16) at (1,3) {};
\node[var-b] (8) at (4,4) {};
\node[var-b] (13) at (3,4) {};
\path[->,>=stealth']
(0) edge (1)
(0) edge (9)
(2) edge (1)
(3) edge (2)
(4) edge (3)
(5) edge (0)
(5) edge (6)
(6) edge (7)
(7) edge (8)
(10) edge (9)
(10) edge (12)
(11) edge (10)
(12) edge (13)
(14) edge (0)
(14) edge (15)
(15) edge (16)
(18) edge (17)
(17) edge (15)
;

\node at (2,-1) {Tree C};
\end{tikzpicture}
\caption{The smallest tree that cannot be solved by Arc Consistency (it has 19 vertices and a majority, but no binary symmetric polymorphism).}
\label{fig:withWNU3withoutWNU2}
\end{figure}
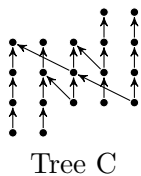

The smallest tree $\T$ that has no binary symmetric polymorphism has 19 vertices and is displayed in Figure~\ref{fig:withWNU3withoutWNU2}. It has a polymorphism satisfying $\WNU3$, and even a majority polymorphism. 
Note that $\csp(\T)$ cannot be solved by the arc-consistency procedure since in this case $\T$ must have a binary symmetric polymorphism~\cite{FederVardi,DalmauPearson}. 
All other trees with at most 19 vertices satisfy $\TS n$ for all $n$. For a tree $\T$ the vertices of the indicator digraph for $\TS{2\cdot |V(\T)|}$ correspond to the nonempty subsets of $T$. Hence, the indicator structure of a tree $\T$ with 19 vertices has $2^{19}-1=524287$ vertices. Using level-wise satisfiability (see Section~\ref{sect:AC-Pol}) the number of vertices of the indicator structure is reduced to something between 19 and 513, depending on the number of vertices on each level. 

\subsection{Open Trees}\label{sec:chaTreesOpenTreesNLConjecture}
\label{sect:opentrees}
In this section we present trees that are interesting test cases, 
in particular regarding the conjectured classification of digraphs in NL (Conjecture~\ref{conj:NL}).

\subsubsection{A Tree not Known to be in NL}

We found a tree with polymorphisms that form a Kearnes-Kiss chain of length five and a Hobby-McKenzie chain of length 2, but has no majority polymorphism, and no (level-wise) J\'onsson chain of length 1000 (see Figure~\ref{fig:KKtree}). This tree is neither known to be P-hard or Mod$_p$L-hard, nor is it known to be in NL.  
It is the smallest tree without a majority polymorphism. 
Note that the existence of a J\'onsson chain of \emph{some} length is decidable because for a given digraph there are only finitely many operations of arity three. Moreover, by the discussion from Section~\ref{sect:level-wise} we know that we may narrow down the set of operations that have to be considered; the resulting number of operations is $12^{36}$. Even if we could show that the tree has no J\'onsson chain
we would not know that the tree is not in NL. We believe that the tree is in NL, but new ideas are needed to prove that (e.g., ideas to prove Conjecture~\ref{conj:NL}).  We mention that it can pp-construct $\Ord$, so it is NL-hard.

\def\scale{0.4}
\def\hdist{1cm}
\def\vdist{1cm}
\begin{figure}
\centering\begin{tikzpicture}[scale=\scale]
\node[bullet,gray] (12) at (1.5,0) {};
\node[bullet] (1) at (0.5,1) {};
\node[bullet] (10) at (3,1) {};
\node[bullet] (11) at (1.5,1) {};
\node[bullet] (2) at (0,2) {};
\node[bullet] (0) at (1,2) {};
\node[bullet] (7) at (3,2) {};
\node[bullet] (13) at (2,2) {};
\node[bullet] (3) at (0,3) {};
\node[bullet] (6) at (1,3) {};
\node[bullet] (8) at (3,3) {};
\node[bullet] (14) at (2,3) {};
\node[bullet] (4) at (0,4) {};
\node[bullet] (9) at (3,4) {};
\node[bullet] (15) at (2,4) {};
\node[bullet,gray] (5) at (0,5) {};
\path[->,>=stealth']
(1) edge (2)
(1) edge (0)
(2) edge (3)
(3) edge (4)
(4) edge[gray] (5)
(0) edge (6)
(7) edge (6)
(7) edge (8)
(10) edge (7)
(8) edge (9)
(11) edge (0)
(11) edge (13)
(12) edge[gray] (11)
(13) edge (14)
(14) edge (15)
;
\node at (1.5,-1) {Tree D};
\end{tikzpicture}
\caption{The smallest tree without a majority polymorphism (unique up to edge reversal; it has 16 vertices). 
It satisfies $\KK 5$ but does not satisfy $\J{1000}$. Therefore, Conjecture~\ref{conj:NL} puts it in NL but we cannot prove this fact; it is an interesting open case. 
}
\label{fig:KKtree}
\end{figure}

\subsubsection{Trees that might be P-hard}
There are 28 trees with 18 vertices that satisfy neither $\HMcK{1000}$ nor $\KK{1000}$, not even level-wise (see Figure~\ref{fig:noKKtrees}). 
They satisfy $\TS n$ for all $n$, so they are in P and cannot pp-construct $\TLinP$ for any $p$. Hence, this is in accordance with Conjecture~\ref{conj:bulin}. All other trees with up to 18 vertices satisfy $\KK 5$ and are in NL assuming Conjecture~\ref{conj:NL}. 
Hence, if this conjecture is true, and if $\operatorname{NL} \neq \operatorname{P}$, and if indeed these 28 trees do not have $\HMcK n$ for any $n$, then they are the smallest trees that are P-hard.

\input{treesNoKK}

\subsection{Majority Polymorphisms}
\label{sect:other} 
Majority polymorphisms play a central role in the early theory of the constraint satisfaction problem~\cite{FederVardi,FederCycles,JeavonsClosure},
in graph theory~\cite{Kazda,HellRafiey-list-homomorphism-digraphs},
and in the algebraic theory of CSPs~\cite{Bulatov-Conservative-Revisited,BulatovFVConjecture}. Dalmau and Krokhin~\cite{DalmauKrokhin08} proved that structures with a majority polymorphism are  in NL already before the mentioned result of Barto, Kozik, and Willard for near-unanimity polymorphisms~\cite{BartoKozikWillard}.
We have therefore also computed a
smallest tree without a majority polymorphism (see Figure~\ref{fig:KKtree}). Interestingly, when solving the indicator problem for the existence of a majority  polymorphism of $\Hb$ for graphs with at most 15 vertices (which all have a majority polymorphism), 
no backtracking was needed: pruning with the arc-consistency procedure sufficed to avoid all dead-ends in the search. Theoretical results only guarantee this behavior for establishing $(2,3)$-consistency (since $\Hb$ has a majority polymorphism). 
So one might ask: can every tree with a majority polymorphism be solved by arc consistency? 
This is not the case; see Lemmata 4.1 and 4.2 in~\cite{SpecialTriads}. 
In our experiments we found the smallest such tree: 
Figure~\ref{fig:withWNU3withoutWNU2} shows a tree with a majority polymorphism which does
not even have a binary symmetric polymorphism and therefore, in particular, cannot be solved by arc consistency.

\subsection{Code Availability and Reproducibility}

All of the code we used can be found at 
\[\text{\codelink}.\] There you also find a list with all core trees with at most 20 vertices and a list with all trees occurring in figures in this section. The trees are represented as lists of edges.


If you implement your own algorithm to generate core trees you can compare the numbers with Table~\ref{table:otrees} (comparing the actual trees in the lists will not work as there are too many).
If you want to verify the satisfiability of various polymorphism conditions independently of our implementation, we recommend the software package PCSP Tools by Opršal \cite{PCSPTools}.

\section{The Smallest NP-hard Cycle}
\begin{figure}
    \centering
    \begin{tikzpicture}[scale=0.5]
\node[bullet] (9) at (-1,0) {};
\node[bullet] (11) at (1,0) {};
\node[bullet] (23) at (4,0) {};
\node[bullet] (0) at (1,1) {};
\node[bullet] (10) at (0,1) {};
\node[bullet] (12) at (2,1) {};
\node[bullet] (14) at (4,1) {};
\node[bullet] (24) at (3,1) {};
\node[bullet] (8) at (-1,1) {};
\node[bullet] (22) at (5,1) {};
\node[bullet] (1) at (1,2) {};
\node[bullet] (5) at (0,2) {};
\node[bullet] (13) at (3,2) {};
\node[bullet] (15) at (4,2) {};
\node[bullet] (19) at (5,2) {};
\node[bullet] (25) at (2,2) {};
\node[bullet] (7) at (-1,2) {};
\node[bullet] (21) at (6,2) {};
\node[bullet] (2) at (1,3) {};
\node[bullet] (6) at (-1,3) {};
\node[bullet] (16) at (4,3) {};
\node[bullet] (20) at (6,3) {};
\node[bullet] (4) at (0,3) {};
\node[bullet] (18) at (5,3) {};
\node[bullet] (3) at (1,4) {};
\node[bullet] (17) at (4,4) {};
\path[->,>=stealth']
(0) edge (1)
(0) edge (25)
(1) edge (2)
(2) edge (3)
(5) edge (6)
(5) edge (4)
(9) edge (10)
(9) edge (8)
(11) edge (12)
(11) edge (10)
(12) edge (13)
(14) edge (15)
(14) edge (13)
(15) edge (16)
(16) edge (17)
(19) edge (20)
(19) edge (18)
(23) edge (24)
(23) edge (22)
(24) edge (25)
(4) edge (3)
(7) edge (6)
(8) edge (7)
(18) edge (17)
(21) edge (20)
(22) edge (21)
;
\node at (2.5,-1) {Cycle A};
\end{tikzpicture}
    \caption{The smallest NP-hard orientation of a cycle (up to edge reversal and assuming $\operatorname{P} \neq\operatorname{NP}$); it has 26 vertices.}
    \label{fig:smallestNPHardCycle}
\end{figure}
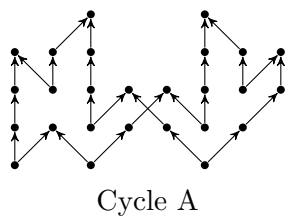

Most of the code we wrote is not specific for trees. Hence, we tried to use it to answer more questions about finite digraphs. We used the program to find the smallest NP-hard orientation of a cycle. It is displayed in Figure~\ref{fig:smallestNPHardCycle}; it is unique up to edge reversal and has 26 vertices. 
To find this cycle we used the following well-known theorem.
\begin{theorem}
Let $\C$ be an orientation of a cycle. If $\C$ is not balanced, then $\cocsp(\C)$ is in linear Datalog, in particular $\csp(\C)$ is in NL.
\end{theorem}
\begin{proof}
If $\C$ is unbalanced, then Theorem~3.1 in \cite{HellOrientedCycles} by Hell and Zhu implies that $\C$ has bounded pathwidth duality. Hence, Theorem~\ref{thm:DalmauBoundedPathwidthDualityIffLinearDL} implies that $\cocsp(\C)$ is in linear Datalog. 
\end{proof}
Using this theorem we can restrict our search to balanced cycles. A balanced cycle with $2n$ vertices can be represented by a string consisting of $n$ many 1s and $n$ many 0s. For each such cycle there are at most $4n$ different string representations: for each of the $2n$ vertices we can start reading off edge directions either clockwise or counterclockwise. Hence, there are roughly 
\[\frac{1}{4n}\cdot{\binom{2n}{n}}\]
balanced cycles with $2n$ vertices up to isomorphism. For $n=13$ we obtain that there are about 200000 balanced cycles with 26 vertices. Note that all of these cycles can easily be enumerated by a computer and that we can use level-wise satisfiability again. Recall that for trees we had to check about 800000 trees with 20 vertices for the existence of a $\KMM$ polymorphism. Hence, our program could compute the smallest balanced cycle without $\KMM$ polymorphisms shown in Figure~\ref{fig:smallestNPHardCycle} by brute force.

\section[Concluding Remarks]{Concluding Remarks and Directions for Further Research}\label{sect:discussion}
The following conjecture is implied by 
Conjecture~\ref{conj:bulin}, but might be easier to answer.

\begin{conjecture}
A tree has Kearnes-Markovi\'{c}-McKenzie polymorphisms if and only if it has a 3-wnu polymorphism.
\end{conjecture}

\begin{question}
Is it true that the probability that a tree drawn uniformly at random from the set of all trees with vertex set $\{1,\dots,n\}$ is NP-hard tends to 1 as $n$ tends to infinity? The answer is yes if we ask the question for random labelled digraphs instead of random labelled trees~\cite{LuczakNesetril}. 
\end{question}

Figure~\ref{fig:HMtrees} suggests that the following conjecture is true.

\begin{conjecture}
The fraction of core trees with $n$ vertices that are $\NL$-hard goes to 1 as $n$ goes to infinity.
\end{conjecture}

One possible step to answer this conjecture is to solve the following problem.
\begin{question}
Let $\S,\T$ be core orientations of trees such that $\S\ppleq\Ord$ and $\S\hookrightarrow\T$. Is it true that $\T\ppleq\Ord$?
\end{question}

\begin{question}
Determine the smallest trees that are P-hard (assuming that $\operatorname{NL}\neq\operatorname{P}$). We know from Section~\ref{sect:other} that they must have at least 16 vertices, since all smaller trees have a majority polymorphism and thus are in NL. 
\end{question}

 \begin{question}
Is our algorithm from Section~\ref{sect:gen} to generate unlabeled core trees a polynomial-delay enumeration algorithm (in the sense of~\cite{JohnYannaPapaGen})?
\end{question}
 
\begin{question}
Characterize minor conditions that can be tested level-wise (in the sense of Section~\ref{sec:levelWiseSatisfyability}) for balanced digraphs, and more specifically, for trees.
\end{question}

It would be interesting to perform experiments similar to the experiments presented here for trees that are equipped with a singleton unary relation $\{a\}$ for each vertex $a$ of the tree; in this case, if $\T_c$ is the resulting expanded tree structure, $\csp(\T_c)$ models the so-called \emph{$\T$-precoloring extension problem}. This setting is particularly nice from the algebraic perspective because 
 then all the polymorphisms of $\T_c$ 
 are idempotent. Note, however, that all these structure $\T_c$ are cores, so there are far more structures to consider, and hardness will dominate more rapidly.

Taking this one step further, it would also be interesting to study the so-called \emph{list homomorphism problem} for trees $\T$ from an experimental perspective. Here, the input contains, besides the graph $\G$, a set (also commonly referred to as a list) of vertices from $\Hb$ for every vertex of $\G$. And we are looking for a homomorphism from $\G$ to $\Hb$ that maps each vertex to an element from its set. This can be seen as a special case of a CSP for a relational structure, which contains besides the edge relation also a unary relation for each subset of the vertices of $\Hb$. On the algebraic side, we are therefore interested in polymorphisms that preserve all subsets of $\Hb$; such polymorphisms (and consequently the respective CSPs) are also called \emph{conservative}. The algorithms and complexities for conservative CSPs are better understood than the general case~\cite{Conservative, Barto-Conservative,Bulatov-Conservative-Revisited,Kazda-binary-conservative}, which will help to determine the complexity of the list homomorphism for trees. On the other hand, as in the case of the precoloring extension problem we have much larger numbers of trees to consider since all the structures that we study are already cores.


\clearpage

\backmatter

\bibliographystyle{plainurl} 
\bibliography{global.bib}

\def\cprime{$'$} \def\cprime{$'$} \def\cprime{$'$}
\begin{thebibliography}{10}

\bibitem{AfratiCosmadakis}
Foto~N. Afrati and Stavros~S. Cosmadakis.
\newblock Expressiveness of restricted recursive queries (extended abstract).
\newblock In {\em Proceedings of the 21st Annual {ACM} Symposium on Theory of
  Computing ({STOC})}, pages 113--126, 1989.
\newblock \href {https://doi.org/10.1145/73007.73018}
  {\path{doi:10.1145/73007.73018}}.

\bibitem{Aigner1967}
Martin Aigner.
\newblock On the linegraph of a directed graph.
\newblock {\em Mathematische Zeitschrift}, 102:56--61, 1967.
\newblock URL: \url{http://eudml.org/doc/170859}.

\bibitem{BartoKozikWillard}
L.~Barto, M.~Kozik, and R.~Willard.
\newblock Near unanimity constraints have bounded pathwidth duality.
\newblock In {\em Proceedings of the 27th ACM/IEEE Symposium on Logic in
  Computer Science (LICS)}, pages 125--134, 2012.

\bibitem{Barto-Conservative}
Libor Barto.
\newblock The dichotomy for conservative constraint satisfaction problems
  revisited.
\newblock In {\em Proceedings of the Symposium on Logic in Computer Science
  (LICS)}, Toronto, Canada, 2011.

\bibitem{Barto-cd}
Libor Barto.
\newblock Finitely related algebras in congruence distributive varieties have
  near unanimity terms.
\newblock {\em Canadian Journal of Mathematics}, 65(1):3--21, 2013.

\bibitem{BartoB13}
Libor Barto and Jakub Bul\'in.
\newblock {CSP} dichotomy for special polyads.
\newblock {\em Int. J. Algebra Comput.}, 23(5):1151--1174, 2013.
\newblock \href {https://doi.org/10.1142/S0218196713500215}
  {\path{doi:10.1142/S0218196713500215}}.

\bibitem{Barto_2021FreeStructures}
Libor Barto, Jakub Bul{\'{\i}}n, Andrei Krokhin, and Jakub Opr{\v{s}}al.
\newblock Algebraic approach to promise constraint satisfaction.
\newblock {\em Journal of the {ACM}}, 68(4):1--66, jul 2021.
\newblock \href {https://doi.org/10.1145/3457606} {\path{doi:10.1145/3457606}}.

\bibitem{BoundedWidthJournal}
Libor Barto and Marcin Kozik.
\newblock Constraint satisfaction problems solvable by local consistency
  methods.
\newblock {\em Journal of the {ACM}}, 61(1):3:1--3:19, 2014.

\bibitem{absorption}
Libor Barto and Marcin Kozik.
\newblock Absorption in universal algebra and {CSP}.
\newblock In {\em The Constraint Satisfaction Problem: Complexity and
  Approximability}, volume~7 of {\em Dagstuhl Follow-Ups}, pages 45--77, 2017.

\bibitem{SpecialTriads}
Libor Barto, Marcin Kozik, Mikl\'os Mar\'oti, and Todd Niven.
\newblock {CSP} dichotomy for special triads.
\newblock {\em Proceedings of the American Mathematical Society},
  137(9):2921--2934, 2009.

\bibitem{SpecialTriadsErratum}
Libor Barto, Marcin Kozik, Mikl\'os Mar\'oti, and Todd Niven.
\newblock Erratum to: {CSP} dichotomy for special triads.
\newblock Available from the website of the first author, 2009.

\bibitem{Barto_modularity}
Libor Barto, Marcin Kozik, Ralph Mckenzie, and Todd Niven.
\newblock Congruence modularity implies cyclic terms for finite algebras.
\newblock {\em Algebra {U}niversalis}, 61(3):365–380, 2009.

\bibitem{BartoKozikNiven}
Libor Barto, Marcin Kozik, and Todd Niven.
\newblock The {CSP} dichotomy holds for digraphs with no sources and no sinks
  (a positive answer to a conjecture of {B}ang-{J}ensen and {H}ell).
\newblock {\em SIAM Journal on Computing}, 38(5), 2009.

\bibitem{Pol}
Libor Barto, Andrei~A. Krokhin, and Ross Willard.
\newblock Polymorphisms, and how to use them.
\newblock In {\em The Constraint Satisfaction Problem: Complexity and
  Approximability}, pages 1--44. Schloss Dagstuhl - Leibniz-Zentrum fuer
  Informatik, 2017.

\bibitem{wonderland}
Libor Barto, Jakub Opr\v{s}al, and Michael Pinsker.
\newblock The wonderland of reflections.
\newblock {\em Israel Journal of Mathematics}, 223(1):363--398, 2018.

\bibitem{birkhoff1940lattice}
Garrett Birkhoff.
\newblock {\em Lattice theory}, volume~25.
\newblock American Mathematical Soc., 1940.

\bibitem{theBodirsky}
Manuel Bodirsky.
\newblock {\em Complexity of Infinite-Domain Constraint Satisfaction}.
\newblock Lecture Notes in Logic. Cambridge University Press, 2021.
\newblock \href {https://doi.org/10.1017/9781107337534}
  {\path{doi:10.1017/9781107337534}}.

\bibitem{BodirskyBodorUIP}
Manuel Bodirsky and Bertalan Bodor.
\newblock Canonical polymorphisms of ramsey structures and the unique
  interpolation property.
\newblock In {\em 36th Annual ACM/IEEE Symposium on Logic in Computer Science
  (LICS)}, pages 1--13, 2021.
\newblock \href {https://doi.org/10.1109/LICS52264.2021.9470683}
  {\path{doi:10.1109/LICS52264.2021.9470683}}.

\bibitem{BodirskyBulinStarkeWernthalerSmallestHardTrees}
Manuel Bodirsky, Jakub Bulín, Florian Starke, and Michael Wernthaler.
\newblock The smallest hard trees.
\newblock {\em Constraints}, 2023.
\newblock \href {https://doi.org/10.1007/s10601-023-09341-8}
  {\path{doi:10.1007/s10601-023-09341-8}}.

\bibitem{Qualitative-Survey}
Manuel Bodirsky and Peter Jonsson.
\newblock A model-theoretic view on qualitative constraint reasoning.
\newblock {\em Journal of Artificial Intelligence Research}, 58:339--385, 2017.

\bibitem{BodirskyStarke2022}
Manuel Bodirsky and Florian Starke.
\newblock Maximal digraphs with respect to primitive positive constructability.
\newblock {\em Combinatorica}, 42:997 -- 1010, 2021.
\newblock \href {https://doi.org/10.1007/s00493-022-4918-1}
  {\path{doi:10.1007/s00493-022-4918-1}}.

\bibitem{StarkeVucajBodirskySmoothDigraphs}
Manuel Bodirsky, Florian Starke, and Albert Vucaj.
\newblock Smooth digraphs modulo primitive positive constructability and cyclic
  loop conditions.
\newblock {\em International Journal of Algebra and Computation},
  31(05):929--967, 2021.

\bibitem{PPPoset}
Manuel Bodirsky and Albert Vucaj.
\newblock Two-element structures modulo primitive positive constructability.
\newblock {\em Algebra {U}niversalis}, 81(20), 2020.
\newblock Preprint available at ArXiv:1905.12333.

\bibitem{VucajZhukBodirsky21ThreeElements}
Manuel Bodirsky, Albert Vucaj, and Dmitriy Zhuk.
\newblock The lattice of clones of self-dual operations collapsed.
\newblock {\em International Journal of Algebra and Computation},
  33(04):717--749, 2023.
\newblock \href {https://doi.org/10.1142/S0218196723500327}
  {\path{doi:10.1142/S0218196723500327}}.

\bibitem{Bodnarchuk1969GaloisTF}
Victor Bodnarchuk, L.~A. Kaluzhnin, V.~N. Kotov, and Boris~A. Romov.
\newblock Galois theory for post algebras. {I} and {II}.
\newblock {\em Cybernetics}, 5:243--252, 1969.

\bibitem{Conservative}
Andrei~A. Bulatov.
\newblock Tractable conservative constraint satisfaction problems.
\newblock In {\em Proceedings of the Symposium on Logic in Computer Science
  {(LICS)}}, pages 321--330, Ottawa, Canada, 2003.

\bibitem{Bulatov-Conservative-Revisited}
Andrei~A. Bulatov.
\newblock Conservative constraint satisfaction re-revisited.
\newblock {\em Journal Computer and System Sciences}, 82(2):347--356, 2016.
\newblock ArXiv:1408.3690.
\newblock \href {https://doi.org/10.1016/j.jcss.2015.07.004}
  {\path{doi:10.1016/j.jcss.2015.07.004}}.

\bibitem{BulatovFVConjecture}
Andrei~A. Bulatov.
\newblock A dichotomy theorem for nonuniform {CSP}s.
\newblock In {\em 58th {IEEE} Annual Symposium on Foundations of Computer
  Science ({FOCS})}, pages 319--330, 2017.

\bibitem{Bulin18}
Jakub Bul\'in.
\newblock On the complexity of {$H$}-coloring for special oriented trees.
\newblock {\em Eur. J. Comb.}, 69:54--75, 2018.
\newblock \href {https://doi.org/10.1016/j.ejc.2017.10.001}
  {\path{doi:10.1016/j.ejc.2017.10.001}}.

\bibitem{BulinDelicJacksonNiven}
Jakub Bul{\'i}n, Dejan Delic, Marcel Jackson, and Todd Niven.
\newblock A finer reduction of constraint problems to digraphs.
\newblock {\em Log. Methods Comput. Sci.}, 11(4), 2015.

\bibitem{CarvalhoDalmauKrokhin}
Catarina Carvalho, V\'ictor Dalmau, and Andrei Krokhin.
\newblock {CSP} duality and trees of bounded pathwidth.
\newblock {\em Theoretical Computer Science}, 411:3188--3208, 2010.

\bibitem{CarvalhoEgriJacksonNiven}
Catarina Carvalho, L{\'{a}}szl{\'{o}} Egri, Marcel Jackson, and Todd Niven.
\newblock On {M}altsev digraphs.
\newblock {\em Electr. J. Comb.}, 22(1):P1.47, 2015.

\bibitem{MetaChenLarose}
Hubie Chen and Beno{\^{\i}}t Larose.
\newblock Asking the metaquestions in constraint tractability.
\newblock {\em {Association for Computing Machinery}}, 9(3):11:1--11:27, 2017.

\bibitem{ChenMengel}
Hubie Chen and Stefan Mengel.
\newblock A trichotomy in the complexity of counting answers to conjunctive
  queries.
\newblock In {\em 18th International Conference on Database Theory ({ICDT})},
  pages 110--126, 2015.
\newblock \href {https://doi.org/10.4230/LIPIcs.ICDT.2015.110}
  {\path{doi:10.4230/LIPIcs.ICDT.2015.110}}.

\bibitem{Dalmau}
V\'ictor Dalmau.
\newblock Computational complexity of problems over generalized formulas.
\newblock Ph{D}-thesis at the Departament de Llenguatges i Sistemes
  Inform\'{a}tics at the Universitat Polit\'{e}cnica de Catalunya, 2000.

\bibitem{Dalmau_2005LinearDatalog}
V\'ictor Dalmau.
\newblock Linear {D}atalog and bounded path duality of relational structures.
\newblock {\em Logical Methods in Computer Science}, 1(1), 2005.
\newblock \href {https://doi.org/10.2168/lmcs-1(1:5)2005}
  {\path{doi:10.2168/lmcs-1(1:5)2005}}.

\bibitem{DalmauKrokhin08}
V{\'{\i}}ctor Dalmau and Andrei~A. Krokhin.
\newblock Majority constraints have bounded pathwidth duality.
\newblock {\em European Journal of Combinatorics}, 29(4):821--837, 2008.
\newblock \href {https://doi.org/10.1016/j.ejc.2007.11.020}
  {\path{doi:10.1016/j.ejc.2007.11.020}}.

\bibitem{DalmauPearson}
V\'ictor Dalmau and Justin Pearson.
\newblock Closure functions and width 1 problems.
\newblock In {\em Proceedings of the International Conference on Principles and
  Practice of Constraint Programming (CP)}, pages 159--173, 1999.

\bibitem{EgriNLhardPaths}
L{\'a}szl{\'o} Egri.
\newblock On constraint satisfaction problems below {P}.
\newblock {\em Journal of Logic and Computation}, 12, 12 2013.
\newblock \href {https://doi.org/10.4230/LIPIcs.CSL.2011.203}
  {\path{doi:10.4230/LIPIcs.CSL.2011.203}}.

\bibitem{EgriLaroseTessonLogspace}
L{\'a}szl{\'o} Egri, Beno\^it Larose, and Pascal Tesson.
\newblock Symmetric {D}atalog and constraint satisfaction problems in logspace.
\newblock In {\em Proceedings of the Symposium on Logic in Computer Science
  ({LICS})}, pages 193--202, 2007.

\bibitem{EgristConNotinSymLinDL}
L{\'a}szl{\'o} Egri, Beno{\^i}t Larose, and Pascal Tesson.
\newblock Directed st-connectivity is not expressible in symmetric {D}atalog.
\newblock In {\em Automata, Languages and Programming}, pages 172--183.
  Springer Berlin Heidelberg, 2008.

\bibitem{FederCycles}
Tom\'as Feder.
\newblock Classification of homomorphisms to oriented cycles and of $k$-partite
  satisfiability.
\newblock {\em SIAM Journal on Discrete Mathematics}, 14(4):471--480, 2001.

\bibitem{FederVardi}
Tom\'as Feder and Moshe~Y. Vardi.
\newblock The computational structure of monotone monadic {SNP} and constraint
  satisfaction: {a} study through {D}atalog and group theory.
\newblock {\em {SIAM} Journal on Computing}, 28:57--104, 1999.

\bibitem{Jana}
Jana Fischer.
\newblock {CSP}s of orientations of trees.
\newblock Master thesis, TU Dresden, 2015.

\bibitem{FreeseV09}
Ralph Freese and Matthew Valeriote.
\newblock On the complexity of some {M}altsev conditions.
\newblock {\em International Journal of Algebra and Computation}, 19(1):41--77,
  2009.
\newblock \href {https://doi.org/10.1142/S0218196709004956}
  {\path{doi:10.1142/S0218196709004956}}.

\bibitem{GaultJeavonsImplementing}
Richard Gault and Peter Jeavons.
\newblock Implementing a test for tractability.
\newblock {\em Constraints}, 9(2):139--160, 2004.
\newblock \href {https://doi.org/10.1023/B:CONS.0000024049.41091.71}
  {\path{doi:10.1023/B:CONS.0000024049.41091.71}}.

\bibitem{Geiger}
David Geiger.
\newblock Closed systems of functions and predicates.
\newblock {\em Pacific Journal of Mathematics}, 27:95--100, 1968.

\bibitem{Gutjahr}
Wolfgang Gutjahr.
\newblock Graph colourings.
\newblock PhD Thesis, Free University Berlin, 1991.

\bibitem{GutjahrWW92}
Wolfgang Gutjahr, Emo Welzl, and Gerhard~J. Woeginger.
\newblock Polynomial graph-colorings.
\newblock {\em Discrete Applied Mathematics}, 35(1):29--45, 1992.
\newblock \href {https://doi.org/10.1016/0166-218X(92)90294-K}
  {\path{doi:10.1016/0166-218X(92)90294-K}}.

\bibitem{SantiagoDualityPairs}
Santiago Guzm{\'{a}}n{-}Pro and C{\'{e}}sar Hern{\'{a}}ndez{-}Cruz.
\newblock Duality pairs and homomorphisms to oriented and unoriented cycles.
\newblock {\em Electronic Journal of Combinatorics}, 2021.
\newblock \href {https://doi.org/10.37236/9747} {\path{doi:10.37236/9747}}.

\bibitem{HagemannMitschke}
J.~Hagemann and A.~Mitschke.
\newblock On $n$-permutable congruences.
\newblock {\em Algebra {U}niversalis}, pages 8--12, 1973.

\bibitem{HellNesetril}
Pavol Hell and Jaroslav Ne\v{s}et\v{r}il.
\newblock On the complexity of {H}-coloring.
\newblock {\em Journal of Combinatorial Theory, Series B}, 48:92--110, 1990.

\bibitem{cores}
Pavol Hell and Jaroslav Ne\v{s}et\v{r}il.
\newblock The core of a graph.
\newblock {\em Discrete Mathematics}, 109:117--126, 1992.

\bibitem{HellNZ96}
Pavol Hell, Jaroslav Ne\v{s}et\v{r}il, and Xuding Zhu.
\newblock Complexity of tree homomorphisms.
\newblock {\em Discret. Appl. Math.}, 70(1):23--36, 1996.
\newblock \href {https://doi.org/10.1016/0166-218X(96)00099-6}
  {\path{doi:10.1016/0166-218X(96)00099-6}}.

\bibitem{HNZ}
Pavol Hell, Jaroslav Ne\v{s}et\v{r}il, and Xuding Zhu.
\newblock Duality and polynomial testing of tree homomorphisms.
\newblock {\em TAMS}, 348(4):1281--1297, 1996.

\bibitem{HellRafiey-list-homomorphism-digraphs}
Pavol Hell and Arash Rafiey.
\newblock The dichotomy of list homomorphisms for digraphs.
\newblock In {\em Proceedings of the Twenty-Second Annual {ACM-SIAM} Symposium
  on Discrete Algorithms ({SODA})}, pages 1703--1713, 2011.
\newblock \href {https://doi.org/10.1137/1.9781611973082.131}
  {\path{doi:10.1137/1.9781611973082.131}}.

\bibitem{HellOrientedCycles}
Pavol Hell and Xuding Zhu.
\newblock The existence of homomorphisms to oriented cycles.
\newblock {\em SIAM Journal on Discrete Mathematics}, 8(2):208--222, 1995.
\newblock \href {https://doi.org/10.1137/S0895480192239992}
  {\path{doi:10.1137/S0895480192239992}}.

\bibitem{HobbyMcKenzie}
David Hobby and Ralph McKenzie.
\newblock {\em The structure of finite algebras}, volume~76 of {\em
  Contemporary Mathematics}.
\newblock American Mathematical Society, 1988.

\bibitem{Immerman}
N.~Immerman.
\newblock {\em Descriptive Complexity}.
\newblock Graduate Texts in Computer Science, Springer, 1998.

\bibitem{NivenDigraphCSP}
Marcel Jackson, Tomasz Kowalski, and Todd Niven.
\newblock Complexity and polymorphisms for digraph constraint problems under
  some basic constructions.
\newblock {\em International Journal of Algebra and Computation}, 26, 09 2016.
\newblock \href {https://doi.org/10.1142/S0218196716500600}
  {\path{doi:10.1142/S0218196716500600}}.

\bibitem{JeavonsCohenGyssensATest}
Peter Jeavons, David Cohen, and Marc Gyssens.
\newblock A test for tractability.
\newblock In Eugene~C. Freuder, editor, {\em Principles and Practice of
  Constraint Programming --- CP96}, pages 267--281. Springer Berlin Heidelberg,
  1996.

\bibitem{JeavonsClosure}
Peter Jeavons, David Cohen, and Marc Gyssens.
\newblock Closure properties of constraints.
\newblock {\em Journal of the ACM}, 44(4):527--548, 1997.

\bibitem{JohnYannaPapaGen}
D.~S. Johnson, M.~Yannakakis, and C.~H. Papadimitriou.
\newblock On generating all maximal independent sets.
\newblock {\em Information Processing Letters}, 27(3):119--123, 1988.
\newblock \href {https://doi.org/10.1016/0020-0190(88)90065-8}
  {\path{doi:10.1016/0020-0190(88)90065-8}}.

\bibitem{Kazda}
Alexandr Kazda.
\newblock Maltsev digraphs have a majority polymorphism.
\newblock {\em European Journal of Combinatorics}, 32:390--397, 2011.

\bibitem{Kazda-n-permute}
Alexandr Kazda.
\newblock {$n$-permutability and linear {D}atalog implies symmetric {D}atalog}.
\newblock {\em {Logical Methods in Computer Science}}, {Volume 14, Issue 2},
  2018.
\newblock \href {https://doi.org/10.23638/LMCS-14(2:3)2018}
  {\path{doi:10.23638/LMCS-14(2:3)2018}}.

\bibitem{Kazda-binary-conservative}
Alexandr Kazda.
\newblock {CSP} for binary conservative relational structures.
\newblock {\em {Algebra {U}niversalis}}, {Volume 75, Issue 1}:75–84, 2019.
\newblock \href {https://doi.org/10.1007/s00012-015-0358-8}
  {\path{doi:10.1007/s00012-015-0358-8}}.

\bibitem{kazdaValerioteMultibraidedness}
Alexandr Kazda and Matt Valeriote.
\newblock Deciding some {M}altsev conditions in finite idempotent algebras.
\newblock {\em The Journal of Symbolic Logic}, 85:539–562, 04 2020.
\newblock \href {https://doi.org/10.1017/jsl.2019.73}
  {\path{doi:10.1017/jsl.2019.73}}.

\bibitem{KearnesMarkovicMcKenzie}
Keith~A. Kearnes, Petar Markovi\'c, and Ralph McKenzie.
\newblock Optimal strong {M}al'cev conditions for omitting type 1 in locally
  finite varieties.
\newblock {\em Algebra {U}niversalis}, 72(1):91--100, 2015.

\bibitem{Kozik-SLAC}
Marcin Kozik.
\newblock Weak consistency notions for all the {CSP}s of bounded width.
\newblock In {\em Proceedings of the 31st Annual ACM/IEEE Symposium on Logic in
  Computer Science ({LICS '16})}, page 633–641, New York, NY, USA, 2016.
  Association for Computing Machinery.
\newblock \href {https://doi.org/10.1145/2933575.2934510}
  {\path{doi:10.1145/2933575.2934510}}.

\bibitem{Maltsev-Cond}
Marcin Kozik, Andrei Krokhin, Matt Valeriote, and Ross Willard.
\newblock Characterizations of several {M}altsev conditions.
\newblock {\em Algebra {U}niversalis}, 73(3):205--224, 2015.
\newblock \href {https://doi.org/10.1007/s00012-015-0327-2}
  {\path{doi:10.1007/s00012-015-0327-2}}.

\bibitem{LLT}
Beno\^it Larose, Cynthia Loten, and Claude Tardif.
\newblock A characterisation of first-order constraint satisfaction problems.
\newblock {\em Logical Methods in Computer Science}, 3(4:6), 2007.

\bibitem{LaroseTesson}
Beno\^it Larose and Pascal Tesson.
\newblock Universal algebra and hardness results for constraint satisfaction
  problems.
\newblock {\em Theoretical Computer Science}, 410(18):1629--1647, 2009.

\bibitem{AbilityToCount}
Beno\^it Larose, Matt Valeriote, and L{\'a}szl{\'o} Z{\'a}dori.
\newblock Omitting types, bounded width and the ability to count.
\newblock {\em International Journal of Algebra and Computation}, 19(5), 2009.

\bibitem{LuczakNesetril}
T.~{\L}uczak and J.~Ne\v{s}et\v{r}il.
\newblock When is a random graph projective?
\newblock {\em Eur. Journal Comb.}, 27(7), 2006.

\bibitem{Mackworth}
A.~K. Mackworth.
\newblock Consistency in networks of relations.
\newblock {\em Artificial {I}ntelligence}, 8:99--118, 1977.

\bibitem{Markowsky}
George Markowsky.
\newblock Free completely distributive lattices.
\newblock {\em Proceedings of The American Mathematical Society - PROC AMER
  MATH SOC}, 74:227--227, 02 1979.
\newblock \href {https://doi.org/10.2307/2043137} {\path{doi:10.2307/2043137}}.

\bibitem{wnuf}
Mikl\'os Mar\'oti and Ralph McKenzie.
\newblock Existence theorems for weakly symmetric operations.
\newblock {\em Algebra Universalis}, 59:463--489, 12 2008.
\newblock \href {https://doi.org/10.1007/s00012-008-2122-9}
  {\path{doi:10.1007/s00012-008-2122-9}}.

\bibitem{nauty}
Brendan~D. McKay and Adolfo Piperno.
\newblock Practical graph isomorphism, {II}.
\newblock {\em Journal of Symbolic Computation}, 60:94--112, 2014.
\newblock \href {https://doi.org/10.1016/j.jsc.2013.09.003}
  {\path{doi:10.1016/j.jsc.2013.09.003}}.

\bibitem{minizinc}
Nicholas Nethercote, Peter~J. Stuckey, Ralph Becket, Sebastian Brand,
  Gregory~J. Duck, and Guido Tack.
\newblock Minizinc: Towards a standard {CP} modelling language.
\newblock In Christian Bessi{\`e}re, editor, {\em Principles and Practice of
  Constraint Programming -- CP 2007}, pages 529--543. Springer Berlin
  Heidelberg, 2007.

\bibitem{olsak-loop}
Miroslav Ol\v{s}\'{a}k.
\newblock Loop conditions.
\newblock {\em Algebra {U}niversalis}, 81(2), 2020.
\newblock Preprint arXiv:1701.00260.

\bibitem{olsak-strong}
Miroslav Ol\v{s}\'{a}k.
\newblock Loop conditions for strongly connected digraphs.
\newblock {\em International Journal of Algebra and Computation},
  30(03):467--499, 2020.
\newblock \href {https://doi.org/10.1142/S0218196720500083}
  {\path{doi:10.1142/S0218196720500083}}.

\bibitem{oprsal2018taylors}
Jakub Opr\v{s}al.
\newblock Taylor's modularity conjecture and related problems for idempotent
  varieties.
\newblock {\em Order}, 35(3):433--460, 11 2018.
\newblock \href {https://doi.org/10.1007/s11083-017-9441-4}
  {\path{doi:10.1007/s11083-017-9441-4}}.

\bibitem{PCSPTools}
Jakub Opršal.
\newblock {PCSP} {T}ools.
\newblock \url{https://github.com/jakub-oprsal/pcsptools}, 2022.

\bibitem{Sabin1994ContradictingCW}
Daniel Sabin and Eugene~C. Freuder.
\newblock Contradicting conventional wisdom in constraint satisfaction.
\newblock In {\em Principles and Practice of Constraint Programming}, pages
  10--20. Springer Berlin Heidelberg, 1994.

\bibitem{gecode}
Christian Schulte, Mikael~Z. Lagerkvist, and Guido Tack.
\newblock Gecode, generic constraint development environment, 2010.
\newblock URL: \url{http://www.gecode.org/}.

\bibitem{Siggers}
Mark~H. Siggers.
\newblock A strong {M}al'cev condition for varieties omitting the unary type.
\newblock {\em Algebra {U}niversalis}, 64(1):15--20, 2010.

\bibitem{Tatarko}
William Tatarko.
\newblock {CSP} over oriented trees.
\newblock Bachelor thesis, Charles University Prague, 2019.

\bibitem{103287}
Peter Taylor.
\newblock Algorithm for generating all unlabeled trees with $n$ nodes?
\newblock Computer Science Stack Exchange.
\newblock URL: https://cs.stackexchange.com/q/103287 (version: 2019-01-23).

\bibitem{vucajzhuk2023submaximal}
Albert Vucaj and Dmitriy Zhuk.
\newblock Submaximal clones over a three-element set up to minor-equivalence,
  2023.
\newblock \href {http://arxiv.org/abs/2304.12807} {\path{arXiv:2304.12807}}.

\bibitem{ZhukFVConjecture}
Dmitriy Zhuk.
\newblock A proof of csp dichotomy conjecture.
\newblock In {\em 2017 IEEE 58th Annual Symposium on Foundations of Computer
  Science (FOCS)}, pages 331--342, 2017.
\newblock \href {https://doi.org/10.1109/FOCS.2017.38}
  {\path{doi:10.1109/FOCS.2017.38}}.

\end{thebibliography}


\listoffigures 
\listofalgorithms

\cleardoublepage
\thispagestyle{empty}

I herewith declare that I have produced this thesis without the prohibited assistance of third parties and without making use of aids other than those specified; notions taken over directly or indirectly from other sources have been identified as such. This thesis has not previously been presented in identical or similar form to any other German or foreign examination board.

\flushright{Florian Starke \\
	Dresden, 22 January 2024\\
}

\end{document}